\documentclass[11pt, DIV10,a4paper]{article}

\usepackage{natbib}
\usepackage{comment}
\newcommand {\ctn}{\citet} 
\usepackage{graphicx,subfigure,amsmath,latexsym,amssymb}
\usepackage{float,epsfig,multirow,rotating,times}
\usepackage{upgreek,wrapfig}

\newcommand{\bb}{\boldsymbol{b}}

\newcommand{\btheta}{\boldsymbol{\theta}}

\newcommand{\bbeta}{\boldsymbol{\beta}}

\newcommand{\bphi}{\boldsymbol{\phi}}

\newcommand{\bmu}{\boldsymbol{\mu}}

\newcommand{\bOmega}{\boldsymbol{\Omega}}

\newcommand{\bV}{\boldsymbol{V}}

\newcommand{\bR}{\boldsymbol{R}}

\newcommand{\bU}{\boldsymbol{U}}

\newcommand{\bzero}{\boldsymbol{0}}

\newtheorem{theorem}{Theorem}

\newtheorem{lemma}[theorem]{Lemma}

\newtheorem{remark}[theorem]{Remark}

\newenvironment{proof}[1][Proof]{\textbf{#1.} }{\ \rule{0.5em}{0.5em}}

\numberwithin{equation}{section}
\numberwithin{algo}{section}
\numberwithin{table}{section}
\numberwithin{figure}{section}

\usepackage{setspace}
\usepackage[mathscr]{euscript}
\usepackage[margin=1in]{geometry}
\begin{document}


\title{\vspace{-0.8in}
On Classical and Bayesian Asymptotics in Stochastic Differential Equations with 
Random Effects having Mixture Normal Distributions}
\author{Trisha Maitra and Sourabh Bhattacharya\thanks{
Trisha Maitra is a postdoctoral researcher and Sourabh Bhattacharya 
is an Associate Professor in Interdisciplinary Statistical Research Unit, Indian Statistical
Institute, 203, B. T. Road, Kolkata 700108.
Corresponding e-mail: sourabh@isical.ac.in.}}
\date{\vspace{-0.5in}}
\maketitle%

\begin{abstract}
\ctn{Maud12} considered a system of stochastic differential equations ($SDE$s) in a random effects setup.
Under the independent and identical ($iid$) situation, and assuming normal distribution of the random effects, 
they established weak consistency of the maximum likelihood estimators ($MLE$s) of the population parameters
of the random effects.

In this article, respecting the increasing importance and versatility of normal mixtures
and their ability to approximate any standard distribution, we consider the 
random effects having mixture of normal distributions and prove asymptotic results 
associated with the $MLE$s in both  independent and identical ($iid$) and 
independent but not identical (non-$iid$) situations. Besides, we consider $iid$ and non-$iid$ setups 
under the Bayesian paradigm and establish posterior consistency and asymptotic normality 
of the posterior distribution of the population parameters, even when the number of mixture components is unknown
and treated as a random variable.

Although ours is an independent work, we later noted that \ctn{Maud16} also assumed the $SDE$ setup with normal mixture distribution
of the random effect parameters but considered only the $iid$ case and proved only weak consistency of the $MLE$ 
under an extra, strong assumption as opposed to strong consistency that we are able to prove without the extra assumption. 
Furthermore, they did not deal with asymptotic normality 
of $MLE$ or the Bayesian asymptotics counterpart which we investigate in details.  

Ample simulation experiments and application to a real, stock market data set reveal the importance and usefulness of our methods even for 
small samples.
\\[2mm]
{\it {\bf Keywords:} Asymptotic normality; Finite mixture of normals; Maximum likelihood estimator; 
Posterior consistency; Random effects; Stochastic differential equations.}
 
\end{abstract}

\section{Introduction}
\label{sec:intro}

Data pertaining to inter-individual variability and intra-individual variability with respect to continuous time can 
be modeled through systems of stochastic differential equations ($SDE$s) consisting of random effects. 
In this regard, \ctn{Maud12}, \ctn{Maitra14b}, \ctn{Maitra14a} investigate asymptotic inference in the context of systems 
of $SDE$s of the following form:
\begin{equation}
d X_i(t)=\phi_ib(X_i(t))dt+\sigma(X_i(t))dW_i(t),\quad\mbox{with}\quad X_i(0)=x^i,~i=1,\ldots,n,
\label{eq:sde1}
\end{equation}
where, for $i=1,\ldots,n$, the stochastic process $X_i(t)$ is assumed to be continuously observed 
on the time interval $[0,T_i]$ with $T_i>0$ known, and corresponding to the $i$-th process initial 
values $\{x^i;~i=1,\ldots,n\}$ are also assumed to be known. Here $\{\phi_i;~i=1,\ldots,n\}$ are random effect 
parameters independent of the Brownian motions $\{W_i(\cdot);~i=1,\ldots,n\}$. 
The above authors assume that $\phi_i$ are independently and identically distributed ($iid$) with common 
distribution $g(\varphi,\theta)d\nu(\varphi)$ where  $g(\varphi,\theta)$ is a density with respect to 
a dominating measure $\nu$ on $\mathbb R^d$ for all $\theta$ ($\mathbb R$ is the real line and $d$ is the dimension). 
Here the unknown parameter to be estimated is $\theta\in\Omega\subset\mathbb R^p$ ($p\geq 2d$). 
In particular, the above authors assume that $g(\varphi,\theta)$ is the Gaussian density with unknown
means and covariance matrix, which are to be learned from the data and the model for classical inference
(see \ctn{Maud12}, \ctn{Maitra14a}),
and from a combination of the data, model and the prior for Bayesian inference (see \ctn{Maitra14b}).
%
%
Statistically, the $i$-th process $X_i(\cdot)$ corresponds to the $i$-th individual 
and the corresponding random effect is $\phi_i$. 
The following conditions (see \ctn{Maud12}, \ctn{Maitra14b} and \ctn{Maitra14a}) that we assume 
ensure existence of solutions of (\ref{eq:sde1}):
\begin{itemize}
\item[(H1)] 
\begin{enumerate}
\item[(i)] $b(\cdot)$ and
$\sigma(\cdot)$ are $C^1$ (differentiable with continuous first derivative) on $\mathbb R$ 
satisfying $b^2(x)\leq K(1+x^2)$  
and $\sigma^2(x)\leq K(1+x^2)$
for all $x\in\mathbb R$, for some $K>0$. 
\item[(ii)] Almost surely for each $i\geq 1$, 
\[
\int_0^{T_i}\frac{b^2(X_i(s))}{\sigma^2(X_i(s))}ds<\infty.
\]
\end{enumerate}
\end{itemize}
\ctn{Maud12} show that the likelihood, depending upon $\theta$, admits a relatively simple form involving
the following sufficient statistics:
\begin{align}
&U_i=\int_0^{T_i}\frac{b(X_i(s))}{\sigma^2(X_i(s))}dX_i(s),\quad V_i=\int_0^{T_i}\frac{b^2(X_i(s))}{\sigma^2(X_i(s))}ds,
\quad i=1,\ldots,n.
\label{eq:sufficient}
\end{align}  
The exact likelihood is given by
\begin{equation}
L(\theta)=\prod_{i=1}^n\lambda_i(X_i,\theta),
\label{eq:likelihood1}
\end{equation}
where
\begin{equation}
\lambda_i(X_i,\theta)=\int_{\mathbb R}g(\varphi,\theta)\exp\left(\varphi U_i-\frac{\varphi^2}{2}V_i\right)d\nu(\varphi).
\label{eq:likelihood2}
\end{equation}

For the Gaussian distribution of $\phi_i$ with mean $\mu$ and variance $\omega^2$, that is, with 
$g(\varphi,\theta)d\nu(\varphi)\equiv N\left(\mu,\omega^2\right)$, 
it is easy to obtain the following form
of $\lambda_i(X_i,\theta)$ (see \ctn{Maud12}):
\begin{equation}
\lambda_i(X_i,\theta)=\frac{1}{\left(1+\omega^2V_i\right)^{1/2}}\exp\left[-\frac{V_i}{2\left(1+\omega^2V_i\right)}
\left(\mu-\frac{U_i}{V_i}\right)^2\right]\exp\left(\frac{U^2_i}{2V_i}\right),
\label{eq:likelihood3}
\end{equation}
where $\theta=(\mu,\omega^2)\in\Omega\subset\mathbb R\times\mathbb R^+$. 
As in \ctn{Maud12}, \ctn{Maitra14a} and \ctn{Maitra14b} here also we assume that
\begin{itemize}
\item[(H2)] $\Omega$ is compact. 
\end{itemize}
\ctn{Maud12} consider $x^i=x$ and $T_i=T$ for $i=1,\ldots,n$, so that the setup boils down to the $iid$ situation, 
and investigate asymptotic properties of the $MLE$ of $\theta$, providing proofs of consistency and asymptotic
normality. As an alternative, \ctn{Maitra14a} verify the regularity conditions of existing results in general 
setups provided in \ctn{Schervish95} and \ctn{Hoadley71} to prove asymptotic properties of the $MLE$ 
in this $SDE$ setup in both $iid$ and non-$iid$ cases. 
Here, by the non-$iid$ setup, we mean that the processes $X_i(\cdot)$ are independent, but not identical, which ensues
when we allow for unequal initial values $x^i$ and unequal time points $T_i$.  

Interestingly, the alternative way of verification
of existing general results allowed \ctn{Maitra14a} to come up with stronger results under weaker assumptions,
compared to \ctn{Maud12}.

\ctn{Maitra14b}, for the first time in the literature, established Bayesian asymptotic results 
$SDE$-based random effects model, for both $iid$ and non-$iid$ setups. Specifically, considering
prior distributions $\pi(\theta)$ of $\theta$, they established asymptotic properties of the corresponding posterior
\begin{equation}
\pi_n(\theta|X_1,\ldots,X_n)=
\frac{\pi(\theta)\prod_{i=1}^n\lambda_i(X_i|\theta)}
{\int_{\psi\in\Omega}\pi(\psi)\prod_{i=1}^n\lambda_i(X_i|\psi)d\psi}
\label{eq:posterior1}
\end{equation}
as the sample size $n$ tends to infinity, through the verification of 
regularity conditions existing in \ctn{Choi04} and \ctn{Schervish95}.

In this article, we extend the asymptotic works of \ctn{Maitra14a} and \ctn{Maitra14b} assuming that 
the random effects are modeled by mixtures of normal distributions. The importance and generality
of such mixture models are briefly discussed in Section \ref{subsec:mixture_random_effects}. 





\subsection{Need for mixture distribution for the random effects parameters}
\label{subsec:mixture_random_effects}
The need for validation of our asymptotic results established in \ctn{Maitra14a} and \ctn{Maitra14b} 
for a bigger class of distributions corresponding to the random effect parameter $\phi$ leads us to 
consider mixture of normal distributions as the distribution of $\phi$, since any continuous 
density can be approximated arbitrarily accurately by an appropriate normal mixture.

Indeed, as is well-known, the unknown data generating process can be flexibly modeled by 
a mixture of parametric distributions. In fact, 
when a single parametric family can not provide a satisfactory model due to local variations in the observed data, 
mixture models can handle quite complex distributions by the appropriate choice of its components. 
So mixture distributions can be thought as the basis approximation to the unknown distributions. 
That mixture of normals with enough components can approximate any multivariate density, is established by \ctn{Norets}.
For instance, Theorems 1--3 of their article show under mild conditions that
when the number of mixture
components $M$ is fixed, for a given $\epsilon>0$ there exists $M$ such that the $L_1$
distance between the predictive density and the data generating
process density is less than $\epsilon$ almost surely. 

As infinite degree of smoothness and wide range of flexibility of mixture of normal densities allow us 
to model any unknown smooth density, it has been used for various inference problems including cluster analysis, 
density estimation, and robust estimation (see, for example, \ctn{Banfield}, \ctn{Lindsay}, \ctn{Roeder}). 
Mixture models have also been adopted in hazards models (\ctn{Louzada}), 
structural equation models (\ctn{Zhu}), analysis of proportions (\ctn{Brooks}), 
disease mapping (\ctn{Green}), neural networks (\ctn{Bishop}), and in many more areas. 
\ctn{Tony} provide a survey on recent developments and applications of normal mixture models in empirical finance.

Apart from the applicabilities of mixture models to various well-known statistical problems, it is worth noting 
that mixtures have important roles to play even in random effects models.  
Indeed, since mixture models are appropriate for data that are expected to arise 
from different sub-populations, these provide interesting extension of ordinary random effects models
to random effects models associated with multiple sub-populations. From the authors' perspective, 
this provides a particularly sound motivation for considering mixture distributions for the random effects. 

It is important to note that while we independently pursued our investigation on system of $SDE$s where the random effect parameters are samples from Gaussian
mixtures, we later came to know that it has also been considered by \ctn{Maud16}. However, asymptotically, their contribution is confined to only proving weak
consistency of the $MLE$ in the $iid$ setup. In contrast, we are able to prove strong and weak consistency of the $MLE$
in the $iid$ and non-$iid$ cases, respectively, asymptotic normality of the $MLE$ in both $iid$ and non-$iid$ setups,
and consistency and asymptotic normality of the Bayesian posterior distribution.  

For finite samples \ctn{Maud16} recommend
the standard $EM$ algorithm (\ctn{Dempster77}) for computation of $MLE$ and the standard Bayes Information Criterion ($BIC$) 
(see, for example, \ctn{Kass95}) for selecting the number of mixture components (see \ctn{Leroux92} and \ctn{Keribin00} for
$BIC$ related to mixture models). In the Bayesian setup it is natural to consider the number of mixture components to be
unknown and place a prior on this unknown quantity. This of course renders the problem variable-dimensional.
We develop the Bayesian asymptotic theory for this variable-dimensional setup as well. 

Bayesian model implementation in the variable-dimensional framework is a non-trivial exercise and it is now well-known that 
the traditional reversible jump Markov Chain Monte Carlo
method (\ctn{Green95}, \ctn{Richardson97}) can be quite inefficient for the purpose. 
Recently, \ctn{Das17} have proposed a novel and efficient 
transformation based variable-dimensional 
Markov chain Monte Carlo method to simulate from variable-dimensional distributions; they refer to this method as
Transdimensional Transformation based Markov chain Monte Carlo (TTMCMC) and demonstrate in particular the success of this
method in mixture problems with unknown number of components. Indeed, as we demonstrate in this article with simulation studies
and a real stock market data application, even in the $SDE$ setup, TTMCMC yields excellent performance. 

There is another issue to be remarked about, namely, the well-known label-switching problem associated with mixtures, which is
also persistent in our $SDE$ setup. Indeed, as is evident from (\ref{eq:likelihood5}) and (\ref{eq:likelihood6}),
the likelihood remains the same even if the labels of $\left\{(a_k,\beta_k):~k=1,\ldots,M\right\}$ are permuted, showing
non-identifiability of the likelihood with respect to label-switching of the parameter components. Under somewhat
restrictive assumption (see (H4) of \ctn{Maud16}), \ctn{Maud16}
show that if two such $SDE$ mixtures with two sets of parameters are the same, then the two sets of parameters are also the same.
However, that simply does not even touch the label-switching problem, and the weak consistency of $MLE$ proved
by \ctn{Maud16} in the $iid$ setup holds up to label-switching. In contrast, none of our results require 
assumption (H4) of \ctn{Maud16}, which are still unique up to label-switching. Specifically, in the $iid$ case we  
have been able to prove strong consistency of the $MLE$ under weaker assumption compared to weak consistency proved
by \ctn{Maud16} under stronger assumption.

The rest of our paper is structured as follows.
In Section \ref{sec:Like} we describe the likelihood corresponding to normal mixture distribution of 
the random effects, and in 
Sections \ref{sec:consistency_iid} and \ref{sec:consistency_non_iid} we investigate asymptotic properties of $MLE$ 
in the $iid$ and non-$iid$ contexts, respectively. In Sections \ref{sec:b_consistency_iid} and 
\ref{sec:b_consistency_non_iid} we investigate Bayesian
asymptotics in the $iid$  and non-$iid$ setups, respectively. In Section \ref{sec:vardim} we address the Bayesian asymptotic theory when the number of mixture
components is unknown and considered random.
Brief discussions on the validity of our asymptotic theory for non-multiplicative and multidimensional linear random effects are provided in Sections \ref{sec:nonmult}
and \ref{sec:multidim}, respectively, of the supplement, whose sections and figures have the prefix ``S-" when referred to in this paper.
In Section \ref{sec:simstudy} we provide a brief overview of our simulation studies, conducted in nine separate cases. 
The complete details of the simulation experiments, illustrating 
the usefulness of our Bayesian methodology with TTMCMC implementation, are presented in Section \ref{sec:simstudy_details} of the supplement.
In Section \ref{sec:realdata} we consider application to a real stock market data set.
We summarize our work and provide concluding remarks in Section \ref{sec:conclusion}. 
Additional technical details are provided in the supplement. 

Notationally, ``$\stackrel{a.s.}{\rightarrow}$", ``$\stackrel{P}{\rightarrow}$" and ``$\stackrel{\mathcal L}{\rightarrow}$"
denote convergence ``almost surely", ``in probability" and ``in distribution", respectively.

\section{Likelihood with respect to normal mixture distribution of the random effects}
\label{sec:Like}
We consider the following normal mixture form of the distribution of $\phi$: 
$$g(\varphi,\theta)d\nu(\varphi)\equiv
\sum_{k=1}^Ma_kN(\mu_k,\omega_k^2),$$ such that $a_k\geq 0$ for $k=1,\ldots,M$
and $\sum_{k=1}^Ma_k=1$. We assume that $\theta=(\gamma,\beta)$ where
$\gamma=(a_1,a_2,\ldots,a_M)$ and 
$\beta=(\beta_1,\ldots,\beta_M)$, such that $\beta_i=(\mu_i,\omega_i^2)$; $i=1,\ldots,M$. 
Here, for all $k=1,\ldots,M$, $(\mu_k,\omega^2_k)\in\Omega_{\beta}\subset\mathbb R\times\mathbb R^+$.
We assume that $\Omega_{\beta}$ is compact, so that, denoting the simplex by the compact space $\Omega_{\gamma}$ on which $\gamma$
lie, it follows that $\Omega=\Omega_{\gamma}\times\Omega_{\beta}$ is compact, in accordance with (H2). 
Our likelihood corresponding to the $i$-th individual is
\begin{equation}
\lambda_i(X_i,\theta)=\sum_{k=1}^M a_k f(X_i|\beta_k),
\label{eq:likelihood5}
\end{equation}
where
\begin{equation}
f(X_i|\beta_k)=\frac{1}{\left(1+\omega_k^2V_i\right)^{1/2}}\exp\left[-\frac{V_i}{2\left(1+\omega_k^2V_i\right)}
\left(\mu_k-\frac{U_i}{V_i}\right)^2\right]\exp\left(\frac{U^2_i}{2V_i}\right).
\label{eq:likelihood6}
\end{equation}
Assuming independence of the individuals conditional on the parameters, it then follows that the complete
likelihood is the product of (\ref{eq:likelihood5}) over the $n$ individuals, having the form
(\ref{eq:likelihood1}).

Note that, due to non-identifiability of the mixture form of the subject-wise likelihood (\ref{eq:likelihood5}), 
our asymptotic results will be unique up to label switching (see, for example, \ctn{Redner}).

\section{Consistency and asymptotic normality of $MLE$ in the $iid$ setup}
\label{sec:consistency_iid}

\subsection{Strong consistency of $MLE$}
\label{subsec:MLE_consistency_iid}

Consistency of the $MLE$ under the $iid$ setup can be verified through the verification of the regularity
conditions of the following theorem (Theorems 7.49 and 7.54 of \ctn{Schervish95}); for our purpose
we present the version for compact $\Omega$.
\begin{theorem}[\ctn{Schervish95}]
\label{theorem:theorem1}
Let $\{X_n\}_{n=1}^{\infty}$ be conditionally $iid$ given $\theta$ with density $\lambda_1(x|\theta)$
with respect to a measure $\nu$ on a space $\left(\mathcal X^1,\mathcal B^1\right)$. Fix $\theta_0\in\Omega$, and define,
for each $S\subseteq\Omega$ and $x\in\mathcal X^1$,
\[
Z(S,x)=\inf_{\psi\in S}\log\frac{\lambda_1(x|\theta_0)}{\lambda_1(x|\psi)}.
\]
Assume that for each $\theta\neq\theta_0$, there is an open set $N_{\theta}$ such that $\theta\in N_{\theta}$ and
that $E_{\theta_0}Z(N_{\theta},X_i)> -\infty$. 
Also assume that $\lambda_1(x|\cdot)$ is continuous at $\theta$ 
for every $\theta$, a.s. $[P_{\theta_0}]$. Then, if $\hat\theta_n$ is the $MLE$ of $\theta$ corresponding to $n$ observations, 
it holds that $\underset{n\rightarrow\infty}{\lim}~\hat\theta_n=\theta_0$, a.s. $[P_{\theta_0}]$.
\end{theorem}

\subsubsection{Verification of strong consistency of $MLE$ in our $SDE$ setup}
\label{subsubsec:MLE_consistency_iid}

Let 
 the true set of parameter be 
$\theta_0=(\gamma_0,\beta_0)$ where $\gamma_0=(a_{0,1},a_{0,2},\ldots,a_{0,M})$ and 
$\beta_0=(\beta_{0,1},\beta_{0,2},\ldots,\beta_{0,M})$, such that $\beta_{0,i}=(\mu_{0,i},\omega_{0,i}^2)$ 
for $i=1,\ldots,M$.

To verify the conditions of Theorem \ref{theorem:theorem1} in our case, we note that for any $x$,
$\lambda_1(x|\theta)=\sum_{k=1}^M a_k f(x|\beta_k)$ where $f(x|\beta_k)$ given by (\ref{eq:likelihood6}) 
is clearly continuous in $\beta_k$, implying continuity of $\lambda_1(x|\theta)$ in $\theta$. 
Also, it follows from the proof of Proposition 7 of \ctn{Maud12} that
\begin{align}
\log\frac{f(x|\beta_{k})}{f(x|\beta_{0,1})}
	&=\frac{1}{2}\log\left(\frac{1+\omega_{0,1}^2V}{1+\omega^2_{k}V}\right)
	+\frac{1}{2}\frac{(\omega^2_{k}-\omega_{0,1}^2)U^2}{(1+\omega_{0,1}^2V)(1+\omega^2_{k}V)}\notag\\
	&\quad+\frac{\mu_{0,1}^2V}{2(1+\omega_{0,1}^2V)}-\frac{\mu_{0,1} U}{1+\omega_{0,1}^2V}
-\left(\frac{\mu^2_{k}V}{2(1+\omega^2_{k}V)}-\frac{\mu_{k}U}{1+\omega^2_{k}V}\right)\notag\\
	&\leq \frac{1}{2}\left\{\log\left(1+\frac{\omega_{0,1}^2}{\omega^2_{k}}\right)+\frac{|\omega_{0,1}^2-\omega^2_{k}|}{\omega_{0,1}^2}\right\}
	+\frac{1}{2}|\omega^2_{k}-\omega_{0,1}^2|\left(\frac{U}{1+\omega^2_{k}V}\right)^2\left(1+\frac{\omega^2_{k}}{\omega_{0,1}^2}\right)\notag\\
	&\quad +|\mu_{0,1}|\left\vert\frac{U}{1+\omega^2_{k}V}\right\vert\left(1+\frac{|\omega^2_{k}-\omega_{0,1}^2|}{\omega_{0,1}^2}\right)
+\left|\frac{\mu^2_{k}V}{2(1+\omega^2_{k}V)}\right|+\left|\frac{\mu_{k}U}{1+\omega^2_{k}V}\right|\notag\\
	&\quad=C_1(U,V,\beta_{0,1},\beta_{k}),~\mbox{(say)}.
\label{eq:upper_bound1}
\end{align}
Hence,
\begin{equation}
	f(x|\beta_{k})\leq\exp(C_1(U,V,\beta_{0,1},\beta_{k}))f(x|\beta_{0,1}).
\label{eq:relation_1}
\end{equation}
Now, since $\log (x)=-\log\left(\frac{1}{x}\right)$, for $x>0$, we have
\begin{align}
\log\frac{\lambda_1(x|\theta_0)}{\lambda_1(x|\theta)}
&=-\log\left(\frac{\sum_{k=1}^M a_{k}f(x|\beta_{k})}{\sum_{k=1}^M a_{0,k} f(x|\beta_{0,k})}\right)\notag\\
&\geq -\log\left(\frac{\sum_{k=1}^Ma_{k}f(x|\beta_{k})}{a_{0,1} f(x|\beta_{0,1})}\right)\notag\\
	&\geq -\log\left(\frac{\sum_{k=1}^Ma_{k}\exp(C_1(U,V,\beta_{0,1},\beta_{k}))f(x|\beta_{0,1}))}{a_{0,1} f(x|\beta_{0,1})}\right)
\quad\mbox{(by (\ref{eq:relation_1}))}\notag\\
	&\geq -\left|\log\left(\sum_{k=1}^Ma_{k}\exp(C_1(U,V,\beta_{0,1},\beta_{k}))\right)\right|-\left|\log a_{0,1}\right|\notag\\
	&\geq -\sum_{k=1}^MC_1(U,V,\beta_{0,1},\beta_{k})-|\log a_{0,1}|.
	\label{eq:upper_bound2}
\end{align}

The last inequality is due to a result of \ctn{Natienza}, which we furnish below along with its proof: 
\begin{lemma}{\ctn{Natienza}}
\label{lemma:atienza}
$\sum_{k=1}^M a_k=1$, with $a_k\geq 0$ implies
$$\quad \left|\log\sum_{k=1}^M a_k f_k\right|\leq\sum_{k=1}^M|\log f_k|,$$
for $f_k>0$; $k=1,\ldots,M$.
\end{lemma}
\begin{proof}
	Letting $\alpha=\min\{a_k:k=1,\ldots,M\}$ and $\beta=\max\{a_k:k=1,\ldots,M\}$, it follows that
	$\alpha\leq \sum_{k=1}^M a_k f_k\leq\beta$, which again leads to 
	$\log \alpha\leq \log\left(\sum_{k=1}^M a_k f_k\right)\leq\log \beta$, by monotonicity of logarithm.
	Hence, $\left|\log\left(\sum_{k=1}^M a_k f_k\right)\right|\leq \max\{|\log \alpha|,|\log\beta|\}\leq \sum_{k=1}^M|\log f_k|$.
\end{proof}

Now taking $N_{\theta}$ to be any open subset of the relevant compact parameter space containing $\theta$, and
noting that 
$$E_{\theta_0}\left[\inf_{\theta\in N_{\theta}}\log\frac{\lambda_1(X|\theta_0)}{\lambda_1(X|\theta)}\right]\geq
-\sum_{k=1}^M\sup_{\theta\in N_{\theta}}C_1(U,V,\beta_{0,1},\beta_{k})-|\log a_1|,$$
it is sufficient to establish that 
$E_{\theta_0}\left[\underset{\theta\in N_{\theta}}{\sup}~C_1(U,V,\beta_{0,1},\beta_{k})\right] <+\infty$ for $k=1,\ldots,M$,
in order to conclude that $E_{\theta_0}Z(N_{\theta},X_i)>-\infty.$

Observation of the terms of $C_1(U,V,\beta_{0,1},\beta_{k})$ makes it clear that it is sufficient to show finiteness of 
$E_{\beta_{0,k}}\left(\frac{U}{1+\omega^2_{k}V}\right)^2$, 
$E_{\beta_{0,k}}\left\vert\frac{U}{1+\omega^2_{k}V}\right\vert$ and 
$E_{\beta_{0,k}}\left(\frac{V}{1+\omega^2_{k}V}\right)$.

Let us first prove finiteness of $E_{\beta_{0,k}}\left|\frac{U}{1+\omega^2_{k}V}\right|$. For $\omega^2_{k}>\omega^2_{0,k}$, note that
\begin{equation*}
E_{\beta_{0,k}}\left|\frac{U}{1+\omega^2_{k}V}\right|=E_{\beta_{0,k}}\left(\frac{|U|}{1+\omega^2_{0,k}\left(\frac{\omega^2_{k}}{\omega^2_{0,k}}\right)V}\right)
<E_{\beta_{0,k}}\left(\frac{|U|}{1+\omega^2_{0,k}V}\right)<\infty.
\end{equation*}
The last inequality holds due to Lemma 1 of \ctn{Maud12}.
For $\omega^2_{k}<\omega^2_{0,k}$,
\begin{equation*}
	E_{\beta_{0,k}}\left|\frac{U}{1+\omega^2_{k}V}\right|=\left(\frac{\omega^2_{0,k}}{\omega^2_k}\right)
	E_{\beta_{0,k}}\left(\frac{|U|}{\left(\frac{\omega^2_{0,k}}{\omega^2_k}\right)+\omega^2_{0,k}V}\right)
	<\left(\frac{\omega^2_{0,k}}{\omega^2_k}\right)E_{\beta_{0,k}}\left(\frac{|U|}{1+\omega^2_{0,k}V}\right)<\infty.
\end{equation*}
In the same way, for $r\geq 1$, 
\begin{equation}
E_{\beta_{0,k}}\left(\frac{|U|}{1+\omega^2_{k}V}\right)^r<\infty. 
\label{eq:maud1}
\end{equation}
Also, since $\frac{V}{1+\omega^2_{k}V}<\frac{1}{\omega^2_k}$
almost surely, it follows that for $r\geq 1$,
\begin{equation}
E_{\beta_{0,k}}\left(\frac{V}{1+\omega^2_{k}V}\right)^r<\infty.
\label{eq:maud2}
\end{equation}	

We have thus shown that $E_{\theta_0}Z(N_{\theta},X_i)> -\infty$.
Hence, $\hat\theta_n\stackrel{a.s.}{\rightarrow}\theta_0$ $[P_{\theta_0}]$.
We summarize the result in the form of the following theorem:
\begin{theorem}
\label{theorem:consistency_iid}
Assume the $iid$ setup and conditions (H1) and (H2). 
Then the $MLE$ is strongly consistent
in the sense that
$\hat\theta_n\stackrel{a.s.}{\rightarrow}\theta_0$~ $[P_{\theta_0}]$.
\end{theorem}

The above theorem requires only a lower bound for $\log \frac{\lambda_i(X_i|\theta_0)}{\lambda_i(X_i|\theta)}$ with finite expectation.
But some of our later results will require existence of expectations of an upper bound for $\left|\log \frac{\lambda_i(X_i|\theta_0)}{\lambda_i(X_i|\theta)}\right|$. 
We state and prove this result as Lemma \ref{lemma:upper_bound}.
\begin{lemma}
\label{lemma:upper_bound}
For any given $\theta\in\Omega$, $\left|\log \frac{\lambda_i(X_i|\theta_0)}{\lambda_i(X_i|\theta)}\right|$ has an upper bound which has finite moments of all orders.
\end{lemma}
\begin{proof}
Since (\ref{eq:upper_bound2}) provides a lower bound $\log \frac{\lambda_i(X_i|\theta_0)}{\lambda_i(X_i|\theta)}$ with finite moments of all orders,	
let us now provide an upper bound for $\log \frac{\lambda_i(X_i|\theta_0)}{\lambda_i(X_i|\theta)}$, which also has finite moments of all orders.
Since
\begin{align}
\log\frac{f(x|\beta_{0,k})}{f(x|\beta_{1})}
&=\frac{1}{2}\log\left(\frac{1+\omega_1^2V}{1+\omega^2_{0,k}V}\right)
+\frac{1}{2}\frac{(\omega^2_{0,k}-\omega_1^2)U^2}{(1+\omega_1^2V)(1+\omega^2_{0,k}V)}\notag\\
&\quad+\frac{\mu_1^2V}{2(1+\omega_1^2V)}-\frac{\mu_1 U}{1+\omega_1^2V}
-\left(\frac{\mu^2_{0,k}V}{2(1+\omega^2_{0,k}V)}-\frac{\mu_{0,k}U}{1+\omega^2_{0,k}V}\right)\notag\\
&\leq \frac{1}{2}\left\{\log\left(1+\frac{\omega_1^2}{\omega^2_{0,k}}\right)+\frac{|\omega_1^2-\omega^2_{0,k}|}{\omega_1^2}\right\}
+\frac{1}{2}|\omega^2_{0,k}-\omega_1^2|\left(\frac{U}{1+\omega^2_{0,k}V}\right)^2\left(1+\frac{\omega^2_{0,k}}{\omega_1^2}\right)\notag\\
&\quad +|\mu_1|\left\vert\frac{U}{1+\omega^2_{0,k}V}\right\vert\left(1+\frac{|\omega^2_{0,k}-\omega_1^2|}{\omega_1^2}\right)
+\left|\frac{\mu^2_{0,k}V}{2(1+\omega^2_{0,k}V)}\right|+\left|\frac{\mu_{0,k}U}{1+\omega^2_{0,k}V}\right|\notag\\
&\quad=C_1(U,V,\beta_1,\beta_{0,k}),\notag
\end{align}
it follows that
\begin{equation}
f(x|\beta_{0,k})\leq\exp(C_1(U,V,\beta_1,\beta_{0,k}))f(x|\beta_{1}).
\label{eq:relation_1_2}
\end{equation}
Now,
\begin{align}
\log\frac{\lambda_1(x|\theta_0)}{\lambda_1(x|\theta)}
&=\log\left(\frac{\sum_{k=1}^M a_{0,k}f(x|\beta_{0,k})}{\sum_{k=1}^M a_k f(x|\beta_k)}\right)\notag\\
&\leq\log\left(\frac{\sum_{k=1}^Ma_{0,k}f(x|\beta_{0,k})}{a_1 f(x|\beta_1)}\right)\notag\\
&\leq\left|\log\frac{\sum_{k=1}^Ma_{0,k}\exp(C_1(U,V,\beta_1,\beta_{0,k}))f(x|\beta_{1}))}{a_1 f(x|\beta_1)}\right|
\quad\mbox{(by (\ref{eq:relation_1_2}))}\notag\\
&\leq\left|\log\sum_{k=1}^Ma_{0,k}\exp(C_1(U,V,\beta_1,\beta_{0,k}))\right|+|\log a_1|\notag\\
&=\sum_{k=1}^MC_1(U,V,\beta_1,\beta_{0,k})+|\log a_1|\quad\mbox{(by Lemma \ref{lemma:atienza})}.
\label{eq:upper_bound2_2}
\end{align}
From (\ref{eq:upper_bound2}) and (\ref{eq:upper_bound2_2}) and existence of moments as shown in \ctn{Maud12}, as well as by the results (\ref{eq:maud1}) 
and (\ref{eq:maud2}),
it follows that $\left|\log \frac{\lambda_i(X_i|\theta_0)}{\lambda_i(X_i|\theta)}\right|$ 
has an upper bound which has finite moments of all orders under $\theta_0$, for any given $\theta\in\Omega$. 
\end{proof}

\subsection{Asymptotic normality of $MLE$}
\label{subsec:MLE_iid_normality}

To verify asymptotic normality of $MLE$ we invoke the following theorem provided in \ctn{Schervish95} (Theorem 7.63):
\begin{theorem}[\ctn{Schervish95}]
\label{theorem:theorem2}
Let $\Omega$ be a subset of $\mathbb R^{3M}$, and let $\{X_n\}_{n=1}^{\infty}$ be conditionally $iid$ given $\theta$
each with density $\lambda_1(\cdot|\theta)$. Let $\hat\theta_n$ be an $MLE$. Assume that
$\hat\theta_n\stackrel{P}{\rightarrow}\theta$ under $P_{\theta}$ for all $\theta$. Assume that $\lambda_1(x|\theta)$
has continuous second partial derivatives with respect to $\theta$ and that differentiation can be passed under the
integral sign in the sense that $0=\frac{\partial}{\partial\theta_j}\int\lambda_1(x|\theta)d\nu(x)
=\int\frac{\partial}{\partial\theta_j}\lambda_1(x|\theta)d\nu(x)$, for $j\geq 1$, where $\nu$ is the relevant dominating measure. 
Assume that there exists $H_r(x,\theta)$ such that, for each $\theta_0\in int(\Omega)$ and each
$k,j$,
\begin{align}
\sup_{\|\theta-\theta_0\|\leq r}\left\vert\frac{\partial^2}{\partial\theta_k\partial\theta_j}\log \lambda_{X_1|\Theta}(x|\theta_0)
-\frac{\partial^2}{\partial\theta_k\partial\theta_j}\log \lambda_{X_1|\Theta}(x|\theta)\right\vert\leq H_r(x,\theta_0),
\label{eq:h1}
\end{align}
with
\begin{equation}
\lim_{r\rightarrow 0}E_{\theta_0}H_r\left(X,\theta_0\right)=0.
\label{eq:h2}
\end{equation}
Assume that the Fisher information matrix $\mathcal I(\theta)$ is finite and non-singular. Then, under $P_{\theta_0}$,
\begin{equation}
\sqrt{n}\left(\hat\theta_n-\theta_0\right)\stackrel{\mathcal L}{\rightarrow}N\left(\bzero,\mathcal I^{-1}(\theta_0)\right). 
\label{eq:MLE_normality_iid}
\end{equation}
\end{theorem}

\subsubsection{Verification of the above regularity conditions for asymptotic normality in our $SDE$ setup}
\label{subsubsec:MLE_normality_iid}
In Section \ref{subsubsec:MLE_consistency_iid} we proved almost sure consistency of the $MLE$ $\hat\theta_n$
in the $SDE$ setup. Hence, $\hat\theta_n\stackrel{P}{\rightarrow}\theta$ under $P_{\theta}$ for all $\theta$.

We assume,
\begin{itemize}
\item[(H3)] For some real constant $K>0$, $$\frac{b^2(x)}{\sigma^2(x)}<K(1+x^{\tau}),\quad\mbox{for some}\quad \tau\geq 1$$

\end{itemize}
Note that (H3) implies moments of all orders of $V_i$ (given by \ref{eq:sufficient}) 
for all $i=1,\ldots,n$ are finite. 
Also note that, since for $k\geq 1$, $E[\phi_i]^{2k}<\infty$, for all $i=1,\ldots,n$, because of normal mixture
distribution, Proposition 1 of \ctn{Maud12} implies that for all 
$T>0,\underset{t\in[0,T]}{\sup}~E[X_i(t)]^{2k}<\infty$ for all $k\geq 1$ and $i=1,\ldots,n$. 

That differentiation can be passed under the integral sign in our case, is proved in Section \ref{sec:diff_int} of the supplement.
For the third order derivatives of $\log\lambda$ note that
\begin{equation}
\frac{\partial^3\log \lambda}{\partial\theta_r\partial\theta_s\partial\theta_t}=\frac{1}{\lambda}\frac{\partial^3\lambda}{\partial\theta_r\partial\theta_s\partial\theta_t}-\frac{1}{\lambda^2}\left[\frac{\partial^2\lambda}{\partial\theta_r\partial\theta_s}\frac{\partial\lambda}{\partial\theta_t}+\frac{\partial^2\lambda}{\partial\theta_s\partial\theta_t}\frac{\partial\lambda}{\partial\theta_r}+\frac{\partial^2\lambda}{\partial\theta_t\partial\theta_r}\frac{\partial\lambda}{\partial\theta_s}\right]
+\frac{2}{\lambda^3}\frac{\partial\lambda}{\partial\theta_r}\frac{\partial\lambda}{\partial\theta_s}\frac{\partial\lambda}{\partial\theta_t}.
\end{equation}

Denoting $\psi=(\mu_1,\cdots,\mu_M,\omega_1^2,\cdots,\omega_M^2),$ the absolute values of the terms on the right hand side of the above equation
are bounded by a sum of terms of the forms
$$ \left|\frac{1}{\lambda}\frac{\partial^3 f}{\partial\psi_r\partial\psi_s\partial\psi_t}\right|,\quad\left|\frac{1}{\lambda}\frac{\partial^2 f}{\partial\psi_r\partial\psi_s}\right|,\quad\left|\frac{1}{\lambda}\frac{\partial f}{\partial\psi_r}\right|\quad\mbox{and}\quad \frac{f}{\lambda}.$$

That expectation of the upper bound of each term is finite, is proved in Section \ref{sec:ub} of the supplement, 
which implies that (\ref{eq:h1}) and (\ref{eq:h2}) clearly hold.

Now, following \ctn{Maud12}, \ctn{Maitra14a}, \ctn{Maitra14b} we assume:
\begin{itemize}
\item[(H4)] The true value $\theta_0\in int\left(\Omega\right)$.
\end{itemize}
Now let us define $\mathcal I(\theta)$ to be the information matrix, whose $(r,s)$-th element is given by
$$\mathcal I_{rs}(\theta)=E_{\theta}\left[\frac{\partial \log \lambda}{\partial \theta_r}\frac{\partial \log \lambda}{\partial \theta_s}\right]=-E_{\theta}\left[\frac{\partial^2 \log \lambda}{\partial \theta_r\partial\theta_s}\right]$$
That $\mathcal I(\theta)$ is well-defined is clear due to existence of derivatives up to the second order.
Now, let $(b_1,b_2,\ldots,b_M)$ be a real row vector where $\theta$ is $M$-dimensional. Then
$$E\left[b_1\frac{\partial\log\lambda}{\partial\theta_1}+b_2\frac{\partial\log\lambda}{\partial\theta_2}+\ldots+b_M\frac{\partial\log\lambda}{\partial\theta_M}\right]^2\geq 0$$
implying $$\sum_{r=1}^M\sum_{s=1}^M b_rb_s\mathcal I_{rs}(\theta)\geq 0.$$
This shows that $\mathcal I(\theta)$ is positive semi-definite. We assume
 \begin{itemize}
\item[(H5)] $\mathcal I(\theta_0)$ is positive definite.
\end{itemize}

Hence, asymptotic normality of the $MLE$, of the form (\ref{eq:MLE_normality_iid}), holds in our case.
Formally,
\begin{theorem}
\label{theorem:asymp_normal_iid}
Assume the $iid$ setup and conditions (H1) -- (H5). 
Then the $MLE$ is asymptotically normally distributed as
(\ref{eq:MLE_normality_iid}).
\end{theorem}

\section{Consistency and asymptotic normality of $MLE$ in the non-$iid$ setup}
\label{sec:consistency_non_iid}

In this section, as in \ctn{Maitra14a} and \ctn{Maitra14b} we allow $T_i \neq T$
and $x^i \neq x$ for each $1\leq i\leq n$. Consequently, here we deal with the setup
where the processes $X_i(\cdot);~i=1,\ldots,n$, are independently,
but not identically distributed.
Following \ctn{Maitra14a} and \ctn{Maitra14b} we assume the following:
\begin{itemize}
\item[(H6)] The sequences $\{T_1,T_2,\ldots\}$ and 
$\{x^1,x^2,\ldots,\}$ are sequences in compact sets $\mathfrak T$ and $\mathfrak X$, respectively, so that
there exist convergent subsequences with limits in $\mathfrak T$ and $\mathfrak X$.
For notational convenience, we continue to denote the convergent subsequences as  $\{T_1,T_2,\ldots\}$
and $\{x^1,x^2,\ldots\}$. Let us denote the limits by $T^{\infty}$ and $x^{\infty}$, where $T^{\infty}\in\mathfrak T$
and $x^{\infty}\in\mathfrak X$. 
\end{itemize}

Following \ctn{Maitra14a} and \ctn{Maitra14b} we denote the process associated with the initial value $x$ and time point $t$ as $X(t,x)$, so that
$X(t,x^i)=X_i(t)$, and $X_i=\left\{X_i(t);~t\in[0,T_i]\right\}$. We also denote by
$\phi(x)$ the random effect
parameter associated with the initial value $x$ such that $\phi(x^i)=\phi_i$. 
We assume 
\begin{itemize}
\item[(H7)] $\phi(x)$ is a real-valued, continuous function of $x$, and that 
for $k\geq 1$, 
$\underset{x\in \mathfrak X}{\sup}~E\left[\phi(x)\right]^{2k}<\infty$. 
\end{itemize}

As in Proposition 1 of \ctn{Maud12}, assumption (H7) implies that for any $T>0$,
\begin{equation}
\underset{t\in [0,T],x\in \mathfrak X}{\sup}~E\left[X(t,x)\right]^{2k}<\infty.
\label{eq:sup_X}
\end{equation}

For $x\in \mathfrak X$ and $T\in \mathfrak T$, let
\begin{align}
U(x,T)&=\int_0^T\frac{b(X(s,x))}{\sigma^2(X(s,x))}d X(s,x)\label{eq:u_x_T};\\
V(x,T)&=\int_0^T\frac{b^2(X(s,x))}{\sigma^2(X(s,x))}ds.\label{eq:v_x_T}
\end{align}
Clearly, $U(x^i,T_i)=U_i$ and $V(x^i,T_i)=V_i$, where $U_i$ and $V_i$ are given by
(\ref{eq:sufficient}).

Even in this non-$iid$ case (H3) ensures that moments of all orders of $V(x,T)$ are finite.
Then, by Theorem 5 of \ctn{Maitra14a}, the moments of uniformly integrable continuous 
functions of $U(x,T)$, $V(x,T)$ and $\theta$ are continuous in 
$x$, $T$ and $\theta$.
In particular, 
the Kullback-Leibler distance and the information matrix, which we denote by $\mathcal K_{x,T}(\theta_0,\theta)$
(or, $\mathcal K_{x,T}(\theta,\theta_0)$)
and $\mathcal I_{x,T}(\theta)$ respectively to emphasize
dependence on the initial values $x$ and $T$, are continuous in $x$, $T$ and $\theta$.
For $x=x^k$ and $T=T_k$, if we denote the Kullback-Leibler 
distance and the Fisher's information
as $\mathcal K_k(\theta_0,\theta)$ ($\mathcal K_k(\theta,\theta_0)$) and $\mathcal I_k(\theta)$, respectively,
then continuity of $\mathcal K_{x,T}(\theta_0,\theta)$ (or $\mathcal K_{x,T}(\theta,\theta_0)$) and $\mathcal I_{x,T}(\theta_0)$
with respect to $x$ and $T$ ensures that as $x^k\rightarrow x^{\infty}$ and $T_k\rightarrow T^{\infty}$,
$\mathcal K_{x^k,T_k}(\theta_0,\theta)\rightarrow \mathcal K_{x^{\infty},T^{\infty}}(\theta_0,\theta)=\mathcal K(\theta_0,\theta)$, say.
Similarly, $\mathcal K_{x^k,T_k}(\theta,\theta_0)\rightarrow \mathcal K(\theta,\theta_0)$ and 
$\mathcal I_{x^k,T_k}(\theta)\rightarrow \mathcal I_{x^{\infty},T^{\infty}}(\theta)=
\mathcal I(\theta)$, say. 
Thanks to compactness, the limits $\mathcal K(\theta_0,\theta)$, $\mathcal K(\theta,\theta_0)$ and $\mathcal I(\theta)$ are well-defined
Kullback-Leibler divergences and Fisher's information, respectively.
Consequently (see \ctn{Maitra14a}, \ctn{Maitra14b}), the following hold for any $\theta\in\Omega$,
\begin{align}
\underset{n\rightarrow\infty}{\lim}~\frac{\sum_{k=1}^n\mathcal K_k(\theta_0,\theta)}{n}&=\mathcal K(\theta_0,\theta);
\label{eq:kl_limit_1}\\
\underset{n\rightarrow\infty}{\lim}~\frac{\sum_{k=1}^n\mathcal K_k(\theta,\theta_0)}{n}&=\mathcal K(\theta,\theta_0);
\label{eq:kl_limit_2}\\
\underset{n\rightarrow\infty}{\lim}~\frac{\sum_{k=1}^n\mathcal I_k(\theta)}{n}&=\mathcal I(\theta).
\label{eq:fisher_limit_1}
\end{align}
We assume that
\begin{itemize}
\item[(H8)] For any $\theta\in\Omega$, $\mathcal I(\theta)$ is positive definite.
\end{itemize}

\subsection{Consistency of $MLE$ in the non-$iid$ setup}
\label{subsec:consistency_non_iid}

Following \ctn{Hoadley71} we define the following:
\begin{align}
R_i(\theta)&=\log\frac{\lambda_i(X_i|\theta)}{\lambda_i(X_i|\theta_0)}\quad\mbox{if}\ \ \lambda_i(X_i|\theta_0)>0\notag\\
&=0 \quad\quad\quad\quad\quad\quad\mbox{otherwise}.
\label{eq:R1}
\end{align}

\begin{align}
R_i(\theta,\rho)&=\sup\left\{R_i(\xi):\|\xi-\theta\|\leq\rho\right\}\label{eq:R2}\\
{\mathcal V}_i(r)&=\sup\left\{R_i(\theta):\|\theta\|>r\right\}.\label{eq:V}
\end{align}
Following \ctn{Hoadley71} we denote by $r_i(\theta)$, $r_i(\theta,\rho)$ and $v_i(r)$
to be expectations of $R_i(\theta)$, $R_i(\theta,\rho)$ and ${\mathcal V}_i(r)$ under $\theta_0$; 
for any sequence $\{A_i;i=1,2,\ldots\}$ we denote
$\sum_{i=1}^nA_i/n$ by $\bar A_n$.

\ctn{Hoadley71} proved that if the following regularity conditions are satisfied, then 
the $MLE$ $\hat\theta_n\stackrel{P}{\rightarrow}\theta_0$:
\begin{itemize}
\item[(1)] $\Omega$ is a closed subset of $\mathbb R^{3M}$.
\item[(2)] $\lambda_i(X_i|\theta)$ is an upper semicontinuous 
function of $\theta$, uniformly in $i$, 
a.s. $[P_{\theta_0}]$.
\item[(3)] There exist $\rho^*=\rho^*(\theta)>0$, $r>0$ and $0<K^*<\infty$ for which 
\begin{enumerate}
\item[(i)] $E_{\theta_0}\left[R_i(\theta,\rho)\right]^2\leq K^*,\quad 0\leq\rho\leq\rho^*$;
\item[(ii)] $E_{\theta_0}\left[{\mathcal V}_i(r)\right]^2\leq K^*$.
\end{enumerate}
\item[(4)]
\begin{enumerate}
\item[(i)]$\underset{n\rightarrow\infty}{\lim}~\bar r_n(\theta)<0,\quad\theta\neq\theta_0$;
\item[(ii)]$\underset{n\rightarrow\infty}{\lim}~\bar v_n(r)<0$.
\end{enumerate}
\item[(5)] $R_i(\theta,\rho)$ and ${\mathcal V}_i(r)$ are measurable functions of $X_i$.
\end{itemize}
Actually, conditions (3) and (4) can be weakened but these are more easily applicable (see \ctn{Hoadley71} for details).

\subsubsection{Verification of the regularity conditions}
\label{subsubsec:consistency_non_iid}

Since $\Omega$ is compact in our case, the first regularity condition
clearly holds. 

For the second regularity condition, note that given $X_i$, 
$\lambda_i(X_i|\theta)$ is continuous (as $\lambda_i(X_i|\theta)=\sum_{k=1}^Ma_kf_k(X_i|\beta_k)$ where each $f_k(X_i|\beta_k)$ is continuous), in fact, uniformly continuous 
  in $\theta$ in our case, 
since $\Omega$ is compact. Hence, for any given $\epsilon>0$, there exists $\delta_i(\epsilon)>0$, independent
of $\theta$,
such that $\|\theta_1-\theta_2\|<\delta_i(\epsilon)$ implies $\left|\lambda(X_i|\theta_1)-\lambda(X_i|\theta_2)\right|<\epsilon$.
Now consider a strictly positive function $\delta_{x,T}(\epsilon)$, continuous in $x\in\mathfrak X$ and $T\in\mathfrak T$,
such that $\delta_{x^i,T_i}(\epsilon)=\delta_i(\epsilon)$. Let 
$\delta(\epsilon)=\underset{x\in\mathfrak X,T\in\mathfrak T}{\inf}\delta_{x,T}(\epsilon)$. Since
$\mathfrak X$ and $\mathfrak T$ are compact, it follows that $\delta(\epsilon)>0$. Now it holds that
$\|\theta_1-\theta_2\|<\delta(\epsilon)$ implies $\left|\lambda(X_i|\theta_1)-\lambda(X_i|\theta_2)\right|<\epsilon$,
for all $i$. Hence, the second regularity condition is satisfied.

Let us now focus attention on condition (3)(i).
It follows from (\ref{eq:upper_bound2}) that
\begin{equation}
	R_i(\theta)\leq \sum_{k=1}^MC_1(U_i,V_i,\beta_{0,1},\beta_{k})+|\log a_{0,1}|,
\label{eq:upper_r}
\end{equation}
 where $C_1(U_i,V_i,\beta_{0,1},\beta_{k})$ is given by (\ref{eq:upper_bound1}).
Let us denote $\left\{\xi\in\mathbb R\times\mathbb R^+:\|\xi-\beta_1\|\leq\rho\right\}$ by 
$B(\rho,\beta_1)$. Here $0<\rho<\rho^*(\beta_1)$, and $\rho^*(\beta_1)$ 
is so small that
$B(\rho,\beta_1)\subset\Omega_{\beta}$ for all $\rho\in (0,\rho^*(\beta_1))$. It then follows from (\ref{eq:upper_bound1}) that

\begin{align}
	&\underset{\xi\in B(\rho,\beta_1)}{\sup}~C_1(U_i,V_i,\beta_{0,1},\xi)\notag\\
	&\leq \underset{(\mu_1,\omega_1^2)\in B(\rho,\beta_1)}{\sup}~\frac{1}{2}\left\{\log\left(1+\frac{\omega_{0,1}^2}{\omega^2_{1}}\right)
	+\frac{|\omega_{0,1}^2-\omega^2_{1}|}{\omega_{0,1}^2}\right\}\notag\\
	&\quad+\underset{(\mu_1,\omega_1^2)\in B(\rho,\beta_1)}{\sup}~\left(\frac{U_i}{1+\omega^2_{1}V_i}\right)^2\times
	\underset{(\mu_1,\omega_1^2)\in B(\rho,\beta_1)}{\sup}~\left[\frac{1}{2}\left|\omega^2_{1}-\omega_{0,1}^2\right|
	\left(1+\frac{\omega^2_{1}}{\omega_{0,1}^2}\right)\right]\notag\\
	&\quad+\underset{(\mu_1,\omega_1^2)\in B(\rho,\beta_1)}{\sup}~\left\vert\frac{U_i}{1+\omega^2_{1}V_i}\right\vert
	\times\underset{(\mu_1,\omega_1^2)\in B(\rho,\beta_1)}{\sup}~
	\left[|\mu_{0,1}|\left(1+\frac{|\omega^2_{1}-\omega_{0,1}^2|}{\omega_{0,1}^2}\right)\right]\notag\\
	&\quad +\underset{(\mu_1,\omega_1^2)\in B(\rho,\beta_1)}{\sup}~\left|\frac{\mu^2_{1}V_i}{2(1+\omega^2_{1}V_i)}\right|
	+\underset{(\mu_1,\omega_1^2)\in B(\rho,\beta_1)}{\sup}~\left|\frac{\mu_{1}U_i}{1+\omega^2_{1}V_i}\right|.
\label{eq:sup_R1}
\end{align}
The supremums in (\ref{eq:sup_R1}) are finite due to compactness of $B(\rho,\beta_1)$ and since $\left\vert\frac{U_i}{1+\omega^2_{1}V_i}\right\vert$
and its square are decreasing in $\omega^2_1$, given $U_i$ and $V_i$. 
Since under $P_{\theta_0}$, $|U_i|/(1+\omega^2_{1}V_i)$ admits moments of all orders and 
$\frac{V_i}{1+\omega^2_{1}V_i}<\frac{1}{\omega^2_{1}}$ almost surely (see Section \ref{subsubsec:MLE_consistency_iid}), 
it follows from (\ref{eq:sup_R1}) and (\ref{eq:upper_r}) that $\underset{\psi\in S(\rho,\theta)}{\sup}~R_i(\psi)$ 
is finite due to  compactness of $S(\rho,\theta)=\left\{\psi\in\Omega:\|\psi-\theta\|\leq\rho\right\}$. 
Then it follows that
\begin{equation}
E_{\theta_0}\left[R_i(\theta,\rho)\right]^2\leq K_i(\theta),
\label{eq:upper_bound3}
\end{equation}
where $K_i(\theta)=K(x^i,T_i,\theta)$, with $K(x,T,\theta)$ being a continuous function of 
$(x,T,\theta)$, continuity being a consequence
of Theorem 5 of \ctn{Maitra14a}. 
Since 
because of compactness of $\mathfrak X$, $\mathfrak T$ and $\Omega$,
$$K_i(\theta)\leq \underset{x\in\mathfrak X,T\in\mathfrak T,\theta\in\Omega}{\sup}~K(x,T,\theta)<\infty,$$
regularity condition (3)(i) follows.

To verify condition (3)(ii), first note that we can choose $r>0$ such that $\|\theta_0\|<r$ and
$\{\theta\in\Omega:\|\theta\|>r\}\neq\emptyset$. 
It then follows that 
$\underset{\left\{\theta\in\Omega:\|\theta\|>r\right\}}{\sup}~R_i(\theta)\leq \underset{\theta\in\Omega}{\sup}~R_i(\theta)$
for every $i\geq 1$. The right hand side is bounded by the finite sum of the same expression as the right hand side of (\ref{eq:sup_R1}) and $\left|\log a_{0,1}\right|$
with only $S(\rho,\theta)$ replaced with $\Omega$. 
The rest of the verification follows in the same way as verification of (3)(i). 

To verify condition (4)(i) note that by (\ref{eq:kl_limit_1}) 
\begin{equation}
\underset{n\rightarrow\infty}{\lim}~\bar r_n=-\underset{n\rightarrow\infty}{\lim}~\frac{\sum_{i=1}^n\mathcal K_i(\theta_0,\theta)}{n}
=-\mathcal K(\theta_0,\theta)<0\quad\mbox{for}~\theta\neq\theta_0.
\label{eq:lim1}
\end{equation}
In other words, (4)(i) is satisfied. 

The verification of (4)(ii) will be in a similar way as it is in \ctn{Maitra14a} except in each case $f_i(X_i|\theta_0)$ will be replaced by $\lambda_i(X_i|\theta_0)$.

Regularity condition (5) holds because for any $\theta\in\Omega$, $R_i(\theta)$ is an almost surely
continuous function of $X_i$ rendering it measurable for all
$\theta\in\Omega$, and due to the fact that supremums of measurable functions
are measurable.

In other words, in the non-$iid$ $SDE$ framework, the following theorem holds:
\begin{theorem}
\label{theorem:consistency_non_iid}
Assume the non-$iid$ $SDE$ setup and conditions (H1) -- (H7). 
Then it holds that $\hat\theta_n\stackrel{P}{\rightarrow}\theta_0$.
\end{theorem}

\subsection{Asymptotic normality of $MLE$ in the non-$iid$ setup}
\label{subsec:normality_non_iid}

Let $\zeta_i(x,\theta)=\log \lambda_i(x|\theta)$; also, let $\zeta'_i(x,\theta)$ be the $3M\times 1$ 
vector with $j$-th component $\zeta'_{i,j}(x,\theta)=\frac{\partial}{\partial\theta_j}\zeta_i(x,\theta)$, and
let $\zeta''_i(x,\theta)$ be the $3M\times 3M$ matrix with $(j,k)$-th element
$\zeta''_{i,jk}(x,\theta)=\frac{\partial^2}{\partial\theta_j\partial\theta_k}\zeta_i(x,\theta)$.

For proving asymptotic normality in the non-$iid$ framework, \ctn{Hoadley71} assumed the
following regularity conditions:
\begin{itemize}
\item[(1)] $\Omega$ is an open subset of $\mathbb R^{3M}$.
\item[(2)] $\hat\theta_n\stackrel{P}{\rightarrow}\theta_0$.
\item[(3)] $\zeta'_i(X_i,\theta)$ and $\zeta''_i(X_i,\theta)$ exist a.s. $[P_{\theta_0}]$.
\item[(4)] $\zeta''_i(X_i,\theta)$ is a continuous function of $\theta$, uniformly in $i$, a.s. $[P_{\theta_0}]$,
and is a measurable function of $X_i$.
\item[(5)] $E_{\theta}[\zeta'_i(X_i,\theta)]=0$ for $i=1,2,\ldots$.
\item[(6)] $\mathcal I_i(\theta)=E_{\theta}\left[\zeta'_i(X_i,\theta)\zeta'_i(X_i,\theta)^T\right]
=-E_{\theta}\left[\zeta''_i(X_i,\theta)\right]$, where for any vector $y$, $y^T$ denotes
the transpose of $y$.
\item[(7)] $\bar{\mathcal I}_n(\theta)\rightarrow\bar{\mathcal I}(\theta)$ as $n\rightarrow\infty$ and 
$\bar{\mathcal I}(\theta)$ is positive definite.
\item[(8)] $E_{\theta_0}\left\vert\zeta'_{i,j}(X_i,\theta_0)\right\vert^3\leq K_2$, for some $0<K_2<\infty$.
\item[(9)] There exist $\epsilon>0$ and random variables $B_{i,jk}(X_i)$ such that
\begin{enumerate}
\item[(i)] $\sup\left\{\left\vert\zeta''_{i,jk}(X_i,\xi)\right\vert:\|\xi-\theta_0\|\leq\epsilon\right\}
\leq B_{i,jk}(X_i)$.
\item[(ii)] $E_{\theta_0}\left\vert B_{i,jk}(X_i)\right\vert^{1+\delta}\leq K_2$, for some $\delta>0$.
\end{enumerate}
\end{itemize}
Condition (8) can be weakened but is relatively easy to handle.
Under the above regularity conditions, \ctn{Hoadley71} prove that 
\begin{equation}
\sqrt{n}\left(\hat\theta_n-\theta_0\right)\stackrel{\mathcal L}{\rightarrow}
N\left(\bzero,\bar{\mathcal I}^{-1}(\theta_0)\right).
\label{eq:MLE_normality_non_iid}
\end{equation}

\subsubsection{Validation of asymptotic normality of $MLE$ in the non-$iid$ $SDE$ setup}
\label{subsubsec:normality_non_iid}

Note that although condition (1) requires the parameter space $\Omega$ to be an open subset, 
the proof of asymptotic normality presented in \ctn{Hoadley71} continues to hold for compact $\Omega$;
see \ctn{Maitra14a}.

Conditions (2), (3), (5), (6) are clearly valid in our case.
Condition (4) can be verified in exactly the same way as condition (2) of Section \ref{subsec:consistency_non_iid}; 
measurability of $\zeta''_i(X_i,\theta)$ follows due to its continuity with respect to $X_i$.
Condition (7) simply follows from (\ref{eq:fisher_limit_1}). 

For conditions (8), (9)(i) and (9)(ii) note that, by the same arguments as in Section \ref{subsubsec:MLE_normality_iid},
finiteness of moments of all orders of the derivatives are seen to hold for every $x\in\mathfrak X$, $T\in\mathfrak T$.
Then compactness of $\mathfrak X$, $\mathfrak T$ and $\Omega$,
ensures that the conditions (8), (9)(i) and (9)(ii) hold.

In other words, in our non-$iid$ $SDE$ case we have the following theorem on asymptotic normality.
\begin{theorem}
\label{theorem:asymp_normal_non_iid}
Assume the non-$iid$ $SDE$ setup and conditions (H1) -- (H8). 
Then (\ref{eq:MLE_normality_non_iid}) holds. 
\end{theorem}

\section{Consistency and asymptotic normality of the Bayesian posterior in the $iid$ setup}
\label{sec:b_consistency_iid}

\subsection{Consistency of the Bayesian posterior distribution}
\label{subsec:Bayesian_consistency_iid}

To verify posterior consistency we make use of Theorem 7.80 presented in \ctn{Schervish95}; below
we state the general form of the theorem.
\begin{theorem}[\ctn{Schervish95}]
\label{theorem:theorem3}
Let $\{X_n\}_{n=1}^{\infty}$ be conditionally $iid$ given $\theta$ with density $\lambda_1(x|\theta)$
with respect to a measure $\nu$ on a space $\left(\mathcal X^1,\mathcal B^1\right)$. Fix $\theta_0\in\Omega$, and define,
for each $S\subseteq\Omega$ and $x\in\mathcal X^1$,
\[
Z(S,x)=\inf_{\psi\in S}\log\frac{\lambda_1(x|\theta_0)}{\lambda_1(x|\psi)}.
\]
Assume that for each $\theta\neq\theta_0$, there is an open set $N_{\theta}$ such that $\theta\in N_{\theta}$ and
that $E_{\theta_0}Z(N_{\theta},X_i)> -\infty$. 
Also assume that $\lambda_1(x|\cdot)$ is continuous at $\theta$ 
for every $\theta$, a.s. $[P_{\theta_0}]$.
For $\epsilon>0$, define 
$C_{\epsilon}=\{\theta:\mathcal K_1(\theta_0,\theta)<\epsilon\}$, where 
\begin{equation}
\mathcal K_1(\theta_0,\theta)=E_{\theta_0}\left(\log\frac{\lambda_1(X_1|\theta_0)}{\lambda_1(X_1|\theta)}\right)
\label{eq:kl1}
\end{equation}
is the Kullback-Leibler divergence measure associated with observation $X_1$.
Let $\pi$ be a prior distribution such that $\pi(C_{\epsilon})>0$, for every $\epsilon>0$.
Then, for every $\epsilon>0$ and open set $\mathcal N_0$ containing $C_{\epsilon}$, the posterior satisfies
\begin{equation}
\lim_{n\rightarrow\infty}\pi_n\left(\mathcal N_0|X_1,\ldots,X_n\right)=1,\quad a.s.\quad [P_{\theta_0}].
\label{eq:posterior_consistency_iid}
\end{equation}
\end{theorem}

\subsubsection{Verification of posterior consistency}
\label{subsubsec:Bayesian_consistency_iid}
The condition $E_{\theta_0}Z(N_{\theta},X_i)> -\infty$ of the above theorem is 
verified in the context of Theorem \ref{theorem:theorem1} in Section \ref{subsubsec:MLE_consistency_iid}.

Now, all we need to ensure is that there exists a prior $\pi$ which gives positive probability to $C_{\epsilon}$
for every $\epsilon>0$. From the identifiability result given by Proposition 7 (i) of \ctn{Maud12}
it follows that $\mathcal K_1(\theta_0,\theta)=0$ if and only if $\theta=\theta_0$ (up to a label switching). Hence, for any $\epsilon>0$,
the set $C_{\epsilon}$ is non-empty, since it contains at least $\theta_0$.
In fact, continuity of
$\mathcal K_1(\theta_0,\theta)$ in $\theta$ follows from the fact that upper bound of 
$\left|\log\frac{\lambda_1(x|\theta_0)}{\lambda_1(x|\theta)}\right|$ has finite $E_{\theta_0}$-expectation as shown in Lemma \ref{lemma:upper_bound},
and since the parameter space $\Omega$ is compact, it follows that 
$\mathcal K_1(\theta_0,\theta)$ is uniformly continuous on $\Omega$. 
The rest of the verification remains the same as Section 2.1.1 of \ctn{Maitra14b}.

In other words, the following result on posterior consistency holds.

\begin{theorem}
\label{theorem:new_theorem3}
Assume the $iid$ setup and conditions (H1), (H2) and (H4). For $\epsilon>0$, define $C_{\epsilon}=\{\theta:\mathcal K_1(\theta_0,\theta)<\epsilon\}$, where $\mathcal K_1(\theta_0,\theta)$ is the Kullback-Leibler divergence measure associated with observation $X_1$.
Let the prior distribution $\pi$ of the parameter $\theta$ satisfy $\frac{d\pi}{d\nu}=h$ 
almost everywhere on $\Omega$, where $h(\theta)$ is any 
positive, continuous density on $\Omega$ with respect to the Lebesgue
measure $\nu$.  
Then the posterior (\ref{eq:posterior1}) is consistent
in the sense that for every $\epsilon>0$ and open set $\mathcal N_0$ containing $C_{\epsilon}$, the posterior satisfies
\begin{equation}
\lim_{n\rightarrow\infty}\pi_n\left(\mathcal N_0|X_1,\ldots,X_n\right)=1,\quad a.s.\quad [P_{\theta_0}].
\label{eq:posterior_consistency_iid2}
\end{equation}
\end{theorem}

\subsection{Asymptotic normality of the Bayesian posterior distribution}
\label{subsec:Bayesian_normality_iid}

To investigate asymptotic normality of our $SDE$-based posterior distributions we exploit 
Theorem 7.102 in conjunction with 
Theorem 7.89 provided in \ctn{Schervish95}. 
Below we state the four requisite conditions for the $iid$ setup.

\subsubsection{Regularity conditions -- $iid$ case}
\label{subsubsec:regularity_iid}
\begin{itemize}
\item[(1)] The parameter space is $\Omega\subseteq\mathbb R^{3M}$ for some finite $M$.
\item[(2)] $\theta_0$ is a point interior to $\Omega$.
\item[(3)] The prior distribution of $\theta$ has a density with respect to Lebesgue measure
that is positive and continuous at $\theta_0$.
\item[(4)] There exists a neighborhood $\mathcal N_0\subseteq\Omega$ of $\theta_0$ on which
$\ell_n(\theta)= \log \lambda(X_1,\ldots,X_n|\theta)$ is twice continuously differentiable with 
respect to all co-ordinates of $\theta$, 
$a.s.$ $[P_{\theta_0}]$.
\end{itemize}

With the above conditions, the relevant theorem (Theorem 7.102 of \ctn{Schervish95}) is as follows:
\begin{theorem}[\ctn{Schervish95}]
\label{theorem:theorem4}
Let $\{X_n\}_{n=1}^{\infty}$ be conditionally $iid$ given $\theta$.
Assume the above four regularity conditions; 
also assume that there exists $H_r(x,\theta)$ such that, for each $\theta_0\in int(\Omega)$ and each
$k,j$,
\begin{align}
\sup_{\|\theta-\theta_0\|\leq r}\left\vert\frac{\partial^2}{\partial\theta_k\partial\theta_j}\log \lambda_{1}(x|\theta_0)
-\frac{\partial^2}{\partial\theta_k\partial\theta_j}\log \lambda_{1}(x|\theta)\right\vert\leq H_r(x,\theta_0),
\label{eq:H1}
\end{align}
with
\begin{equation}
\lim_{r\rightarrow 0}E_{\theta_0}H_r\left(X,\theta_0\right)=0.
\label{eq:H2}
\end{equation}
Further suppose that
the conditions of Theorem \ref{theorem:theorem3} hold, and that the Fisher's information matrix
$\mathcal I(\theta_0)$ is positive definite. 
Now denoting 
by $\hat\theta_n$ the $MLE$ associated with
$n$ observations, 
let
\begin{equation}
\Sigma^{-1}_n=\left\{\begin{array}{cc}
-\ell''_n(\hat\theta_n) & \mbox{if the inverse and}\ \ \hat\theta_n\ \ \mbox{exist}\\
\mathbb I_{3M} & \mbox{if not},
\end{array}\right.
\label{eq:information1}
\end{equation}
where for any $t$,
\begin{equation}
\ell''_n(t)=\left(\left(\frac{\partial^2}{\partial\theta_i\partial\theta_j}\ell_n(\theta)\bigg\vert_{\theta=t}\right)\right),
\label{eq:information2}
\end{equation}
and $\mathbb I_{3M}$ is the identity matrix of order $3M$.
Thus, $\Sigma^{-1}_n$ is the observed Fisher's information matrix.

Letting $\Psi_n=\Sigma^{-1/2}_n\left(\theta-\hat\theta_n\right)$, for each compact subset 
$B$ of $\mathbb R^{3M}$ and each $\epsilon>0$, the following holds:
\begin{equation}
\lim_{n\rightarrow\infty}P_{\theta_0}
\left(\sup_{\Psi_n\in B}\left\vert\pi_n(\Psi_n\vert X_1,\ldots,X_n)-\xi(\Psi_n)\right\vert>\epsilon\right)=0,
\label{eq:Bayesian_normality_iid}
\end{equation}
where $\xi(\cdot)$ denotes the density of the standard normal distribution.
\end{theorem}

\subsubsection{Verification of posterior normality}
\label{subsubsec:Bayesian_normality_iid}
Firstly, note that (H5) ensures positive definiteness of $\mathcal I(\theta_0)$.
Now observe that the 
four regularity conditions in Section \ref{subsubsec:regularity_iid} trivially hold. The remaining conditions of Theorem \ref{theorem:theorem4} are verified in the context of 
Theorem \ref{theorem:consistency_iid} in Section \ref{subsec:MLE_consistency_iid}. Briefly, $\frac{\partial^2}{\partial\theta_k\partial\theta_j}\log \lambda_{1}(x|\theta)$
is differentiable in $\theta$ and the derivative has finite expectation, which ensure (\ref{eq:H1})
and (\ref{eq:H2}).
Hence, (\ref{eq:Bayesian_normality_iid}) holds in our $SDE$ setup.
Thus, the following theorem holds:
\begin{theorem}
\label{theorem:new_theorem4}
Assume the $iid$ setup and conditions (H1) --  (H5). 
Let the prior distribution $\pi$ of the parameter $\theta$ satisfy $\frac{d\pi}{d\nu}=h$ 
almost everywhere on $\Omega$, where $h(\theta)$ is any 
density with respect to the Lebesgue measure $\nu$ which is positive and continuous at $\theta_0$. 
Then, letting $\Psi_n=\Sigma^{-1/2}_n\left(\theta-\hat\theta_n\right)$, for each compact subset 
$B$ of $\mathbb R^{3M}$ and each $\epsilon>0$, the following holds:
\begin{equation}
\lim_{n\rightarrow\infty}P_{\theta_0}
\left(\sup_{\Psi_n\in B}\left\vert\pi_n(\Psi_n\vert X_1,\ldots,X_n)-\xi(\Psi_n)\right\vert>\epsilon\right)=0.
\label{eq:Bayesian_normality_iid2}
\end{equation}
\end{theorem}

\section{Consistency and asymptotic normality of the Bayesian posterior in the non-$iid$ setup}
\label{sec:b_consistency_non_iid}

In this section, as in Section \ref{sec:consistency_non_iid} we assume (H7) -- (H9). The equations (\ref{eq:kl_limit_1}), (\ref{eq:kl_limit_2}) and (\ref{eq:fisher_limit_1}) as described in Section \ref{sec:consistency_non_iid} 
will also have important roles in our proceedings. For consistency in the Bayesian framework we utilize the theorem 
of \ctn{Choi04}, and for asymptotic normality of the posterior we make use of Theorem 7.89
of \ctn{Schervish95}.

\subsection{Posterior consistency in the non-$iid$ setup}
\label{subsec:posterior_consistency_non_iid}

Analogous to Section 3.1 of \ctn{Maitra14b}, here we need to ensure existence of moments
of the form
$$\underset{x\in\mathfrak X,T\in\mathfrak T}{\sup}~E_{\theta}\left[\exp\left\{\alpha \left|\omega^2_{k}-\omega_{0,1}^2\right| 
\left(\frac{U(x,T)}{1+\omega^2_{k}V(x,T)}\right)^2\left(1+\frac{\omega^2_{k}}{\omega_{0,1}^2}\right)\right\}\right],$$ for some
$0<\alpha<\infty$. Hence, we assume the following assumption analogous to assumption (H10$^\prime$) of \ctn{Maitra14b}.
\begin{itemize}
\item[(H9)] For $k=1,\ldots,M$, there exists a strictly positive function $\alpha^*(x,T,\beta_1)$, 
continuous in $(x,T,\beta_1)$, such that for any $(x,T,\beta_1)$,
\begin{equation*}
E_{\theta}\left[\exp\left\{\alpha^*(x,T,\beta_1)K_1 U^2(x,T)\right\}\right]<\infty, 
\label{eq:moment_alpha1}
\end{equation*}
where
$K_1=\underset{\omega_1:~\beta_1\in\Omega_{\beta},~1\leq k\leq M}{\sup}~
\left|\omega^2_{k}-\omega_{0,1}^2\right|\left(1+\frac{\omega^2_{k}}{\omega_{0,1}^2}\right)$.
\end{itemize}

Now, let 
\begin{equation}
\alpha^*_{\min}=\underset{x\in\mathfrak X,T\in\mathfrak T,\beta_1\in\Omega_{\beta}}{\inf}\alpha^*(x,T,\beta_1),
\label{eq:alpha_star}
\end{equation}
and 
\begin{equation}
\alpha=\min\left\{\alpha^*_{\min},c^*\right\},
\label{eq:alpha2}
\end{equation}
where $0<c^*<1/16$. 

Compactness ensures that $\alpha^*_{\min}>0$, so that $0<\alpha<1/16$.
It also holds due to compactness that for $\beta_1\in\Omega_{\beta}$,
\begin{equation}
\underset{x\in\mathfrak X,T\in\mathfrak T}{\sup}~E_{\theta}\left[\exp\left\{\alpha K_1U^2(x,T)\right\}\right]<\infty. 
\label{eq:moment_alpha2}
\end{equation}
This ensures that 
\begin{align}
	&\underset{x\in\mathfrak X,T\in\mathfrak T}{\sup}~E_{\theta}\left[\exp\left\{\alpha \left|\omega^2_{k}-\omega_{0,1}^2\right| 
	\left(\frac{U(x,T)}{1+\omega^2_{k}V(x,T)}\right)^2\left(1+\frac{\omega^2_{k}}{\omega_{0,1}^2}\right)\right\}\right]\notag\\
&\leq \underset{x\in\mathfrak X,T\in\mathfrak T}{\sup}~E_{\theta}\left[\exp\left\{\alpha K_1U^2(x,T)\right\}\right]\notag\\
&<\infty.
\label{eq:u_square_moment}
\end{align}

This choice of $\alpha$ ensuring (\ref{eq:moment_alpha2}) will be useful in verification
of the conditions of Theorem \ref{theorem:theorem5}, which we next state.
\begin{theorem}[\ctn{Choi04}]
\label{theorem:theorem5}
Let $\{X_i\}_{i=1}^{\infty}$ be independently distributed with densities $\{\lambda_i(\cdot|\theta)\}_{i=1}^{\infty}$,
with respect to a common $\sigma$-finite measure, where $\theta\in\Omega$, a measurable space. The densities
$\lambda_i(\cdot|\theta)$ are assumed to be jointly measurable. Let $\theta_0\in\Omega$ and let $P_{\theta_0}$
be the joint distribution of $\{X_i\}_{i=1}^{\infty}$ when $\theta_0$ is the true value of $\theta$.
Let $\{\Theta_n\}_{n=1}^{\infty}$ be a sequence of subsets of $\Omega$. Let $\theta$ have prior $\pi$ on $\Omega$.
Define the following:
\begin{align}
\Lambda_i(\theta_0,\theta) &=\log\frac{\lambda_i(X_i|\theta_0)}{\lambda_i(X_i|\theta)};\notag\\
\mathcal K_i(\theta_0,\theta) &= E_{\theta_0}\left(\Lambda_i(\theta_0,\theta)\right);\notag\\
\varrho_i(\theta_0,\theta) &= Var_{\theta_0}\left(\Lambda_i(\theta_0,\theta)\right).\notag
\end{align}
Make the following assumptions:
\begin{itemize}
\item[(1)] Suppose that there exists a set $B$ with $\pi(B)>0$ such that
\begin{enumerate}
\item[(i)] $\sum_{i=1}^{\infty}\frac{\varrho_i(\theta_0,\theta)}{i^2}<\infty,\quad\forall~\theta\in B$,
\item[(ii)] For all $\epsilon>0$, $\pi\left(B\cap\left\{\theta:\mathcal K_i(\theta_0,\theta)<\epsilon,~\forall~i\right\}\right)>0$.
\end{enumerate}
\item[(2)] Suppose that there exist test functions $\{\Phi_n\}_{n=1}^{\infty}$, sets $\{\Omega_n\}_{n=1}^{\infty}$
and constants $C_1,C_2,c_1,c_2>0$ such that
\begin{enumerate}
\item[(i)] $\sum_{n=1}^{\infty}E_{\theta_0}\Phi_n<\infty$,
\item[(ii)] $\underset{\theta\in \Theta^c_n\cap\Omega_n}{\sup}~E_{\theta}\left(1-\Phi_n\right)\leq C_1e^{-c_1n}$,
\item[(iii)] $\pi\left(\Omega^c_n\right)\leq C_2e^{-c_2n}$.
\end{enumerate}
\end{itemize}
Then,
\begin{equation}
\pi_n\left(\theta\in \Theta^c_n|X_1,\ldots,X_n\right)\rightarrow 0\quad a.s.~[P_{\theta_0}].
\label{eq:posterior_consistency_non_iid}
\end{equation}
\end{theorem}

\subsubsection{Validation of posterior consistency}
\label{subsubsec:posterior_consistency_non_iid}

From Lemma \ref{lemma:upper_bound} it follows that $\left|\log \frac{\lambda_i(X_i|\theta_0)}{\lambda_i(X_i|\theta)}\right|$ 
has an upper bound which has finite expectation and square of expectation under $\theta_0$, and is uniform for all 
$\theta\in B$, where $B$ is an appropriate compact subset of the relevant parameter space.
The rest of the verification of condition (1)(i) and the verification of (1)(ii) are similar to those of \ctn{Maitra14b}.

In verification of (2)(iii), we let $\Omega_n=\left(\Omega_{1n}\times\mathbb R^{3M-1}\right)$ 
(since our parameter set $\theta$ contains $3M$ parameters),
where $\Omega_{1n}=\left\{a_1:|a_1|<\tilde M_n\right\}$, where $\tilde M_n=O(e^n)$. 
Note that
\begin{equation}
\pi\left(\Omega^c_n\right)=\pi\left(\Omega^c_{1n}\right)=\pi(|a_1|>\tilde M_n)
<E_{\pi}\left(|a_1|\right)\tilde M^{-1}_n
\label{eq:sieve1}
\end{equation}
implies (2)(iii) holds, assuming that the prior $\pi$ is such that the expectation 
$E_{\pi}\left(|a_1|\right)$ is finite (which holds for proper priors on $a_1$).

The verification of (2)(i) will follow in the same way as the verification in \ctn{Maitra14b} except the 
corresponding changes. Hence we will only mention the changes at which the verifications differ. 
Firstly, in our setup $L_n(\theta)=\prod_{i=1}^n\lambda_i(X_i|\theta)$
and $\ell_n(\theta)=\sum_{i=1}^n\log \lambda_i(X_i|\theta)$, that is, 
$f_i$ is now replaced with $\lambda_i$. The existence of the 
third order derivative of $\log\lambda_i$ is already established in Section \ref{subsubsec:MLE_normality_iid}.

Here also the continuity of the moments of $V(x, T)$ and $\left|\frac{U(x,T)}{V(x,T)}\right|$ with respect to $x$ and
$T$ holds (which follows from Theorem 5 of \ctn{Maitra14a} where uniform integrability
is ensured by finiteness of the moments of the aforementioned functions for every $x, T$ belonging to
compact sets $\mathfrak X$ and $\mathfrak T$). Moreover, Kolmogorov's strong law of large numbers for 
non-$iid$ cases holds due to finiteness of the moments of $V$ and $\left|\frac{U}{V}\right|$ 
for every $x$ and $T$ belonging to the compact spaces $\mathfrak X$ and $\mathfrak T$.

Now, assuming $\hat\theta_n=\zeta=(\gamma,\beta)$ we obtain (\ref{eq:upper_bound2}) from Section \ref{subsubsec:MLE_consistency_iid}, the following:
$$\log\lambda_i(x|\theta_0)-\log\lambda_i(x|\hat\theta_n)\geq -\sum_{k=1}^MC_1(U_i,V_i,\beta_{0,1},\beta_{k})-|\log a_{0,1}|,$$
where $C_1(U_i,V_i,\beta_{0,1},\beta_{k})$ is given by (\ref{eq:upper_bound1}). 
The rest of the verification of (2)(i) is the same as in \ctn{Maitra14b}.

For the verification of (2)(ii), we define 
$\Theta_n=\Theta_{\delta}=\left\{(\gamma,\beta):\mathcal K(\theta,\theta_0)<\delta\right\}$ 
where $\mathcal K(\theta,\theta_0)$, defined as in (\ref{eq:kl_limit_2}), is the proper Kullback-Leibler
divergence and the verification will be in a similar manner as in \ctn{Maitra14b}. Hence, 
posterior consistency (\ref{eq:posterior_consistency_non_iid}) holds in our non-$iid$
$SDE$ setup. The result can be summarized in the form of the following theorem.
\begin{theorem}
\label{theorem:new_theorem5}
Assume the non-$iid$ $SDE$ setup. Also assume conditions (H1) -- (H9).  
For any $\delta>0$, let $\Theta_{\delta}=\left\{(\gamma,\beta):\mathcal K(\theta,\theta_0)<\delta\right\}$,
where $\mathcal K(\theta,\theta_0)$, defined as in (\ref{eq:kl_limit_2}), is the proper Kullback-Leibler
divergence. 
Let the prior distribution $\pi$ of the parameter $\theta$ satisfy $\frac{d\pi}{d\nu}=h$ 
almost everywhere on $\Omega$, where $h(\theta)$ is any 
positive, continuous density on $\Omega$ with respect to the Lebesgue
measure $\nu$.  
Then,
\begin{equation}
\pi_n\left(\theta\in \Theta^c_{\delta}|X_1,\ldots,X_n\right)\rightarrow 0\quad a.s.~[P_{\theta_0}].
\label{eq:posterior_consistency_non_iid2}
\end{equation}
\end{theorem}

\subsection{Asymptotic normality of the posterior distribution in the non-$iid$ setup} 
\label{subsec:Bayesian_normality_non_iid}

For asymptotic normality of the posterior in the $iid$ situation, four regularity conditions,
stated in Section \ref{subsubsec:regularity_iid}, were necessary. In the non-$iid$ framework,
three more are necessary, in addition to the already presented four conditions. They are as follows
(see \ctn{Schervish95} for details).

\subsubsection{Extra regularity conditions in the non-$iid$ setup}
\label{subsubsec:regularity_non_iid}

\begin{itemize}
\item[(5)] The largest eigenvalue of $\Sigma_n$ goes to zero in probability.
\item[(6)] For $\delta>0$, define $\mathcal N_0(\delta)$ to be the open ball of radius $\delta$ around $\theta_0$.
Let $\rho_n$ be the smallest eigenvalue of $\Sigma_n$. If $\mathcal N_0(\delta)\subseteq\Omega$, there exists
$K(\delta)>0$ such that
\begin{equation}
\underset{n\rightarrow\infty}{\lim}~P_{\theta_0}\left(\underset{\theta\in\Omega\backslash\mathcal N_0(\delta)}{\sup}~
\rho_n\left[\ell_n(\theta)-\ell_n(\theta_0)\right]<-K(\delta)\right)=1.
\label{eq:extra1}
\end{equation}
\item[(7)] For each $\epsilon>0$, there exists $\delta(\epsilon)>0$ such that
\begin{equation}
\underset{n\rightarrow\infty}{\lim}~P_{\theta_0}\left(\underset{\theta\in\mathcal N_0(\delta(\epsilon)),\|\eta\|=1}{\sup}~
\left\vert 1+\eta^T\Sigma^{\frac{1}{2}}_n\ell''_n(\theta)\Sigma^{\frac{1}{2}}_n\eta\right\vert<\epsilon\right)=1.
\label{eq:extra2}
\end{equation}
\end{itemize}

In the non-$iid$ case, the four regularity conditions presented in Section \ref{subsubsec:regularity_iid}
and additional three provided above, are sufficient to guarantee (\ref{eq:Bayesian_normality_iid}).

\subsubsection{Verification of the regularity conditions}
\label{subsubsec:verify_last}

The verification follows from the verification presented in Section 3.2.2 of \ctn{Maitra14b}, 
as finiteness of the expectation of $\log\lambda$ up to the third order derivative is already 
justified in Section \ref{subsubsec:MLE_normality_iid}. So, we have our result in the form of the following theorem.
\begin{theorem}
\label{theorem:new_theorem6}
Assume the non-$iid$ setup and conditions (H1) -- (H8). 
Let the prior distribution $\pi$ of the parameter $\theta$ satisfy $\frac{d\pi}{d\nu}=h$ 
almost everywhere on $\Omega$, where $h(\theta)$ is any 
density with respect to the Lebesgue measure $\nu$ which is positive and continuous at $\theta_0$. 
Then, letting $\Psi_n=\Sigma^{-1/2}_n\left(\theta-\hat\theta_n\right)$, for each compact subset 
$B$ of $\mathbb R^{3M}$ and each $\epsilon>0$, the following holds:
\begin{equation}
\lim_{n\rightarrow\infty}P_{\theta_0}
\left(\sup_{\Psi_n\in B}\left\vert\pi_n(\Psi_n\vert X_1,\ldots,X_n)-\xi(\Psi_n)\right\vert>\epsilon\right)=0.
\label{eq:Bayesian_normality_non_iid2}
\end{equation}
\end{theorem}

\section{Bayesian asymptotics when the number of mixture components is unknown}
\label{sec:vardim}

In this section we allow the number of mixture components $M$ to be a random variable. 
Now $\theta_M=(\beta_M,\gamma_M)$, where $\beta_M$ and $\gamma_M$ are $M$-dimensional, with $M\geq 1$ being the number of positive components of $\gamma_M$. 
Suppose that $\theta_{0M_0}=(\beta_{0M_0},\gamma_{0M_0})$ is the true set of parameters, where 
$M_0\geq 1$ is the number of positive components of $\gamma_{M_0}$. 
Then 
the true likelihood is
\begin{equation}
\lambda_i(X_i,\theta_{0M})=\sum_{k=1}^{M_0} a_{0k} f(X_i|\beta_{0k}),
\label{eq:likelihood_var_true}
\end{equation}
and
\begin{equation}
\lambda_i(X_i,\theta_M)=\sum_{k=1}^M a_k f(X_i|\beta_k)
\label{eq:likelihood_var}
\end{equation}
is the modeled likelihood.

In order to prove posterior consistency in this setup, we shall apply the theorem from \ctn{Schervish95}. However, the theorem demands that the true and 
modeled likelihoods must have the same dimension, that is, the parameter set for both the models must have the same number of components.
To handle this variable dimensional situation, we consider the following idea.

Whenever $M\neq M_0$  we are left with two possibilities, $M<M_0$ or $M>M_0$. Let $M^\prime=\max\{M_0,M\}$. 
We can then rewrite both the true and the postulated likelihoods in terms of $M^\prime$. 
So, now with the general notation $\theta_{0M^\prime}=(\beta_{0M^\prime},\gamma_{0M^\prime})$ and $\theta_{M^\prime}=(\beta_{M^\prime},\gamma_{M^\prime})$, 
our comparable models are
$$\lambda_i(X_i,\theta_{0M^\prime})=\sum_{k=1}^{M^\prime} a_{0k} f(X_i|\beta_{0k})$$ and $$\lambda_i(X_i,\theta_{M^\prime})=\sum_{k=1}^{M^\prime} a_k f(X_i|\beta_k),$$
If $M<M_0$, then $a_{M+1}=a_{M+2}=\cdots=a_{M^\prime}=a_{M_0}=0$. 
When $M>M_0$, then $a_{0,M_0+1}=a_{0,M_0+2}=\cdots=a_{0M}=0$ and values of $\beta_{0,M_0+1},\beta_{0,M+2},\ldots,\beta_{0M}$ set arbitrarily, 
matching the number of components in both the models.

Hence, with $\theta_{0M^\prime}$ and $\theta_{M^\prime}$ having the same dimension, we can exploit our previous result established in
Section \ref{subsec:Bayesian_consistency_iid} to prove posterior consistency
even in this variable-dimensional setup. Theorem \ref{theorem:new_theorem7} formalizes our result in this regard.


\begin{theorem}
\label{theorem:new_theorem7}
Assume the $iid$ setup and conditions (H1), (H2) and (H4). For $\epsilon>0$, define $C_{\epsilon M}
=\{\theta_{M}:\mathcal K_{1}(\theta_{0M^{\prime}},\theta_{M^\prime})<\epsilon\}$, where 
$\mathcal K_1(\theta_{0M^\prime},\theta_{M^\prime})$ is the Kullback-Leibler divergence measure associated with observation $X_1$. 
Let $M$ have a prior with respect to the counting measure on the set $\mathcal S=\{1,2,\ldots,M_{max}\}$, where $M_0\leq M_{max}\leq\infty$.
For $M\in\mathcal S$, let the prior distribution $\pi$ of the parameter $\theta_{M}$ satisfy $\frac{d\pi}{d\nu}=h$ 
almost everywhere on $\Omega_{M}=\Omega_{\gamma_{M}}\times\Omega_{\beta_{M}}$, where $h(\theta_{M})$ is any positive, continuous density on 
$\Omega_{M}$ with respect to the Lebesgue
measure $\nu$. In addition we assume that {\it a priori}, all the components of $\gamma_M$ are bounded away from zero. 
Then 
for every $\epsilon>0$ and open set $\mathcal N_{0M}(\epsilon)$ containing $C_{\epsilon M}$, the posterior satisfies
\begin{equation}
\lim_{n\rightarrow\infty}\pi_n\left(\mathcal N_{0M}(\epsilon)|M,X_1,\ldots,X_n\right)=1~\mbox{a.s.}~[P_{\theta_{0M_0}}],~\mbox{if and only if}~M=M_0.
\label{eq:posterior_consistency_mixvar}
\end{equation}
\end{theorem}
\begin{proof}[Proof]
We write the posterior as follows:
\begin{align}
&\pi_n(\mathcal N_{0M}(\epsilon),M|X_1,X_2,\ldots,X_n)\notag\\
&\qquad=\pi_n(\mathcal N_{0M}(\epsilon)|M,X_1,X_2,\ldots,X_n)\pi_n(M|X_1,X_2,\ldots,X_n)\notag\\
&\qquad=\pi_n(\mathcal N_{0M^\prime}(\epsilon)|M^\prime,X_1,X_2,\ldots,X_n)\pi_n(M|X_1,X_2,\ldots,X_n).
\label{eq:posterior_var1}
\end{align}
Note that in the right hand side of (\ref{eq:posterior_var1}), $\pi_n(\mathcal N_{0M^\prime}|M^\prime,X_1,X_2,\ldots,X_n)$ is the posterior when 
$M^\prime$ is fixed. Hence, this problem is fixed-dimensional. 

In our case, observe that when $M<M_0$, then since $a_{M+1}=a_{M+2}=\ldots=a_{M_0}=0$ but $a_{01},a_{02},\ldots,a_{0,M_0}$ are all non-zero,
it follows that almost surely $K_1\left(\theta_{0M^\prime},\theta_{M^\prime}\right)>\varepsilon$, for some $\varepsilon>0$.
Hence, for $M<M_0$, if $0<\epsilon\leq\varepsilon$, 
$\underset{n\rightarrow\infty}{\lim}~\pi_n\left(\mathcal N_{0M^\prime}(\epsilon)|M^\prime,X_1,\ldots,X_n\right)=0$, almost surely.

If, on the other hand, $M>M_0$, then since $a_{0,M_0+1}=a_{0,M_0+2}=\cdots=a_{0M}=0$, but all the components of $\gamma_M$ are bounded
away from zero {\it a priori}, again there exists $\varepsilon>0$ such that $K_1\left(\theta_{0M^\prime},\theta_{M^\prime}\right)>\varepsilon$ almost surely. Hence
for $0<\epsilon\leq\varepsilon$, $\underset{n\rightarrow\infty}{\lim}~\pi_n\left(\mathcal N_{0M^\prime}(\epsilon)|M^\prime,X_1,\ldots,X_n\right)=0$, almost surely,
if $M>M_0$.

Only if $M=M_0$, it follows from the proof in Section \ref{subsec:Bayesian_consistency_iid}, that
$\underset{n\rightarrow\infty}{\lim}~\pi_n\left(\mathcal N_{0M^\prime}(\epsilon)|M^\prime,X_1,\ldots,X_n\right)=1$ holds almost surely, for any $\epsilon>0$.  


\end{proof}
\begin{remark}
It follows from Theorem \ref{theorem:new_theorem7} that asymptotic posterior normality holds in the same way as Theorem \ref{theorem:new_theorem4} when $M=M_0$.
Moreover, it is worth mentioning that both posterior consistency and asymptotic posterior normality in the non-$iid$ set up can be handled in the same way as in the $iid$ case.
\end{remark}

\section{A brief overview of the simulation studies}
\label{sec:simstudy}

For our simulation studies we consider the three $SDE$ mixture models and the five normal mixture models for the random effects
illustrated by \ctn{Maud16} for their simulation experiments. 

As the frequentist $MLE$ based approach using the $EM$ algorithm is already demonstrated by \ctn{Maud16}, here we focus on the Bayesian approach. Importantly, in all the cases we
assume that the number of mixture components $M$, is unknown, and consider a prior for $M$. 
Since $M$ is random in our approach, this makes the problem a variable-dimensional one, where the usual MCMC methods fail.
We resort to TTMCMC to obtain samples from the variable-dimensional posterior. 

Our variable-dimensional Bayesian approach is in sharp contrast with the approach of \ctn{Maud16}, where 
$M$ is either considered known or held fixed after selecting a value using $BIC$.  
We demonstrate that even in the cases where the number of mixture components are not well-separated, our Bayesian methodology puts up quite reasonable
performance, unlike the frequentist approach based on the $EM$ algorithm and $BIC$. 
The complete details are presented in Section \ref{sec:simstudy_details} of the supplement.

\section{Application to real stock market data}
\label{sec:realdata} 

With the simulation experiments we have demonstrated the usefulness of our Bayesian approach to $SDE$ mixtures. 
We now consider application of our $SDE$ system to a real, stock market data. The data, available at {\it www.nseindia.com}, consists of
$467$ observations from August $5$, 2013, to June $30$, 2015, for $15$ companies. For our purpose, we consider 
the ``close price" of each company as our data $X(t)$. We are interested in understanding if the companies 
comprise a single cluster with respect to the close price time series.

\subsection{Choice of the $SDE$ random effects model}
\label{subsec:sde_realdata}
To select an appropriate model for the data, we first model the company-wise data sets by the available standard financial $SDE$ models. 
These models are made available and amenable to easy implementation in the ``fitsde" package of $R$.  
We then obtain obtain the best model among such models by $BIC$. Here the minimum value of $BIC$ turned out to correspond to the $CKLS$ 
(\ctn{Chan92}) model. 
Denoting the data for the $i$-th company by $X_i(t)$, the $CKLS$ model is described by
\begin{equation}
dX_i(t)=(\theta_{1i}+\theta_{2i}X(t))dt+\theta_{3i}X_i(t)^{\theta_{4i}}dW_i(t).
\label{eq:sde_realdata}
\end{equation}
We fix the values of $\theta_{3i}$ and $\theta_{4i}$ as estimated by the ``fitsde" function. 

In this application, we consider the bivariate random effects $\bphi_i=(\theta_{1i},\theta_{2i})$; $i=1,\ldots,15$, assuming that for $i=1,\ldots,15$,
\begin{equation}
\bphi_i\stackrel{iid}{\sim}\sum_{k=1}^Ma_kN_2\left(\bmu_k,\bOmega_k\right),
\label{eq:bivariate_phi}
\end{equation}
where $N_2\left(\bmu_k,\bOmega_k\right)$ stands for bivariate normal with mean $\bmu_k$ and covariance matrix $\bOmega_k$, and given $M\geq 1$, for $k=1,\ldots,M$,
$a_k\geq 0$ and $\sum_{k=1}^Ma_k=1$. 
This setup thus is of the form multidimensional linear random effects $SDE$ of the form (\ref{eq:likelihood_mult}) discussed in Section \ref{sec:multidim}.

\subsection{Prior structure}
\label{subsec:prior_realdata}
As in the simulation studies, we consider a uniform prior for $M$ on $\{1,\ldots,30\}$. As in \ctn{Das17}, we assume the following prior structure for $(\bmu_k,\bOmega_k)$:
\begin{align}
[\bmu_k|\bOmega_k]&\sim N_2\left(\bzero,\bOmega_k\right);\notag\\
[\bOmega_k]&\sim W^{-1}\left(3,\tilde\bOmega\right).\notag
\end{align}
In the above, $W^{-1}\left(3,\tilde\bOmega\right)$ denotes the inverse-Wishart distribution with $3$ degrees of freedom and positive definite scale matrix $\tilde\bOmega$.
In our application, we choose $\tilde\bOmega=10\times\mathbb I_2$. Note that such a prior structure is a simple generalization of the one-dimensional
mixtures considered in the simulation studies. Indeed, here we also use the same prior for $(a_1,\ldots,a_M)$ as in the simulation studies.
We also experimented with other choices of the priors as detailed in \ctn{Das17}, but the results remained almost identical, suggesting considerable robustness
of our Bayesian inference with respect to choice of the priors.

\subsection{Implementation}
\label{subsec:implementation_realdata}
Although implementation of TTMCMC in this case is more involved than in the one-dimensional mixture setups, the issue of multivariate mixture implementation
was treated in details in \ctn{Das17}, who illustrate implementation up to $20$-dimensional mixtures. We follow the general strategy used in their paper
for our 2-dimensional setup. As in the simulation studies, we discard the first $15\times 10^5$ TTMCMC iterations as burn-in and store one out of $150$ iterations in the next
$15\times 10^5$ iterations, to obtain $10000$ TTMCMC realizations. The time taken is only one minute for the entire exercise.

\subsection{Results}
\label{subsec:results_realdata}
\begin{figure}
\centering
\subfigure[Trace plot of $\mu_{11}$.]{ \label{fig:sim11_trace_mu1}
\includegraphics[width=7cm,height=5cm]{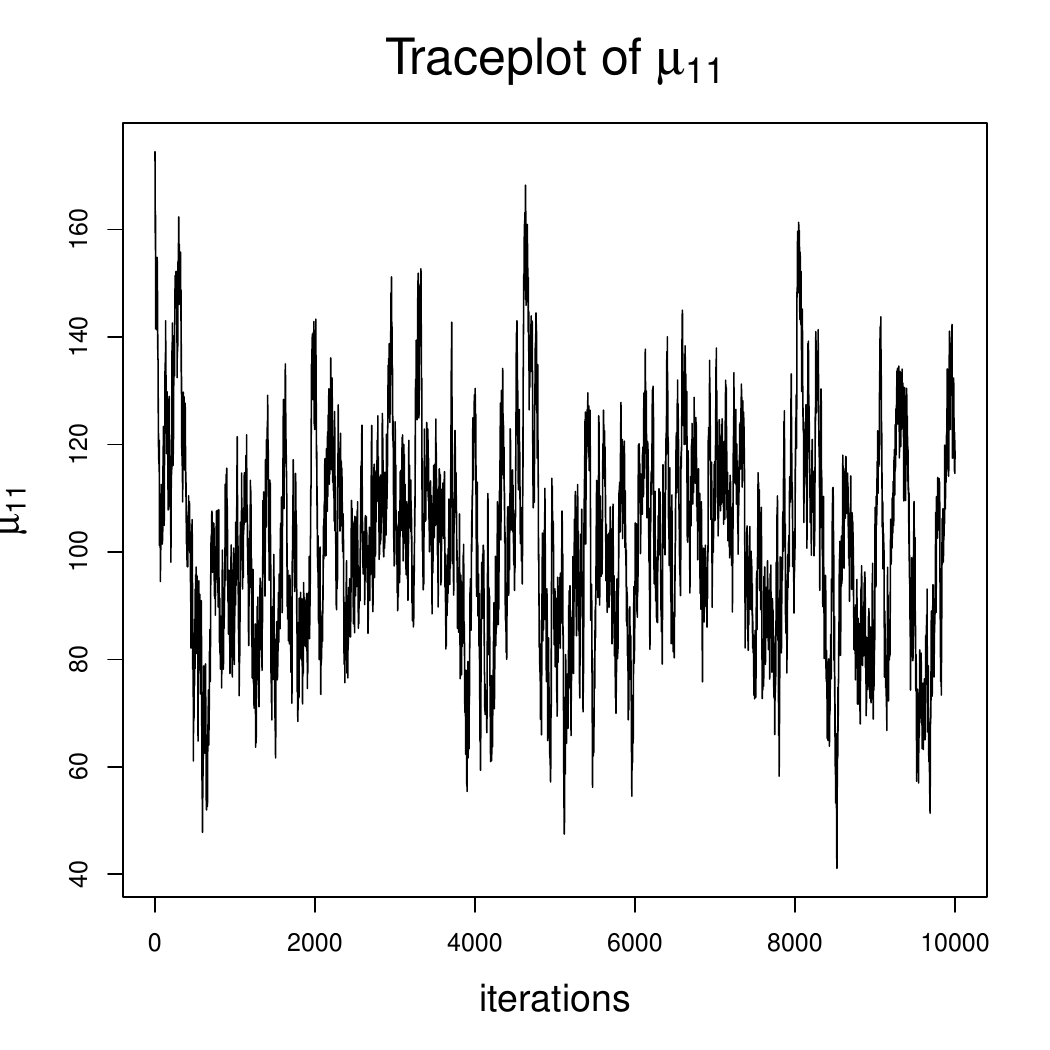}}
\hspace{2mm}
\subfigure[Trace plot of $\mu_{21}$.]{ \label{fig:sim11_trace_mu2}
\includegraphics[width=7cm,height=5cm]{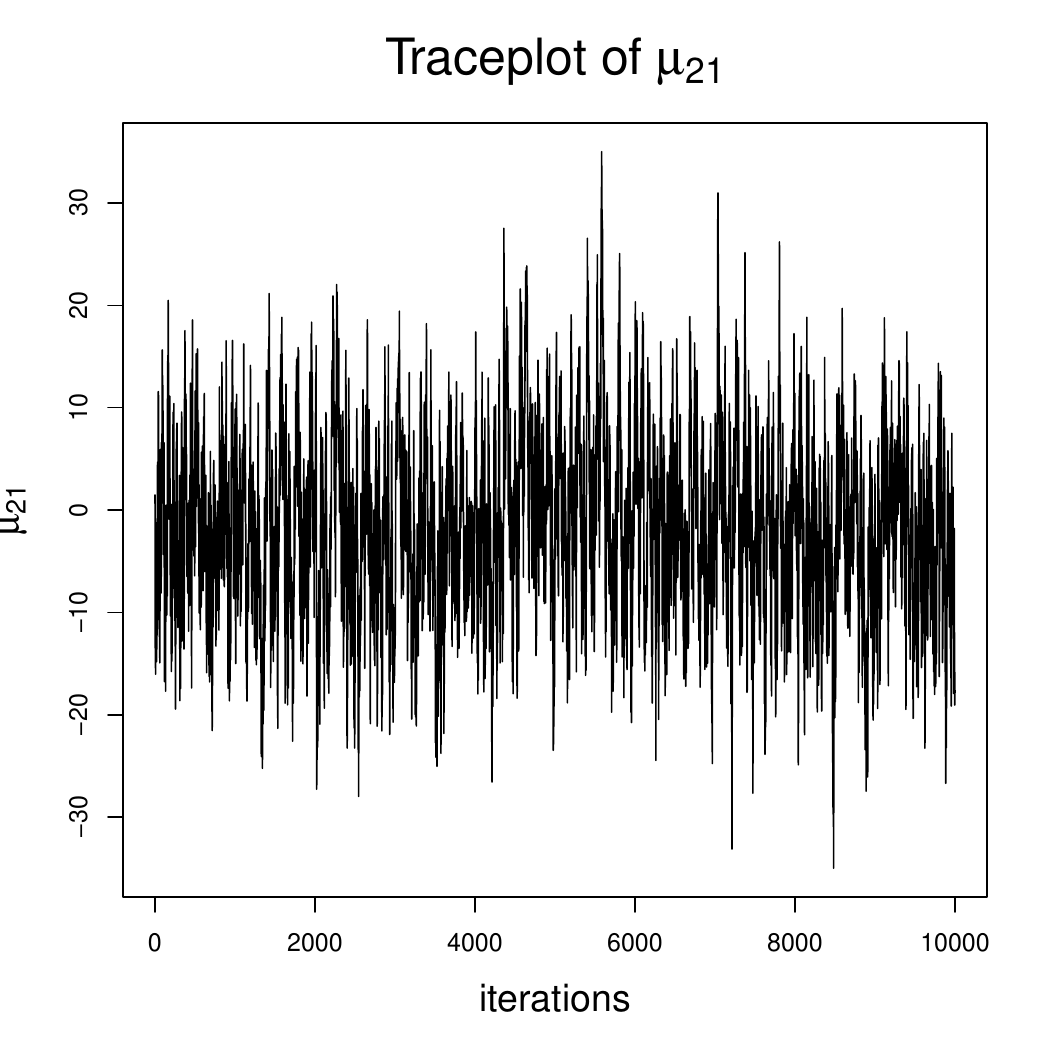}}\\
\vspace{2mm}
\subfigure[Trace plot of $\Omega_{11,1}$.]{ \label{fig:sim11_trace_om1}
\includegraphics[width=7cm,height=5cm]{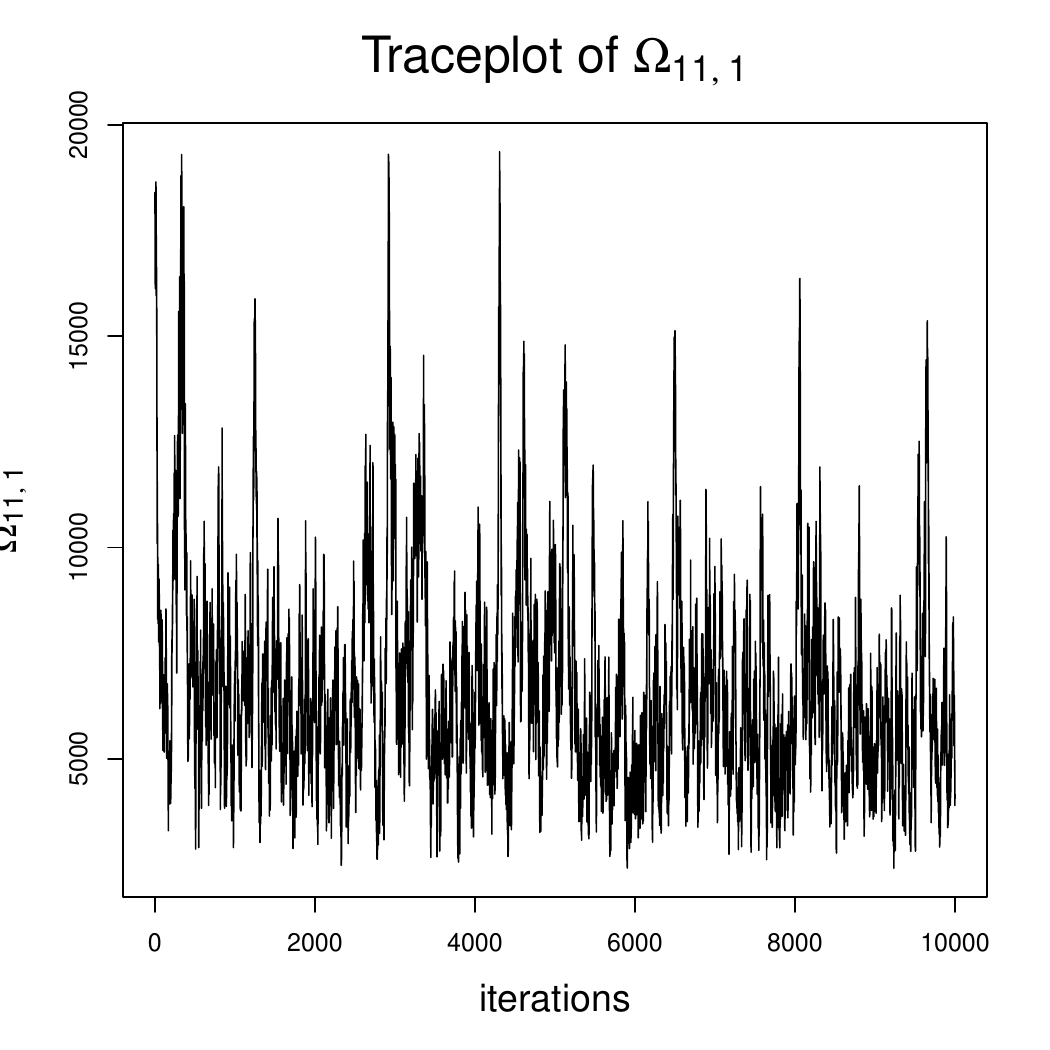}}
\hspace{2mm}
\subfigure[Trace plot of $\Omega_{12,1}$.]{ \label{fig:sim11_trace_om2}
\includegraphics[width=7cm,height=5cm]{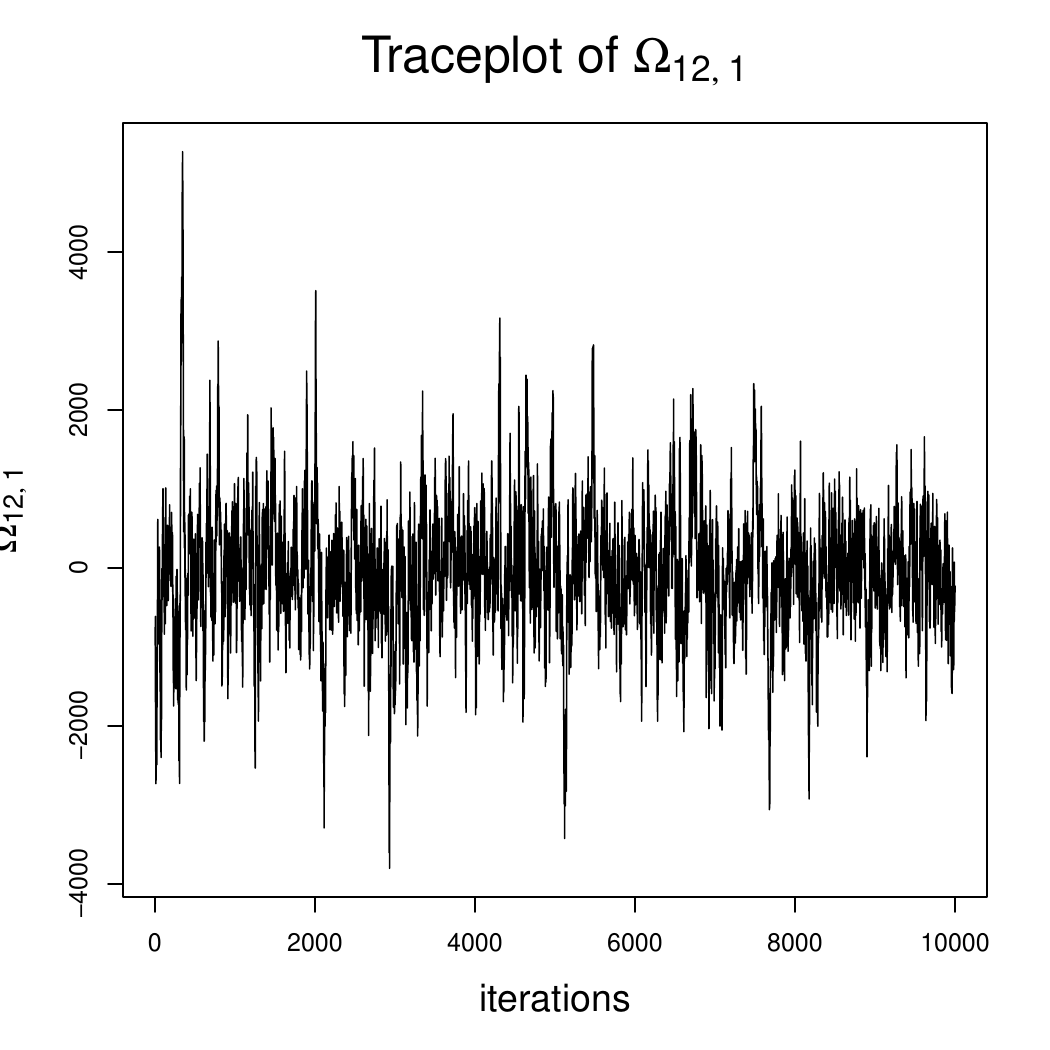}}\\
\vspace{2mm}
\subfigure[Trace plot of $\Omega_{22,1}$.]{ \label{fig:sim11_trace_om4}
\includegraphics[width=7cm,height=5cm]{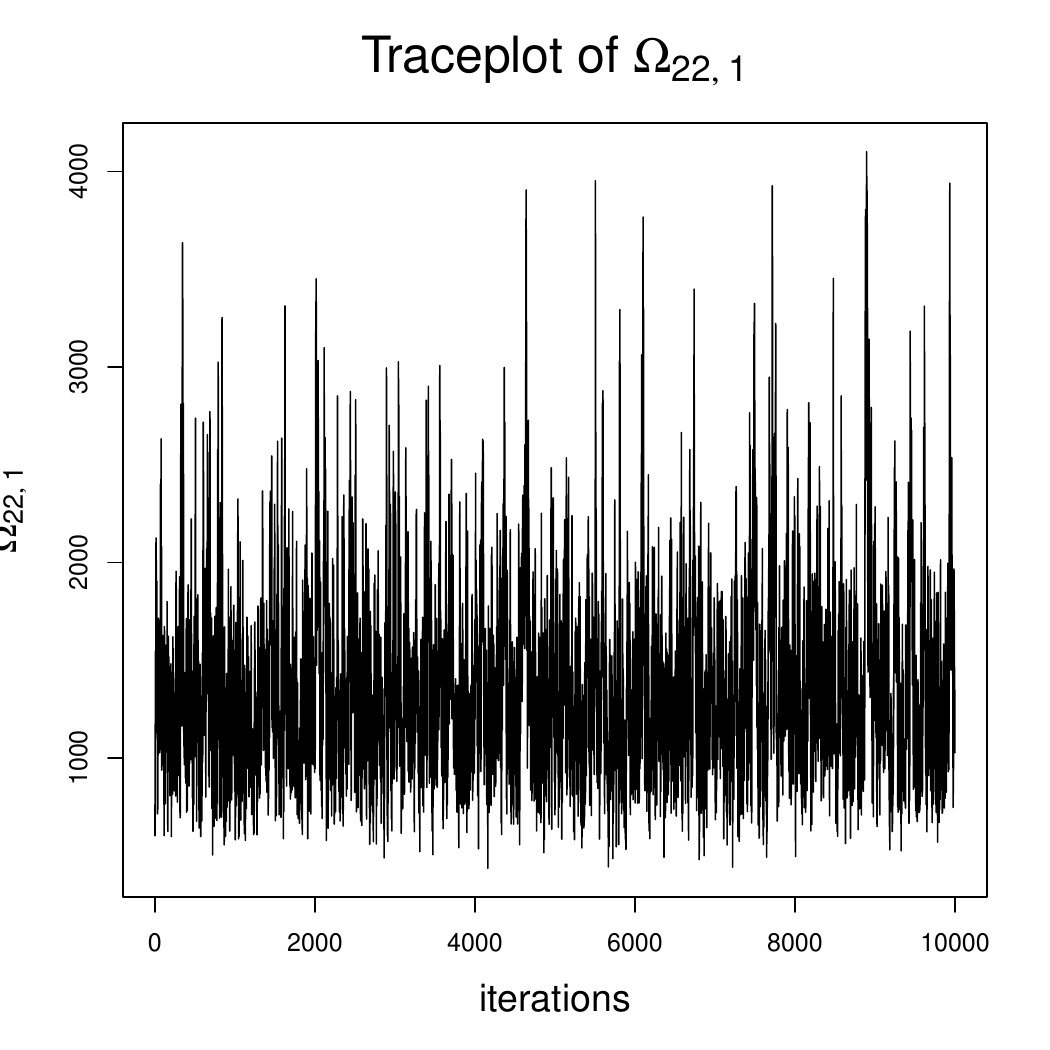}}
\caption{{\bf TTMCMC for the real data case:} Trace plots of $\mu_{11}$, $\mu_{21}$, $\Omega_{11,1}$, $\Omega_{12,1}$ and $\Omega_{22,1}$.} 
\label{fig:sim11_trace_plots}
\end{figure}

\begin{figure}
\centering
\subfigure[Posterior of $\mu_{11}$.]{ \label{fig:sim11_mu1}
\includegraphics[width=7cm,height=5cm]{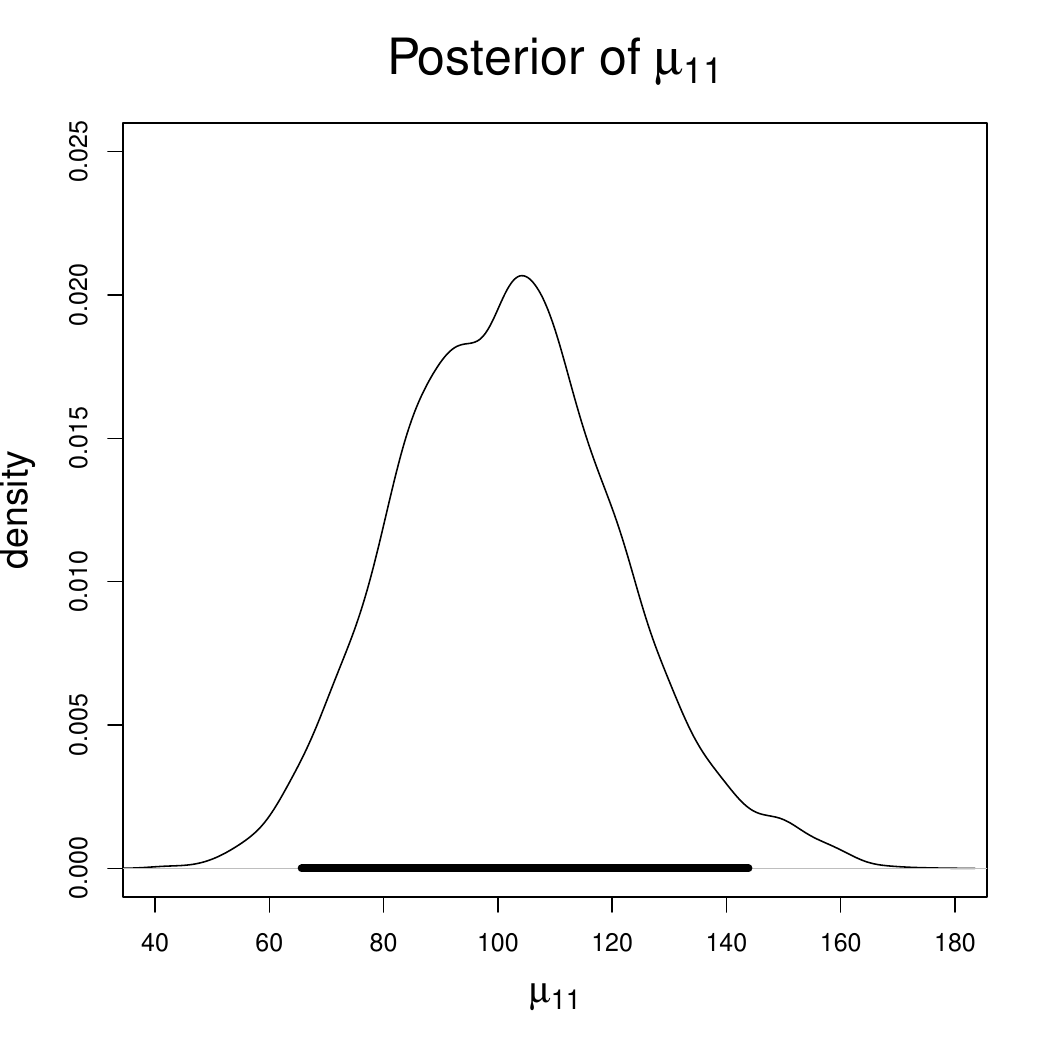}}
\hspace{2mm}
\subfigure[Posterior of $\mu_{21}$.]{ \label{fig:sim11_mu2}
\includegraphics[width=7cm,height=5cm]{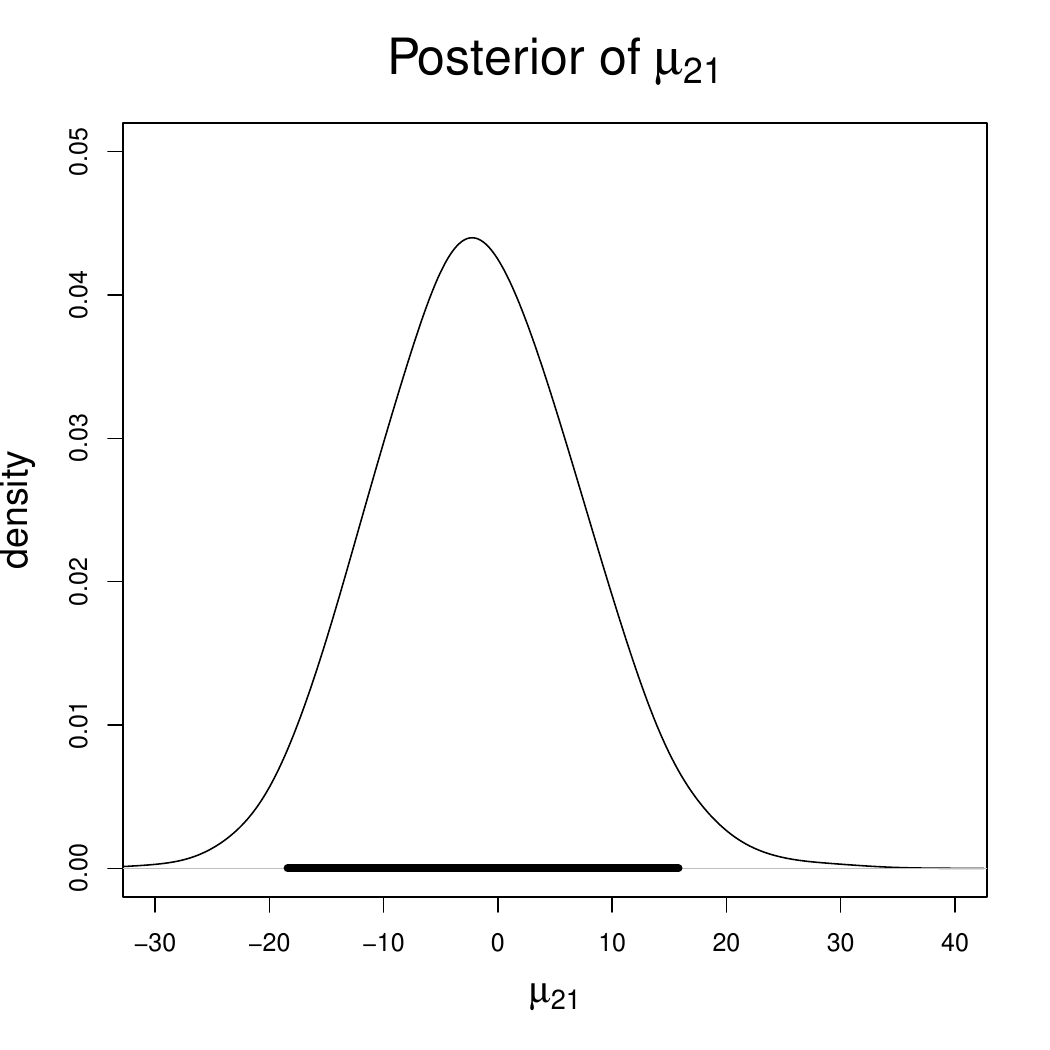}}\\
\vspace{2mm}
\subfigure[Posterior of $\Omega_{11,1}$.]{ \label{fig:sim11_om1}
\includegraphics[width=7cm,height=5cm]{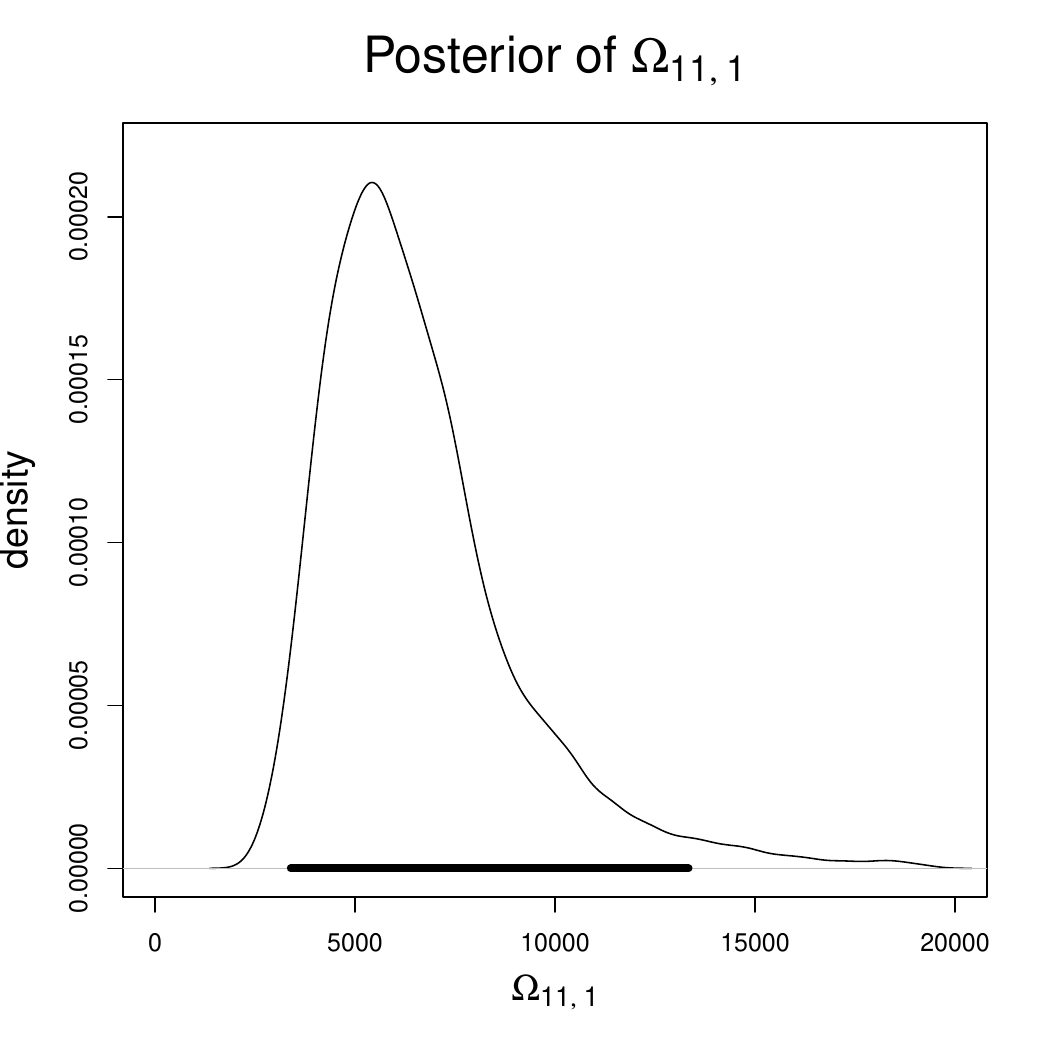}}
\hspace{2mm}
\subfigure[Posterior of $\Omega_{12,1}$.]{ \label{fig:sim11_om2}
\includegraphics[width=7cm,height=5cm]{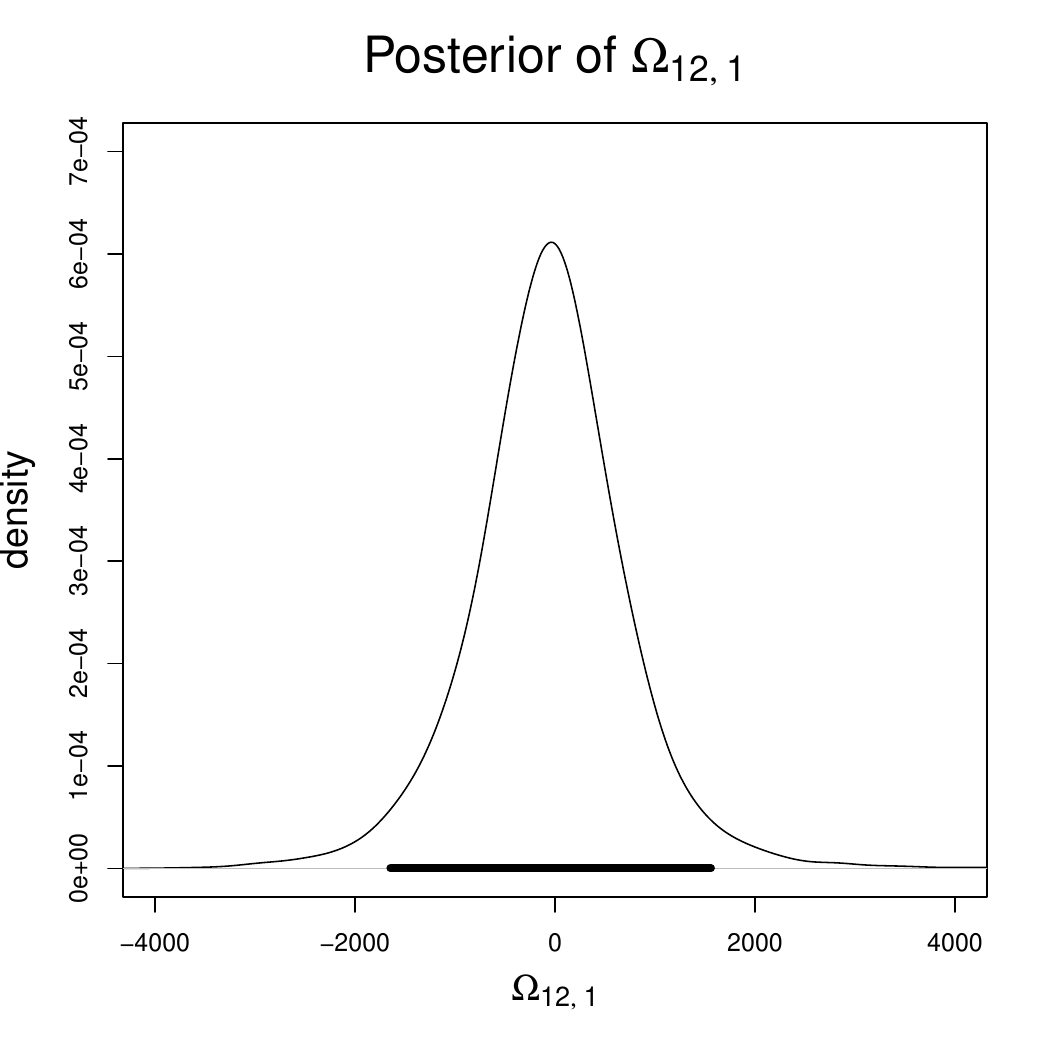}}\\
\vspace{2mm}
\subfigure[Posterior of $\Omega_{22,1}$.]{ \label{fig:sim11_om4}
\includegraphics[width=7cm,height=5cm]{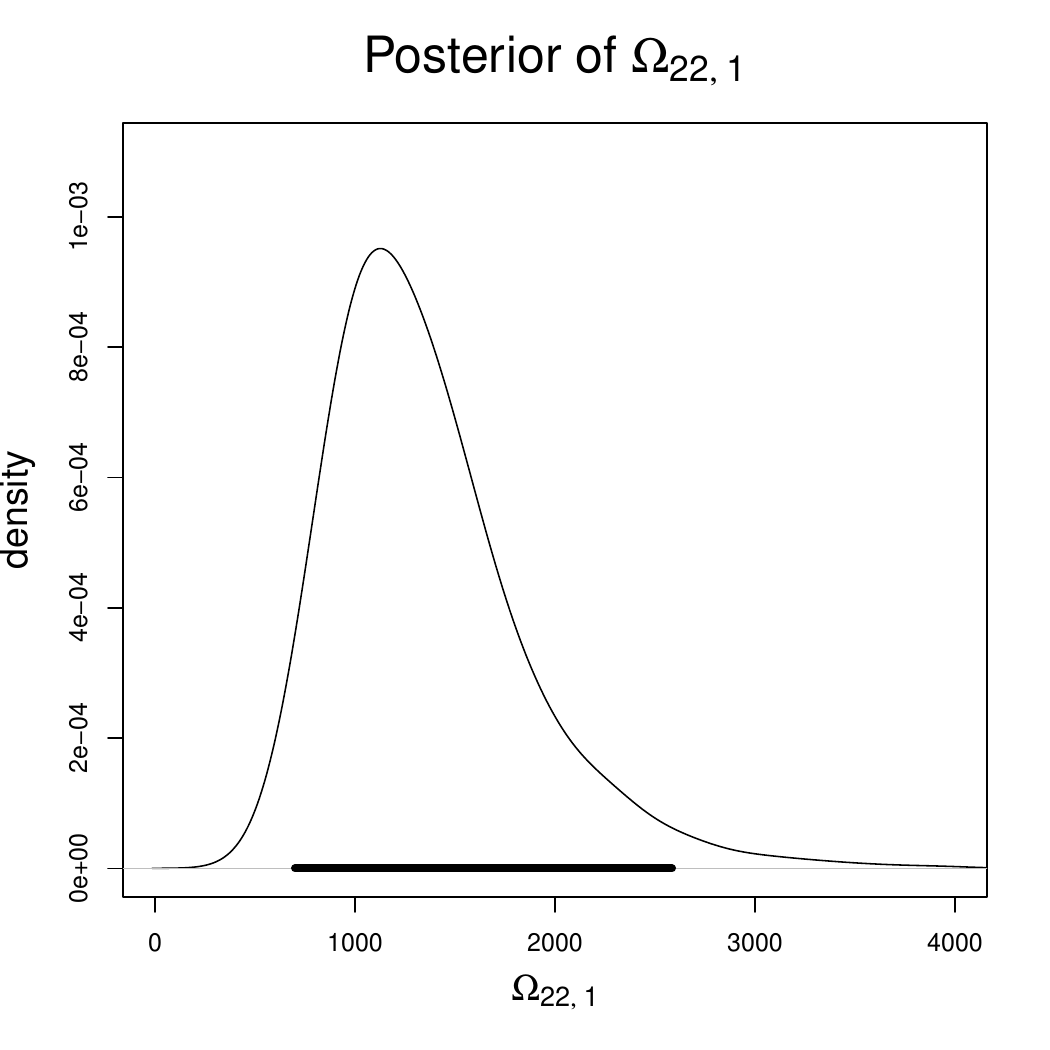}}
\caption{{\bf TTMCMC for the real data case:} Posteriors of $\mu_{11}$, $\mu_{21}$, $\Omega_{11,1}$, $\Omega_{12,1}$ and $\Omega_{22,1}$.} 
\label{fig:sim11_posterior_plots}
\end{figure}

\begin{figure}
\centering
\subfigure[{\bf Realdata analysis:} Histogram of the first co-ordinate of $\hat\phi_i$'s.]{ \label{fig:hist1_realdata}
\includegraphics[width=7cm,height=6cm]{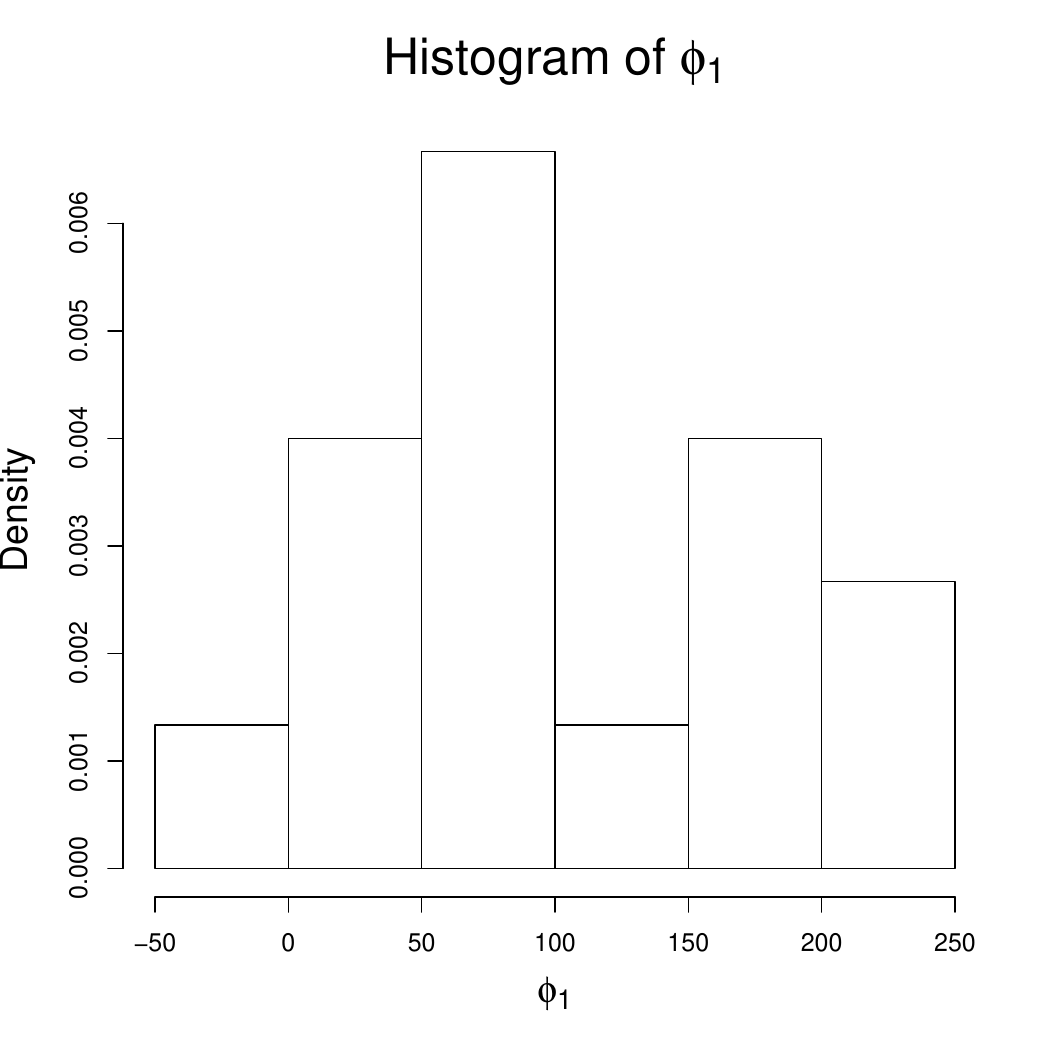}}
\hspace{2mm}
\subfigure[{\bf Realdata analysis:} Histogram of the second co-ordinate of $\hat\phi_i$'s.]{ \label{fig:hist2_realdata}
\includegraphics[width=7cm,height=6cm]{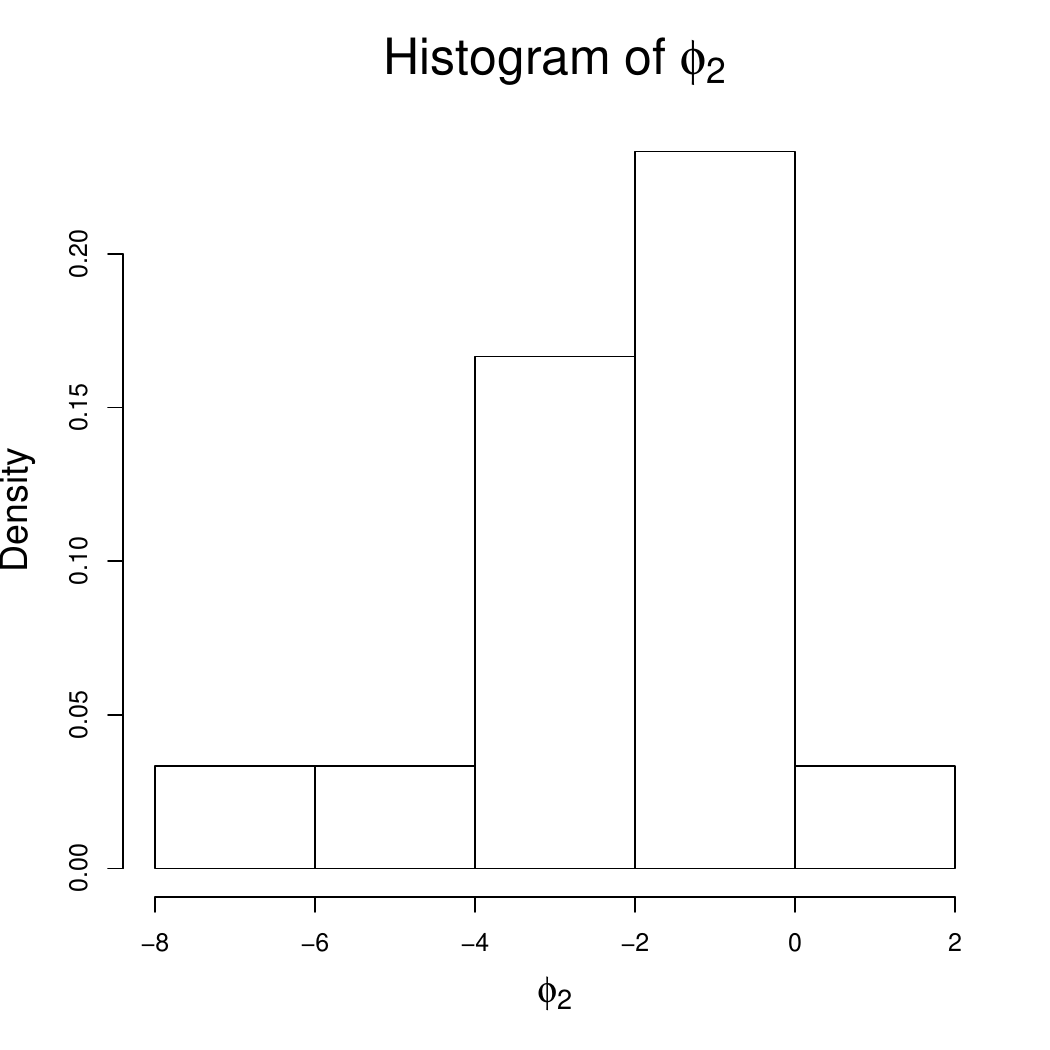}}
\end{figure}

In this applications, the posterior of $M$ turned out to give full mass to $1$; the result
remained the same for the other priors and all starting values. This rules out any possibility of clustering among the companies.
The trace plots of the relevant parameters, shown in Figure \ref{fig:sim11_trace_plots}, signifies quite adequate mixing properties for all the relevant parameters. 
Figure \ref{fig:sim11_posterior_plots} depicts the posterior densities of the relevant parameters, along with their 95\% credible intervals.
In particular, observe that the posterior distributions of the co-ordinates of $\bmu$ are in keeping with the histograms of the co-ordinates
of $\hat\bphi_i=V^{-1}_iU_i$'s, shown in Figures \ref{fig:hist1_realdata} and \ref{fig:hist2_realdata}, respectively, particularly with respect to their modal values.

\begin{figure}
\centering
\subfigure[]{ \label{fig:postpred1}
\includegraphics[width=4.5cm,height=5cm]{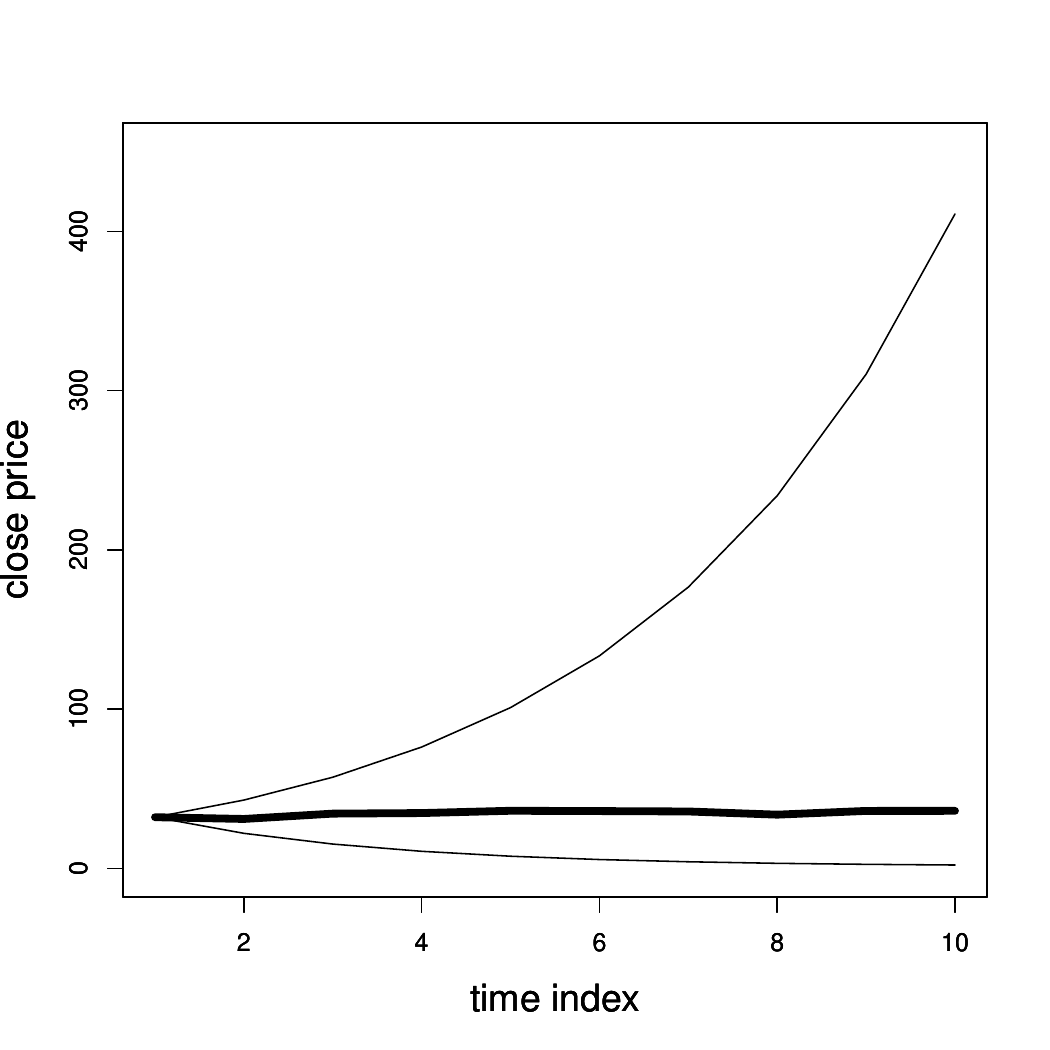}}
\hspace{2mm}
\subfigure[]{ \label{fig:postpred2}
\includegraphics[width=4.5cm,height=5cm]{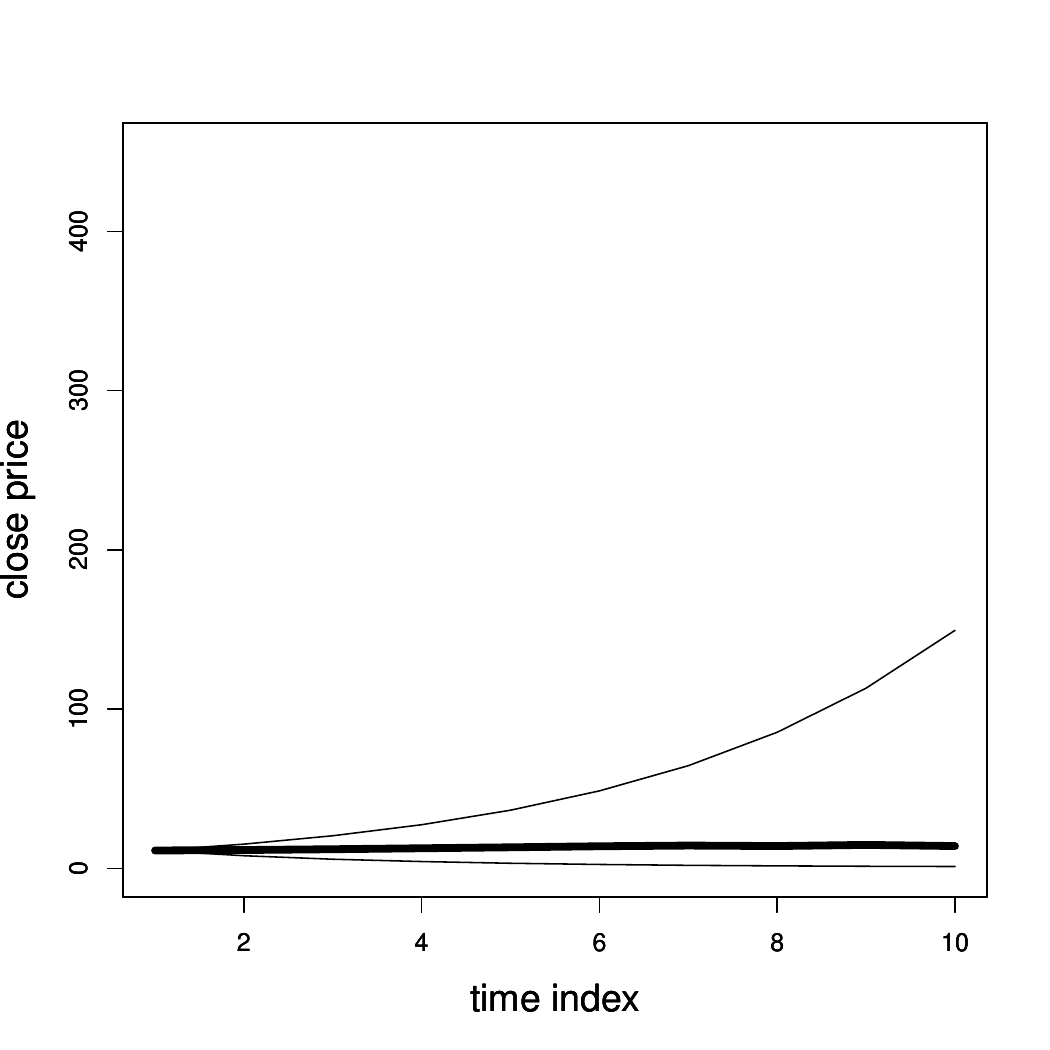}}
\hspace{2mm}
\subfigure[]{ \label{fig:postpred3}
\includegraphics[width=4.5cm,height=5cm]{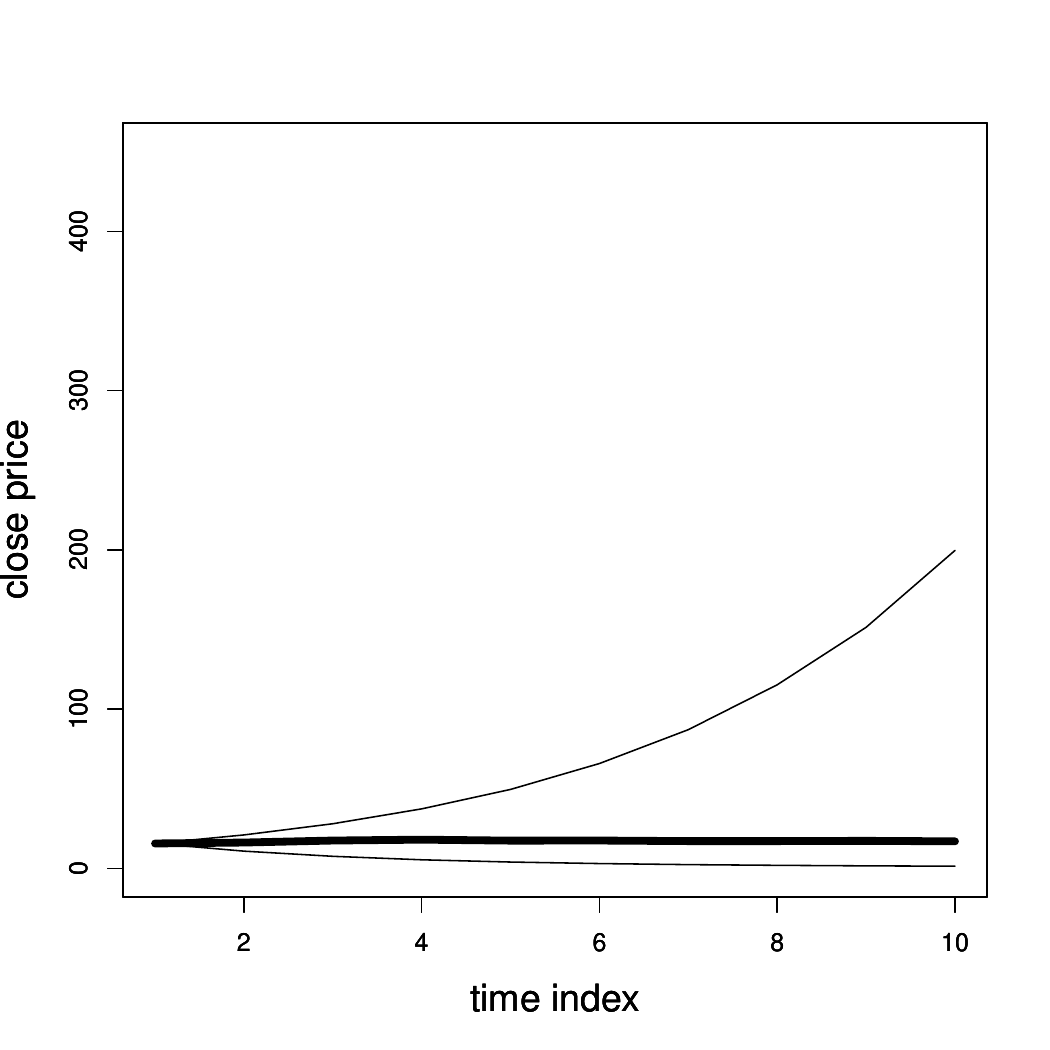}}\\
\vspace{2mm}
\subfigure[]{ \label{fig:postpred4}
\includegraphics[width=4.5cm,height=5cm]{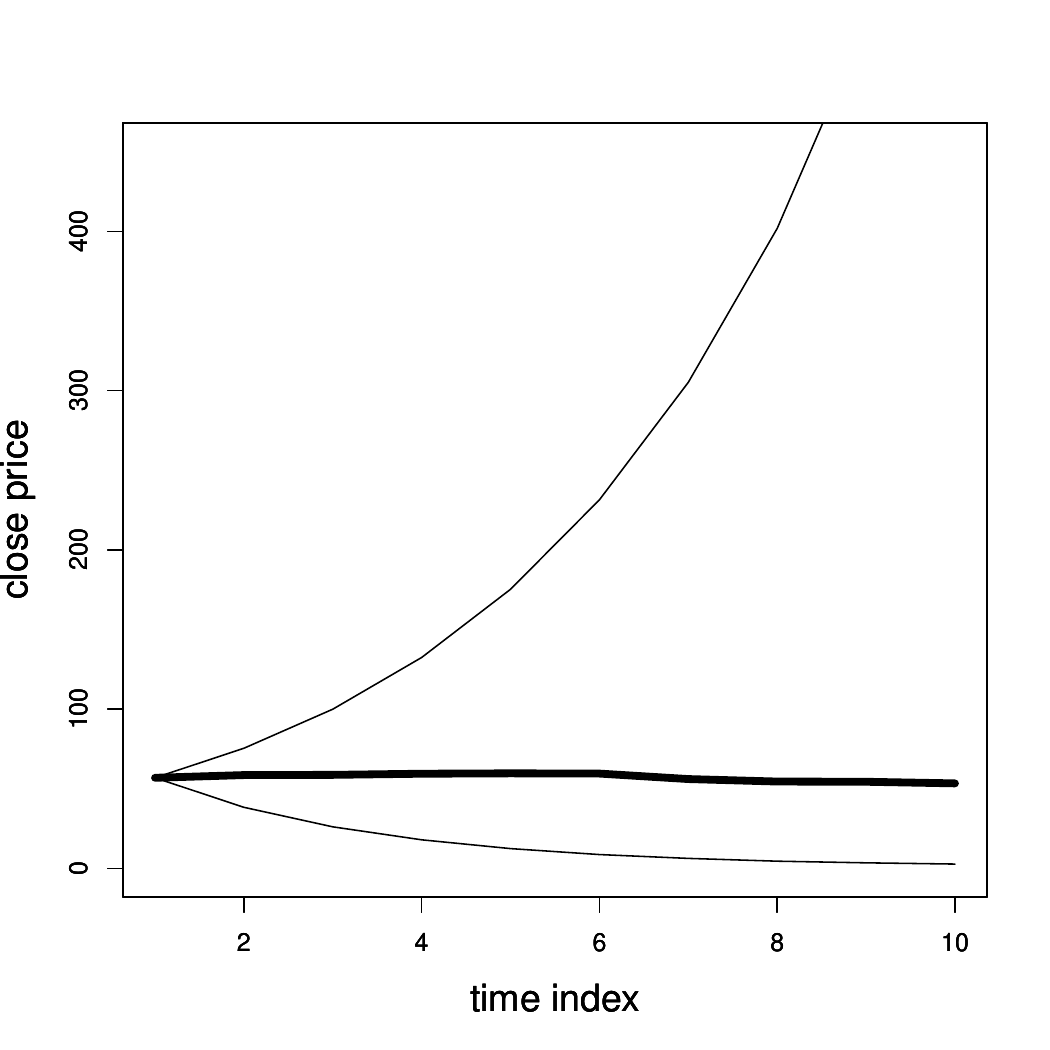}}
\hspace{2mm}
\subfigure[]{ \label{fig:postpred5}
\includegraphics[width=4.5cm,height=5cm]{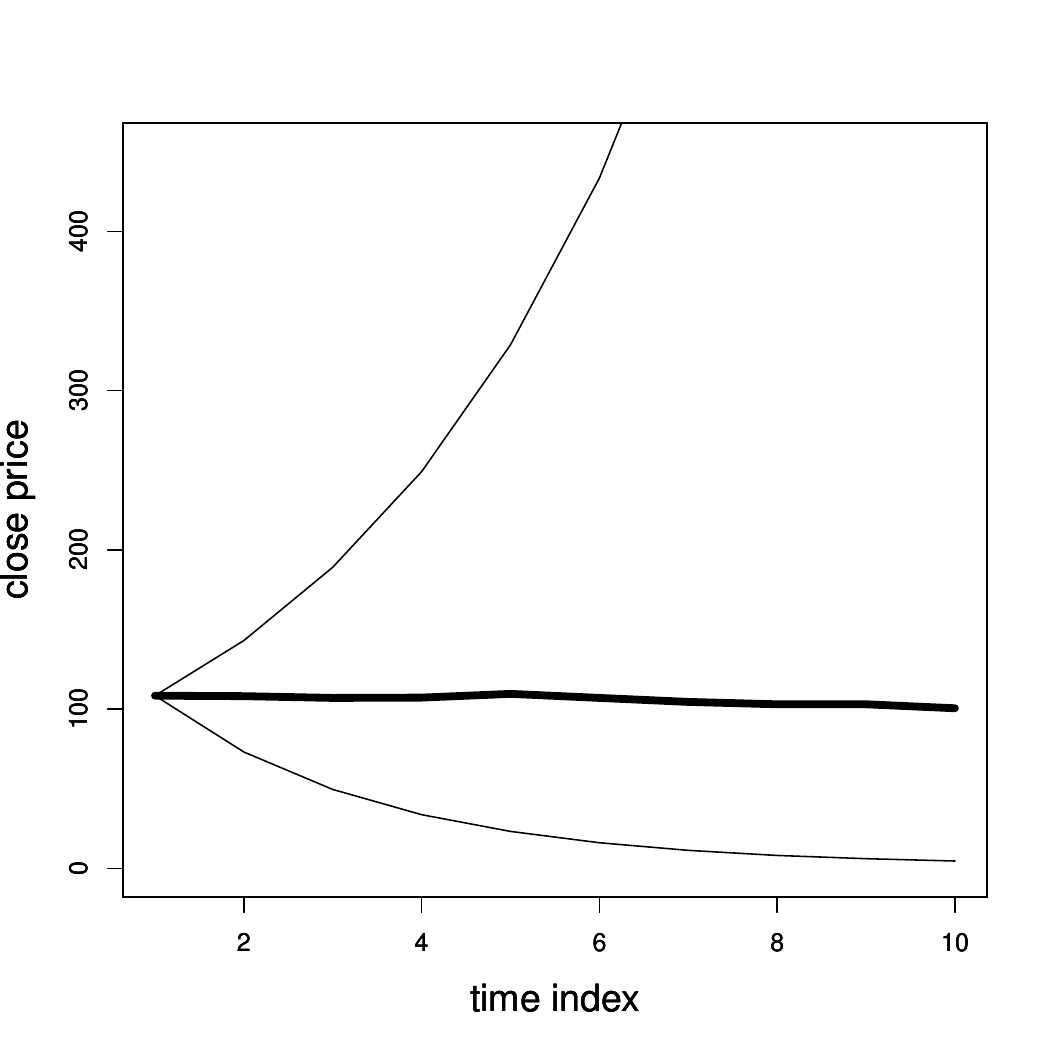}}
\hspace{2mm}
\subfigure[]{ \label{fig:postpred6}
\includegraphics[width=4.5cm,height=5cm]{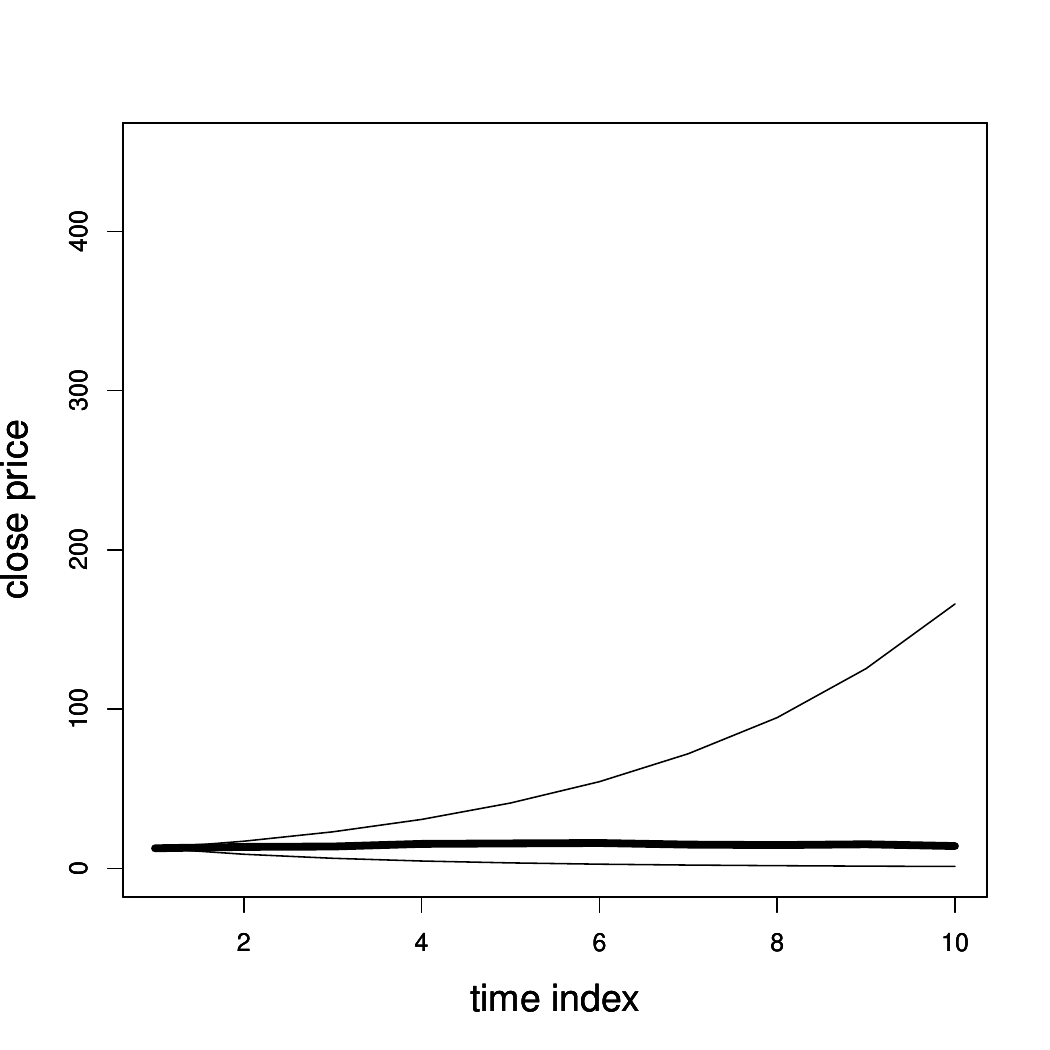}}\\
\vspace{2mm}
\subfigure[]{ \label{fig:postpred7}
\includegraphics[width=4.5cm,height=5cm]{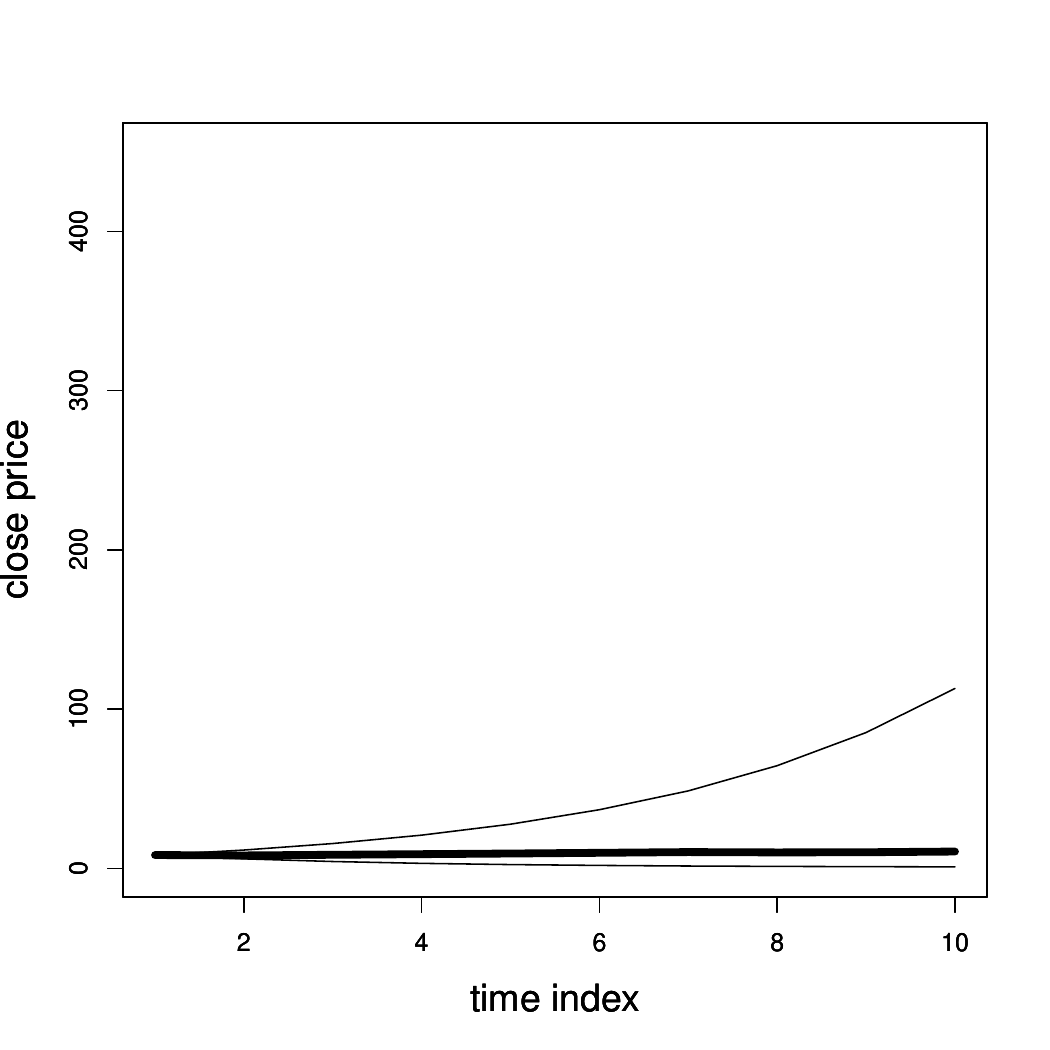}}
\hspace{2mm}
\subfigure[]{ \label{fig:postpred8}
\includegraphics[width=4.5cm,height=5cm]{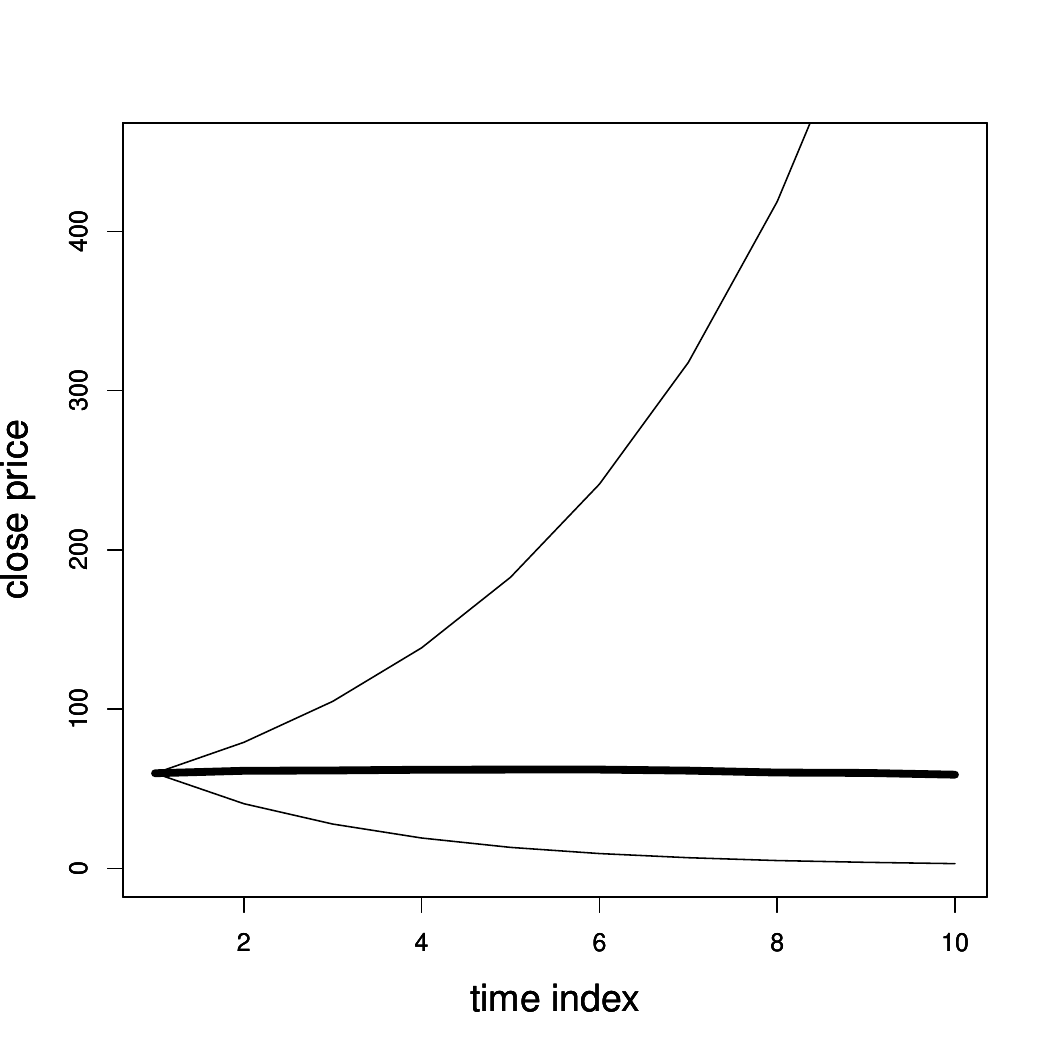}}
\hspace{2mm}
\subfigure[]{ \label{fig:postpred9}
\includegraphics[width=4.5cm,height=5cm]{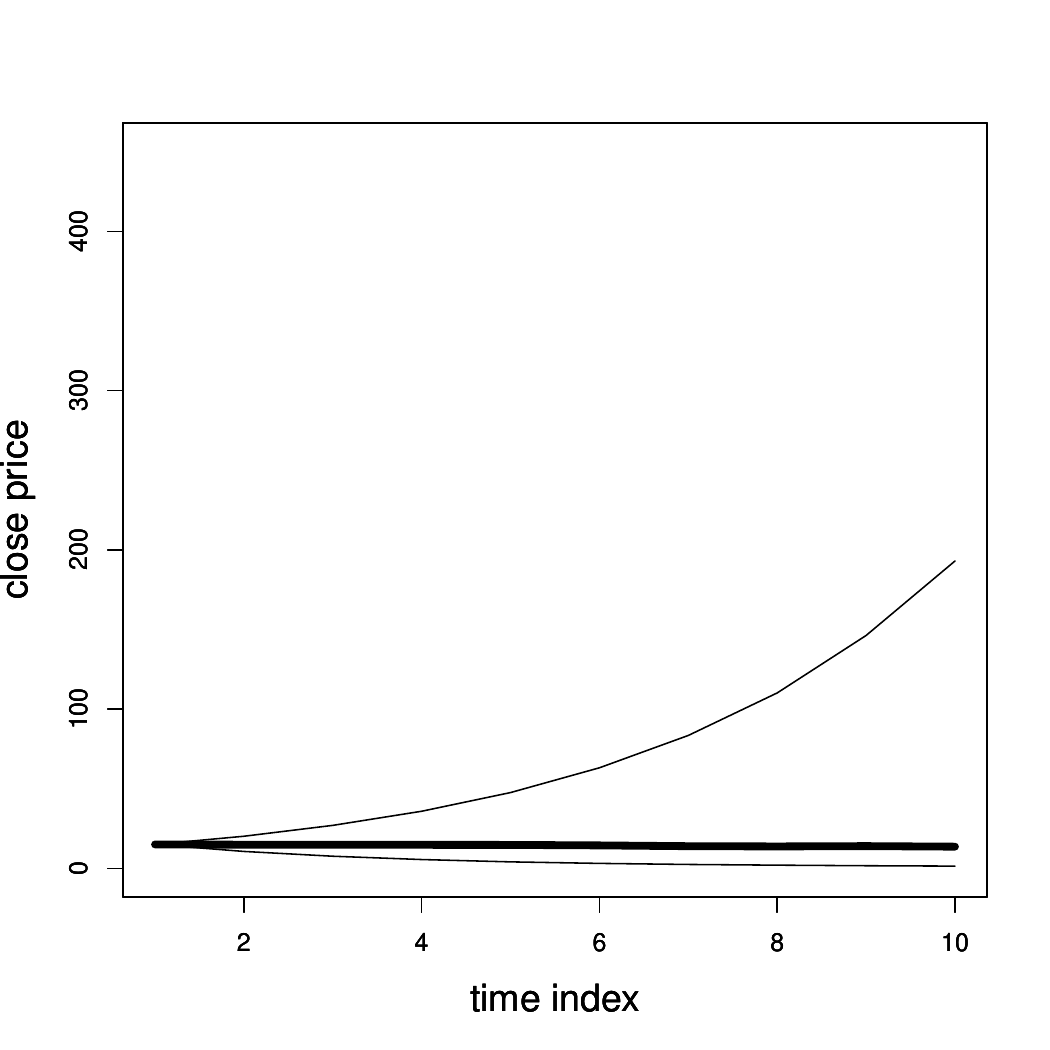}}
\caption{{\bf Realdata analysis:} Posterior predictive pointwise 95\% credible intervals of the times series for the companies. The thick curves
denote the actual data.}
\label{fig:postpred_plots1}
\end{figure}

\begin{figure}
\centering
\subfigure[]{ \label{fig:postpred10}
\includegraphics[width=4.5cm,height=5cm]{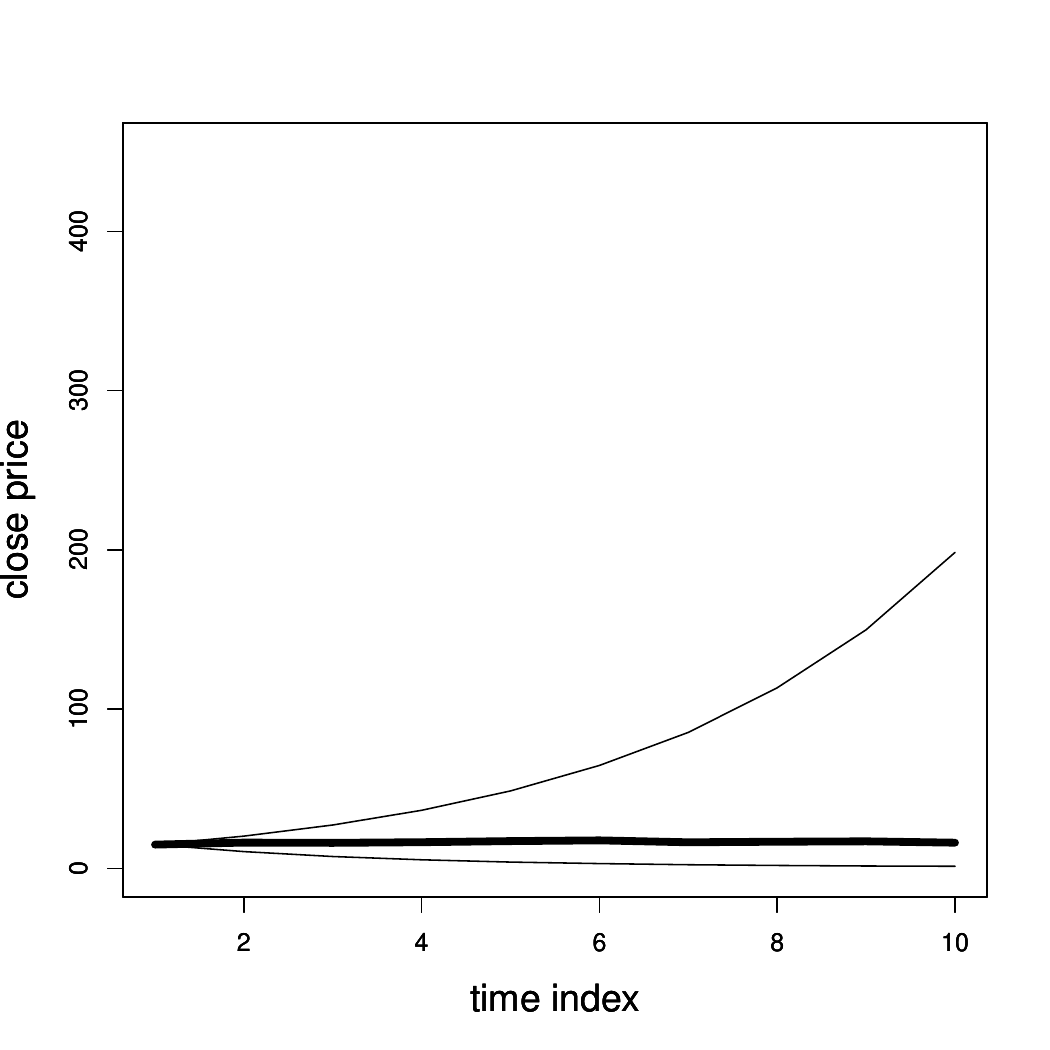}}
\hspace{2mm}
\subfigure[]{ \label{fig:postpred11}
\includegraphics[width=4.5cm,height=5cm]{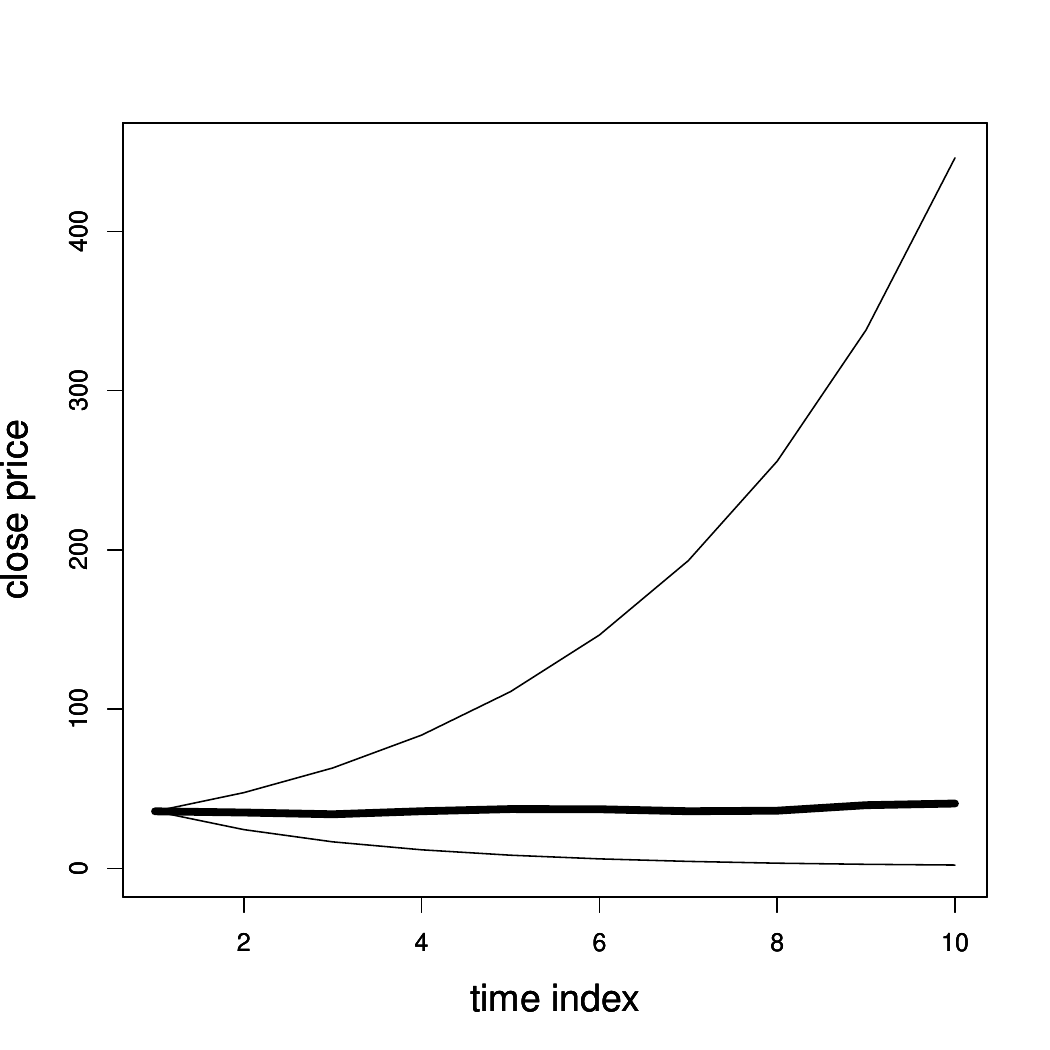}}
\hspace{2mm}
\subfigure[]{ \label{fig:postpred12}
\includegraphics[width=4.5cm,height=5cm]{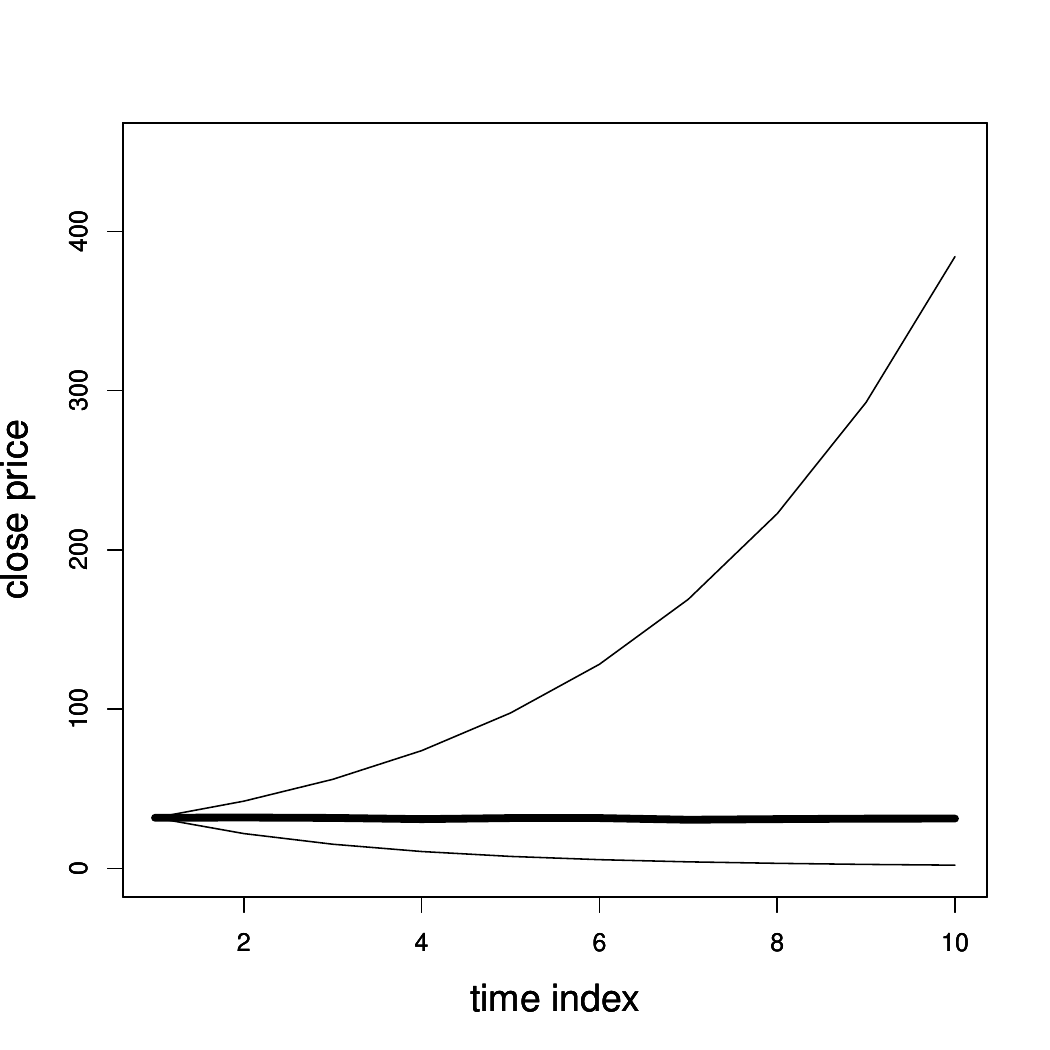}}\\
\vspace{2mm}
\subfigure[]{ \label{fig:postpred13}
\includegraphics[width=4.5cm,height=5cm]{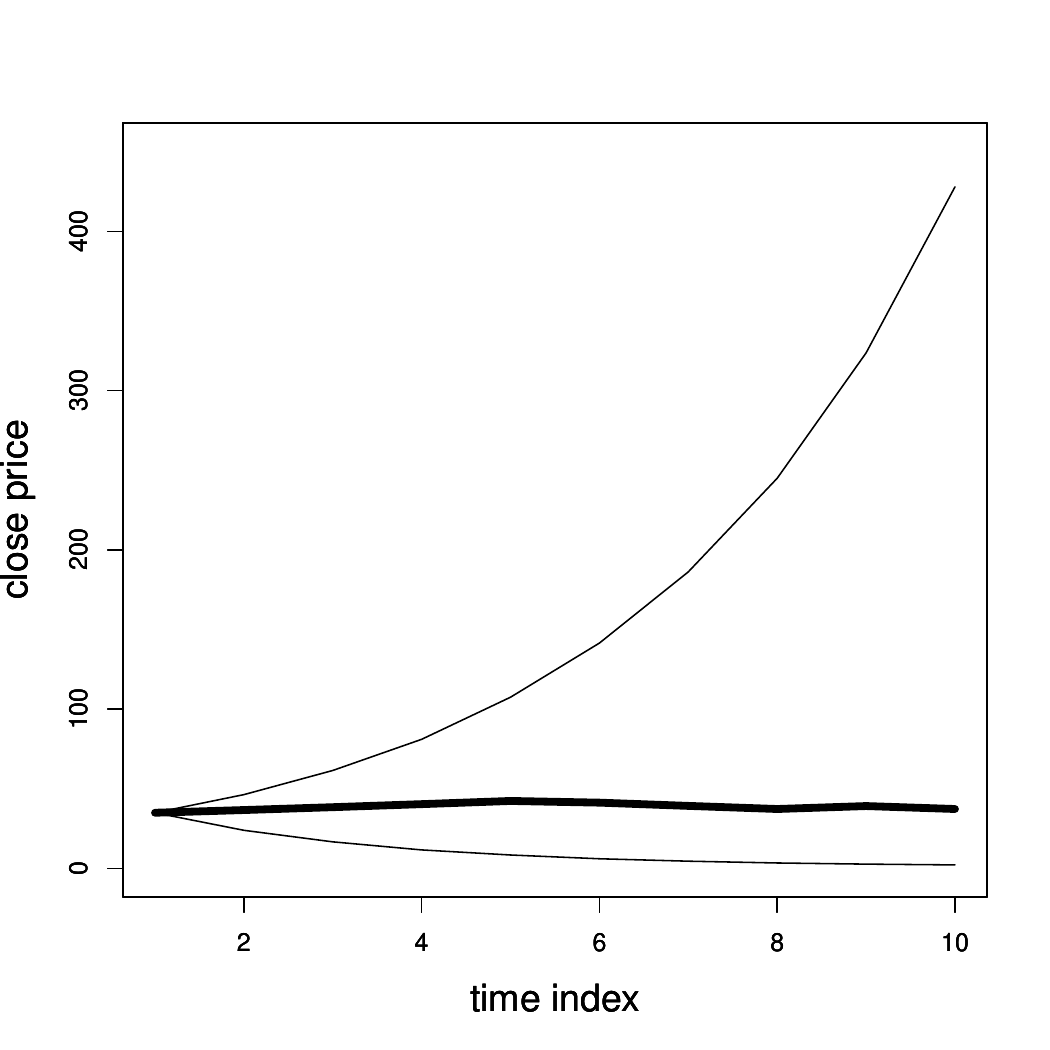}}
\hspace{2mm}
\subfigure[]{ \label{fig:postpred14}
\includegraphics[width=4.5cm,height=5cm]{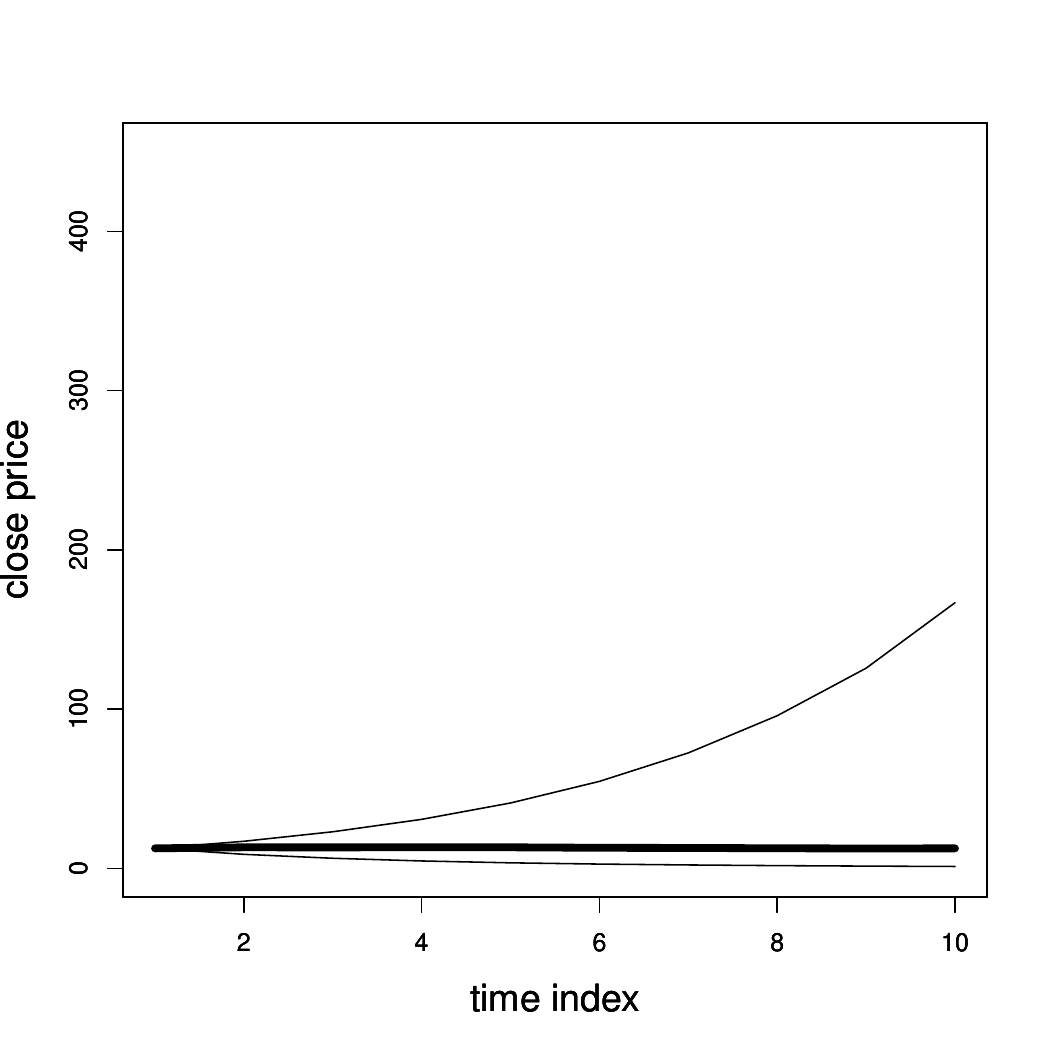}}
\hspace{2mm}
\subfigure[]{ \label{fig:postpred15}
\includegraphics[width=4.5cm,height=5cm]{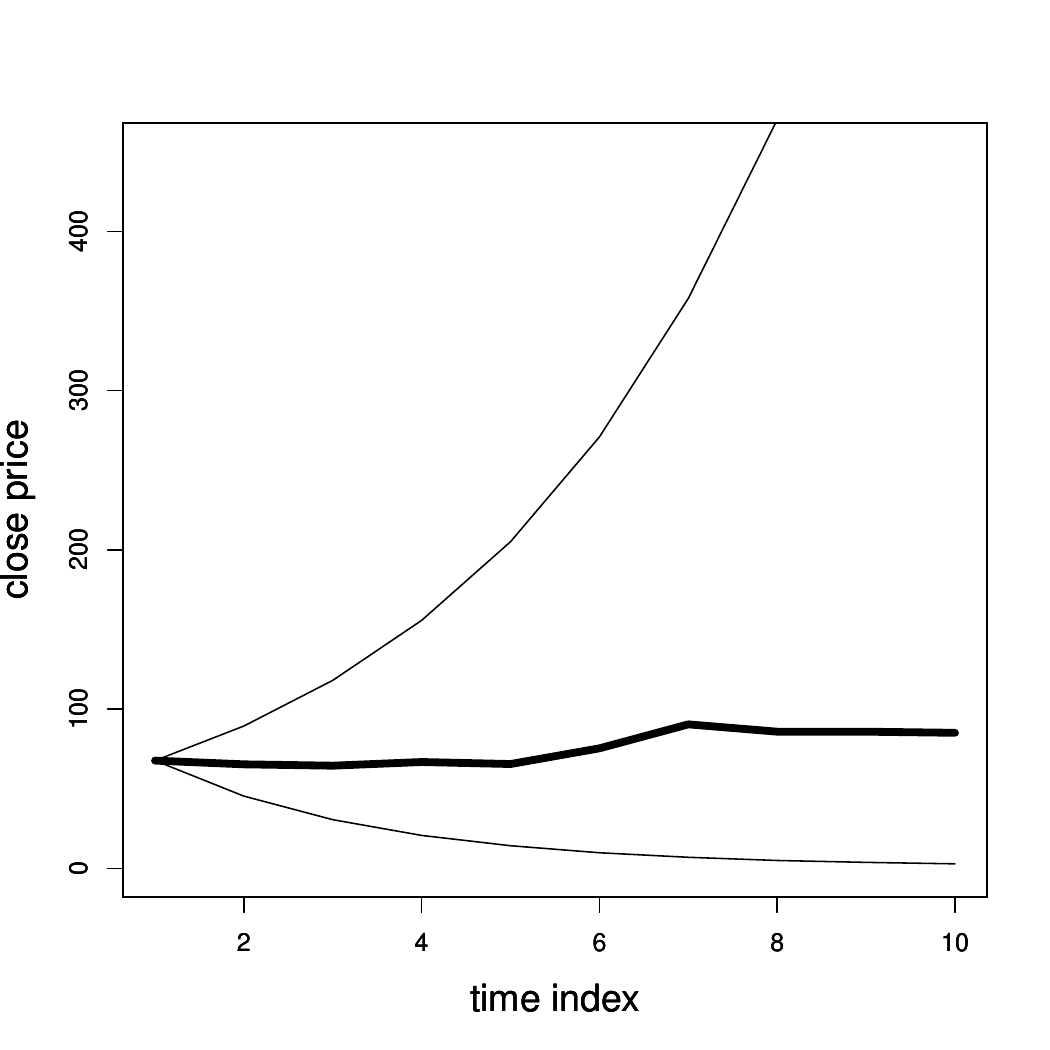}}
\caption{{\bf Realdata analysis:} Posterior predictive pointwise 95\% credible intervals of the times series for the companies. The thick curves
denote the actual data.}
\label{fig:postpred_plots2}
\end{figure}

Additionally, we have obtained the posterior predictive distributions of the $15$ different $SDE$'s. To obtain these, for each company $i=1,\ldots,15$,
for every TTMCMC-generated realization of the parameters, we simulated $\bphi_i$ from the mixture distribution (\ref{eq:bivariate_phi}), using which
we obtained a realization of the $SDE$ (\ref{eq:sde_realdata}). Thus, for each $i$, we obtained $10000$ realizations of (\ref{eq:sde_realdata}), using which
we computed point-wise 95\% credible intervals for each time point. The credible intervals and the actual time series data for the $15$ companies 
are shown in Figures \ref{fig:postpred_plots1} and \ref{fig:postpred_plots2}. We focus on only the first $10$ time points as the upper limits of the credible
intervals increase very rapidly with time. The diagrams show that the actual time series data are wholly contained in the credible intervals, for all the companies,
indicating adequate model fit.

\section{Summary and discussion}
\label{sec:conclusion}

In this work, we consider the random effect parameter having mixture of normal distributions.
Inherent flexibility of the mixture framework offers far greater generality compared to normal
distributions considered so far in this direction of research.
Even in this setup we are able to prove strong consistency and asymptotic normality of the $MLE$  
and in the corresponding Bayesian framework we establish posterior consistency and asymptotic normality of 
the  posterior distribution, our results encompassing both $iid$ and 
non-$iid$ situations, for both the paradigms. 
It is important to note that 
no extra assumptions are needed besides the assumptions of \ctn{Maitra14a} and \ctn{Maitra14b} to conclude the corresponding 
results under the consideration of normal mixture. 
In other words, without any extra assumption 
we could achieve asymptotic results regarding the random effect parameters having a much 
larger class of distributions, which seems to hold importance from theoretical and practical perspectives. 
Importantly, we have also developed the Bayesian asymptotic theory when the number of mixture components
is unknown and considered random. 

Although \ctn{Maud16} also considered the $SDE$ mixture model, they required an extra, strong assumption
to prove even weak consistency of $MLE$ in the $iid$ setup, compared to our stronger, almost sure convergence result, without
the assumption. Moreover, the weak consistency result in the $iid$ situation is the only asymptotic result they provided, 
while we delved into consistency and asymptotic normality of both $MLE$ and the Bayesian posterior distribution, in both
$iid$ and non-$iid$ setups, even when the number mixture components is random.

In all our simulation studies and the real stock market data analysis we fitted the data assuming that the number of mixture components is unknown and 
implemented our Bayesian methodology with TTMCMC in those variable-dimensional models. Excellent performance of our methods in all the cases lead us to recommend the 
variable-dimensional Bayesian paradigm to be implemented with TTMCMC. In contrast, the classical $BIC$ based method of \ctn{Maud16} often yielded poor performance 
when the number of mixtures is unknown.

In this work we have considered one dimensional $SDE$s. The generalization of our asymptotic theories to high dimensions 
will be considered in our future endeavors.


\section*{Acknowledgments}
We are extremely grateful to the referee whose comments have led to improved presentation of our manuscript. 
The first author gratefully acknowledges her NBHM Fellowship, Govt. of India.

\newpage
renewcommand\thefigure{S-\arabic{figure}}
\renewcommand\thetable{S-\arabic{table}}
\renewcommand\thesection{S-\arabic{section}}

\setcounter{section}{0}
\setcounter{figure}{0}
\setcounter{table}{0}

\begin{center}
{\bf \LARGE Supplementary Material}
\end{center}

\section{Proof that differentiation can be passed under the integral sign}
\label{sec:diff_int}
 
Let us denote the marginal distribution of $\{X_i(t): t\in [0,T]\}$ by $Q_{\theta}^i$ 
on $(C_{T},\mathcal C_{T})$, where $C_{T}$ is the space of real continuous functions on $[0,T]$ 
and $\mathcal C_{T}$ is the corresponding $\sigma$-algebra. 
Recall that \\$\theta=(\gamma,\beta)=(a_1,\ldots,a_M,\mu_1,\ldots,\mu_M,\omega_1^2,\ldots,\omega_M^2)$.

Let $\tau=(1,0,\ldots,0,\mu_1^*,\ldots,\mu_M^*,{\omega_1^*}^2,\ldots,{\omega_M^*}^2)$ and set 
\begin{align}
p_1(\theta)&=\frac{\lambda_1(X_1,\theta)}{\lambda_1(X_1,\tau)}\notag\\
&=\frac{\sum_{k=1}^M \frac{a_k}{\left(1+\omega_k^2V_1\right)^{1/2}}\exp\left[-\frac{V_1}{2\left(1+\omega_k^2V_1\right)}
\left(\mu_k-\frac{U_1}{V_1}\right)^2\right]\exp\left(\frac{U^2_1}{2V_1}\right)}{\frac{1}{\left(1+{\omega_1^*}^2V_1\right)^{1/2}}\exp\left[-\frac{V_1}{2\left(1+{\omega_1^*}^2V_1\right)}
\left(\mu_1^*-\frac{U_1}{V_1}\right)^2\right]\exp\left(\frac{U^2_1}{2V_1}\right)}\notag\\
&=\sum_{k=1}^M\frac{a_k(1+{\omega_1^*}^2V_1)^{\frac{1}{2}}}{(1+{\omega_k}^2V_1)^{\frac{1}{2}}}\exp\left[-\frac{V_1}{2\left(1+\omega_k^2V_1\right)}\left(\mu_k-\frac{U_1}{V_1}\right)^2+\frac{V_1}{2\left(1+{\omega_1^*}^2V_1\right)}
\left(\mu_1^*-\frac{U_1}{V_1}\right)^2\right]
\label{eq:p_1}
\end{align}
so that $\int_{C_T} p_1(\theta)dQ_{\tau}^1=1$. The assurance of interchange of integration with respect to 
$Q_{\tau}^1$ and differentiation with respect to $\theta$ 
implies interchange of integration with respect to $Q_{\theta}^1$ and differentiation with respect to $\theta$. 
So, here we will justify interchange of integration with respect to $Q_{\tau}^1$ and differentiation 
with respect to $\theta$. 

Note that
\begin{equation}
\frac{\partial p_1(\theta)}{\partial a_k}\leq (1+{\omega_1^*}^2V_1)^{\frac{1}{2}}\exp\left[-\frac{V_1}{2\left(1+\omega_k^2V_1\right)}\left(\mu_k-\frac{U_1}{V_1}\right)^2+\frac{V_1}{2\left(1+{\omega_1^*}^2V_1\right)}
\left(\mu_1^*-\frac{U_1}{V_1}\right)^2\right].
\label{eq:p_2}
\end{equation}

Let
$$\overline{\mu}^2=\max\{\mu^2_k,{\mu_1^*}^2\};\quad \underline{\mu}^2=\min\{\mu^2_k,{\mu_1^*}^2\};\quad 
\underline{\omega}=\min\{\omega_k,\omega_1^*\}.$$
Now 
\begin{align}
\mu_k^2-2\mu_k\frac{U_1}{V_1}+\frac{U_1^2}{V_1^2}\geq \underline{\mu}^2-2\mu_k\frac{U_1}{V_1}+\frac{U_1^2}{V_1^2}\notag\\
\implies -\left(\mu_k^2-2\mu_k\frac{U_1}{V_1}+\frac{U_1^2}{V_1^2}\right)&\leq -\left(\underline{\mu}^2-2\mu_k\frac{U_1}{V_1}
+\frac{U_1^2}{V_1^2}\right).
\label{eq:p_3}
\end{align}

Also,
\begin{equation}
{\mu_1^*}^2-2\mu_1^*\frac{U_1}{V_1}+\frac{U_1^2}{V_1^2}\leq \overline{\mu}^2-2\mu_1^*\frac{U_1}{V_1}+\frac{U_1^2}{V_1^2},
\label{eq:p_4}
\end{equation}

\begin{equation}
\frac{V_1}{2\left(1+\omega_k^2V_1\right)}\leq\frac{V_1}{2\left(1+\underline{\omega}^2V_1\right)}\quad 
\mbox{and}\quad\frac{V_1}{2\left(1+{\omega_1^*}^2V_1\right)}\leq\frac{V_1}{2\left(1+\underline{\omega}^2V_1\right)}.
\label{eq:p_5}
\end{equation}

Due to (\ref{eq:p_3}), (\ref{eq:p_4}) and (\ref{eq:p_5}) it follows that
\begin{align}
\frac{\partial p_1(\theta)}{\partial a_k} & 
\leq (1+{\omega_1^*}^2V_1)^{\frac{1}{2}}\exp\left[-\frac{V_1}{2\left(1+\underline{\omega}^2V_1\right)}
\left( \underline{\mu}^2-2\mu_k\frac{U_1}{V_1}+\frac{U_1^2}{V_1^2}\right)\right.\notag\\
&\quad\quad\quad\quad\quad\quad\quad\quad\quad\quad\left.+ \frac{V_1}{2\left(1+\underline{\omega}^2V_1\right)}
\left(\overline{\mu}^2-2\mu_1^*\frac{U_1}{V_1}+\frac{U_1^2}{V_1^2}\right)\right]\notag\\
&=(1+{\omega_1^*}^2V_1)^{\frac{1}{2}}\exp\left[\frac{V_1}{2\left(1+\underline{\omega}^2V_1\right)}
(\overline{\mu}^2-\underline{\mu}^2)+\frac{U_1}{\left(1+\underline{\omega}^2V_1\right)}(\mu_k-\mu_1^*)\right]\notag\\
&\leq (1+{\omega_1^*}^2V_1)^{\frac{1}{2}}\exp\left[\frac{1}{2\underline{\omega}^2}(\overline{\mu}^2-\underline{\mu}^2)
+\frac{U_1}{\left(1+\underline{\omega}^2V_1\right)}(\mu_k-\mu_1^*)\right]=K_1(U_1,V_1)~\mbox{(say)},\notag
\end{align}
where the last inequality follows as $\frac{V_1}{2(1+\underline{\omega}^2V_1)}< \frac{1}{2\underline{\omega}^2}$.
Note that $K_1(U_1,V_1)$ has finite expectation, that is, integrable with respect to $Q_{\tau}^1$, 
thanks to existence of all order moments of $V_1$, Lemma 1 of \ctn{Maud12}, and the Cauchy-Schwartz inequality.


Further, with $\mu_{\max}=\max\left\{|\mu_1|,\ldots,|\mu_M|\right\}$, note that,
\begin{align}
\left|\frac{\partial p_1}{\partial\mu_k}\right|&\leq K_1(U_1,V_1)[V_1\mu_{\max}+|U_1|];\\
\left|\frac{\partial p_1}{\partial \omega_k^2}\right|&\leq K_1(U_1,V_1)[V_1+(V_1\mu_{\max}+|U_1|)^2];\\
\left|\frac{\partial ^2 p_1}{\partial\mu_k^2}\right|&\leq K_1(U_1,V_1)[V_1+(V_1\mu_{\max}+|U_1|)^2];\\
\left|\frac{\partial ^2 p_1}{\partial\omega_k^2\partial\omega_k^2}\right|&\leq K_1(U_1,V_1)\left[V_1^2+3\left(V_1^{3/2}\mu_{\max}+V_1^{1/2}|U_1|\right)^2+\left(V_1\mu_{\max}+|U_1|\right)^4\right];\\
\left|\frac{\partial ^2 p_1}{\partial\mu_k\partial\omega_k^2}\right|&\leq K_1(U_1,V_1)\left[2\left(V_1^2\mu_{\max}+V_1|U_1|\right)+(V_1\mu_{\max}+|U_1|)^3\right],
\end{align}
where each upper bound has finite expectation. This easily follows due to existence of all order moments of $V_1$ 
and $|U_1|$ (which follows from Lemma 1 of \ctn{Maud12}), compactness of the parameter space, the fact that the 
expectation of $K_1(U_1,V_1)$ is finite, and the Cauchy-Schwartz inequality. Therefore, the interchange is justified.

\section{Upper bounds of the third order partial derivatives of log-likelihood}
\label{sec:ub}
Likelihood corresponding to the $i$-th individual is
\begin{equation*}
\lambda_i(X_i,\theta)=\sum_{k=1}^M a_k f(X_i|\beta_k),
\label{eq:likelihood5_supp}
\end{equation*}
where
\begin{equation*}
f(X_i|\beta_k)=\frac{1}{\left(1+\omega_k^2V_i\right)^{1/2}}\exp\left[-\frac{V_i}{2\left(1+\omega_k^2V_i\right)}
\left(\mu_k-\frac{U_i}{V_i}\right)^2\right]\exp\left(\frac{U^2_i}{2V_i}\right).
\label{eq:likelihood6_supp}
\end{equation*}
Hence we obtain,
\begin{align}
\left|\frac{f}{\lambda}\right|&\leq 1;\\
\left|\frac{1}{\lambda}\frac{\partial f}{\partial\mu}\right|&\leq \left[V\mu_{\max}+V|U|\right];\\
\left|\frac{1}{\lambda}\frac{\partial f}{\partial \omega^2}\right|&\leq \left[V+\left(V\mu_{\max}+|U|\right)^2\right];\\
\left|\frac{1}{\lambda}\frac{\partial ^2 f}{\partial\mu^2}\right|&\leq\left[V+\left(V\mu_{\max}+|U|\right)^2\right];\\
\left|\frac{1}{\lambda}\frac{\partial ^2 f}{\partial\omega^2\partial\omega^2}\right|&\leq\left[V^2+3\left(V^{3/2}\mu_{\max}+V^{1/2}|U|\right)^2+\left(V\mu_{\max}+|U|\right)^4\right];\\
\left|\frac{1}{\lambda}\frac{\partial ^2 f}{\partial\mu\partial\omega^2}\right|&\leq \left[2\left(V^2\mu_{\max}+V|U|\right)+\left(V\mu_{\max}+|U|\right)^3\right];\\
\left|\frac{1}{\lambda}\frac{\partial^3 f}{\partial\mu^3}\right|&\leq\left[3\left(V^2\mu_{\max}+V|U|\right)+\left(V\mu_{\max}+|U|\right)^3\right];\\
\left|\frac{1}{\lambda}\frac{\partial ^3 f}{\partial\mu^2\partial\omega^2}\right|&\leq\left[2V^2+3\left(V^{3/2}\mu_{\max}+V^{1/2}|U|\right)^2+\left(V\mu_{\max}+|U|\right)^4\right];\\
\left|\frac{1}{\lambda}\frac{\partial ^3 f}{\partial\mu\partial\omega^2\partial\omega^2}\right|&\leq\left[4\left(V^3\mu_{\max}+V^2|U|\right)+4\left(V^{4/3}\mu_{\max}+V^{1/3}|U|\right)^3+\left(V\mu_{\max}+|U|\right)^5\right];\\
\left|\frac{1}{\lambda}\frac{\partial ^3 f}{\partial\omega^2\partial\omega^2\partial\omega^2}\right|&\leq\left[2V^3+11\left(V^2\mu_{\max}+V|U|\right)^2+5\left(V^{5/4}\mu_{\max}+V^{1/4}|U|\right)^4+\left(V\mu_{\max}+|U|\right)^6\right],
\end{align}
where each of the upper bounds has finite expectation. This can be seen from the existence of all order moments of 
$V$ and $|U|$ (which follows from Lemma 1 of \ctn{Maud12}), along with the compactness of the parameter space, and using the 
Cauchy-Schwartz inequality. 

\section{Non-multiplicative random effects}
\label{sec:nonmult}
So far we have considered systems of $SDE$s with multiplicative random effects given by (\ref{eq:sde1}).
However, it is possible to consider models of the following form with non-multiplicative random effects (see, for example, \ctn{Maud16}):
\begin{equation}
d X_i(t)=(\phi_ib(X_i(t))+a(X_i(t)))dt+\sigma(X_i(t))dW_i(t),\quad\mbox{with}\quad X_i(0)=x^i,~i=1,\ldots,n,
\label{eq:sde1_nonmult}
\end{equation}
Under the assumptions that $a(\cdot)$ and $b(\cdot)$ are Lipschitz continuous on $\mathbb R$, $\sigma(\cdot)$ is Holder continuous with exponent
$\alpha\in\left[\frac{1}{2},1\right]$ on $\mathbb R$, and $\int_0^{T_i}\frac{b^2(X_i(s))+a^2(X_i(s))}{\sigma^2(X_i(s))}<\infty$ almost surely for $i\geq 1$,
it follows that the form likelihood of the $i$-th individual remains the same as (\ref{eq:likelihood5}), with $V_i$ remaining the same as before and $U_i$ 
having the form $U_i=\int_0^{T_i}\frac{b(X_i(s))}{\sigma^2(X_i(s))}(dX_i(s)-a(X_i(s))ds)$.
Note that this is the density with respect to the measure associated with (\ref{eq:sde1_nonmult}) where $\phi_i=0$.

It is easy to perceive that all our asymptotic results with respect to (\ref{eq:sde1_nonmult}) remain the same as before. 

\section{Multidimensional linear random effects}
\label{sec:multidim}

Here we consider $d$-dimensional random effect, that is, we consider $SDE$s of the following form:
\begin{equation}
d X_i(t)=\bphi^T_i\bb(X_i(t))dt+\sigma(X_i(t))dW_{i}(t),\quad\mbox{with}\quad X_i(0)=x^i,~i=1,\ldots,n.
\label{eq:sde_mult}
\end{equation}
where $\bphi_i=(\phi_i^1,\phi_i^2,\ldots,\phi_i^d)^T$ is a $d$-dimensional random vector and 
$\bb(x)=(b^1(x),b^2(x),\ldots,b^d(x))^T$ is a function from $\mathbb R$ to $\mathbb R^d$. 
Here $b(x,\varphi)=\sum_{i=1}^d \varphi^i b^i(x)$ satisfies (H1). We consider $\bphi_i$ having, say, $M$ mixture of normal distributions having expectation vectors $\bmu_k$ 
and covariance matrices $\bOmega_k$ for $k=1,\ldots, M$, with density 
$$g(\varphi,\btheta)d\nu(\varphi)\equiv
\sum_{k=1}^Ma_kN(\bmu_k,\bOmega_k)$$ such that $a_k\geq 0$ for $k=1,\ldots,M$
and $\sum_{k=1}^Ma_k=1$. Here the parameter set is
$$\btheta=(a_1,\bmu_1,\bOmega_1, \ldots, a_M,\bmu_M,\bOmega_M)=(\gamma,\bbeta),$$
where, $\gamma=(a_1,\ldots,a_M)$, and $\bbeta=(\bbeta_1,\ldots,\bbeta_M)$, where, 
for $k=1,\ldots,M$, $\bbeta_k=(\bmu_k,\bOmega_k)$.

The sufficient statistics for $i=1,\ldots,n$ are
$$\bU_i=\int_0^{T_i} \frac{\bb(X_i(s))}{\sigma^2(X_i(s))}dX_i(s)$$ and
$$\bV_i=\int_0^{T_i} \frac{\bb(X_i(s))\bb^T(X_i(s))}{\sigma^2(X_i(s))}ds.$$
Note that $\bU_i$ are $d$-dimensional random vectors and $\bV_i$ are random matrices of order $d\times d$. 
As in \ctn{Maud12} we need to assume that $\bV_i$ is positive definite for each $i\geq 1$ and for all $\btheta$.

By Lemma 2 of \ctn{Maud12} it follows, almost surely, for all $i\geq 1$, for $1\leq k\leq M$ and for all $\btheta$, that
$\bV_i+\bOmega_k^{-1}, \mathbb I_d +\bV_i\bOmega_k, \mathbb I_d+\bOmega_k \bV_i$ are invertible.

Setting $\bR_i^{-1}=(\mathbb I_d+\bV_i\bOmega_k)^{-1}\bV_i$ we have 
\begin{equation}
\lambda_i(X_i,\btheta)=\sum_{i=1}^M a_k f(X_i|\bbeta_k),
\end{equation}
where
\begin{equation}
f(X_i|\bbeta_k)=\frac{1}{\sqrt{\det(\mathbb I_d+\bV_i\bOmega_k)}}
\exp\left(-\frac{1}{2}(\bmu_k-\bV_i^{-1}\bU_i)^T\bR_i^{-1}(\bmu_k-\bV_i^{-1}\bU_i)\right)\exp\left(\frac{1}{2}\bU_i^T\bV_i^{-1}\bU_i\right)
\label{eq:likelihood_mult}
\end{equation}
The asymptotic investigation for both classical and Bayesian paradigms in this multidimensional case can be carried out in 
the same way as the one-dimensional problem, using Theorem 5 of \ctn{Maitra14a} with relevant modifications 
and Proposition 10 (i) of \ctn{Maud12} which is valid here for each $(\bmu_k,\bOmega_k)$.

\section{Simulation studies}
\label{sec:simstudy_details}

For our convenience, we denote the $SDE$ models by the following: for $i\geq 1$ and $t\in[0,T]$,
\begin{align}
& SDE_1:~ dX_i(t) = \phi_iX_i(t)dt+\sigma dW_i(t),~X_i(0)=x;\label{eq:SDE_1}\\
& SDE_2:~ dX_i(t) = (\phi_i-X_i(t))dt+\sigma dW_i(t),~X_i(0)=x;\label{eq:SDE_2}\\
& SDE_3:~ dX_i(t) = (\phi_iX_i(t)+2\sigma^2)dt+2\sigma\sqrt{X_i(t)} dW_i(t),~X_i(0)=x.\label{eq:SDE_3}
\end{align}
As in \ctn{Maud16}, we choose $\sigma=0.1$, $x=1$, $T=1$. Note that although $SDE_1$ is of the form (\ref{eq:sde1}) with multiplicative random effect, $SDE_2$ abd $SDE_3$ are
of the form (\ref{eq:sde1_nonmult}). Unless otherwise stated, $n=100$. For sampling the discrete paths, we divide $[0,T]$ into $T/\delta$ intervals, each of length
$\delta$. We set $\delta=0.0002$.

For the distribution of $\phi$, we choose the following 5 normal mixture distributions, as in \ctn{Maud16}:
\begin{align}
& \pi_1:~0.5N\left(-0.5,0.25^2\right)+0.5N\left(-1.8,0.25^2\right);\label{eq:pi_1}\\ 
& \pi_2:~0.7N\left(-0.5,0.5^2\right)+0.3N\left(-1.8,0.5^2\right);\label{eq:pi_2}\\ 
& \pi_3:~0.2N\left(-0.5,0.25^2\right)+0.3N\left(-3.5,0.25^2\right)+0.5N\left(-5.5,0.25^2\right);\label{eq:pi_3}\\ 
& \pi_4:~0.2N\left(-0.5,0.25^2\right)+0.3N\left(-1.8,0.25^2\right)+0.5N\left(-2.5,0.25^2\right);\label{eq:pi_4}\\ 
& \pi_5:~0.2N\left(-0.5,0.5^2\right)+0.3N\left(-1.8,0.5^2\right)+0.5N\left(-2.5,0.5^2\right).\label{eq:pi_5} 
\end{align}

\subsection{Prior structure for our postulated mixture model}
\label{subsec:prior}
In all the simulation studies, we choose the uniform prior for $M$ on $\{1,2,\ldots,M_{\max}\}$, with $M_{\max}=30$.
We reparameterize $\omega^2_k$ as $\omega^2_k=\exp(-\tau_k)$, where $\tau_k\sim\log\left(\mathcal G\left(s/2,S/2\right)\right)$. 
Here, for positive $a$ and $b$, $\mathcal G(a,b)$ denotes the gamma distribution with mean $a/b$ and variance $a/b^2$.
Given $\omega^2_k$,
we assume that $[\mu_k|\omega^2_k]\sim N\left(\mu_0,\psi\omega^2_k\right)$. For the prior on $a_k$'s, we write $a_k=\frac{\exp(w_k)}{\sum_{j=1}^M\exp(w_j)}$,
where $w_j\stackrel{iid}{\sim}N\left(\mu_w,\sigma^2_w\right)$. All these choices are motivated by \ctn{Das17}. Further, we set the standard choices 
$s=S=2\times 0.1$ and $\mu_w=0$
in all our simulation experiments. We choose $\sigma^2_w$ to be $0.5$ and $0.35$ for the different simulation studies. These somewhat small variances reflect our
belief that the $a_k$'s are approximately the same {\it a priori}. 
The values of $\psi$ and $\sigma^2_w$ are particularly chosen to ensure good mixing properties of TTMCMC and
little uncertainty about the posterior of $M$. As it turned out, when $\pi_1$ and $\pi_2$ are the true models, $\psi=10$ and $\sigma^2_w=0.5$ achieved
this purpose, while for $\pi_3$, $\pi_4$ and $\pi_5$, $\psi=100$ and $\sigma^2_w=0.35$ were the most appropriate in terms of ensuring good mixing and small uncertainty
about $M$ {\it a posteriori}. 

\subsection{Implementation details}
\label{subsec:implementation}
For TTMCMC, we implement Algorithm S-3.1 of
the supplement of \ctn{Das17}, with additive transformation driven by the standard normal density on the positive part of the real line. 
We set equal probabilities for birth (emergence of a new mixture component), death (deletion of an existing mixture component) 
and no-change (no change to the number of mixture components) moves, and equal probabilities for the forward (addition of a scaled, positive standard normal variable 
to the current value
of any given parameter co-ordinate) and backward (subtraction of a scaled, positive standard normal variable from the current value of any given parameter co-ordinate) 
transformations. 
We choose all the scales of the additive transformation to be $0.5$. 
As we shall demonstrate, these choices led to excellent mixing properties of our TTMCMC algorithm.
For all the experiments, we discarded the first $15\times 10^5$ TTMCMC iterations as burn-in and stored one in 150 out of a further
$15\times 10^5$ iterations to obtain $10000$ realizations from the variable-dimensional posteriors.
Our codes, all written in $C$, are implemented on a 64-bit 4-core laptop with 7.6 GB memory, 2.3 GHz processor speed. 
When the true models are two-component mixtures, the entire TTMCMC implementation takes 2-3 minutes, while for the three component mixtures, our codes take
3-4 minutes. When $n=200$, as we consider in one of our simulation based illustrations, then it takes about 5 minutes for the two-component mixtures.

\subsection{First simulation study: $SDE_1$ with $\pi_1$} 
\label{subsec:sim1}
\begin{figure}
\centering
\subfigure[Trace plot of $M$.]{ \label{fig:sim1_trace_comp}
\includegraphics[width=7cm,height=5cm]{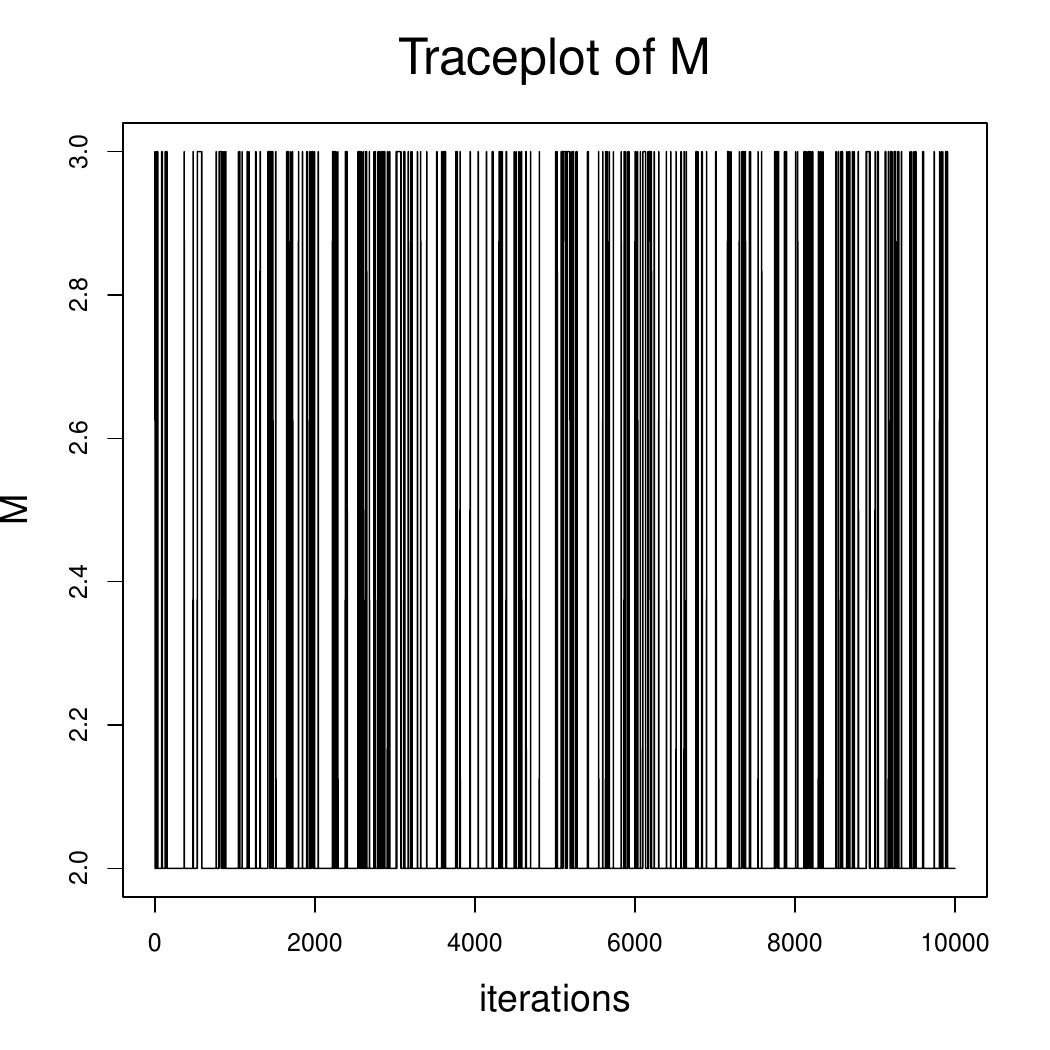}}
\hspace{2mm}
\subfigure[Trace plot of $\mu_1$.]{ \label{fig:sim1_trace_mu1}
\includegraphics[width=7cm,height=5cm]{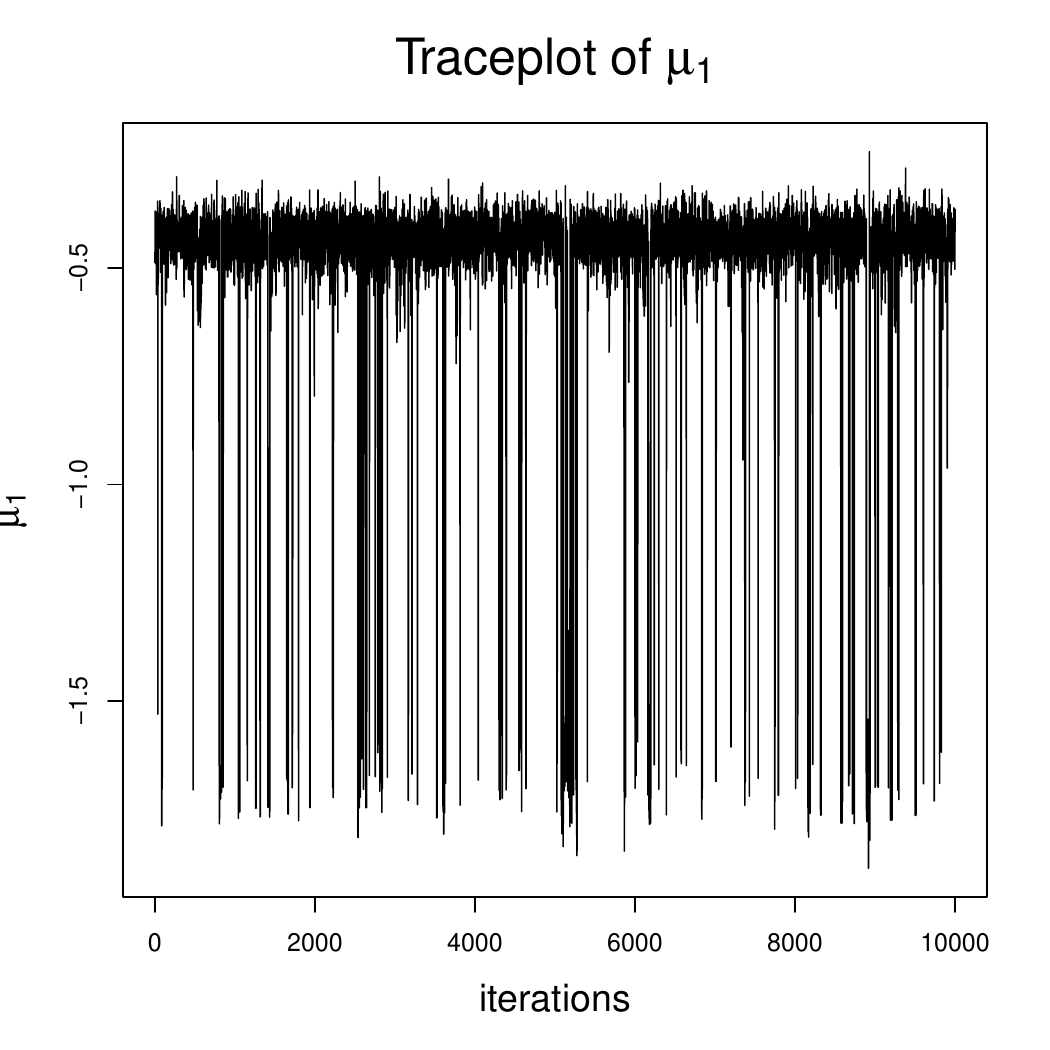}}\\
\vspace{2mm}
\subfigure[Trace plot of $\mu_2$.]{ \label{fig:sim1_trace_mu2}
\includegraphics[width=7cm,height=5cm]{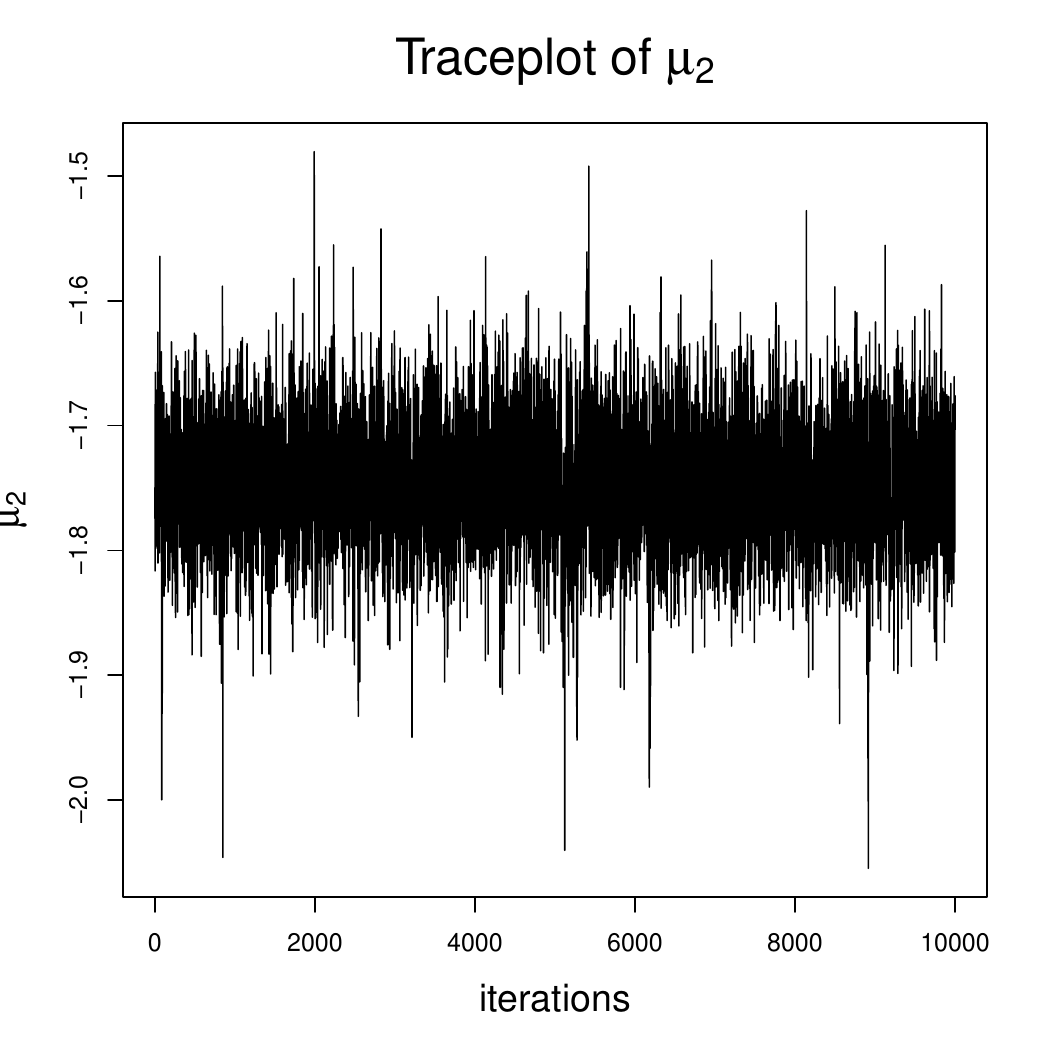}}
\vspace{2mm}
\subfigure[Trace plot of $\omega^2_1$.]{ \label{fig:sim1_trace_omegasq1}
\includegraphics[width=7cm,height=5cm]{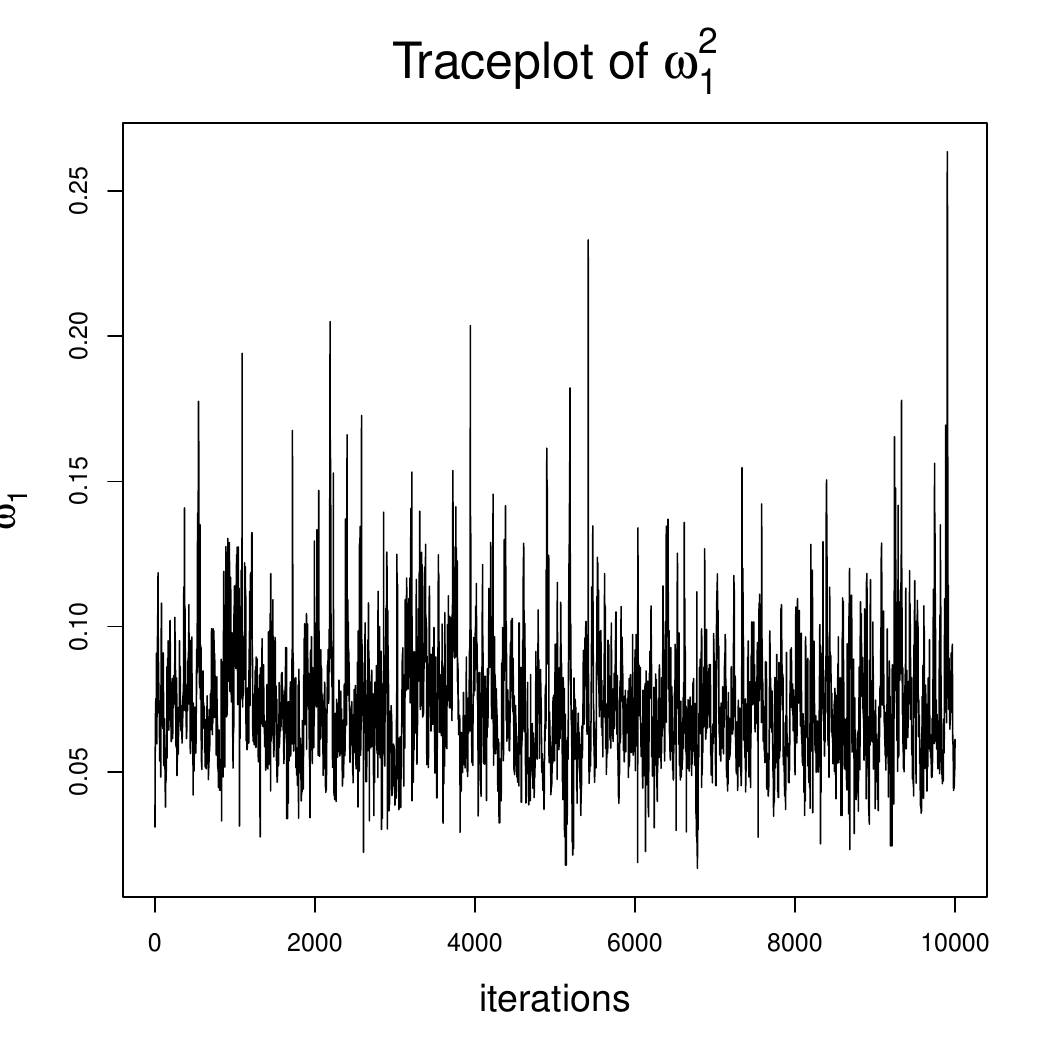}}\\
\vspace{2mm}
\subfigure[Trace plot of $\omega^2_2$.]{ \label{fig:sim1_trace_omegasq2}
\includegraphics[width=7cm,height=5cm]{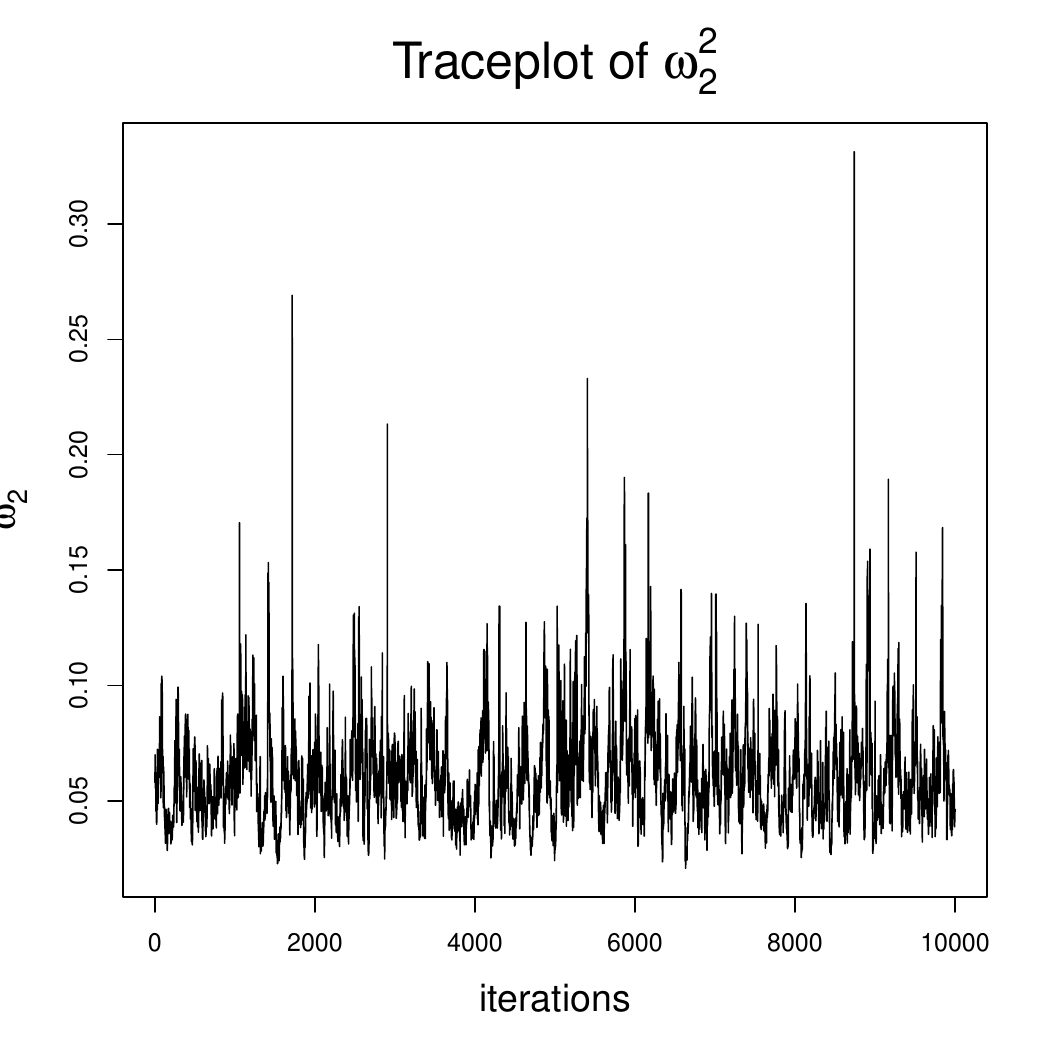}}
\hspace{2mm}
\subfigure[Trace plot of $a_1$.]{ \label{fig:sim1_trace_p1}
\includegraphics[width=7cm,height=5cm]{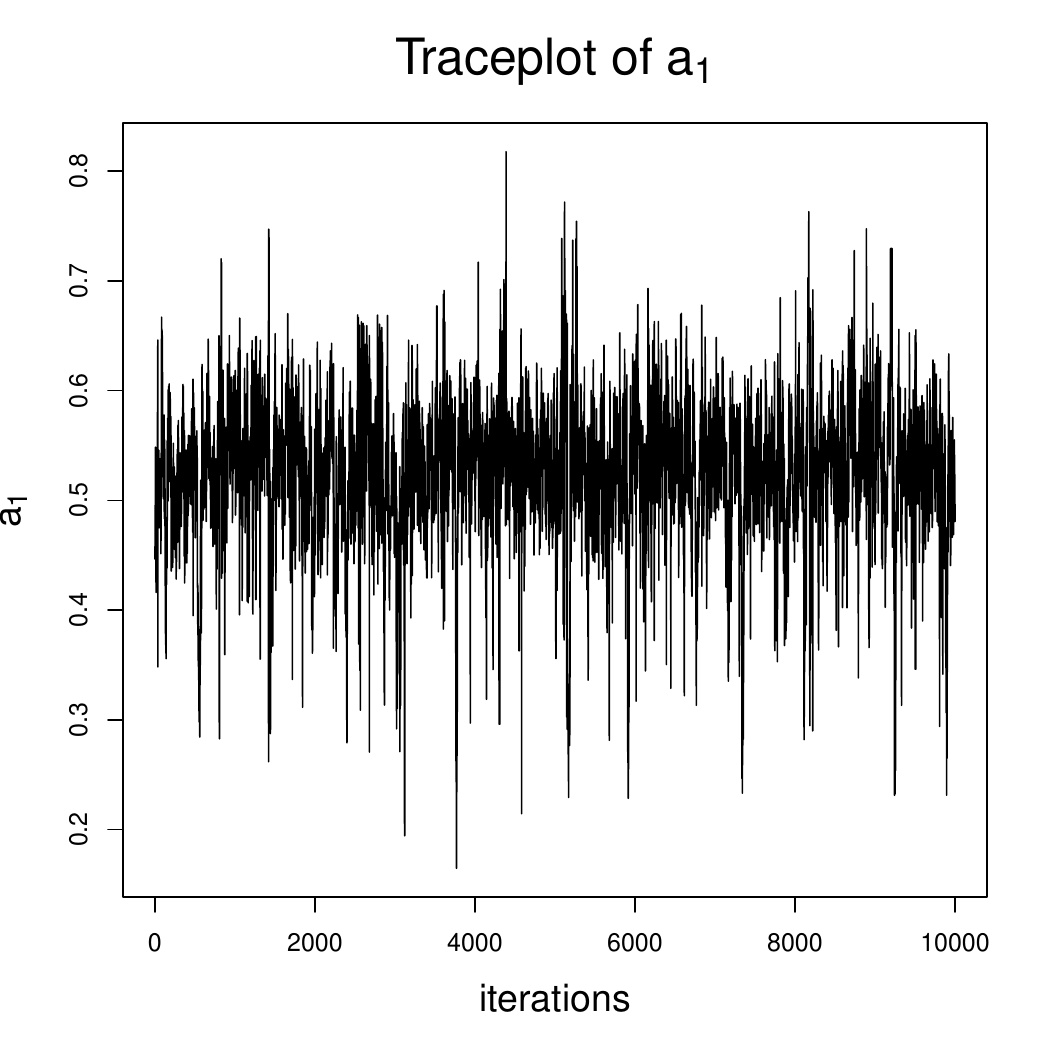}}\\
\vspace{2mm}
\subfigure[Trace plot of $a_2$.]{ \label{fig:sim1_trace_p2}
\includegraphics[width=7cm,height=5cm]{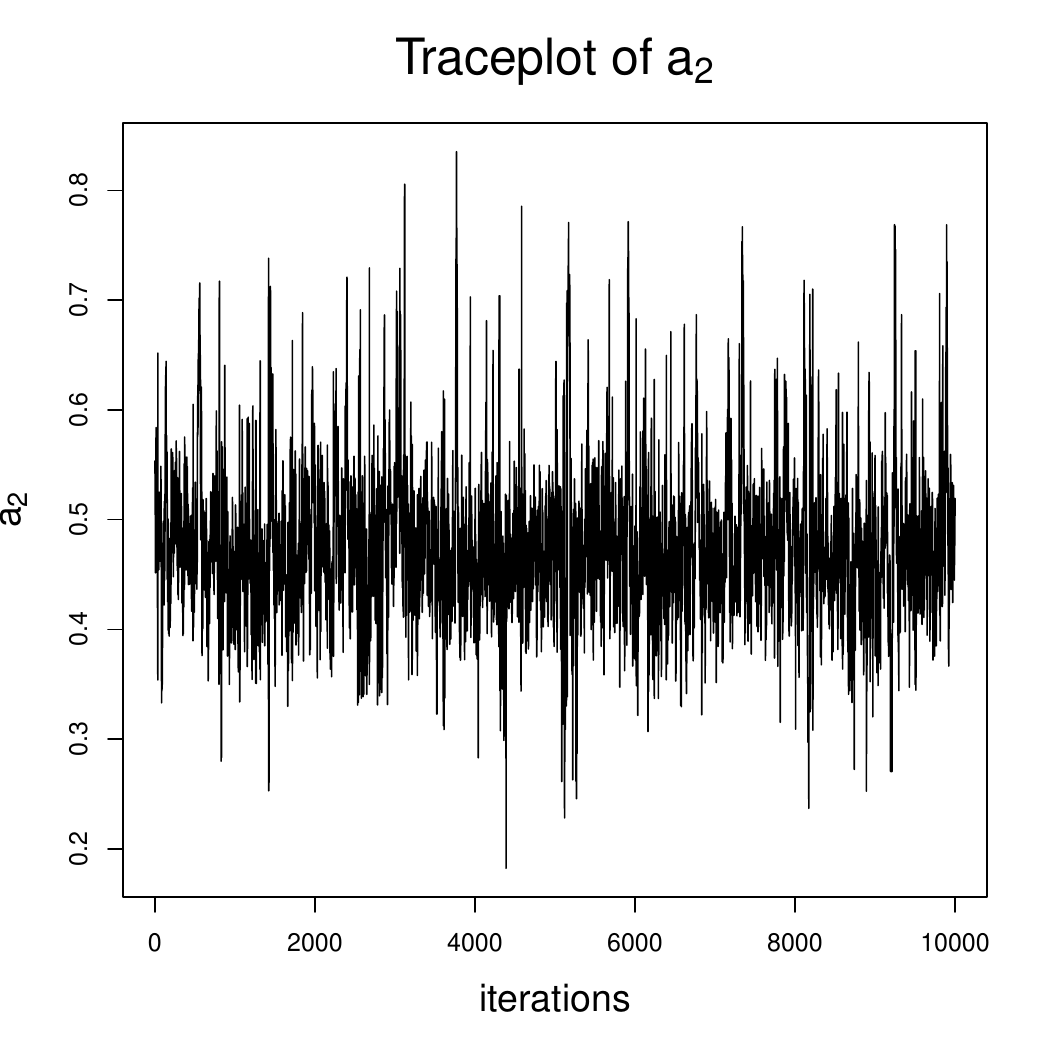}}
\caption{{\bf TTMCMC for $SDE_1$ and $\pi_1$:} Trace plots of $M$, $\mu_1$, $\mu_2$, $\omega^2_1$, $\omega^2_2$, $a_1$ and $a_2$.} 
\label{fig:sim1_trace_plots}
\end{figure}

\begin{figure}
\centering
\subfigure[Posterior of $\mu_1$.]{ \label{fig:sim1_mu1}
\includegraphics[width=7cm,height=6cm]{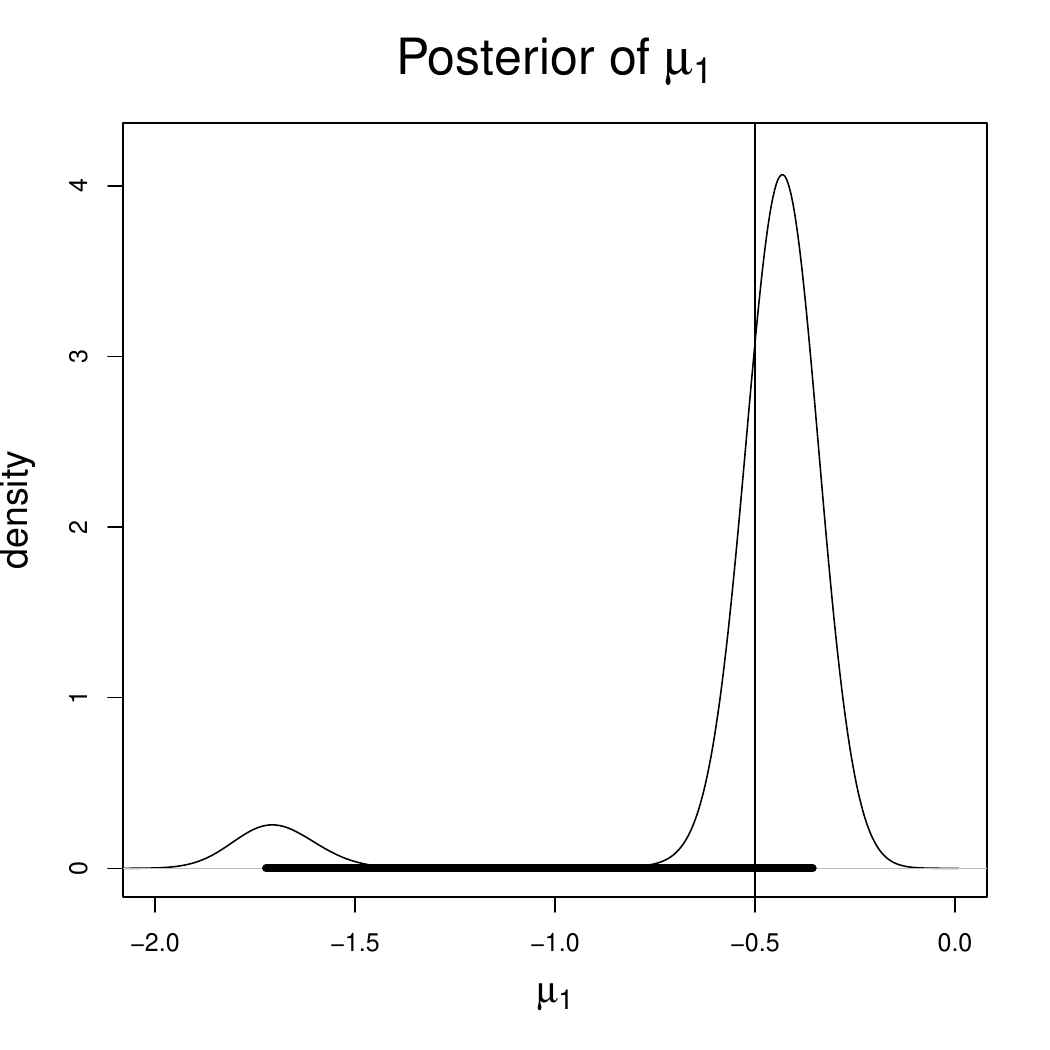}}
\hspace{2mm}
\subfigure[Posterior of $\mu_2$.]{ \label{fig:sim1_mu2}
\includegraphics[width=7cm,height=6cm]{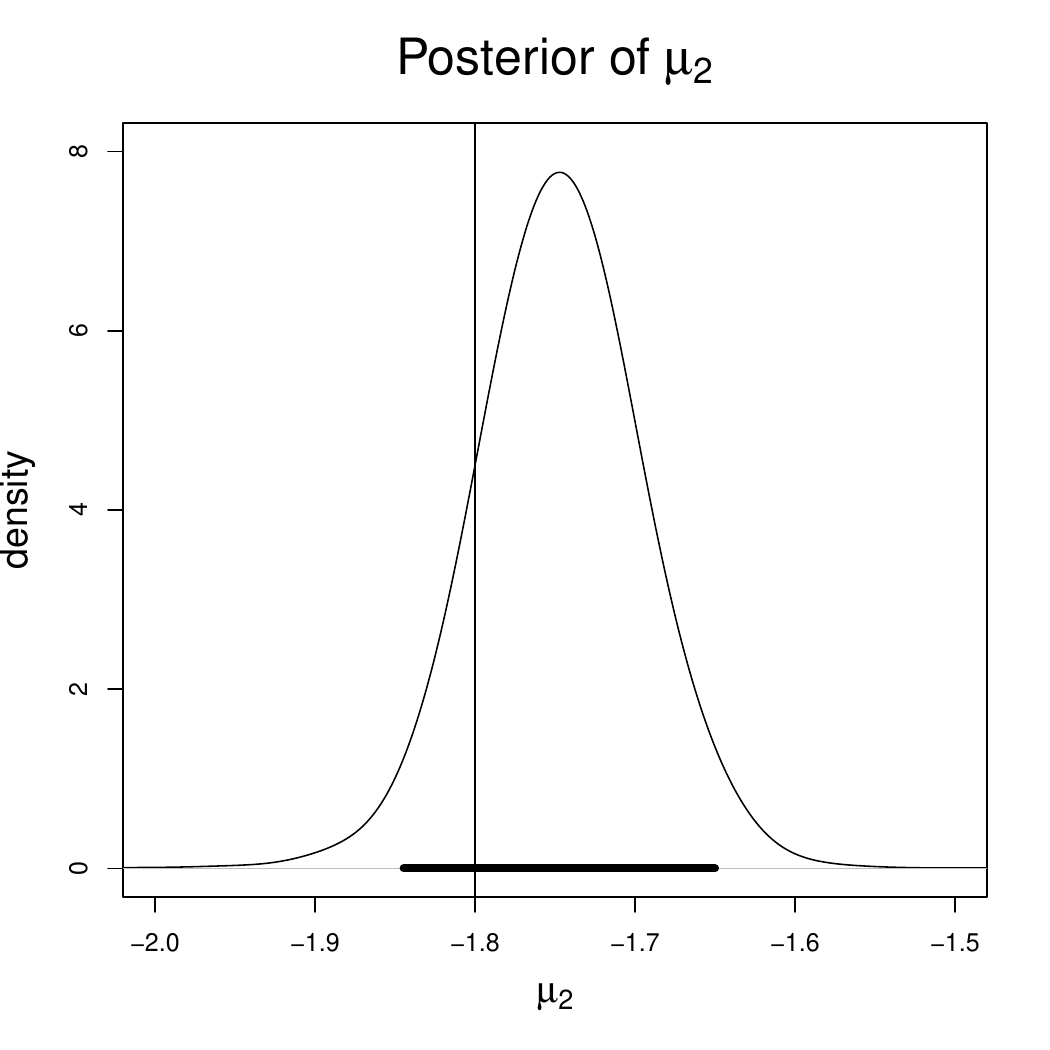}}\\
\vspace{2mm}
\subfigure[Posterior of $\omega^2_1$.]{ \label{fig:sim1_omegasq1}
\includegraphics[width=7cm,height=6cm]{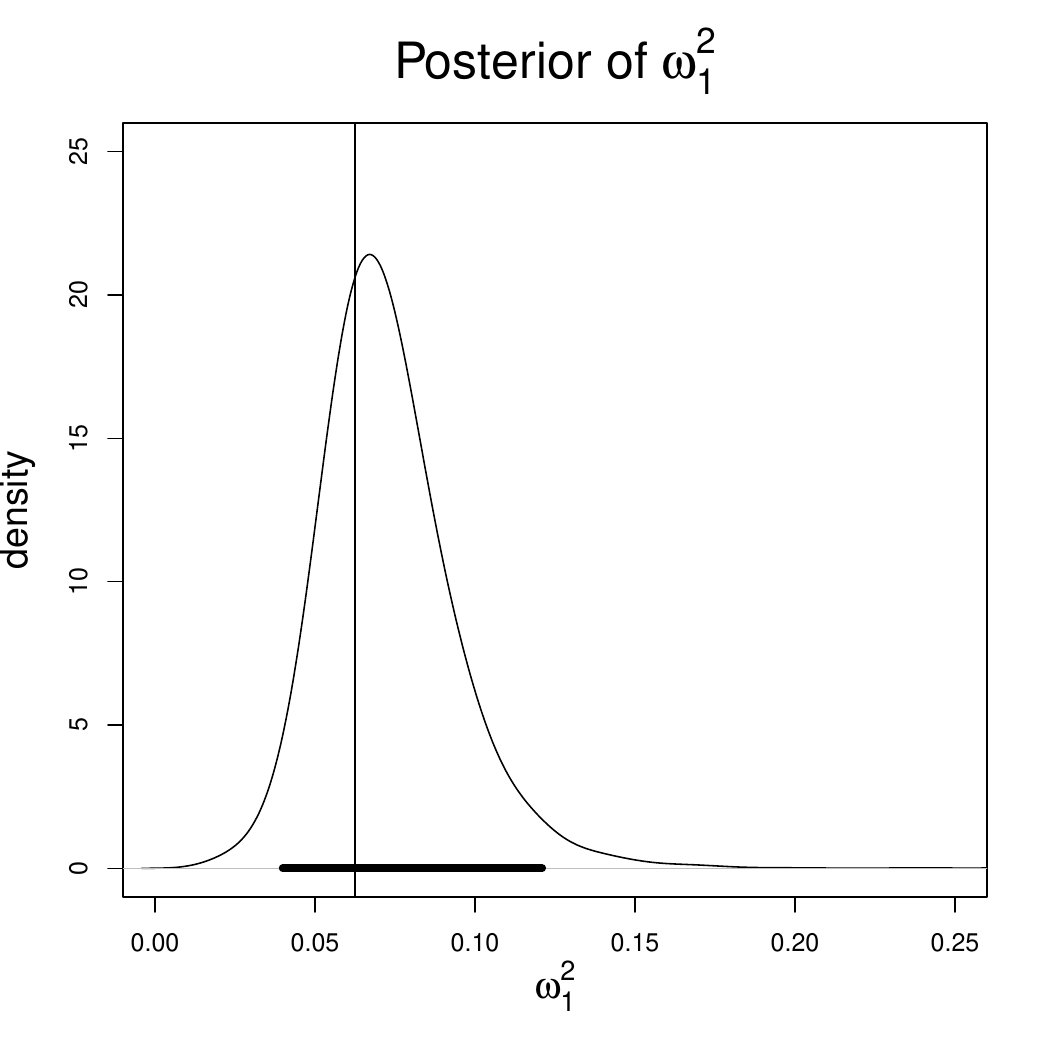}}
\hspace{2mm}
\subfigure[Posterior of $\omega^2_2$.]{ \label{fig:sim1_omegasq2}
\includegraphics[width=7cm,height=6cm]{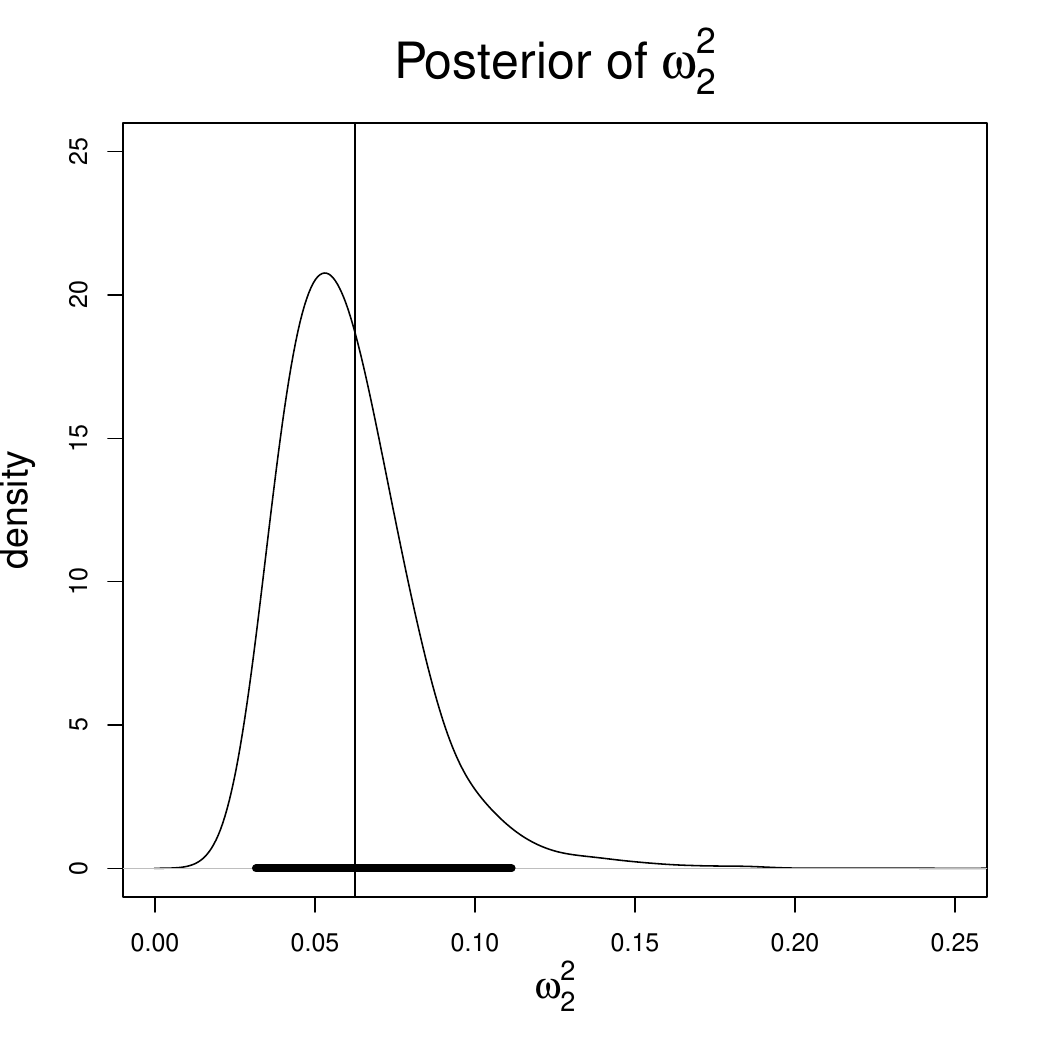}}\\
\vspace{2mm}
\subfigure[Posterior of $a_1$.]{ \label{fig:sim1_p1}
\includegraphics[width=7cm,height=6cm]{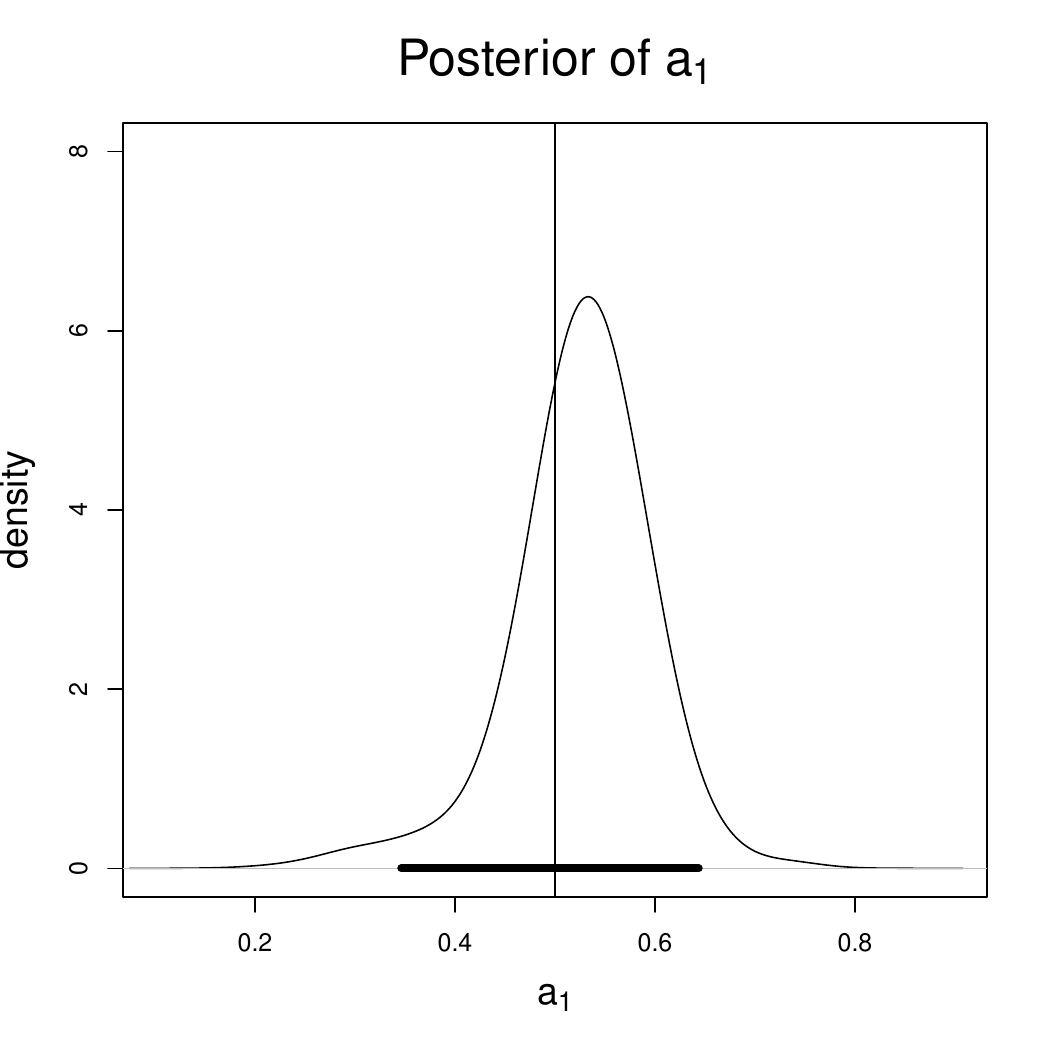}}
\hspace{2mm}
\subfigure[Posterior of $a_2$.]{ \label{fig:sim1_p2}
\includegraphics[width=7cm,height=6cm]{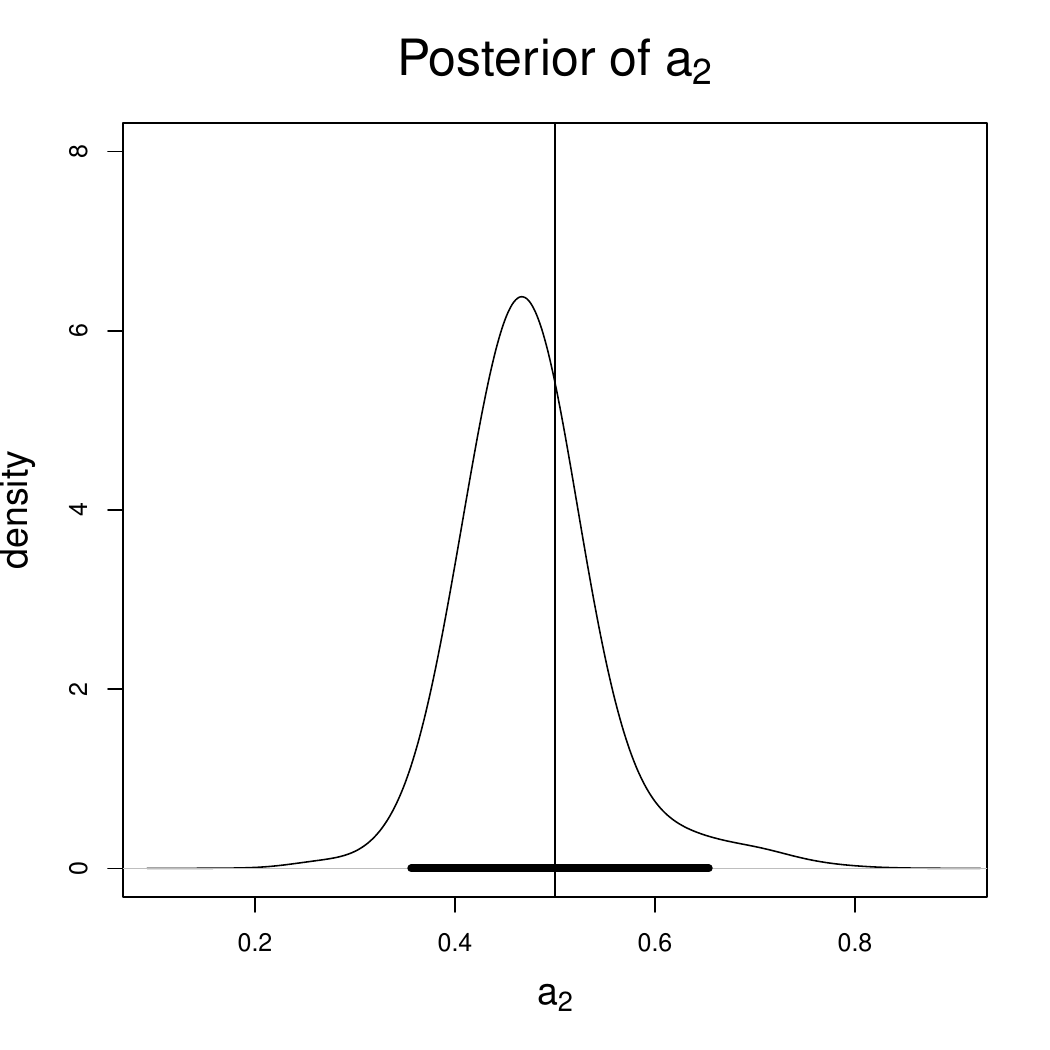}}
\caption{{\bf TTMCMC for $SDE_1$ and $\pi_1$:} Posteriors of $M$, $\mu_1$, $\mu_2$, $\omega^2_1$, $\omega^2_2$, $a_1$ and $a_2$. The vertical lines stand
for the true values, while the thick horizontal lines denote the 95\% credible intervals.} 
\label{fig:sim1_posterior_plots}
\end{figure}

Here we simulate data in accordance with $\pi_1$ and $SDE_1$ and fit our mixture model with a maximum of 30 components. The trace plots and the posterior
density estimates of the parameters are provided in 
Figures \ref{fig:sim1_trace_plots} and \ref{fig:sim1_posterior_plots}, respectively. 

The trace plots indicate excellent mixing of all the unknowns
and the posterior plots in Figure \ref{fig:sim1_posterior_plots} show that all the true values, represented by the vertical lines fall well within the 
respective 95\% credible regions, denoted by the thick horizontal line.
The posterior probabilities of $M=2$ and $M=3$ turned out to be 
$0.8772$ and $0.1228$, respectively, and all other values of $M$ have zero posterior mass. 
Thus, the true value, $M_0=2$ receives significantly larger posterior probability compared to all the other values.

In other words, our Bayesian method puts up excellent performance in this experiment.
Since the mixture components of $\pi_1$ are well-separated, this is to be expected.

The overall acceptance rate of our algorithm in this case is about $2\%$.

\subsection{Second simulation study: $SDE_1$ with $\pi_2$} 
\label{subsec:sim2}
We now simulate the data in accordance with $\pi_2$ and $SDE_1$ and again fit our mixture model with a maximum of 30 components. Since the mixture components
of $\pi_2$ are not as well-separated as in $\pi_1$, this setup is somewhat more challenging than the previous one.

\begin{figure}
\centering
\subfigure[Trace plot of $M$.]{ \label{fig:sim2_trace_comp}
\includegraphics[width=7cm,height=5cm]{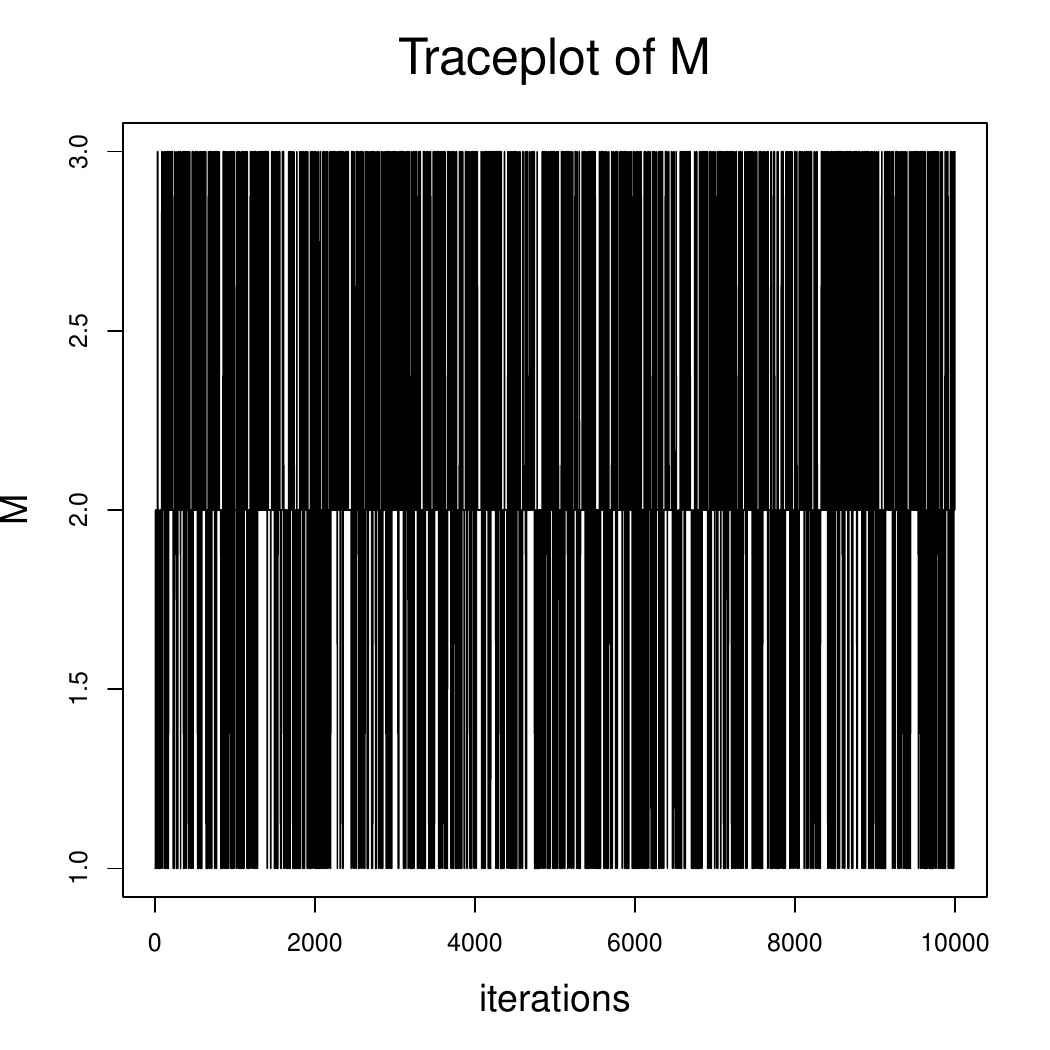}}
\hspace{2mm}
\subfigure[Trace plot of $\mu_1$.]{ \label{fig:sim2_trace_mu1}
\includegraphics[width=7cm,height=5cm]{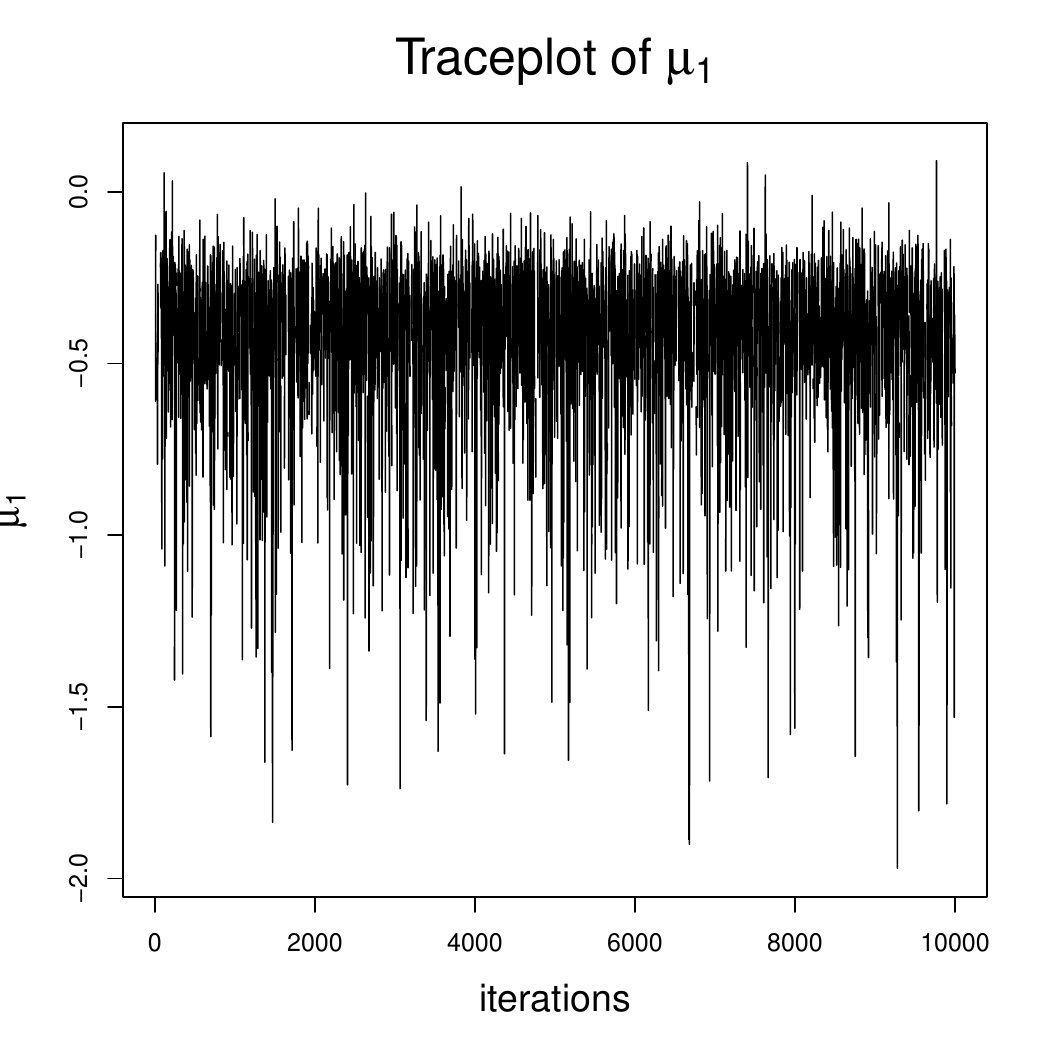}}\\
\vspace{2mm}
\subfigure[Trace plot of $\mu_2$.]{ \label{fig:sim2_trace_mu2}
\includegraphics[width=7cm,height=5cm]{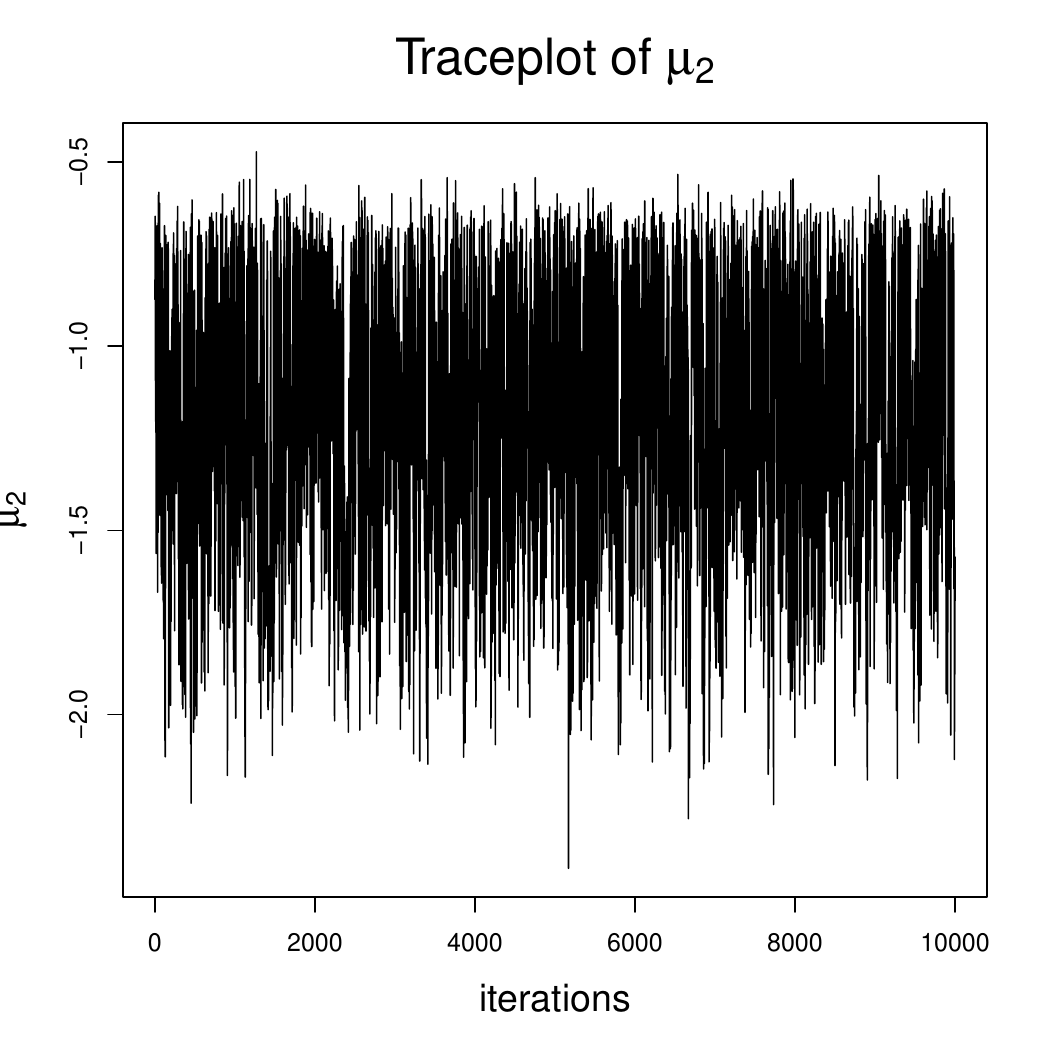}}
\vspace{2mm}
\subfigure[Trace plot of $\omega^2_1$.]{ \label{fig:sim2_trace_omegasq1}
\includegraphics[width=7cm,height=5cm]{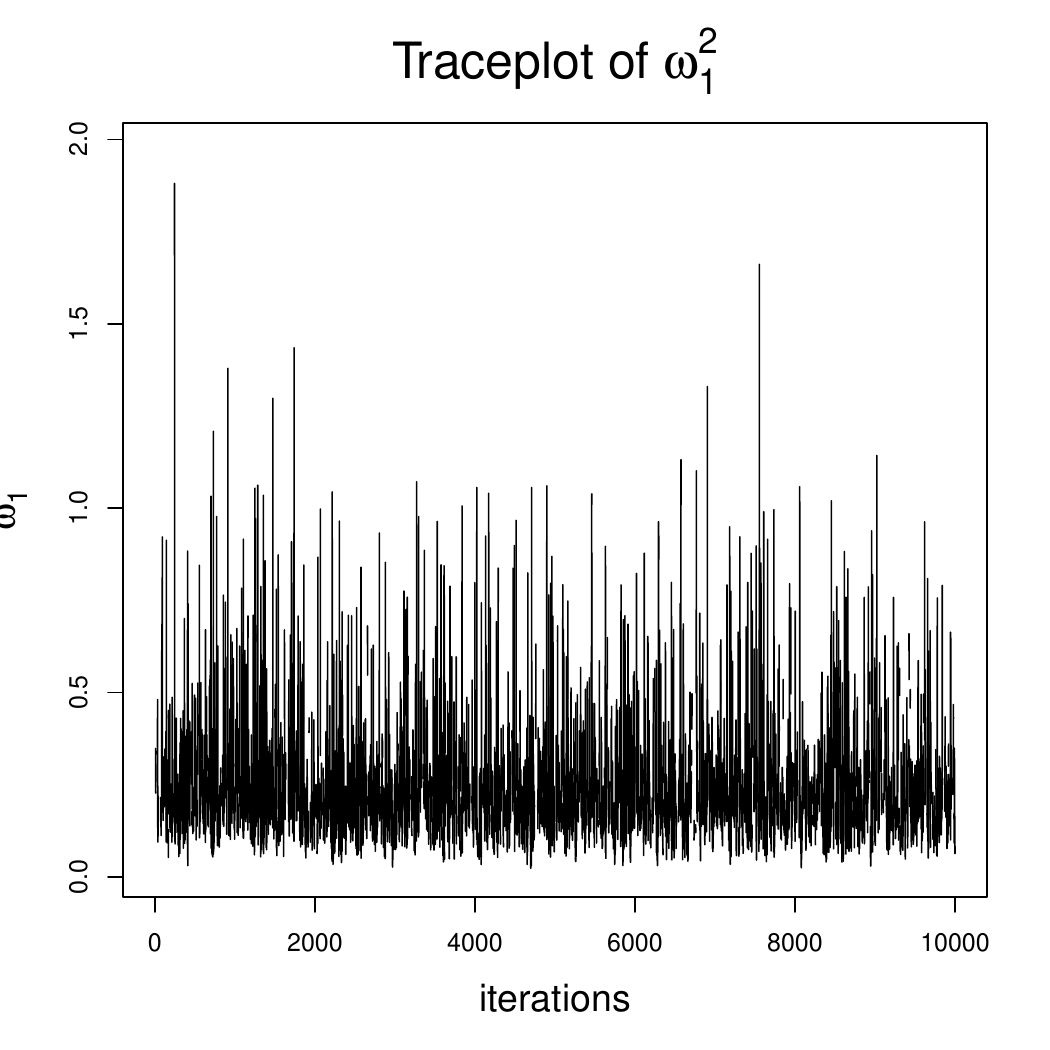}}\\
\vspace{2mm}
\subfigure[Trace plot of $\omega^2_2$.]{ \label{fig:sim2_trace_omegasq2}
\includegraphics[width=7cm,height=5cm]{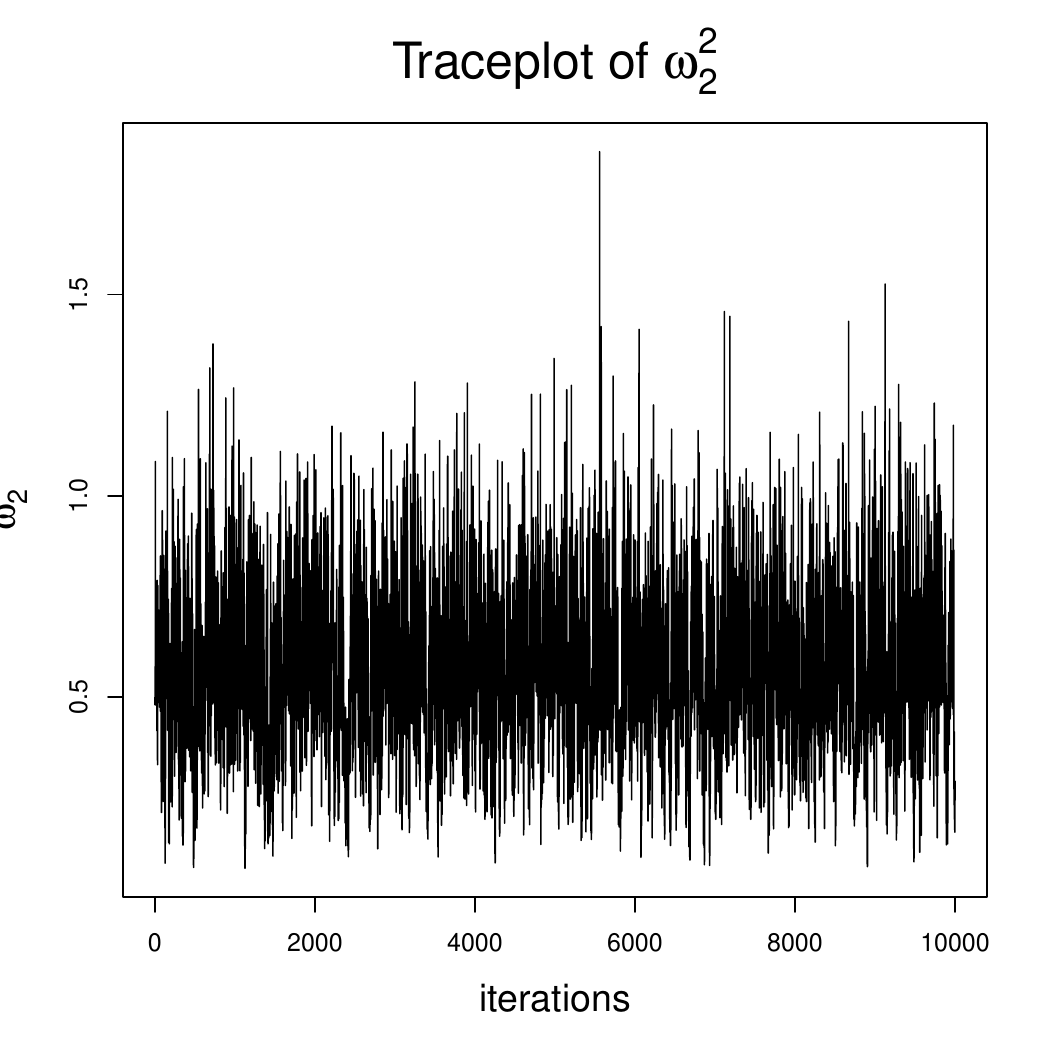}}
\hspace{2mm}
\subfigure[Trace plot of $a_1$.]{ \label{fig:sim2_trace_p1}
\includegraphics[width=7cm,height=5cm]{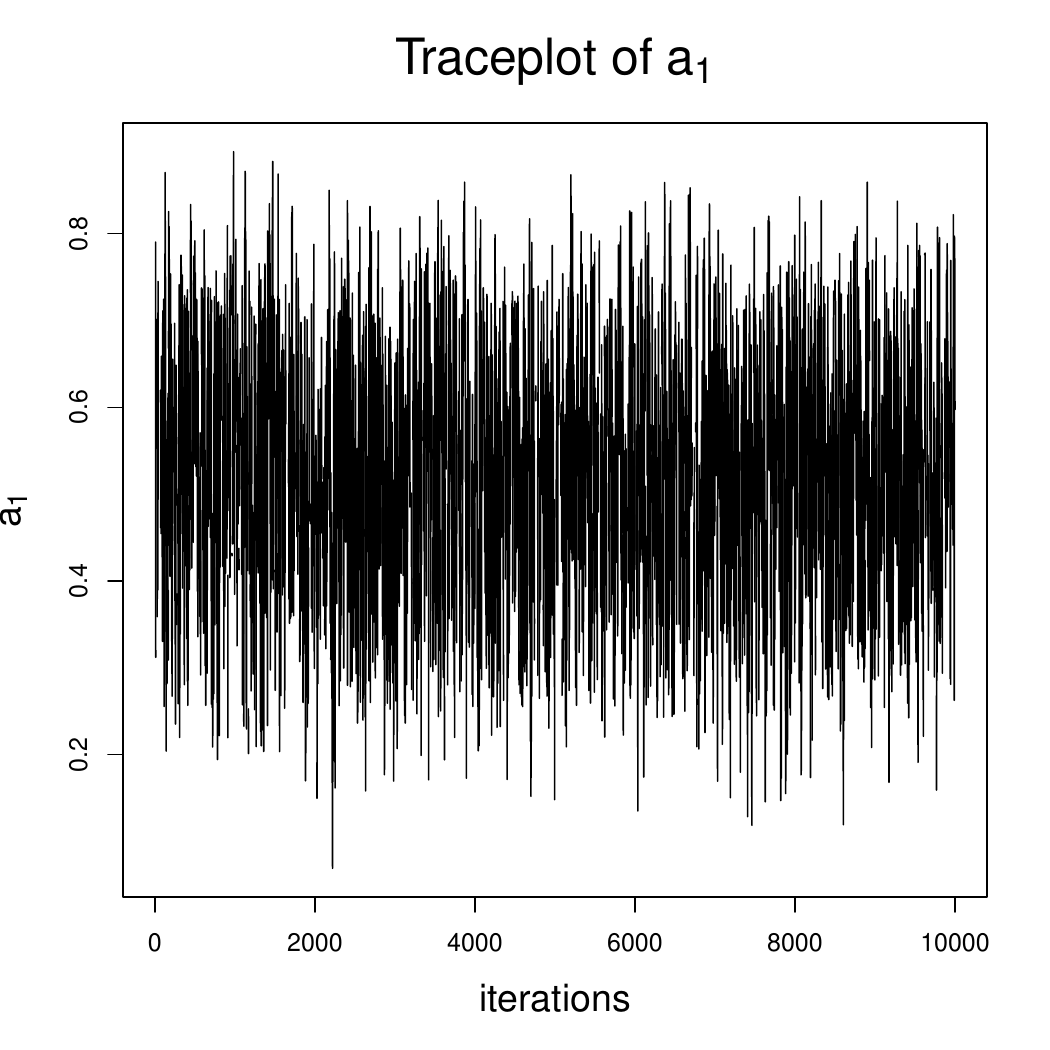}}\\
\vspace{2mm}
\subfigure[Trace plot of $a_2$.]{ \label{fig:sim2_trace_p2}
\includegraphics[width=7cm,height=5cm]{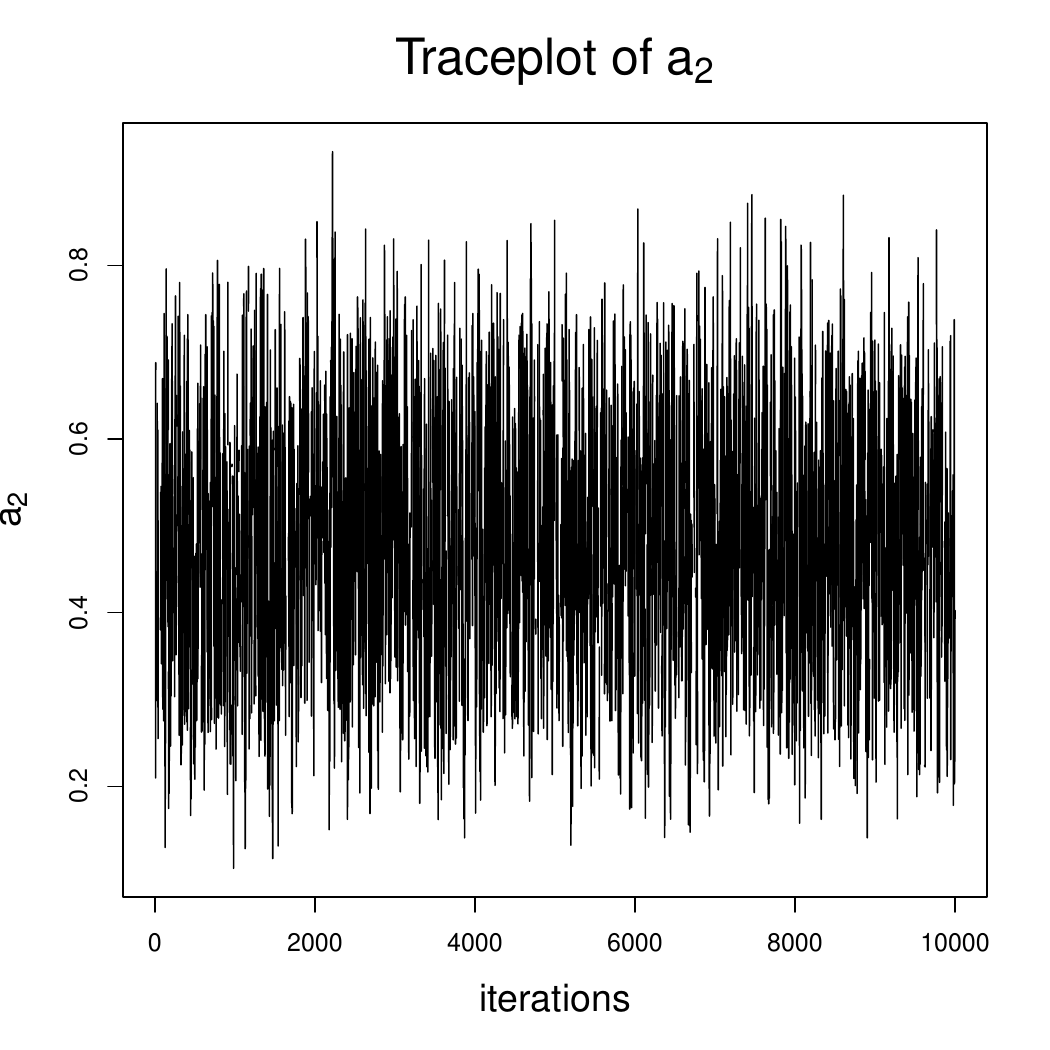}}
\caption{{\bf TTMCMC for $SDE_1$ and $\pi_2$:} Trace plots of $M$, $\mu_1$, $\mu_2$, $\omega^2_1$, $\omega^2_2$, $a_1$ and $a_2$.} 
\label{fig:sim2_trace_plots}
\end{figure}

\begin{figure}
\centering
\subfigure[Posterior of $\mu_1$.]{ \label{fig:sim2_mu1}
\includegraphics[width=7cm,height=6cm]{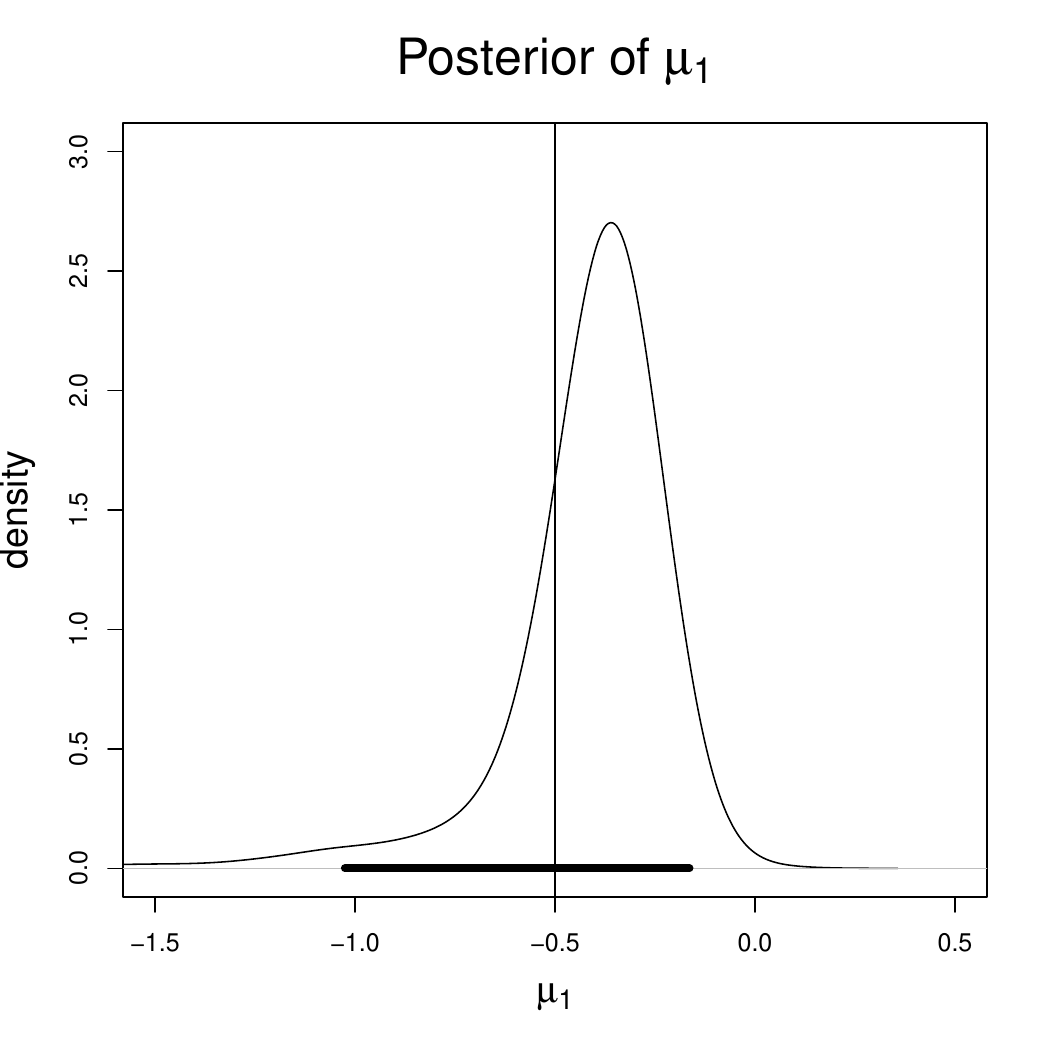}}
\hspace{2mm}
\subfigure[Posterior of $\mu_2$.]{ \label{fig:sim2_mu2}
\includegraphics[width=7cm,height=6cm]{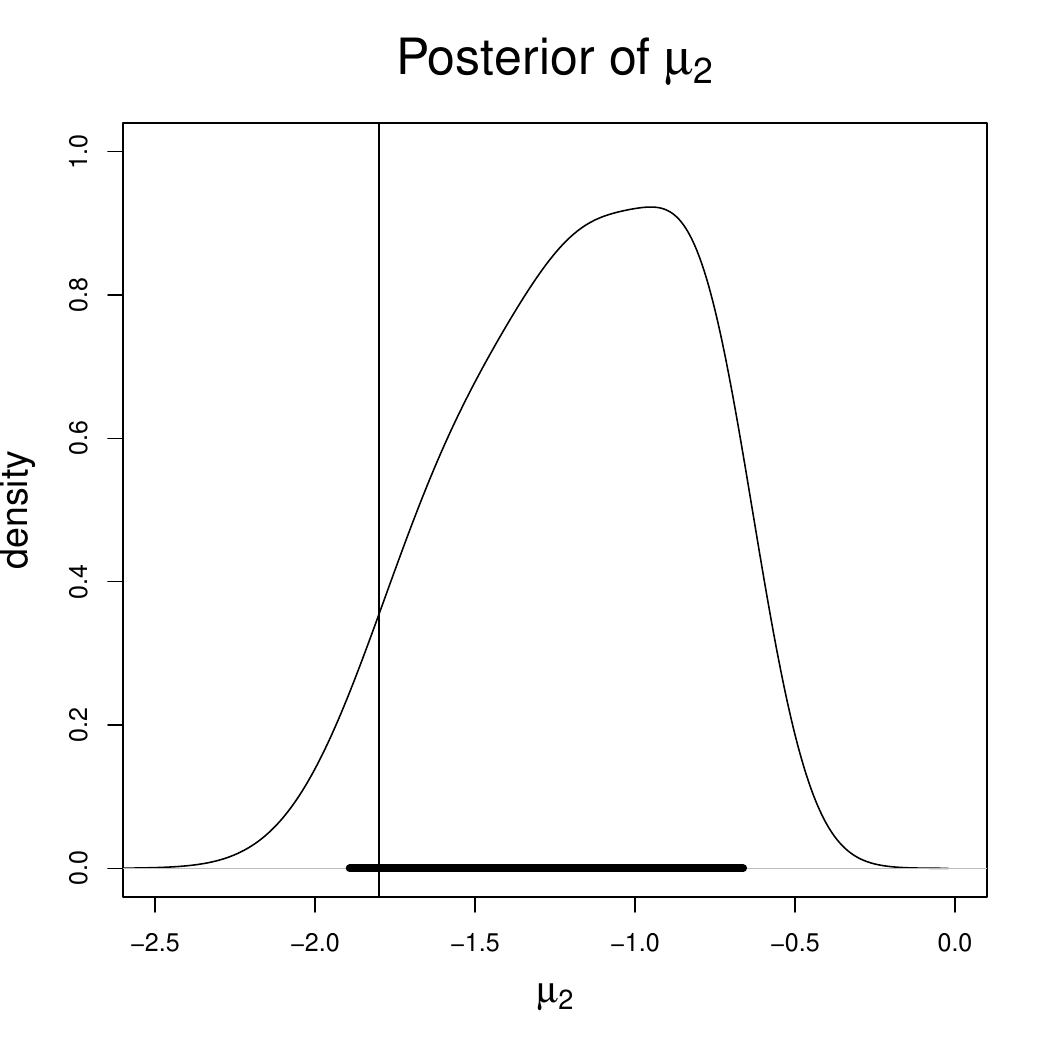}}\\
\vspace{2mm}
\subfigure[Posterior of $\omega^2_1$.]{ \label{fig:sim2_omegasq1}
\includegraphics[width=7cm,height=6cm]{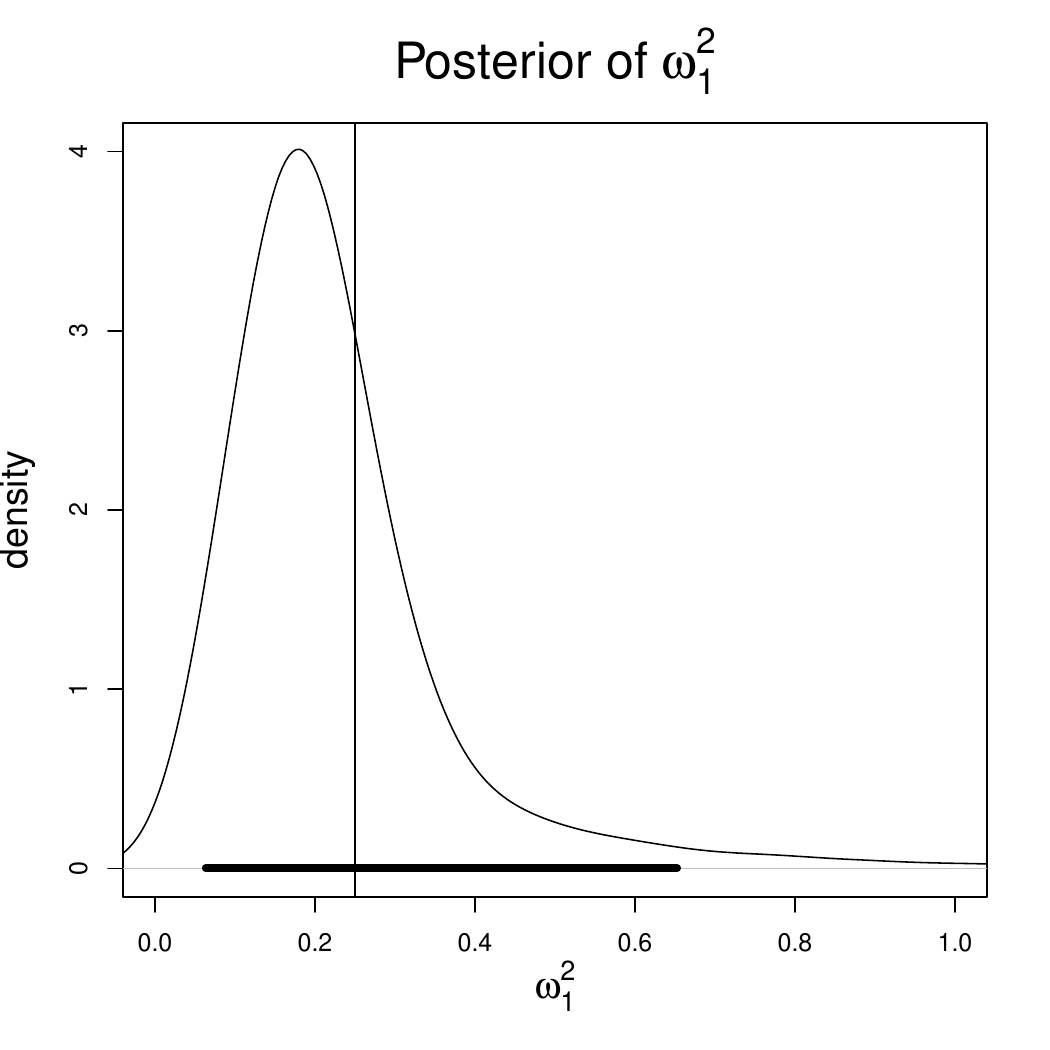}}
\hspace{2mm}
\subfigure[Posterior of $\omega^2_2$.]{ \label{fig:sim2_omegasq2}
\includegraphics[width=7cm,height=6cm]{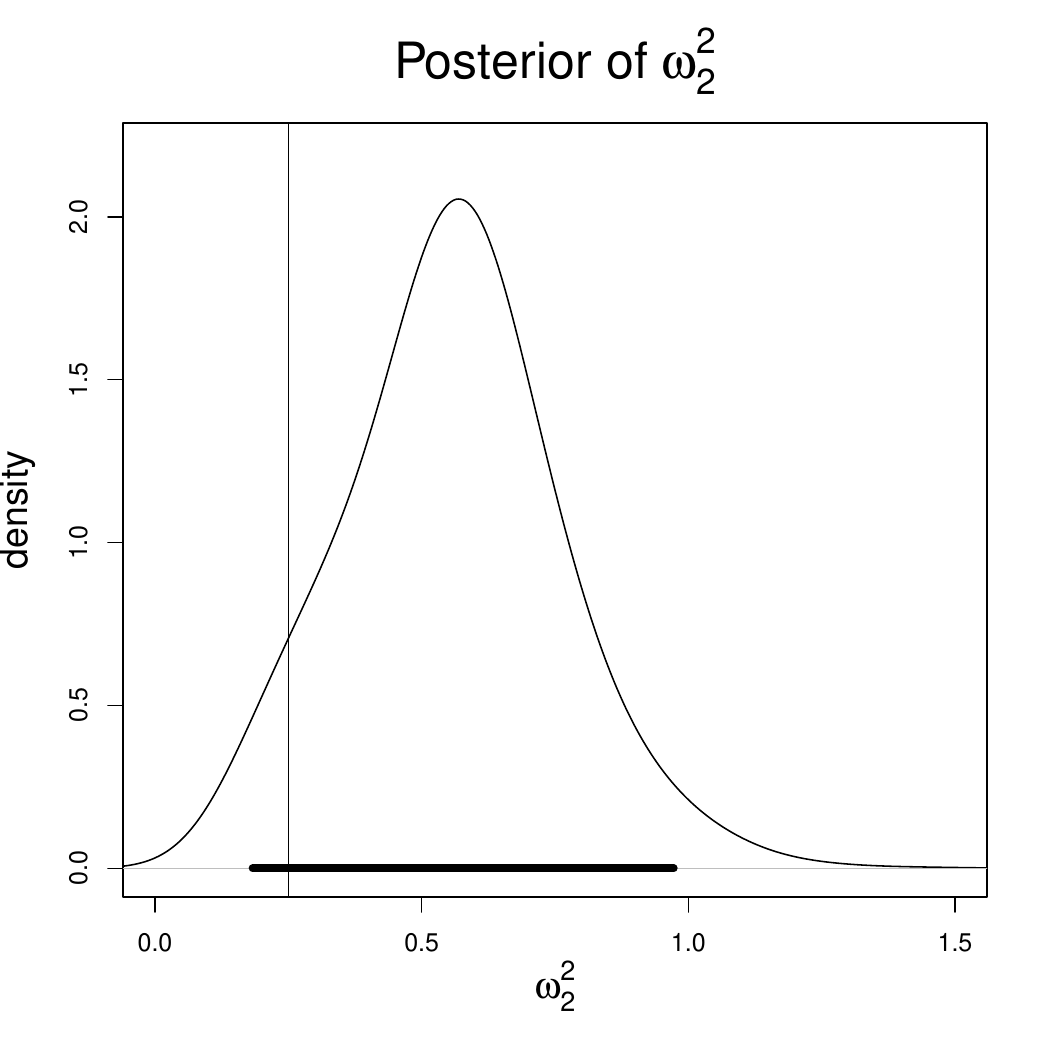}}\\
\vspace{2mm}
\subfigure[Posterior of $a_1$.]{ \label{fig:sim2_p1}
\includegraphics[width=7cm,height=6cm]{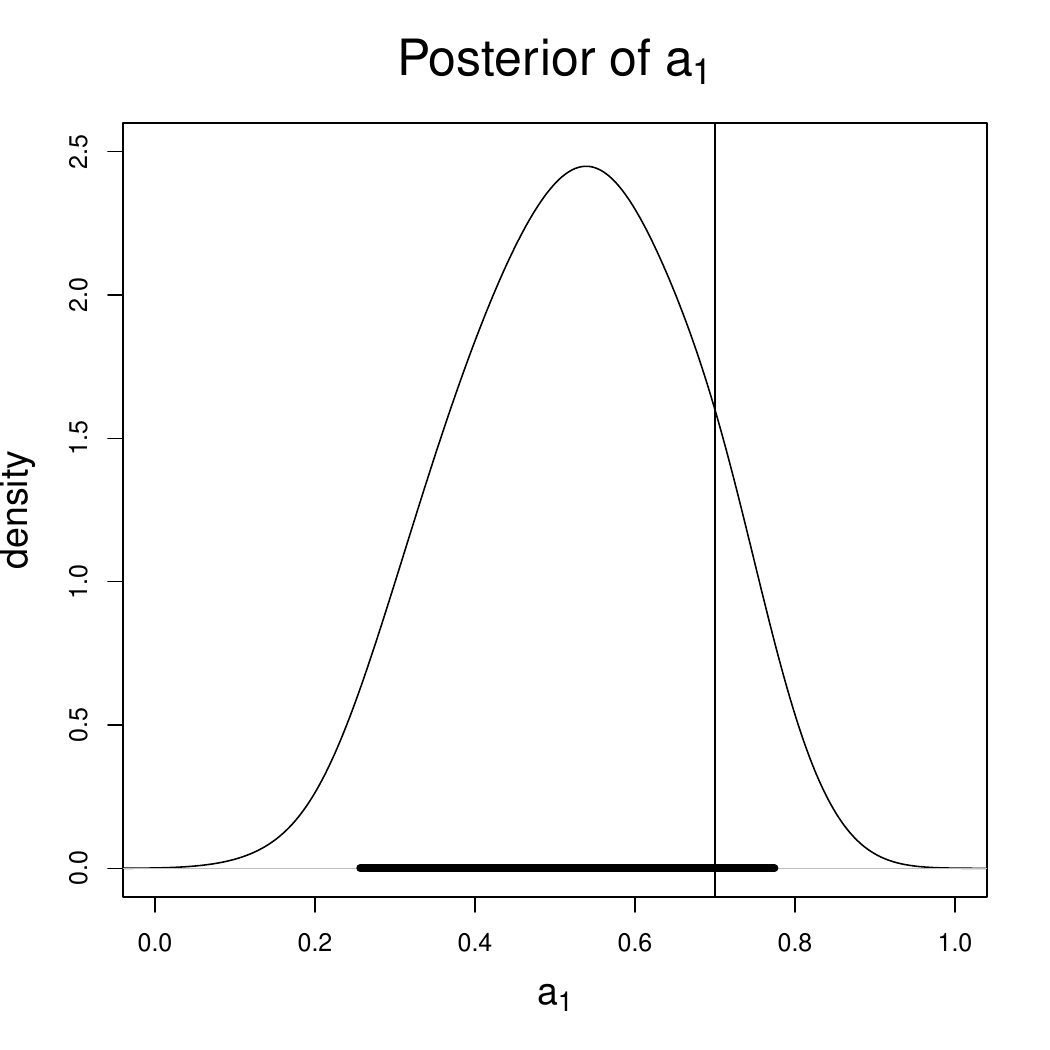}}
\hspace{2mm}
\subfigure[Posterior of $a_2$.]{ \label{fig:sim2_p2}
\includegraphics[width=7cm,height=6cm]{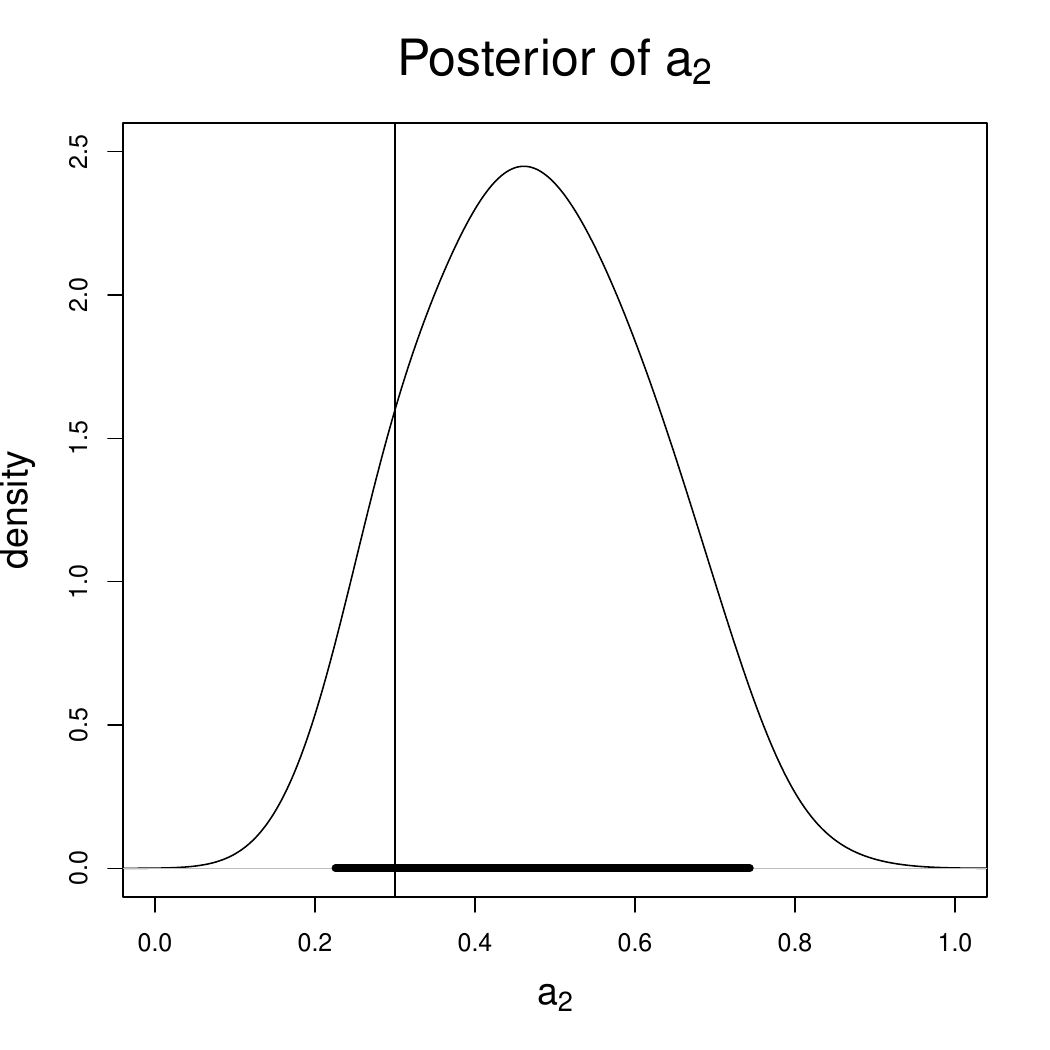}}
\caption{{\bf TTMCMC for $SDE_1$ and $\pi_2$:} Posteriors of $M$, $\mu_1$, $\mu_2$, $\omega^2_1$, $\omega^2_2$, $a_1$ and $a_2$. The vertical lines stand
for the true values, while the thick horizontal lines denote the 95\% credible intervals.} 
\label{fig:sim2_posterior_plots}
\end{figure}

However, the trace plots and the posterior
density estimates of the parameters are provided in 
Figures \ref{fig:sim2_trace_plots} and \ref{fig:sim2_posterior_plots}, demonstrate excellent performance, despite the more challenging nature of the problem.
The overall acceptance rate of our algorithm in this case is about $8\%$.

Here the posterior distribution of $M$ gave positive mass to $M=1,2,3$, with the respective probabilities being $0.2258$, $0.6391$ and $0.1351$.
Although the posterior probability of $M=2$ is less than in the previous case of well-separated components, this is still very significantly large compared
to the posterior probabilities of the other components, demonstrating that the Bayesian method has handled the challenge quite efficiently.

\subsection{Third simulation study: $SDE_2$ with $\pi_1$} 
\label{subsec:sim3}

\begin{figure}
\centering
\subfigure[Trace plot of $M$.]{ \label{fig:sim3_trace_comp}
\includegraphics[width=7cm,height=5cm]{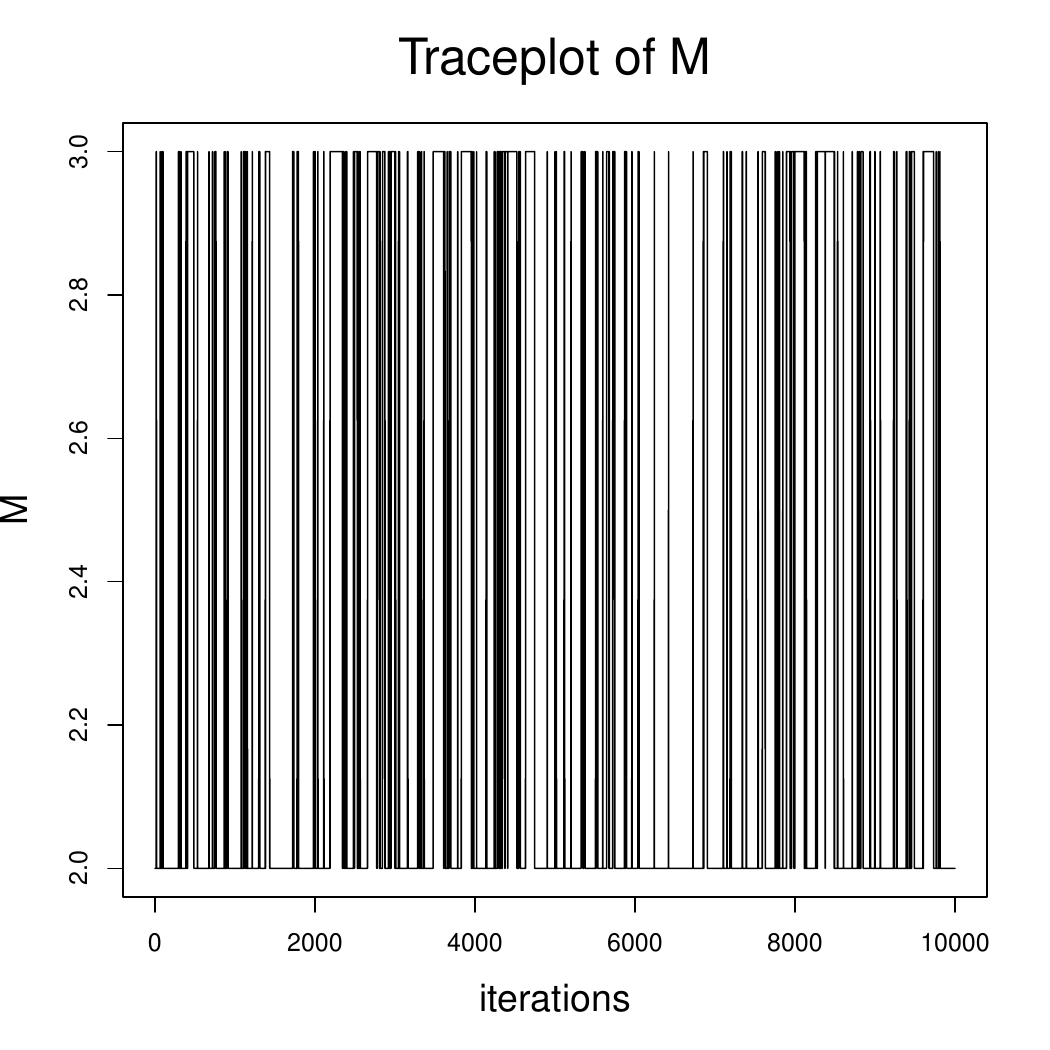}}
\hspace{2mm}
\subfigure[Trace plot of $\mu_1$.]{ \label{fig:sim3_trace_mu1}
\includegraphics[width=7cm,height=5cm]{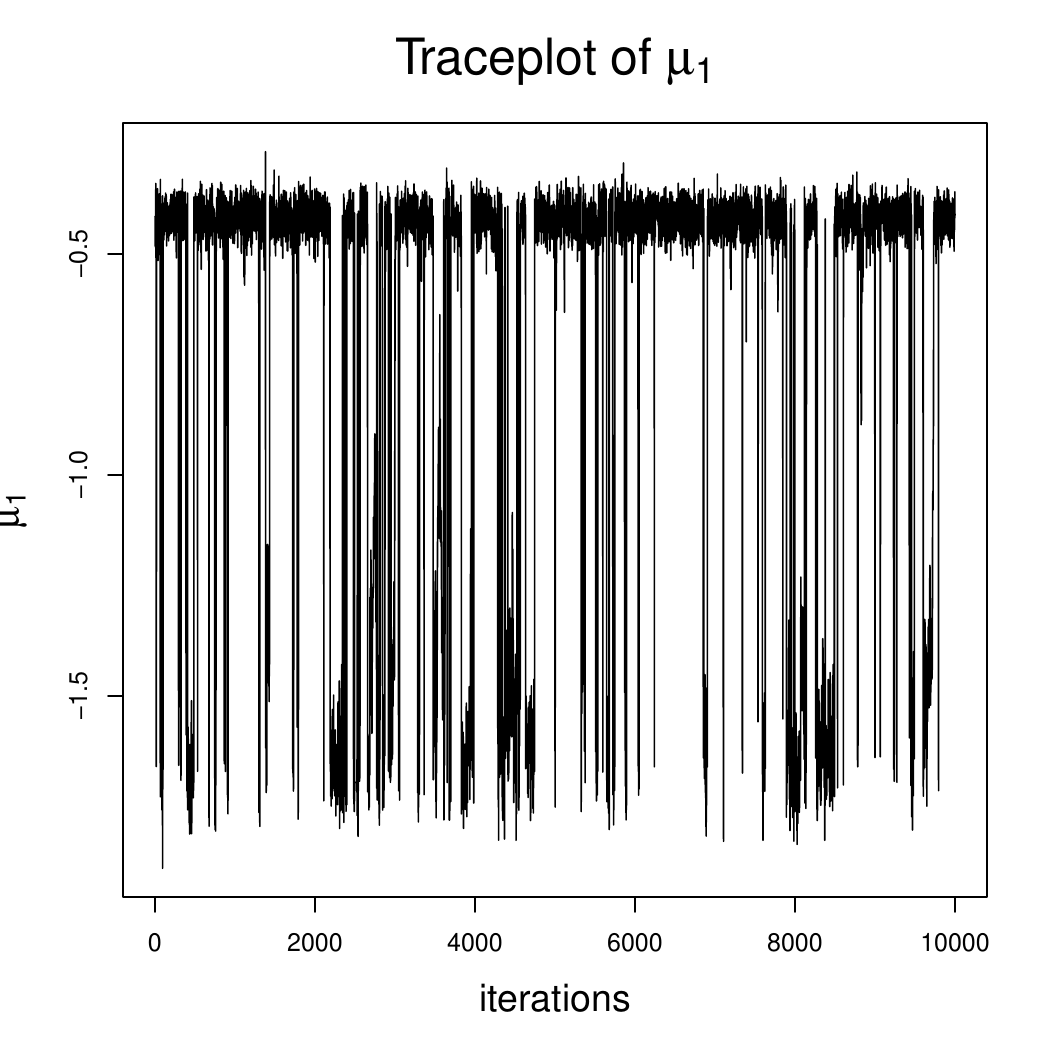}}\\
\vspace{2mm}
\subfigure[Trace plot of $\mu_2$.]{ \label{fig:sim3_trace_mu2}
\includegraphics[width=7cm,height=5cm]{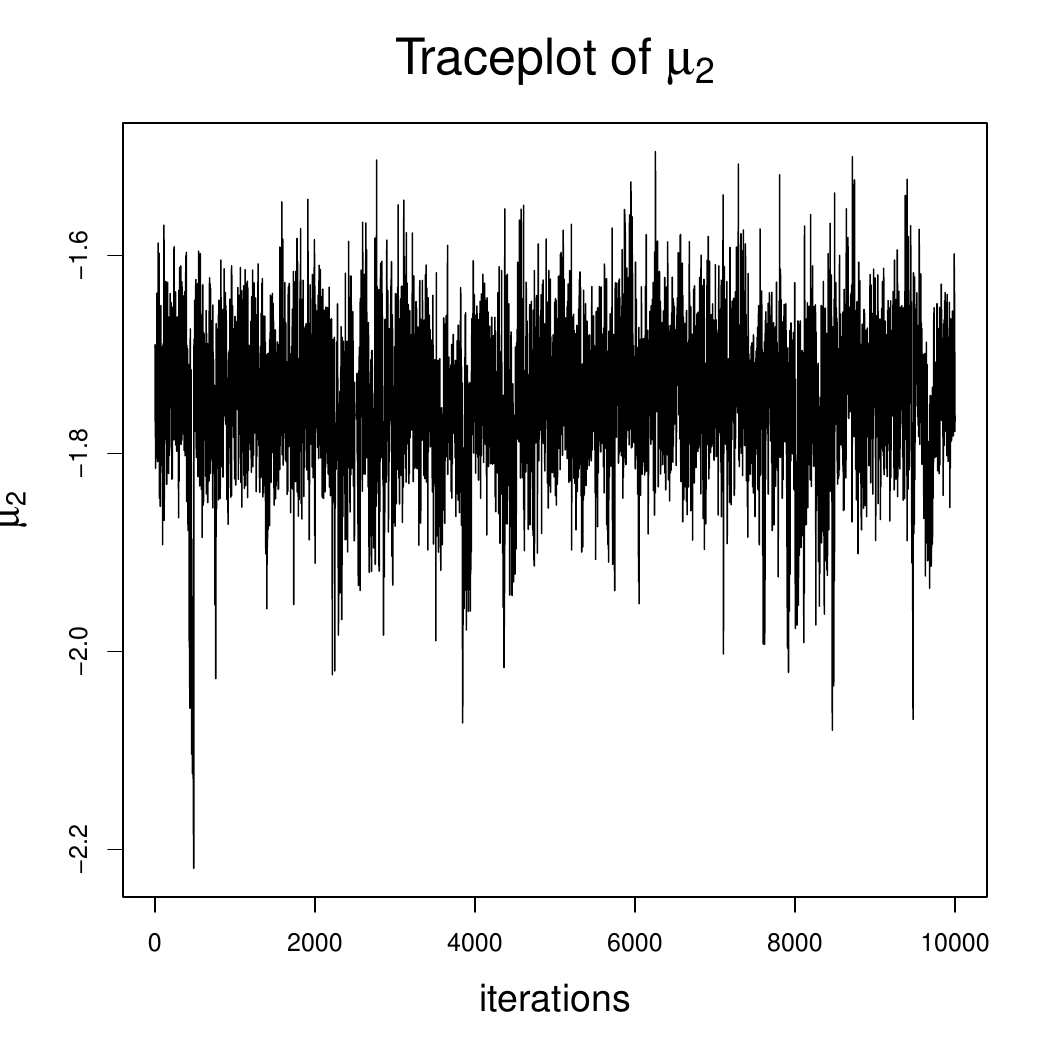}}
\vspace{2mm}
\subfigure[Trace plot of $\omega^2_1$.]{ \label{fig:sim3_trace_omegasq1}
\includegraphics[width=7cm,height=5cm]{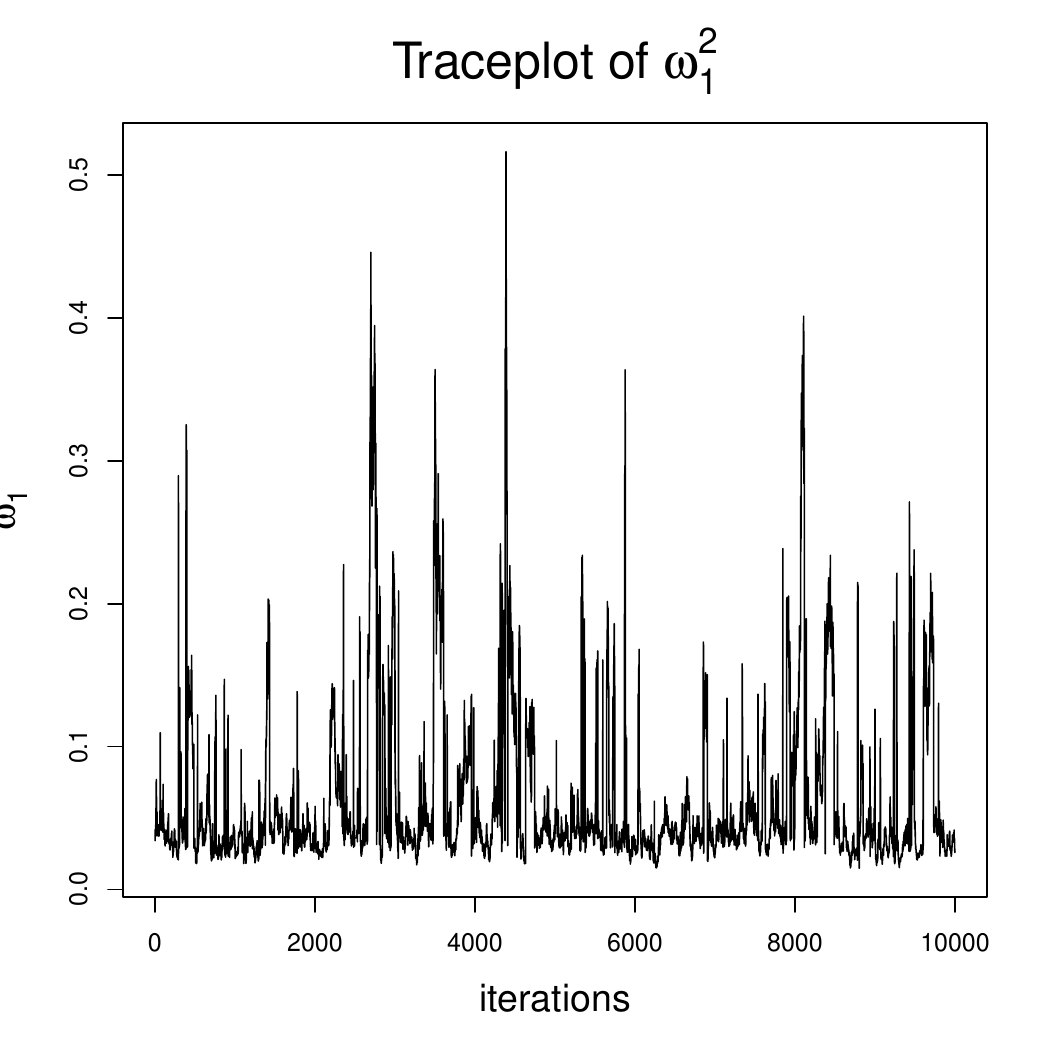}}\\
\vspace{2mm}
\subfigure[Trace plot of $\omega^2_2$.]{ \label{fig:sim3_trace_omegasq2}
\includegraphics[width=7cm,height=5cm]{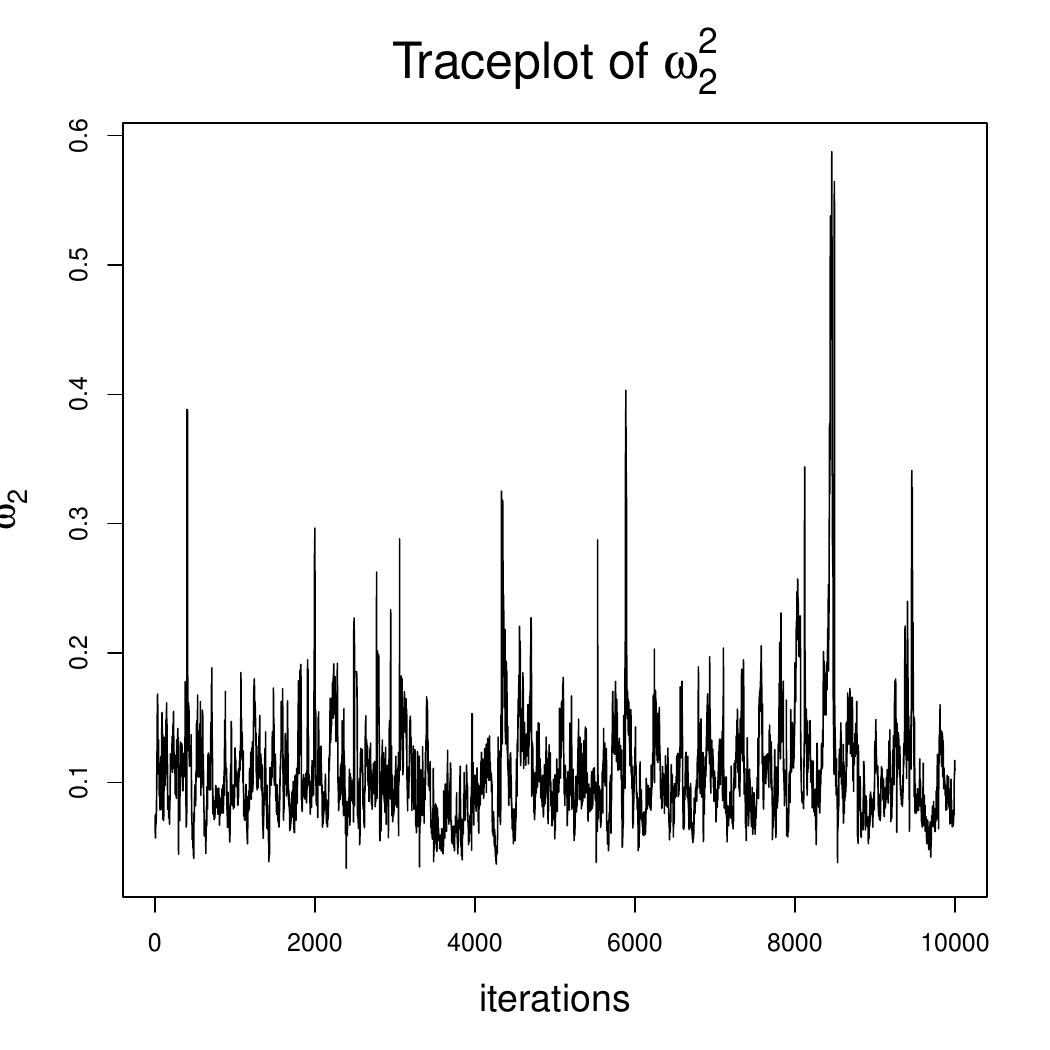}}
\hspace{2mm}
\subfigure[Trace plot of $a_1$.]{ \label{fig:sim3_trace_p1}
\includegraphics[width=7cm,height=5cm]{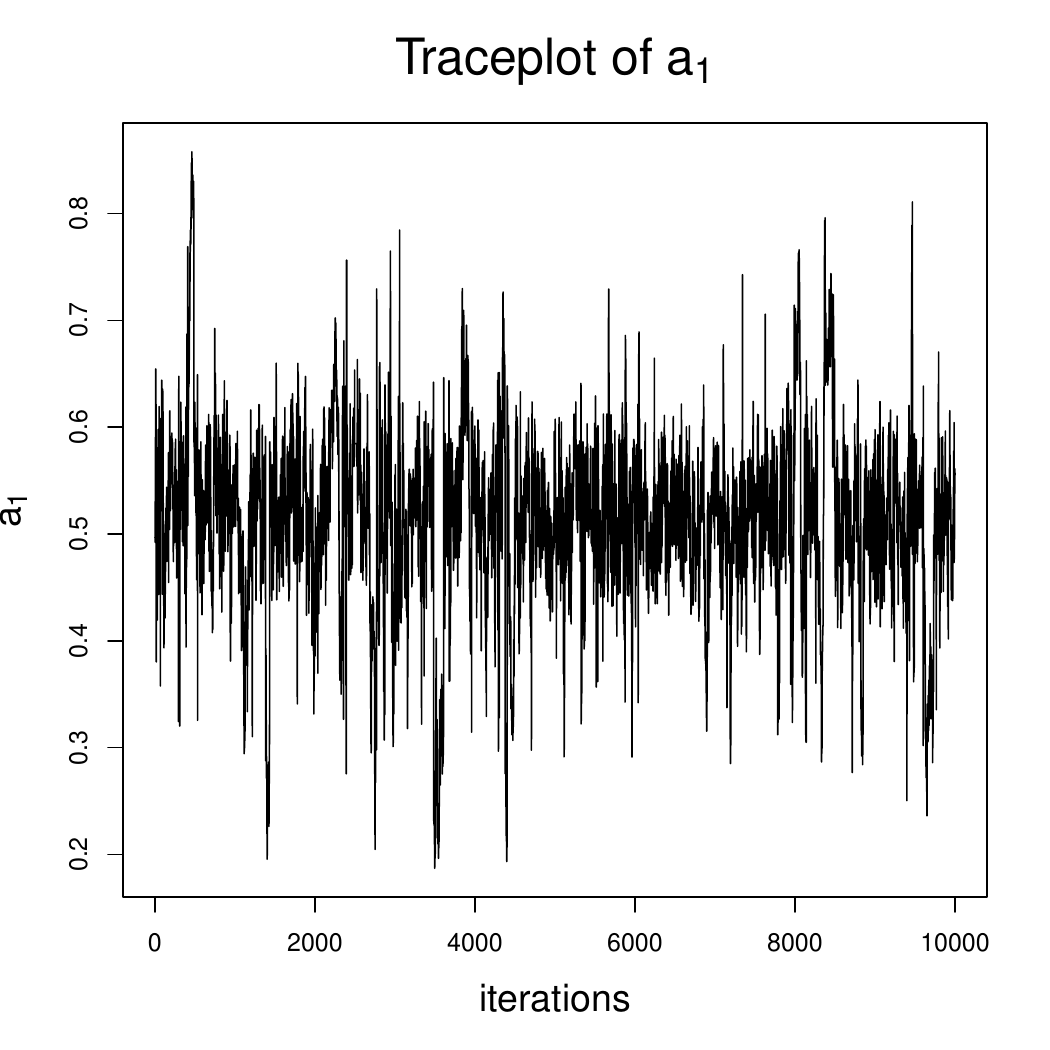}}\\
\vspace{2mm}
\subfigure[Trace plot of $a_2$.]{ \label{fig:sim3_trace_p2}
\includegraphics[width=7cm,height=5cm]{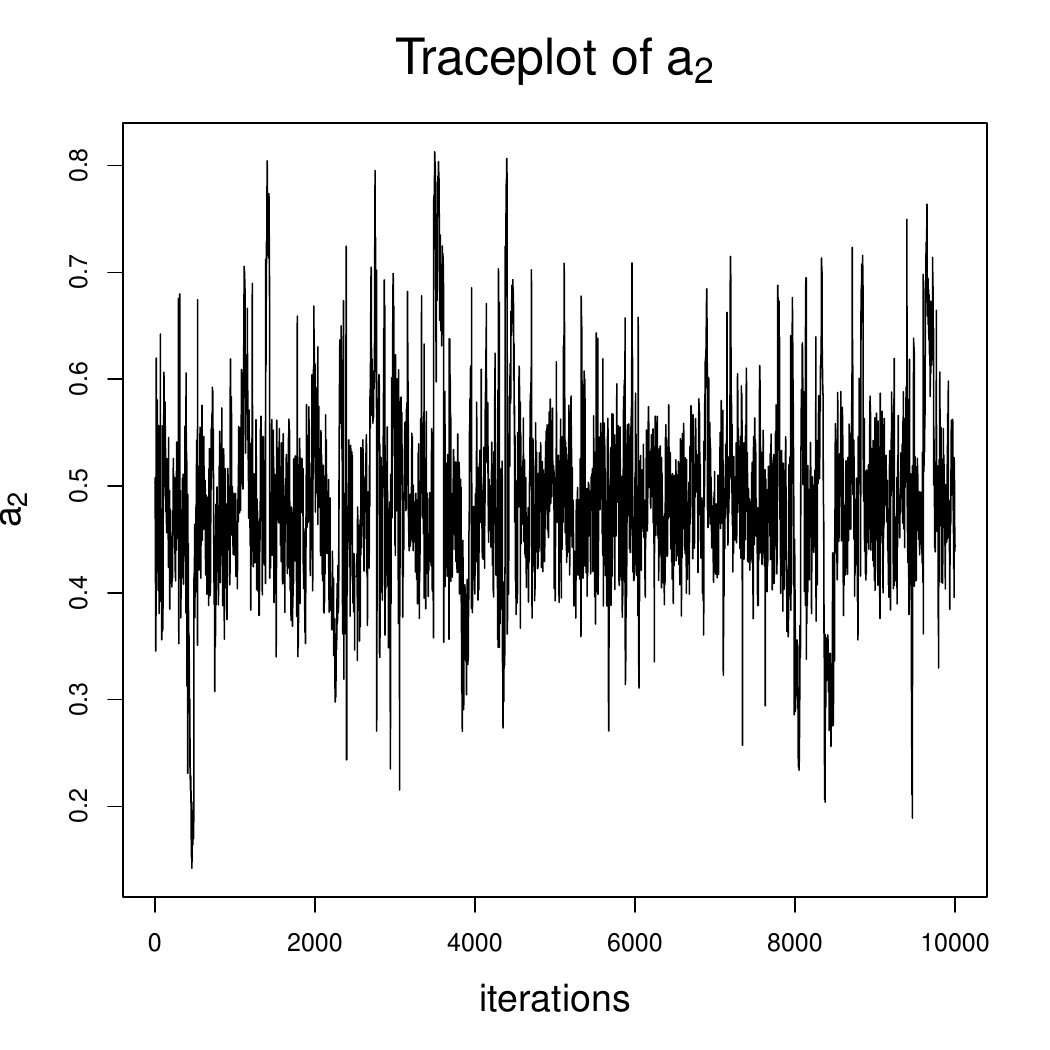}}
\caption{{\bf TTMCMC for $SDE_2$ and $\pi_1$:} Trace plots of $M$, $\mu_1$, $\mu_2$, $\omega^2_1$, $\omega^2_2$, $a_1$ and $a_2$.} 
\label{fig:sim3_trace_plots}
\end{figure}

\begin{figure}
\centering
\subfigure[Posterior of $\mu_1$.]{ \label{fig:sim3_mu1}
\includegraphics[width=7cm,height=6cm]{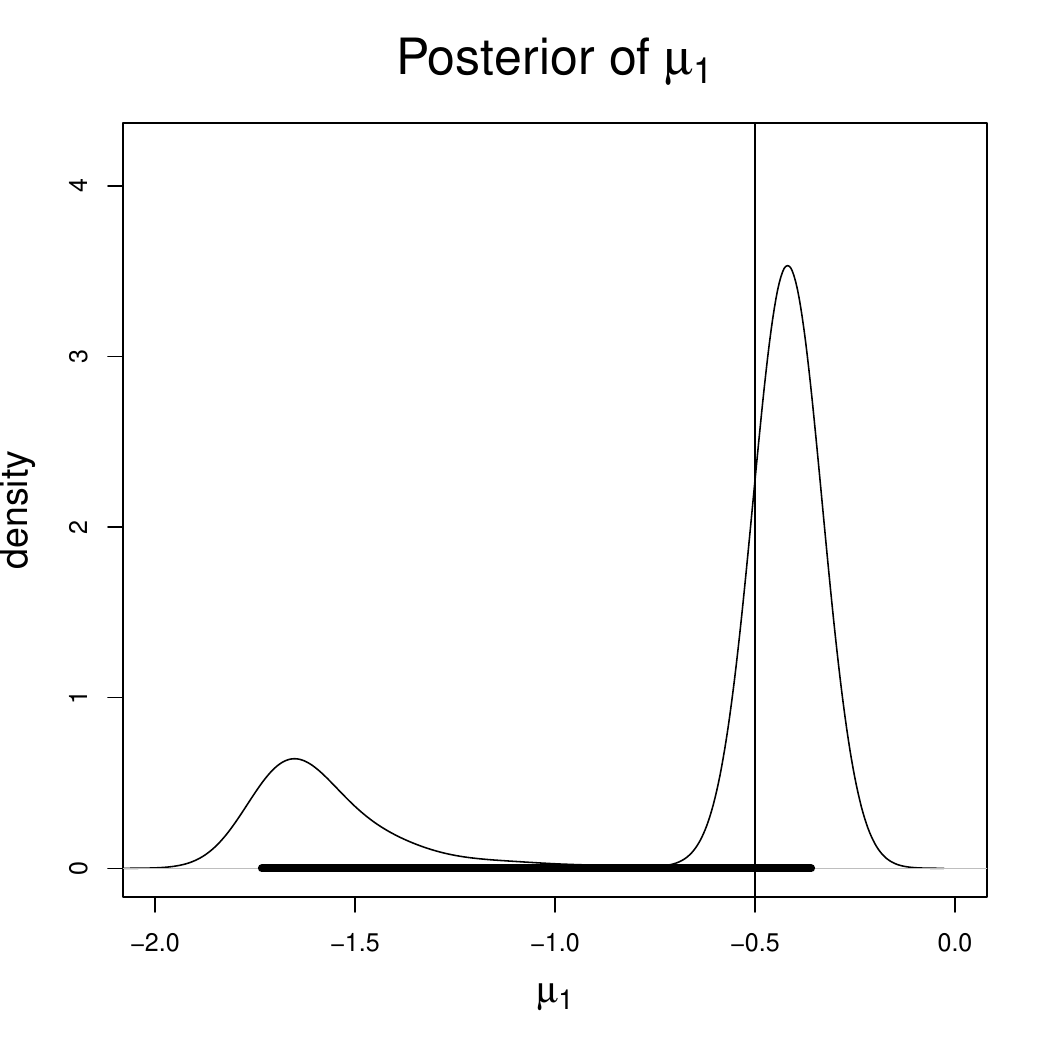}}
\hspace{2mm}
\subfigure[Posterior of $\mu_2$.]{ \label{fig:sim3_mu2}
\includegraphics[width=7cm,height=6cm]{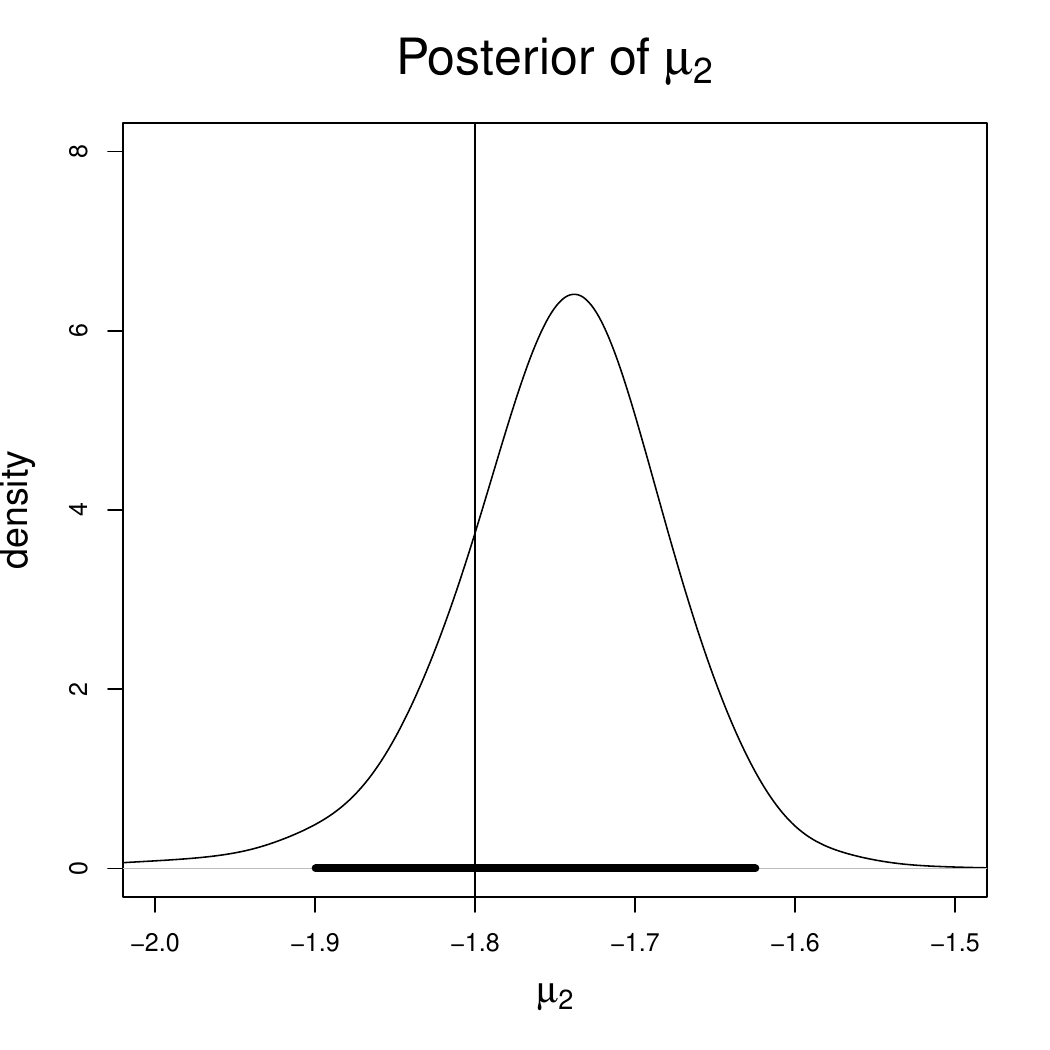}}\\
\vspace{2mm}
\subfigure[Posterior of $\omega^2_1$.]{ \label{fig:sim3_omegasq1}
\includegraphics[width=7cm,height=6cm]{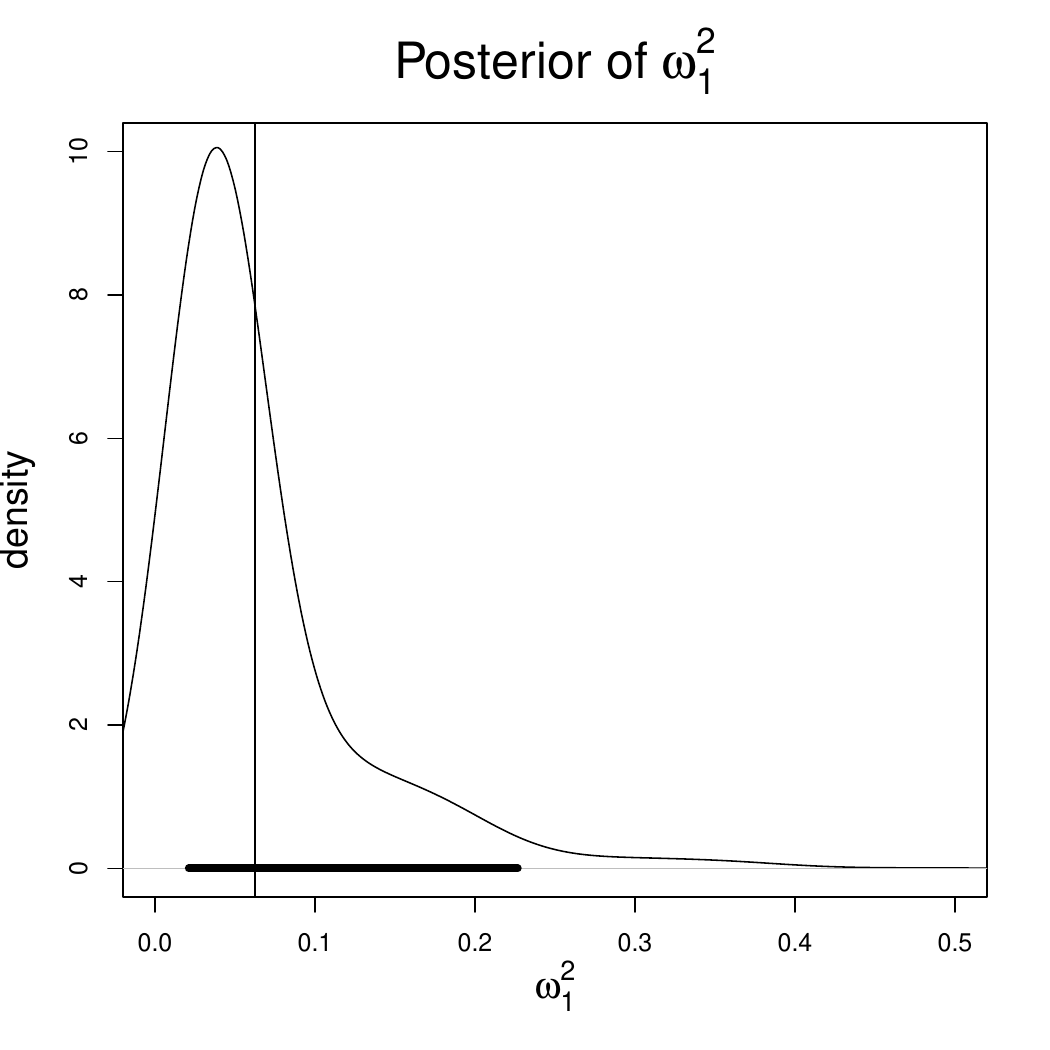}}
\hspace{2mm}
\subfigure[Posterior of $\omega^2_2$.]{ \label{fig:sim3_omegasq2}
\includegraphics[width=7cm,height=6cm]{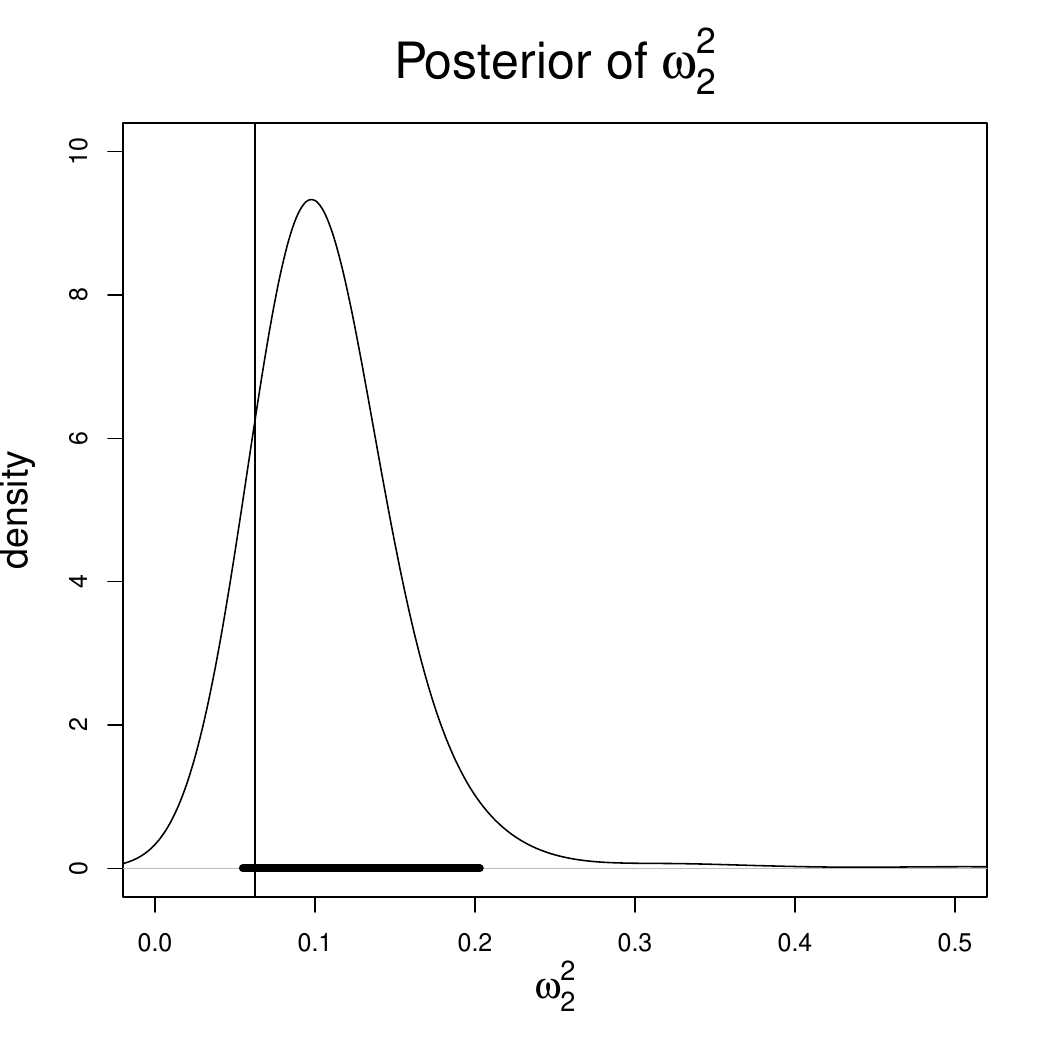}}\\
\vspace{2mm}
\subfigure[Posterior of $a_1$.]{ \label{fig:sim3_p1}
\includegraphics[width=7cm,height=6cm]{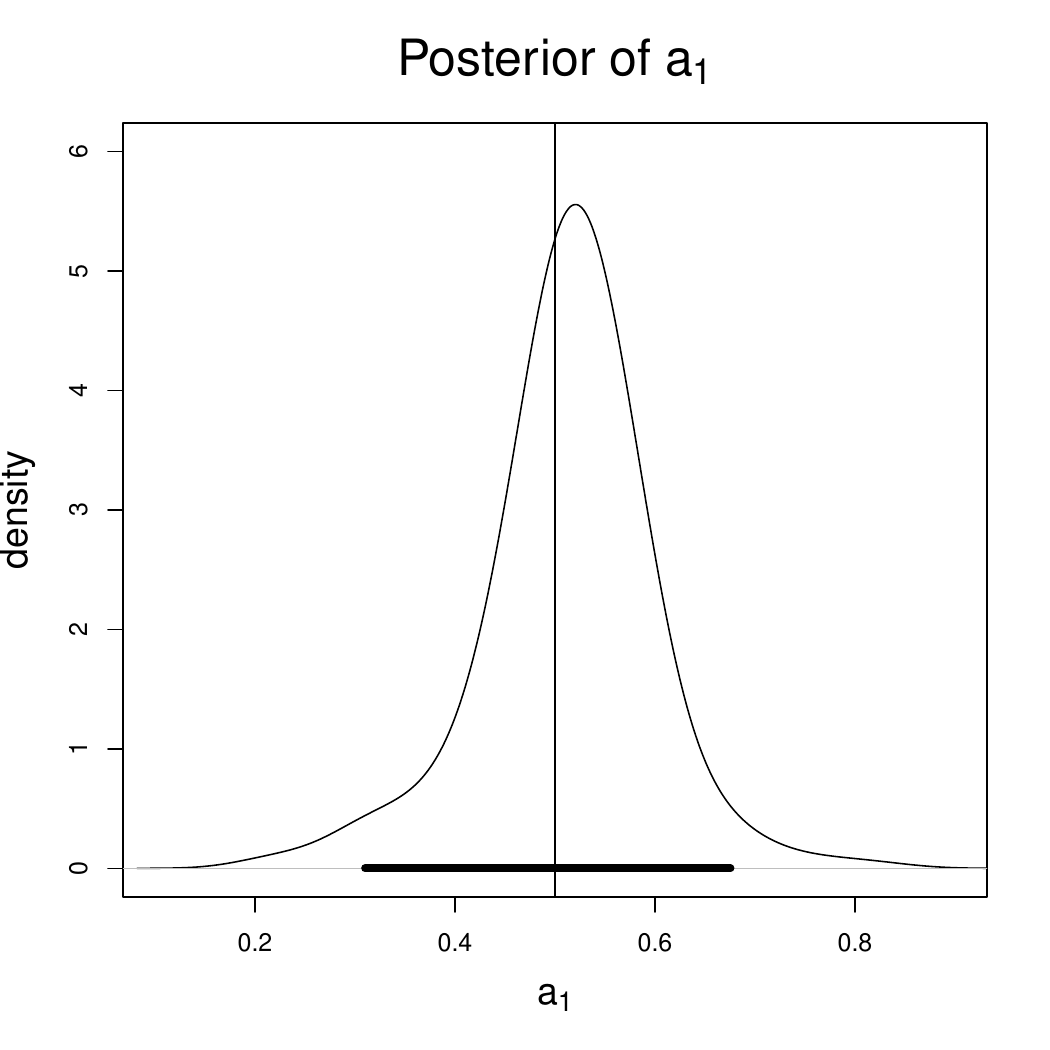}}
\hspace{2mm}
\subfigure[Posterior of $a_2$.]{ \label{fig:sim3_p2}
\includegraphics[width=7cm,height=6cm]{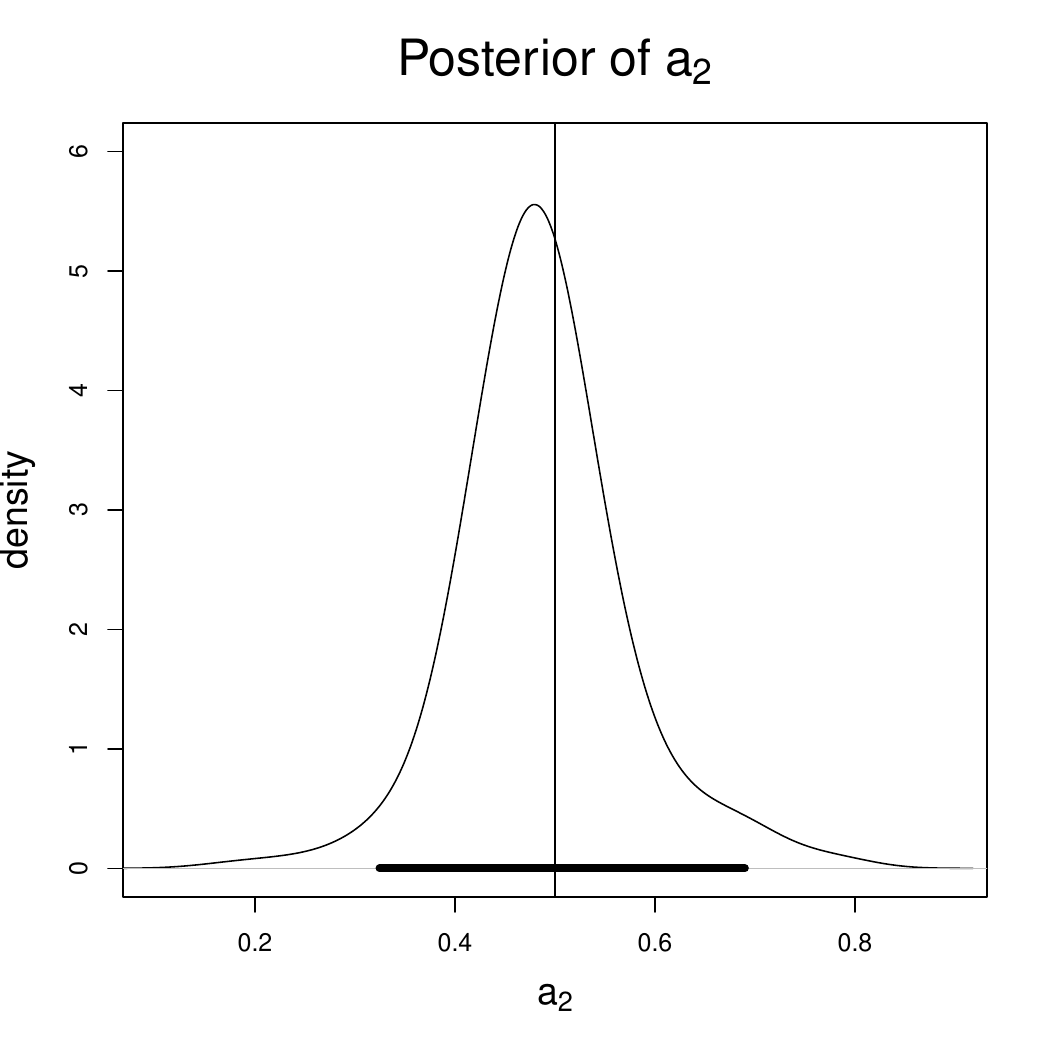}}
\caption{{\bf TTMCMC for $SDE_2$ and $\pi_1$:} Posteriors of $M$, $\mu_1$, $\mu_2$, $\omega^2_1$, $\omega^2_2$, $a_1$ and $a_2$. The vertical lines stand
for the true values, while the thick horizontal lines denote the 95\% credible intervals.} 
\label{fig:sim3_posterior_plots}
\end{figure}

When the true model is given by $SDE_2$ and $\pi_1$, the trace plots and the posterior
density estimates of the parameters are provided in 
Figures \ref{fig:sim3_trace_plots} and \ref{fig:sim3_posterior_plots}, respectively. The performance, as indicated by the plots, is again very satisfactory. 
Here $M=2$ and $M=3$ got the posterior probabilities $0.7549$ and $0.2451$, and all other values received the zero posterior probability.
The overall acceptance rate in this case is about $2\%$.

\subsection{Fourth simulation study: $SDE_2$ with $\pi_2$} 
\label{subsec:sim4}

\begin{figure}
\centering
\subfigure[Trace plot of $M$.]{ \label{fig:sim4_trace_comp}
\includegraphics[width=7cm,height=5cm]{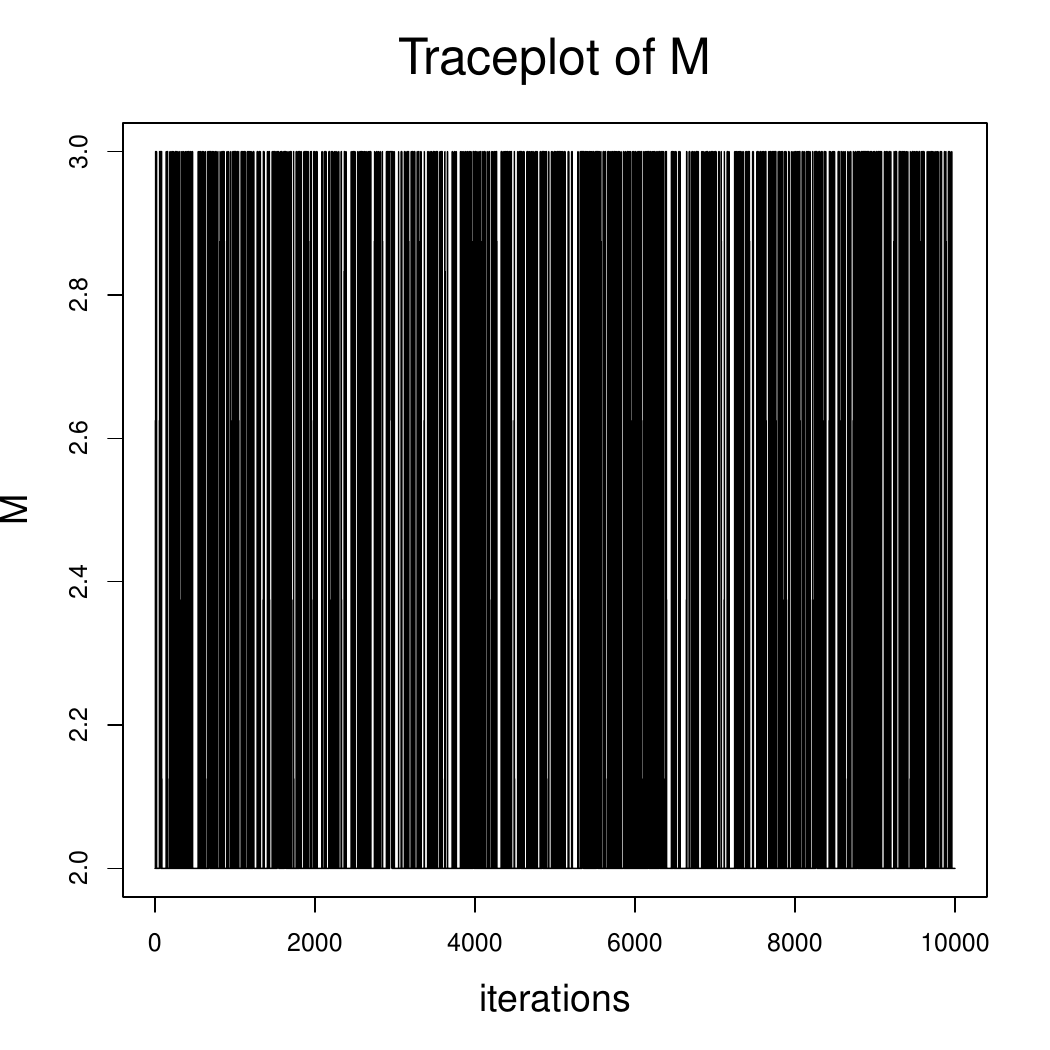}}
\hspace{2mm}
\subfigure[Trace plot of $\mu_1$.]{ \label{fig:sim4_trace_mu1}
\includegraphics[width=7cm,height=5cm]{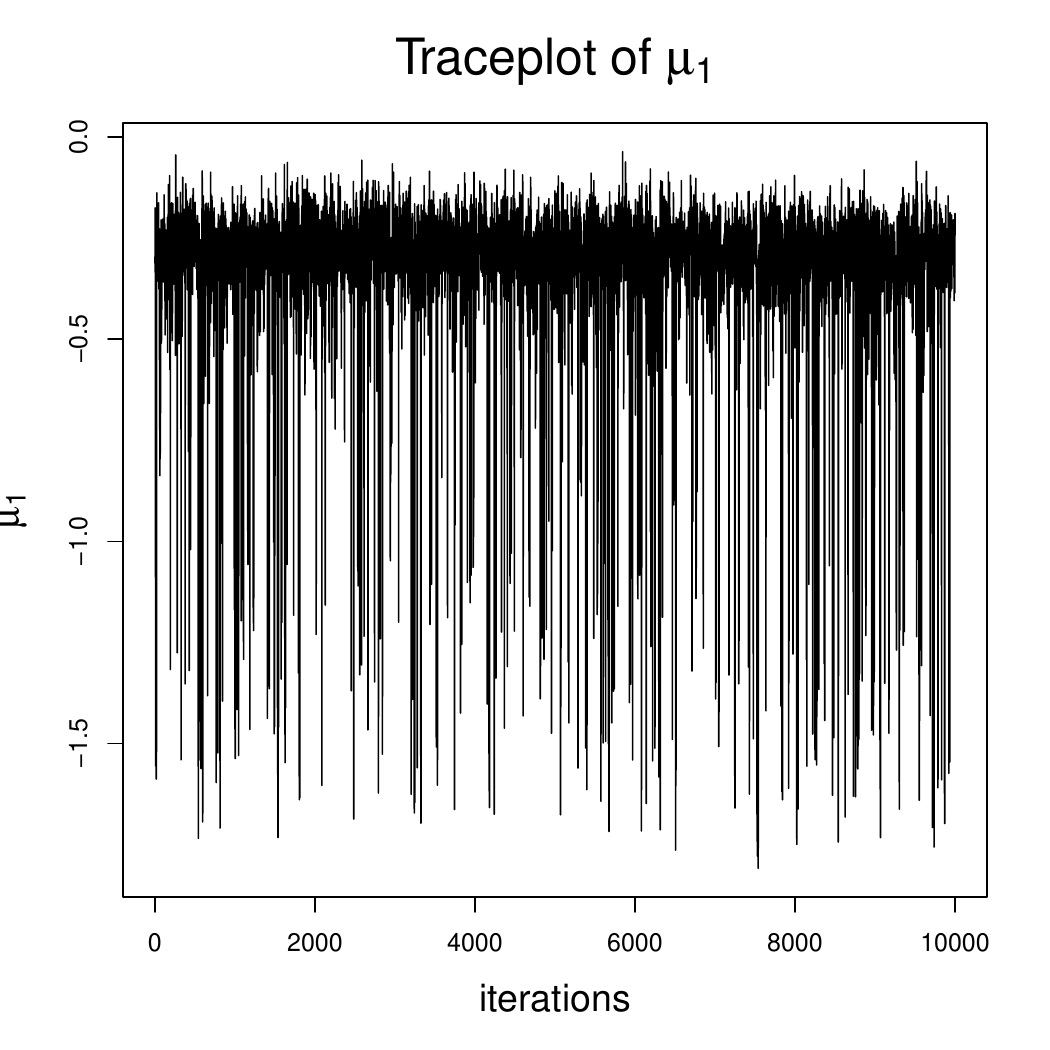}}\\
\vspace{2mm}
\subfigure[Trace plot of $\mu_2$.]{ \label{fig:sim4_trace_mu2}
\includegraphics[width=7cm,height=5cm]{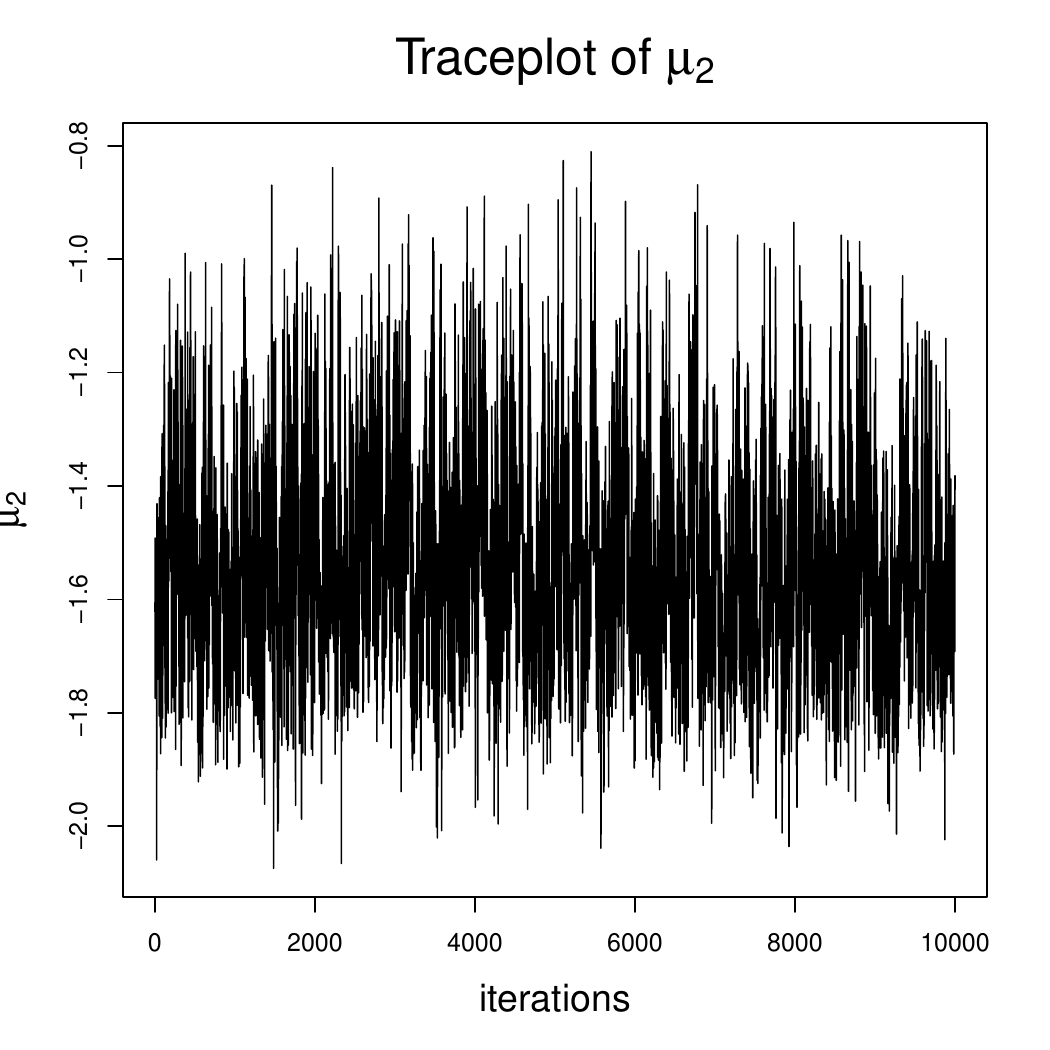}}
\vspace{2mm}
\subfigure[Trace plot of $\omega^2_1$.]{ \label{fig:sim4_trace_omegasq1}
\includegraphics[width=7cm,height=5cm]{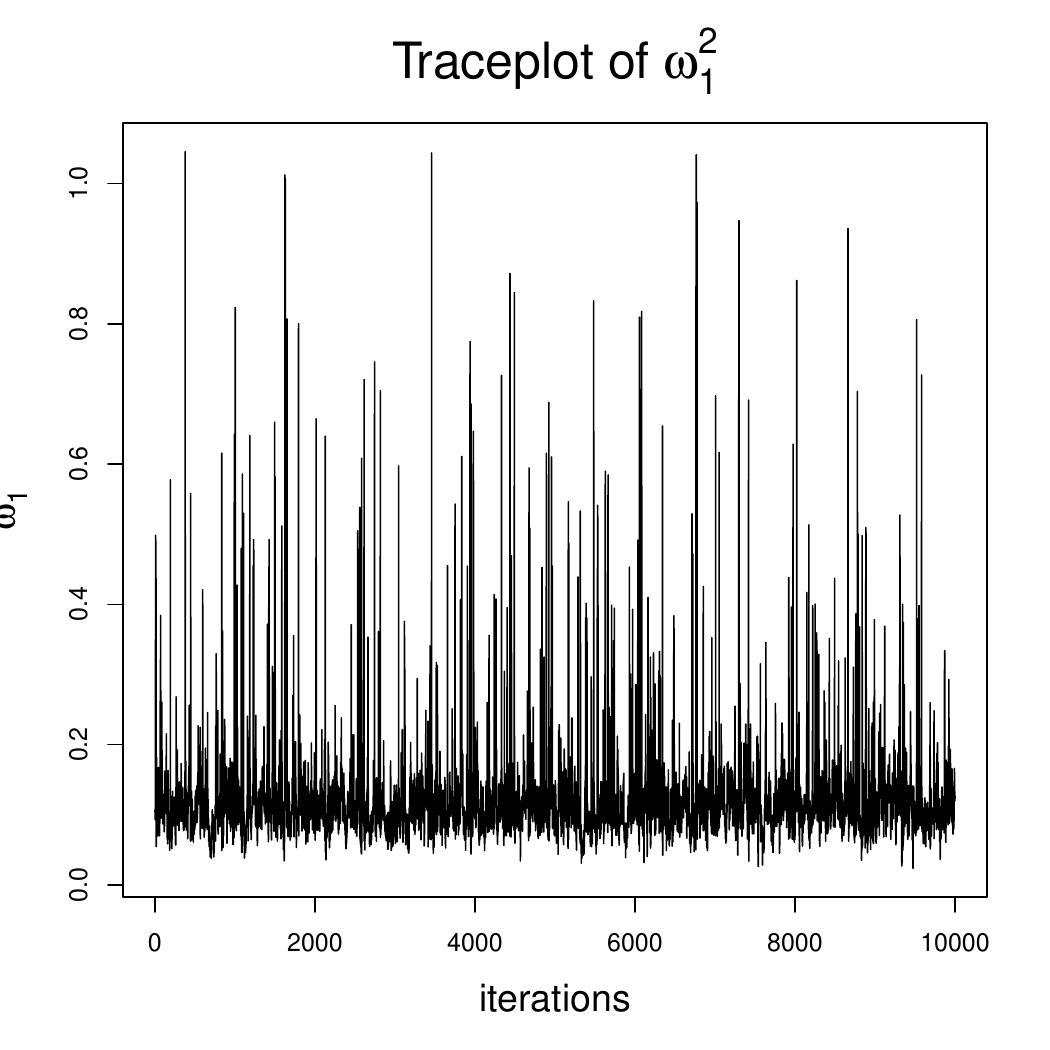}}\\
\vspace{2mm}
\subfigure[Trace plot of $\omega^2_2$.]{ \label{fig:sim4_trace_omegasq2}
\includegraphics[width=7cm,height=5cm]{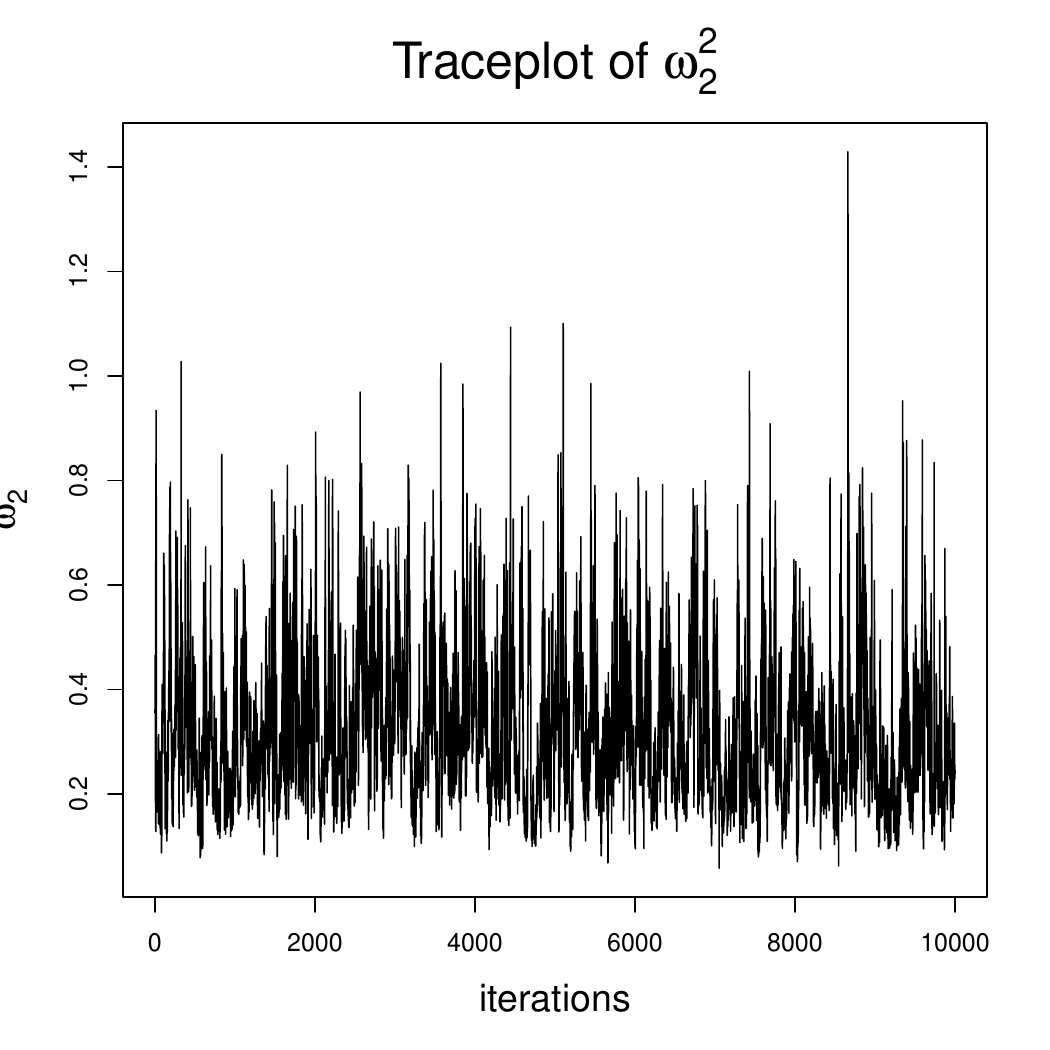}}
\hspace{2mm}
\subfigure[Trace plot of $a_1$.]{ \label{fig:sim4_trace_p1}
\includegraphics[width=7cm,height=5cm]{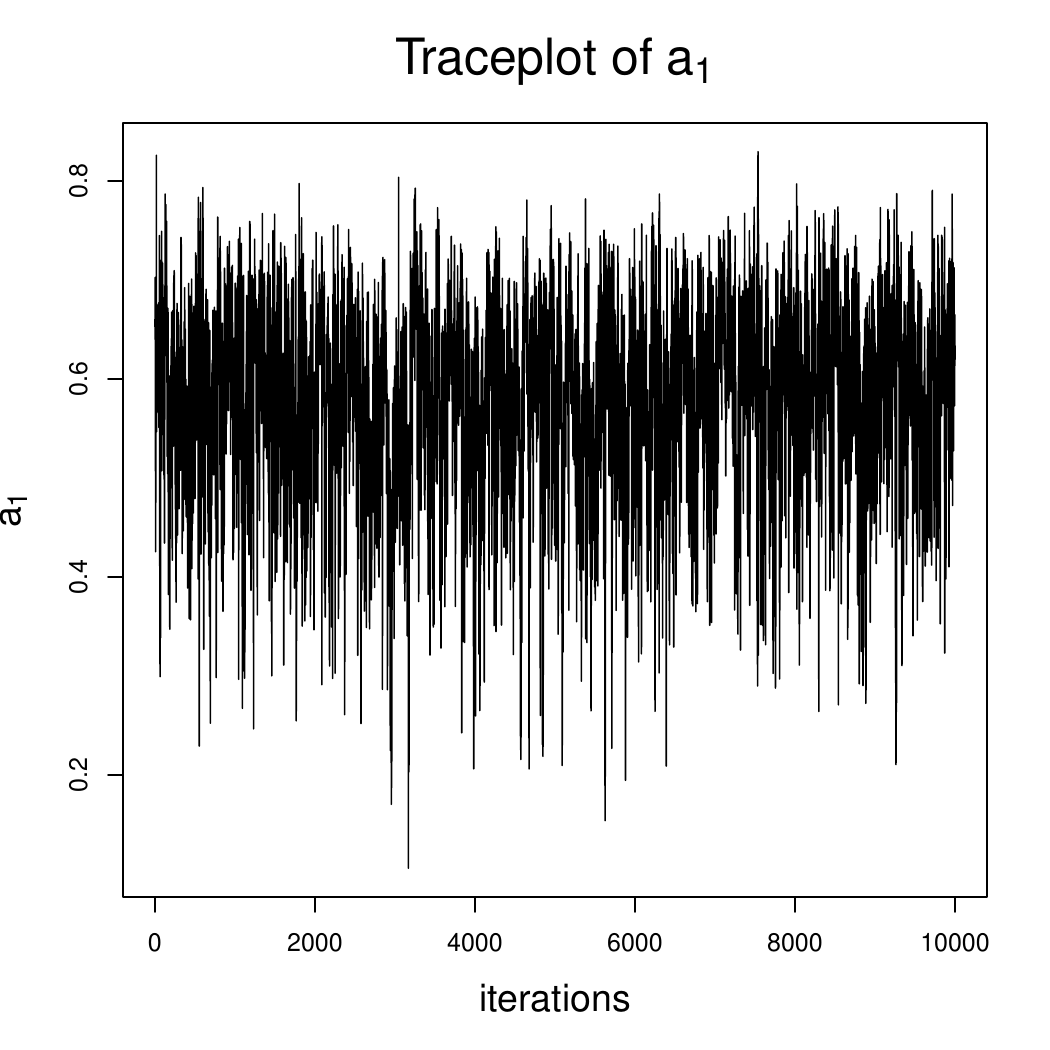}}\\
\vspace{2mm}
\subfigure[Trace plot of $a_2$.]{ \label{fig:sim4_trace_p2}
\includegraphics[width=7cm,height=5cm]{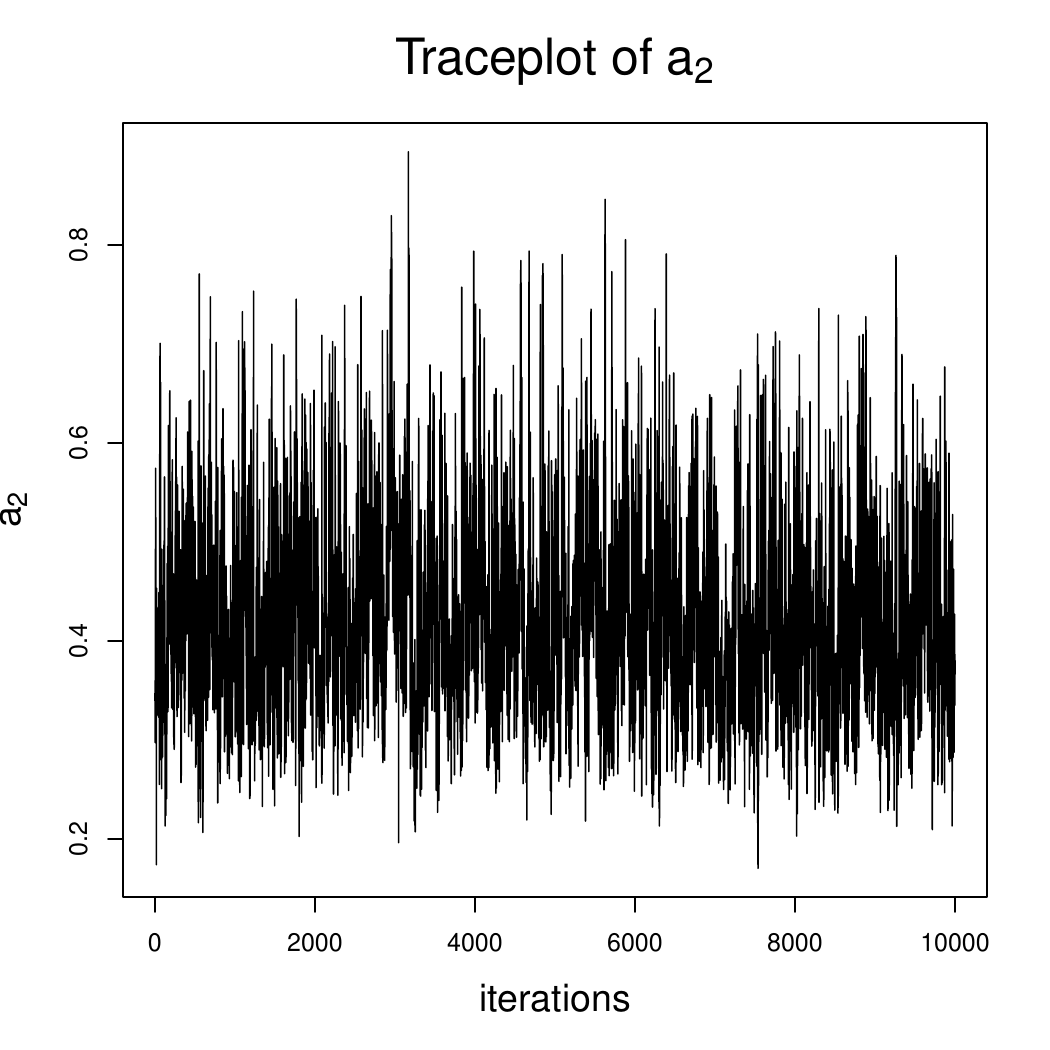}}
\caption{{\bf TTMCMC for $SDE_2$ and $\pi_2$:} Trace plots of $M$, $\mu_1$, $\mu_2$, $\omega^2_1$, $\omega^2_2$, $a_1$ and $a_2$.} 
\label{fig:sim4_trace_plots}
\end{figure}

\begin{figure}
\centering
\subfigure[Posterior of $\mu_1$.]{ \label{fig:sim4_mu1}
\includegraphics[width=7cm,height=6cm]{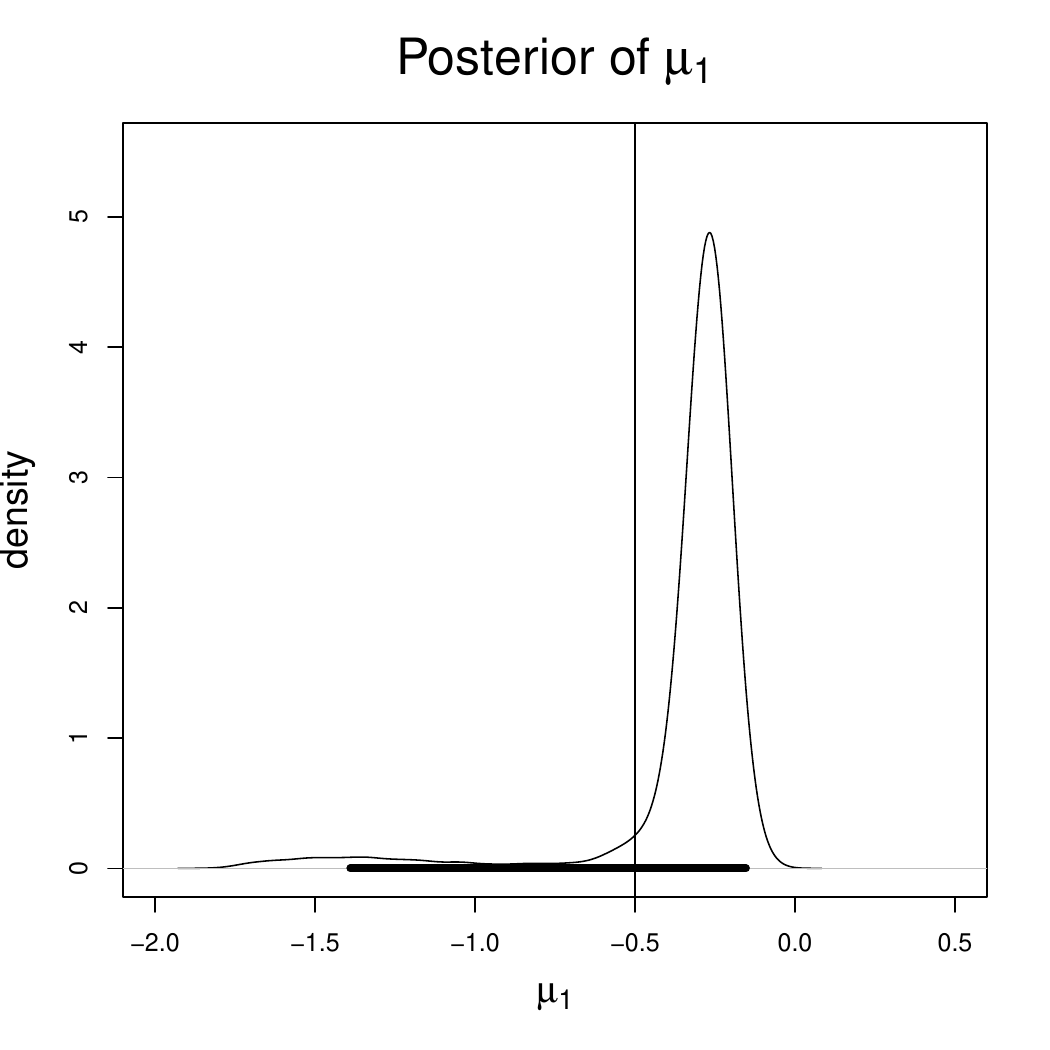}}
\hspace{2mm}
\subfigure[Posterior of $\mu_2$.]{ \label{fig:sim4_mu2}
\includegraphics[width=7cm,height=6cm]{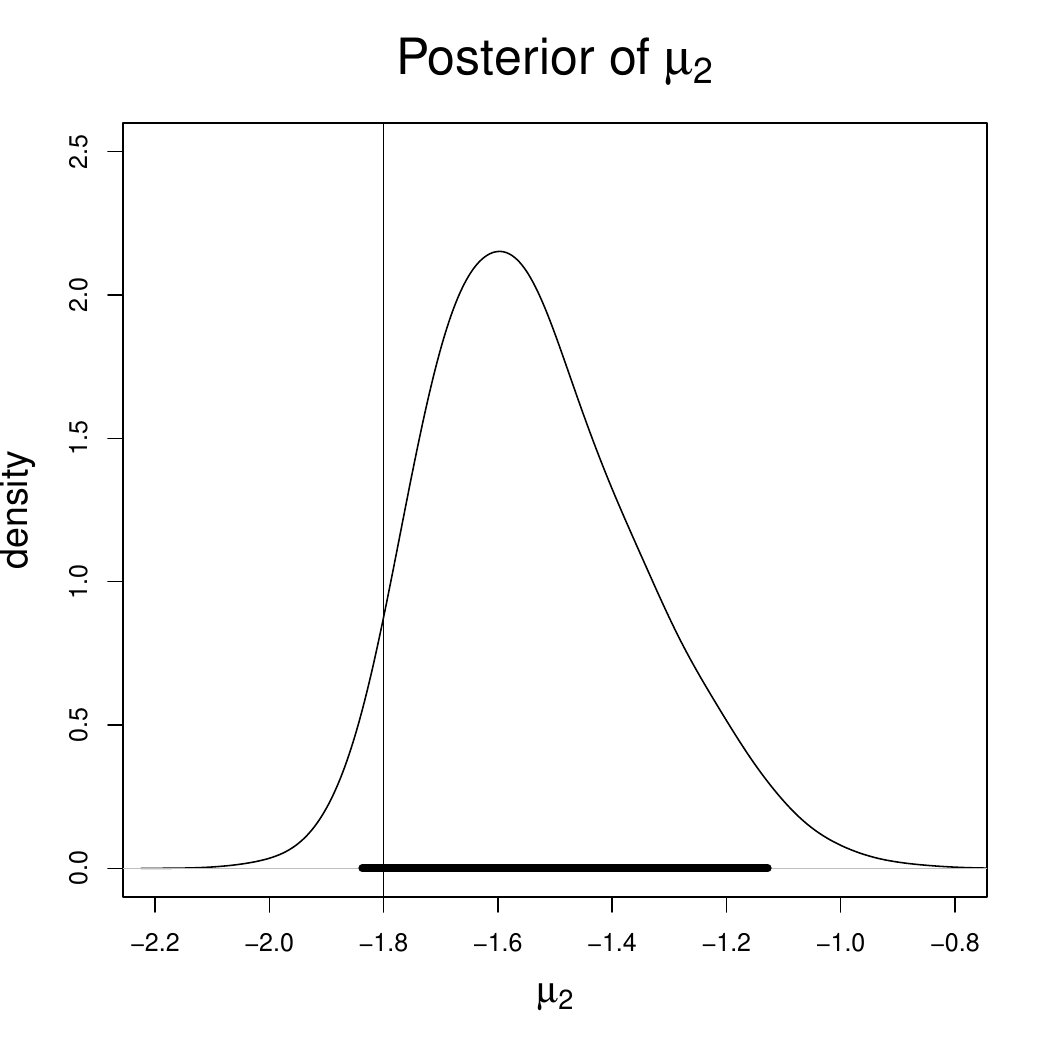}}\\
\vspace{2mm}
\subfigure[Posterior of $\omega^2_1$.]{ \label{fig:sim4_omegasq1}
\includegraphics[width=7cm,height=6cm]{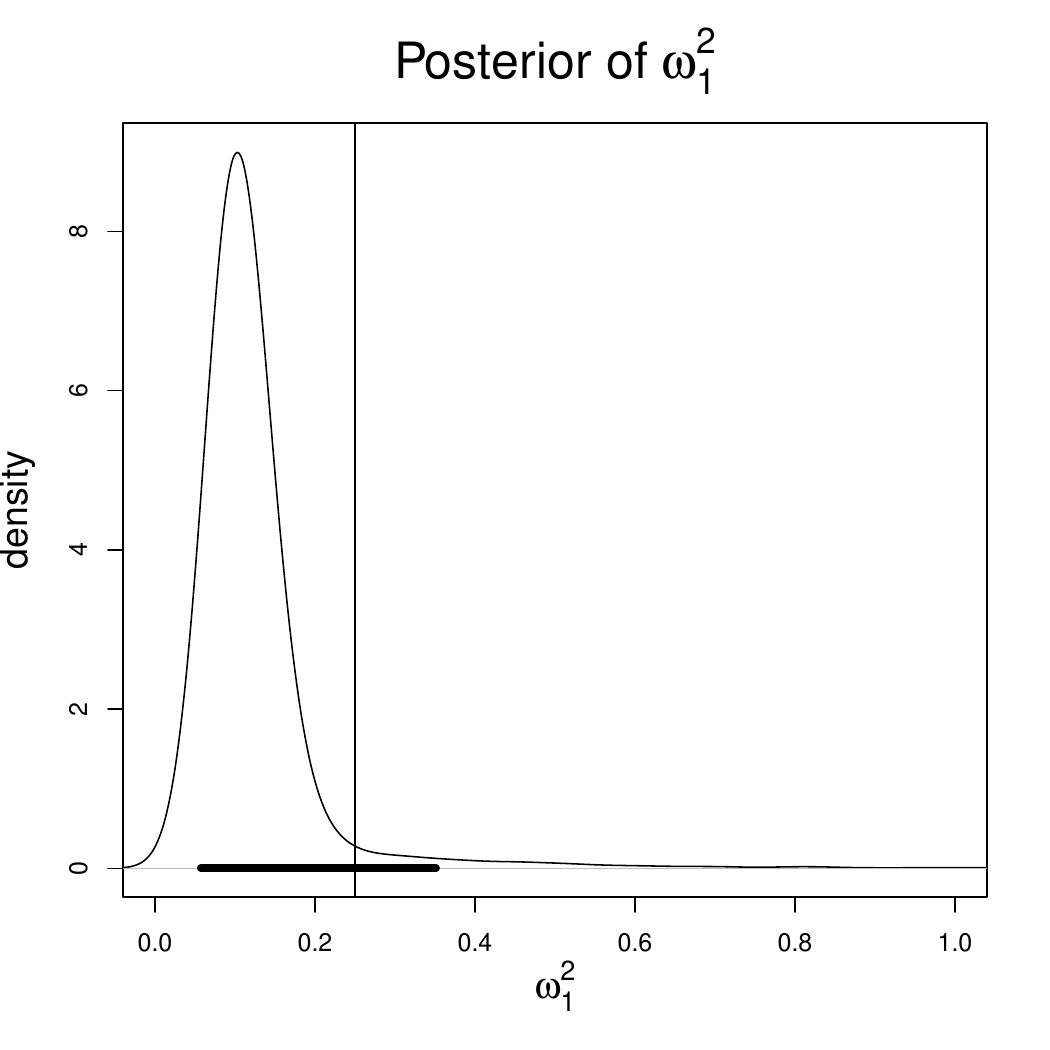}}
\hspace{2mm}
\subfigure[Posterior of $\omega^2_2$.]{ \label{fig:sim4_omegasq2}
\includegraphics[width=7cm,height=6cm]{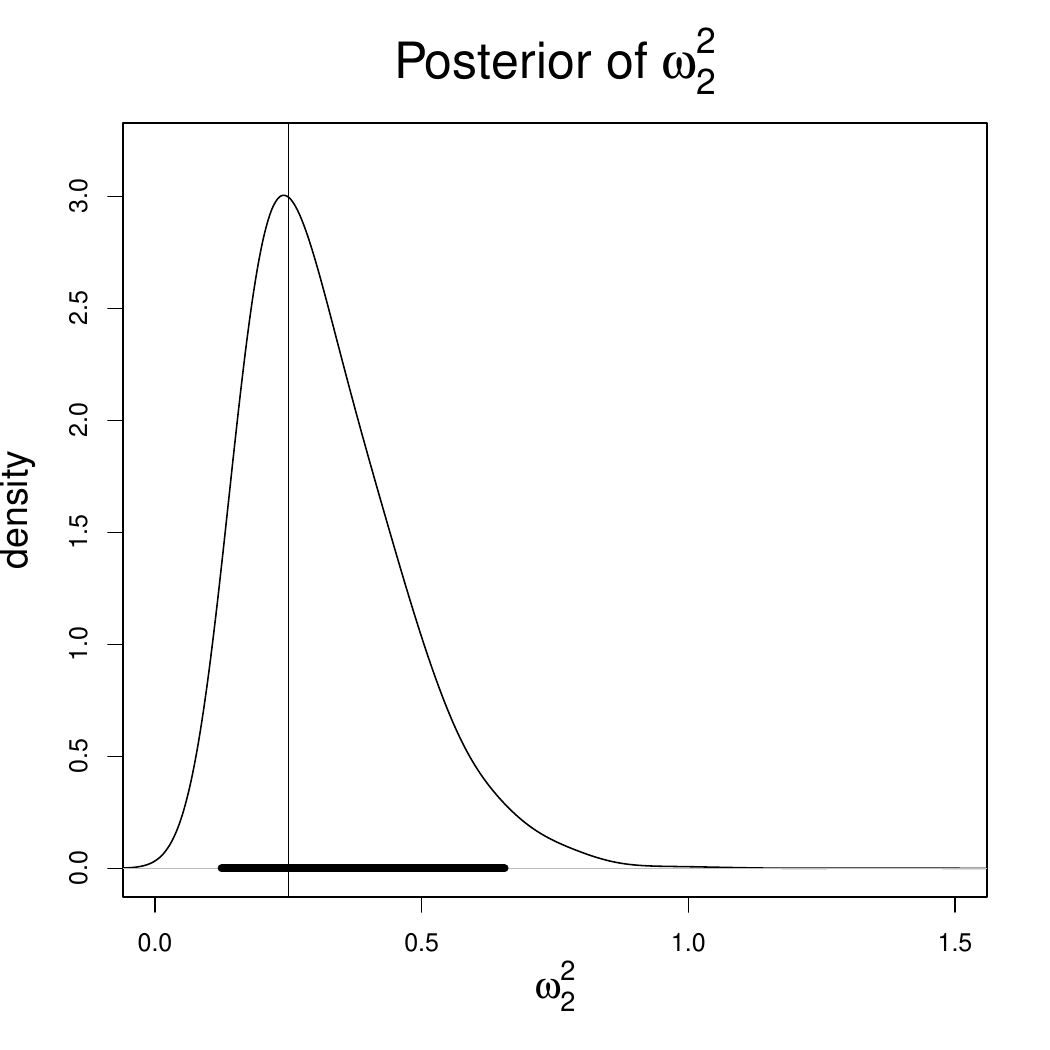}}\\
\vspace{2mm}
\subfigure[Posterior of $a_1$.]{ \label{fig:sim4_p1}
\includegraphics[width=7cm,height=6cm]{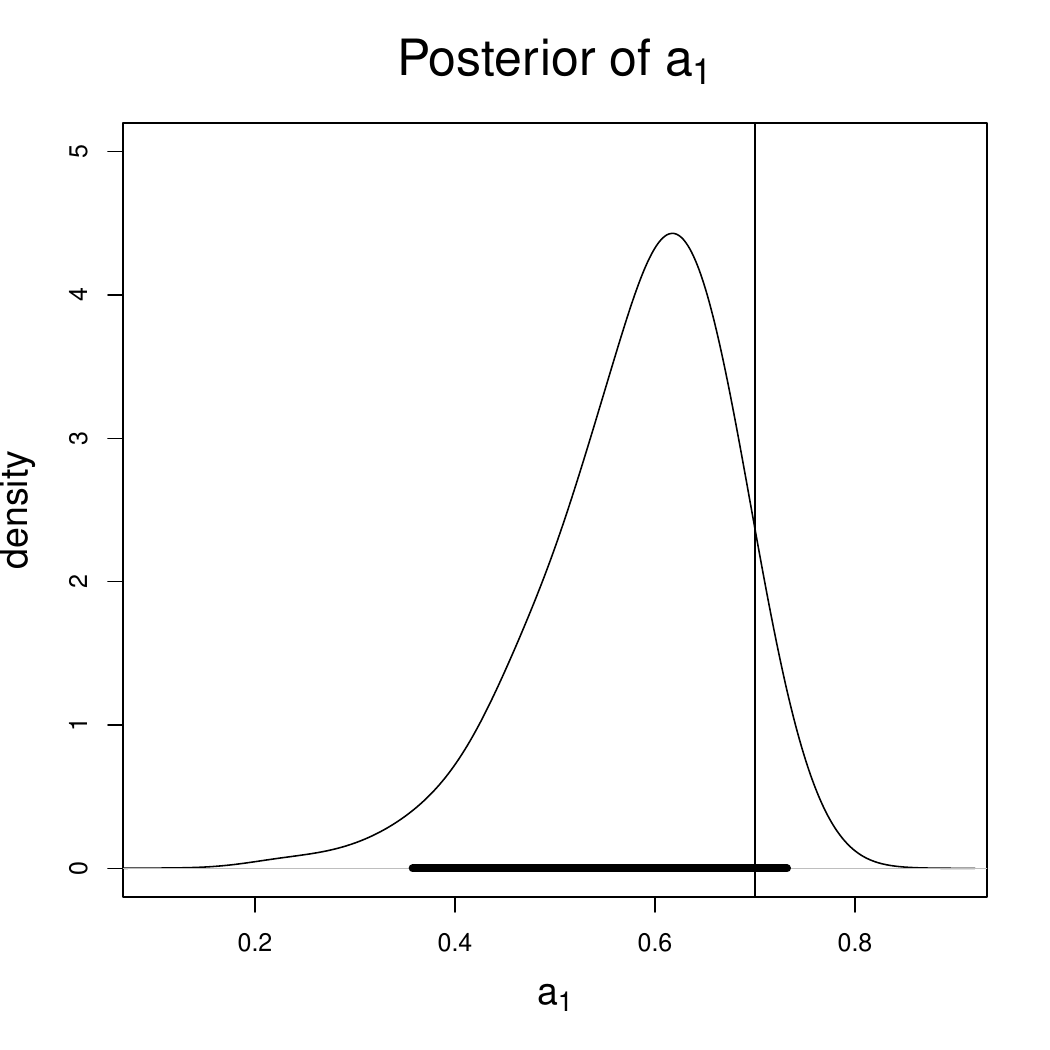}}
\hspace{2mm}
\subfigure[Posterior of $a_2$.]{ \label{fig:sim4_p2}
\includegraphics[width=7cm,height=6cm]{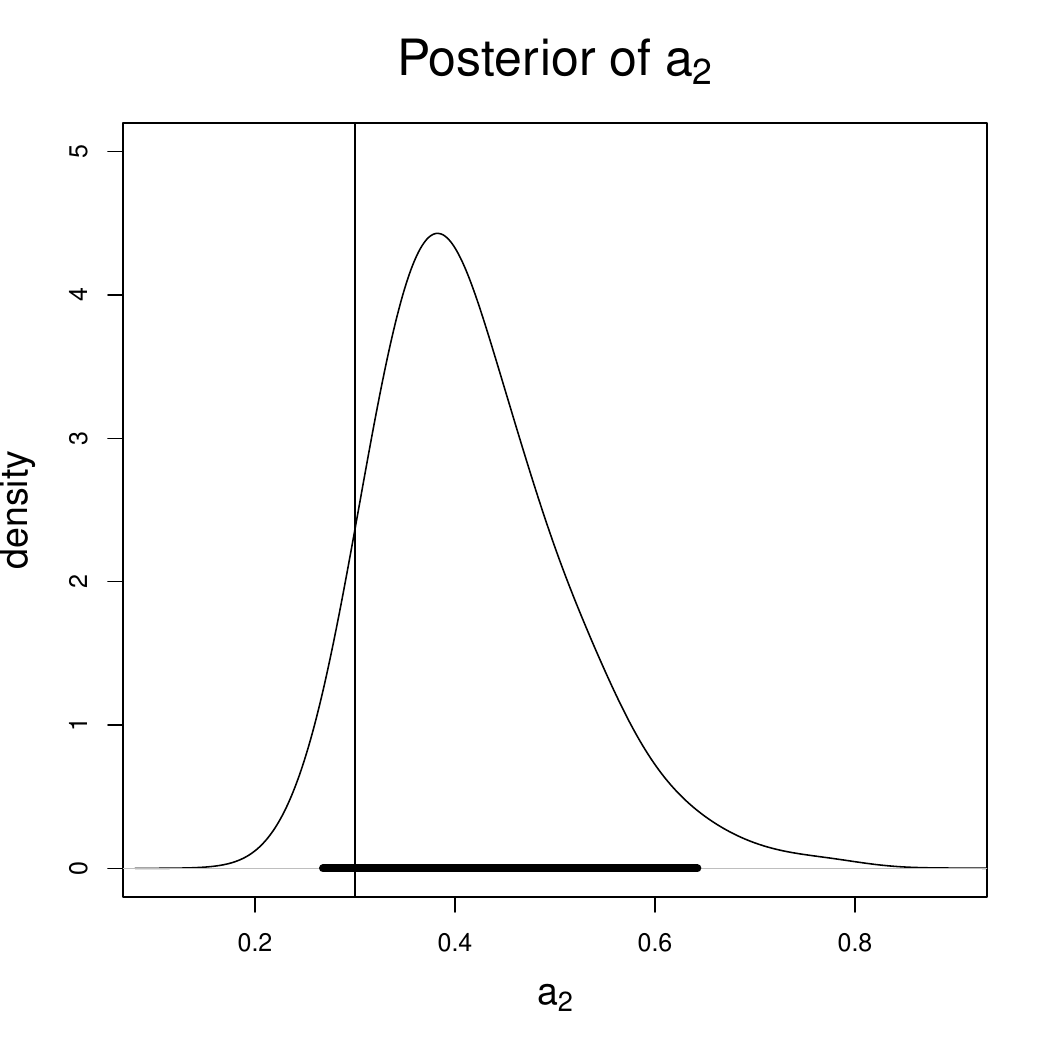}}
\caption{{\bf TTMCMC for $SDE_2$ and $\pi_2$:} Posteriors of $M$, $\mu_1$, $\mu_2$, $\omega^2_1$, $\omega^2_2$, $a_1$ and $a_2$. The vertical lines stand
for the true values, while the thick horizontal lines denote the 95\% credible intervals.} 
\label{fig:sim4_posterior_plots}
\end{figure}

In this case, the relevant diagrams provided in Figures \ref{fig:sim4_trace_plots} and \ref{fig:sim4_posterior_plots}, again indicate quite satisfactory performance. 
Here $M=2$ and $M=3$ received the posterior probabilities $0.8492$ and $0.1508$, respectively.
The overall acceptance rate in this case is about $4\%$.

\subsection{Fifth simulation study: $SDE_2$ with $\pi_3$} 
\label{subsec:sim5}

\begin{figure}
\centering
\subfigure[Trace plot of $M$.]{ \label{fig:sim5_trace_comp}
\includegraphics[width=4.5cm,height=5cm]{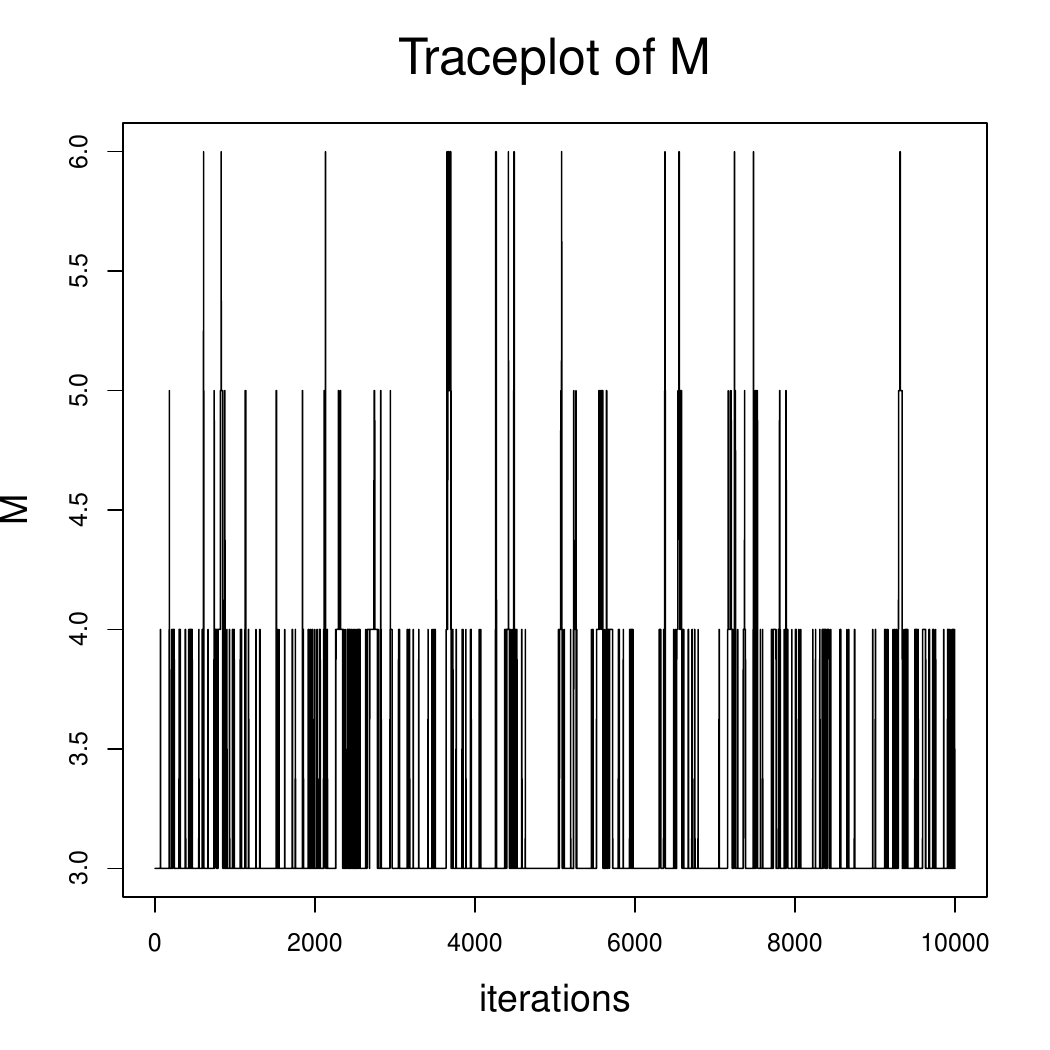}}
\hspace{2mm}
\subfigure[Trace plot of $\mu_1$.]{ \label{fig:sim5_trace_mu1}
\includegraphics[width=4.5cm,height=5cm]{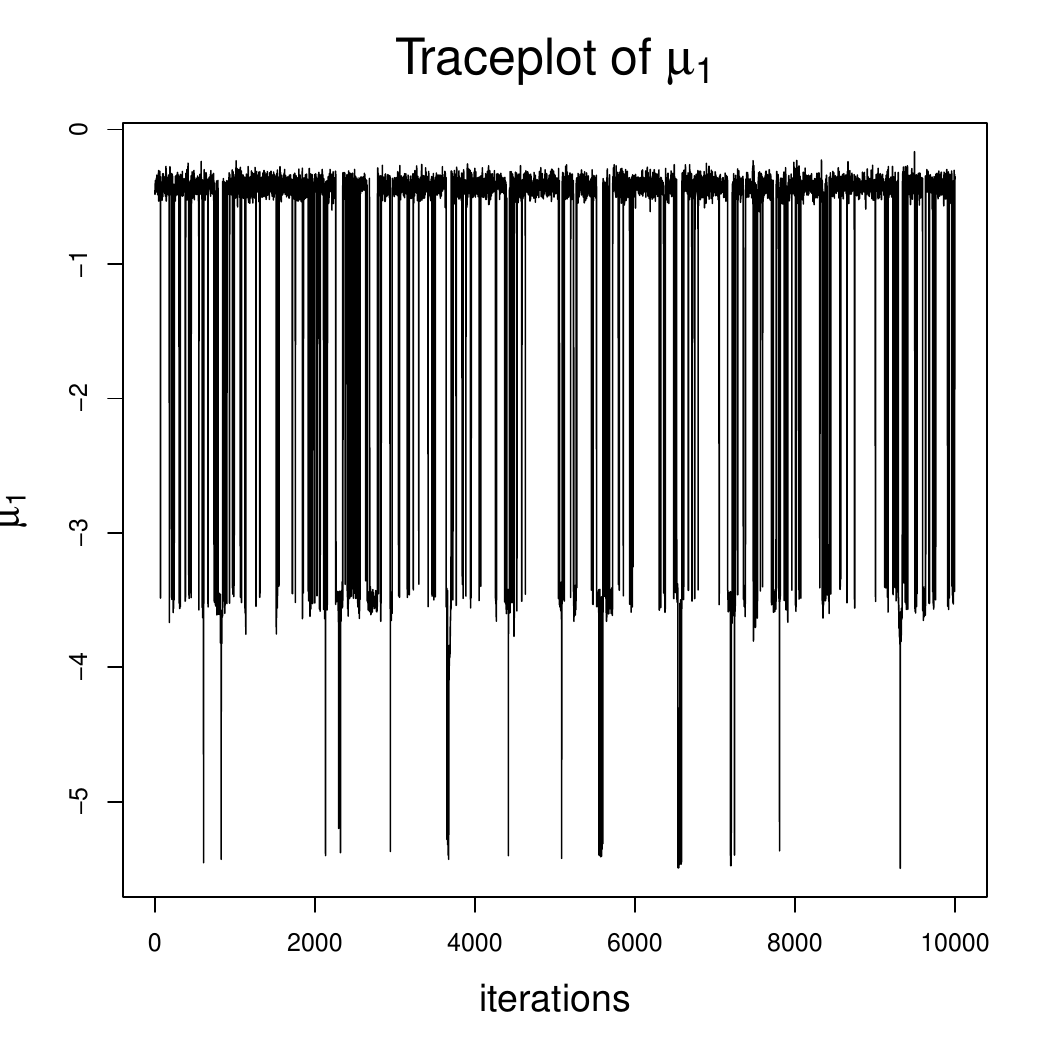}}
\hspace{2mm}
\subfigure[Trace plot of $\mu_2$.]{ \label{fig:sim5_trace_mu2}
\includegraphics[width=4.5cm,height=5cm]{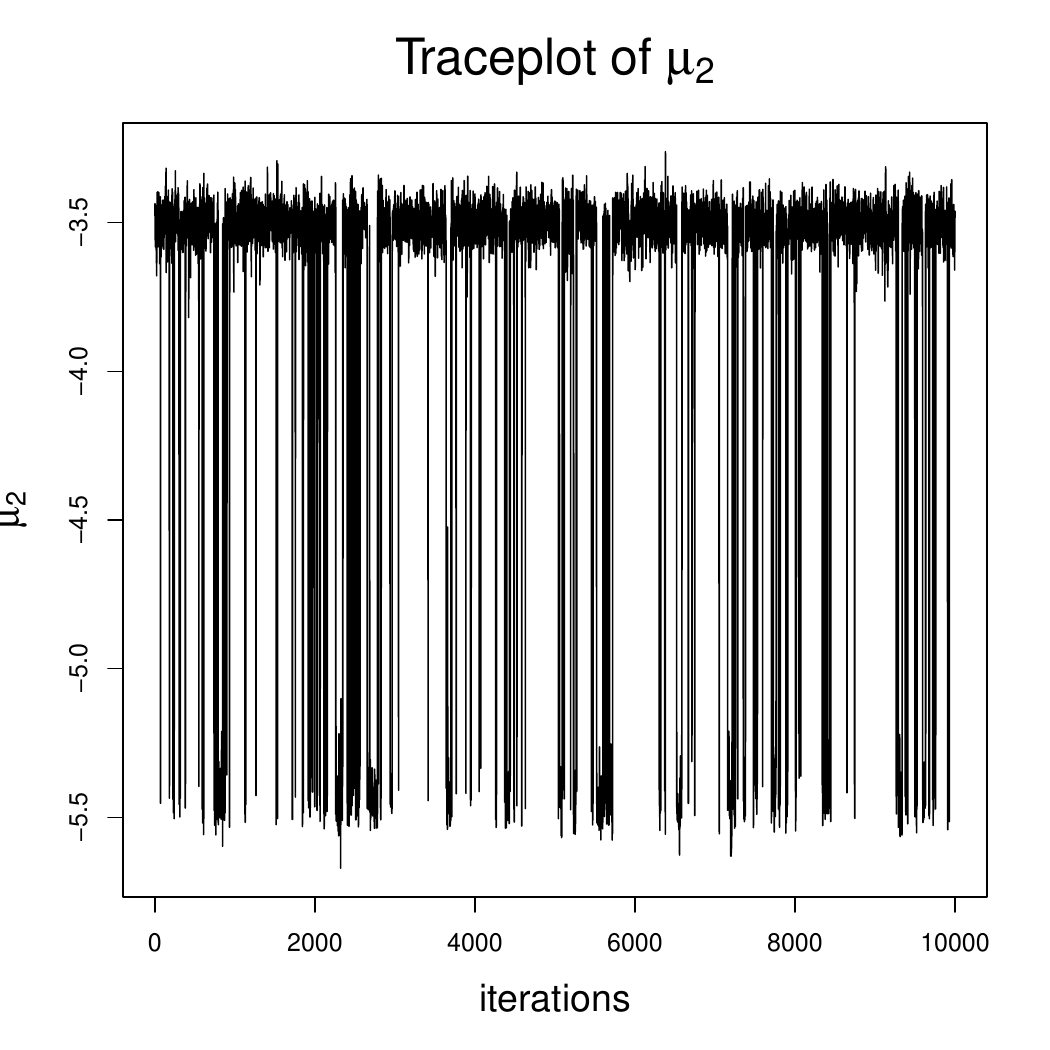}}\\
\vspace{2mm}
\subfigure[Trace plot of $\mu_3$.]{ \label{fig:sim5_trace_mu3}
\includegraphics[width=4.5cm,height=5cm]{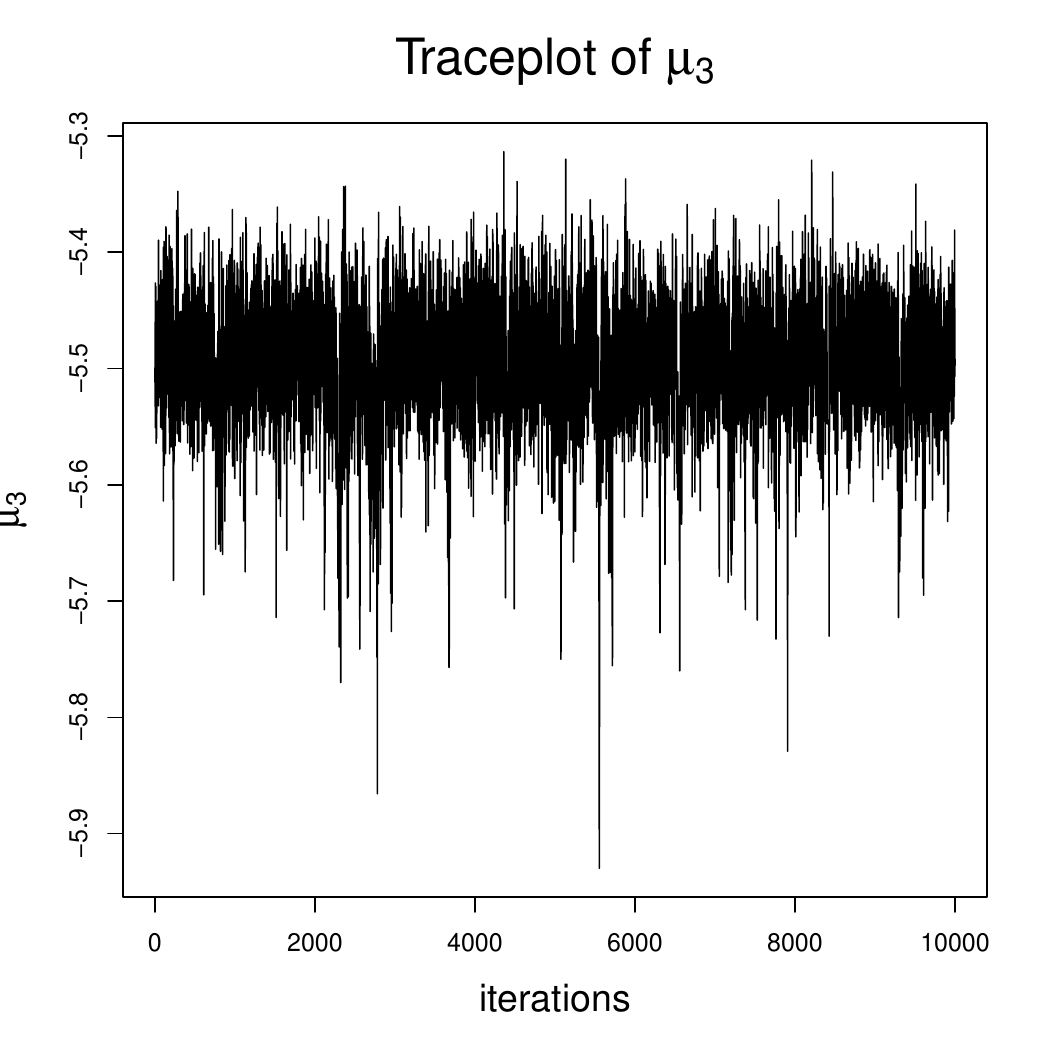}}
\hspace{2mm}
\subfigure[Trace plot of $\omega^2_1$.]{ \label{fig:sim5_trace_omegasq1}
\includegraphics[width=4.5cm,height=5cm]{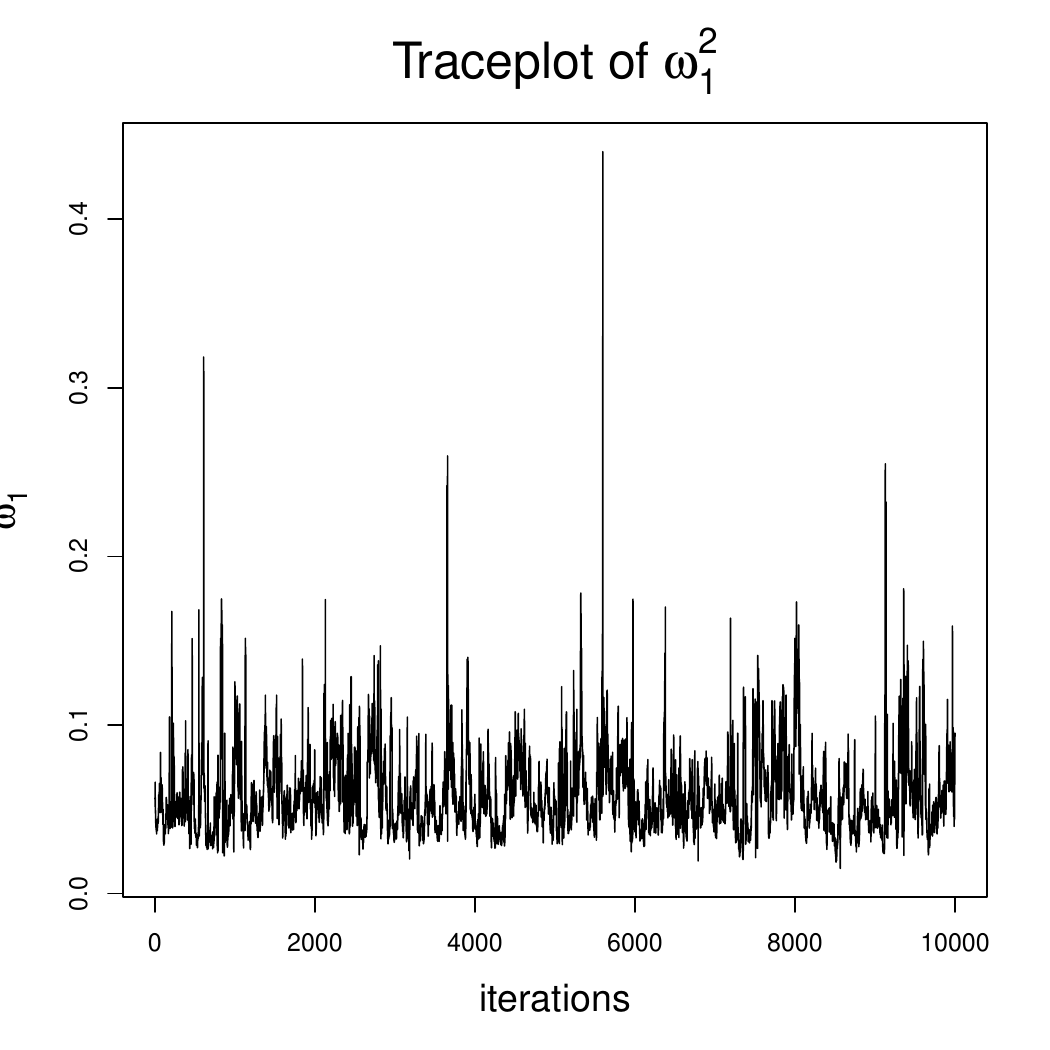}}
\hspace{2mm}
\subfigure[Trace plot of $\omega^2_2$.]{ \label{fig:sim5_trace_omegasq2}
\includegraphics[width=4.5cm,height=5cm]{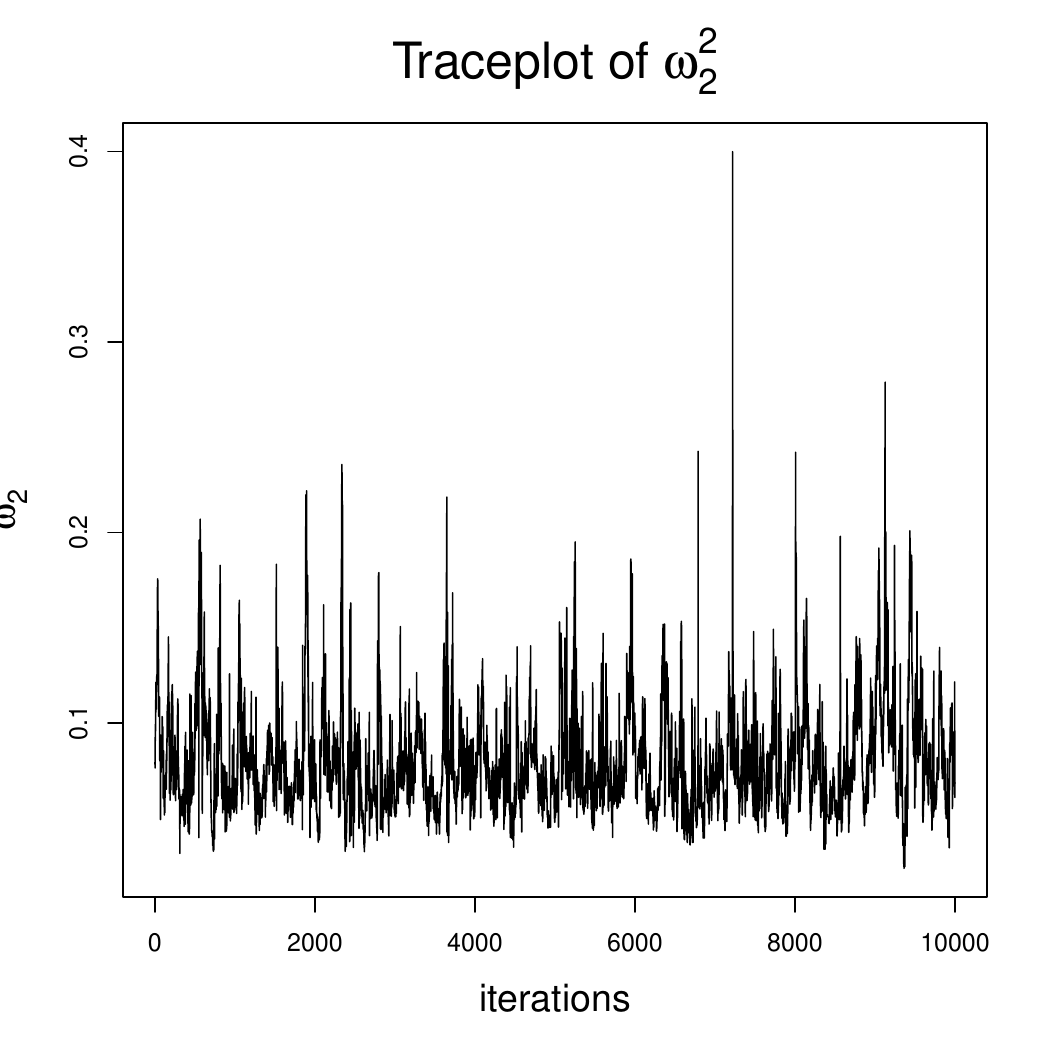}}\\
\vspace{2mm}
\subfigure[Trace plot of $\omega^2_3$.]{ \label{fig:sim5_trace_omegasq3}
\includegraphics[width=4.5cm,height=5cm]{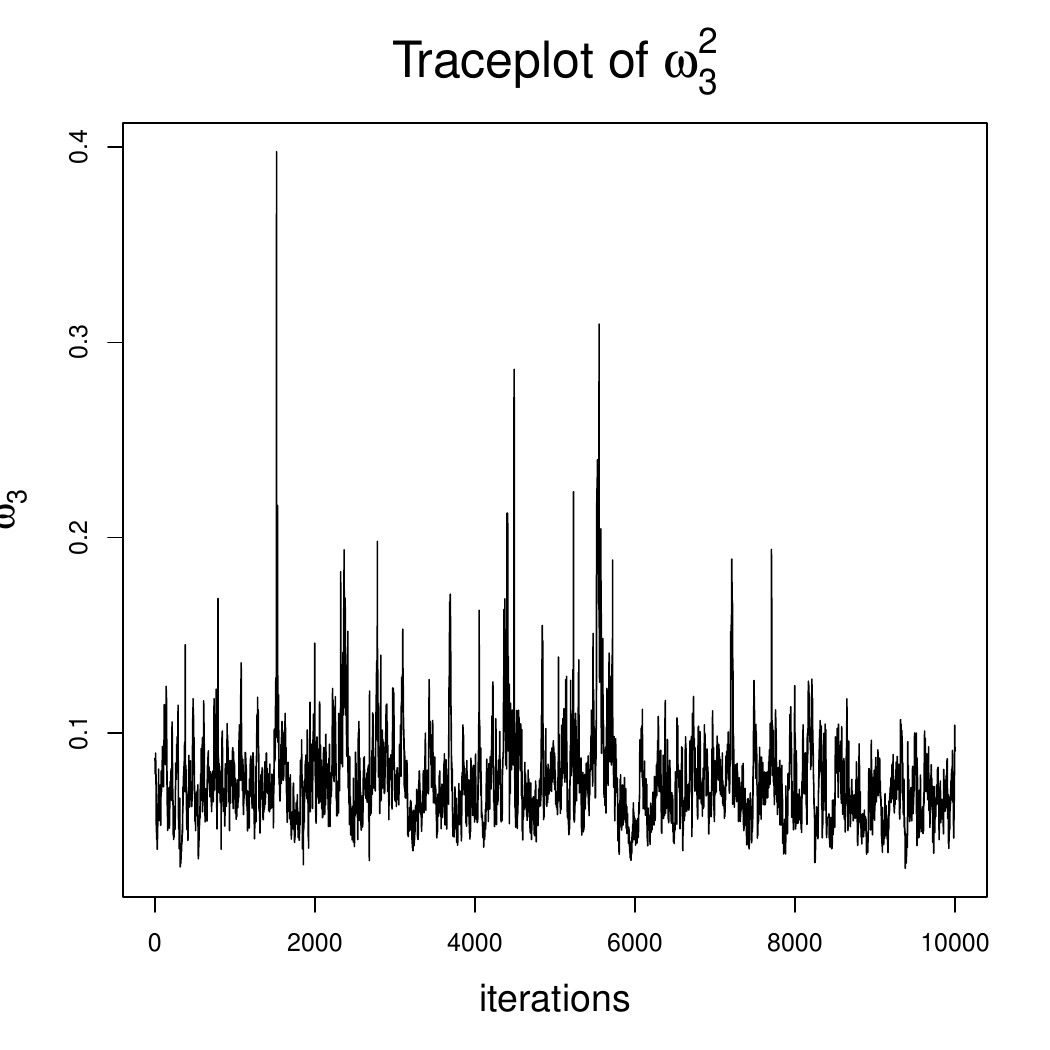}}
\hspace{2mm}
\subfigure[Trace plot of $a_1$.]{ \label{fig:sim5_trace_p1}
\includegraphics[width=4.5cm,height=5cm]{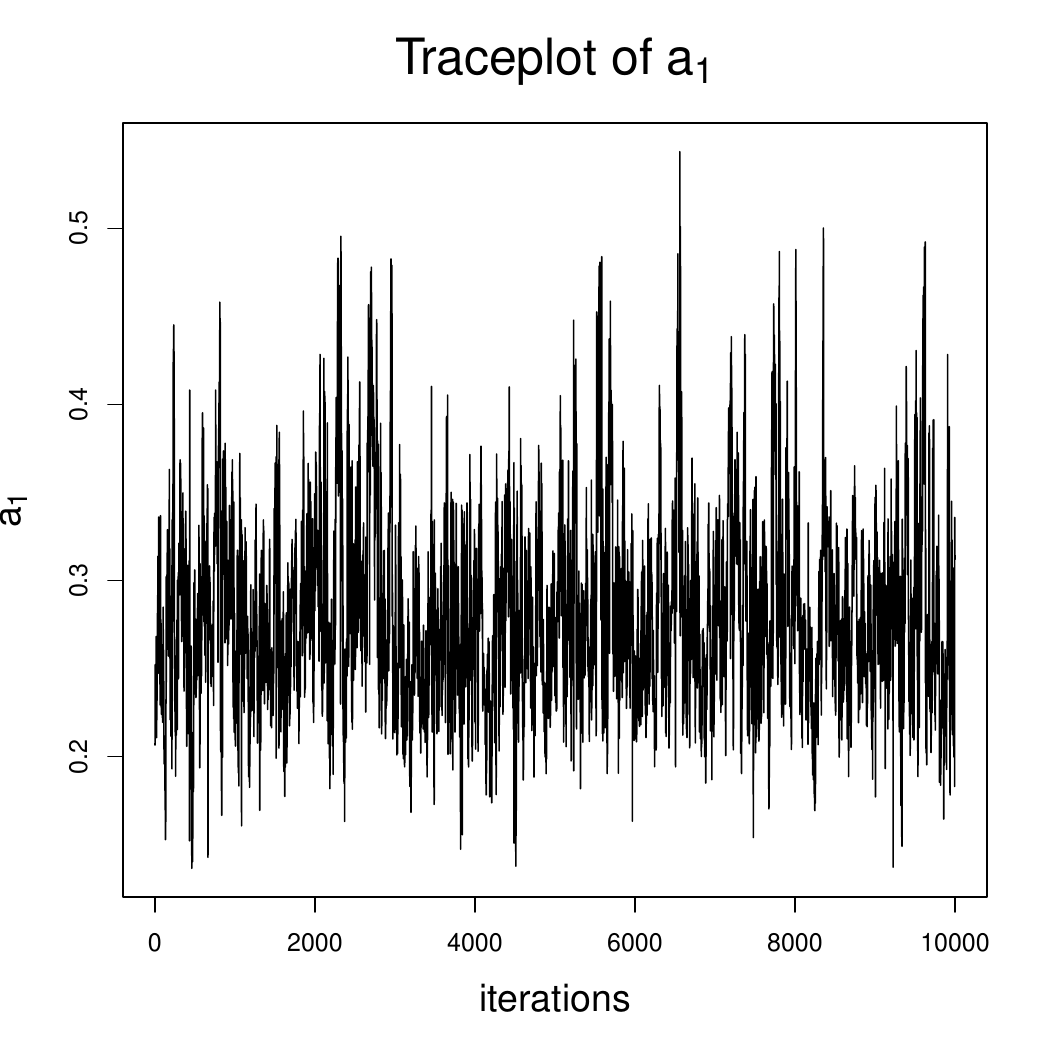}}
\hspace{2mm}
\subfigure[Trace plot of $a_2$.]{ \label{fig:sim5_trace_p2}
\includegraphics[width=4.5cm,height=5cm]{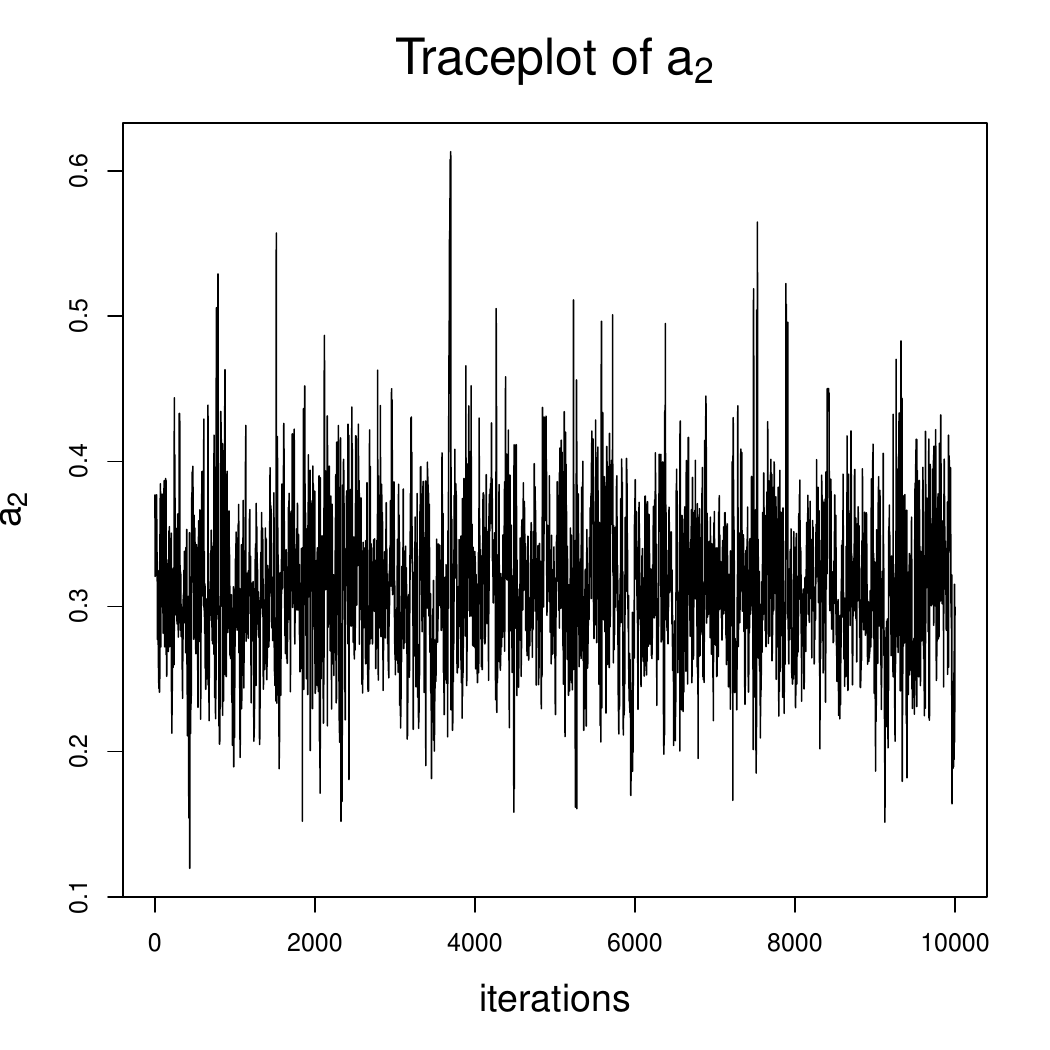}}\\
\vspace{2mm}
\subfigure[Trace plot of $a_3$.]{ \label{fig:sim5_trace_p3}
\includegraphics[width=4.5cm,height=5cm]{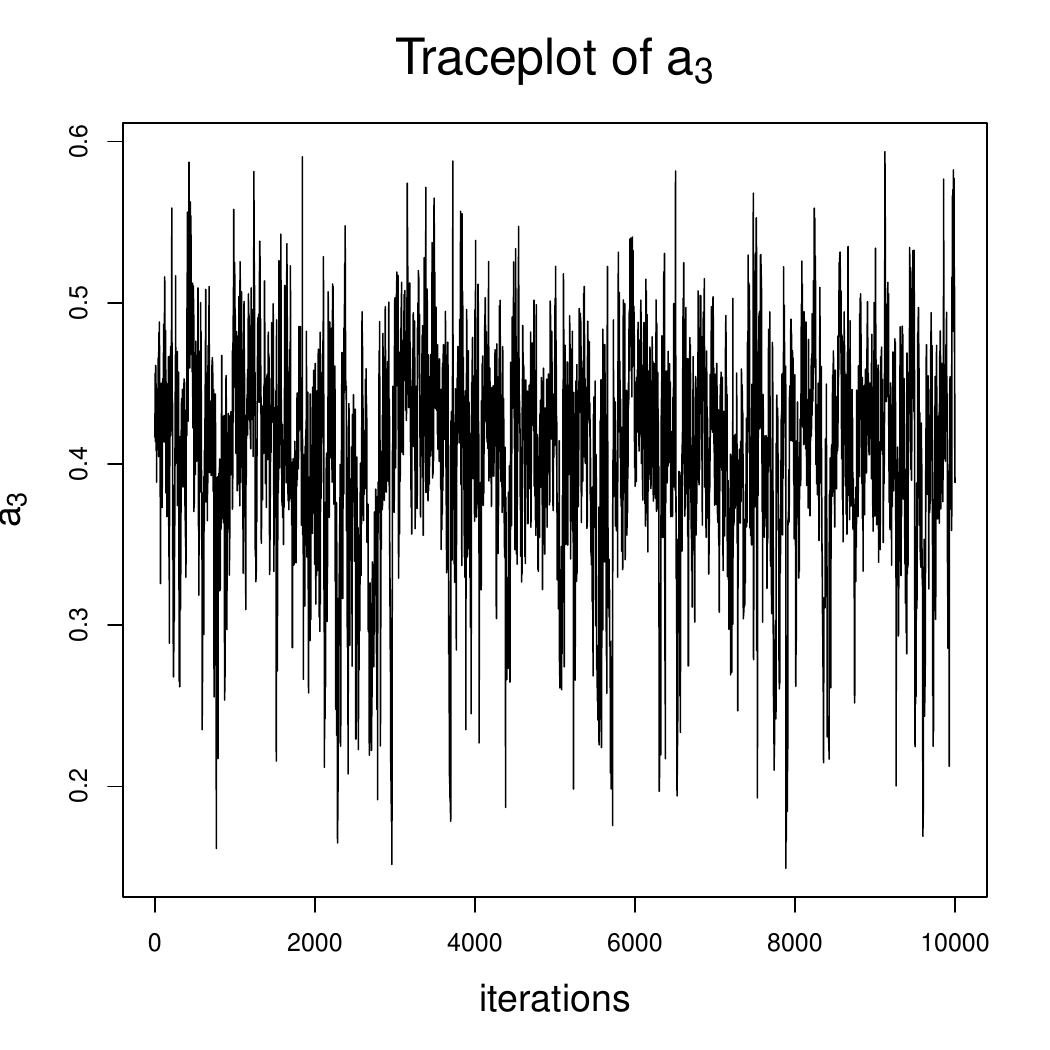}}
\caption{{\bf TTMCMC for $SDE_2$ and $\pi_3$:} Trace plots of $M$, $\mu_1$, $\mu_2$, $\nu_3$, $\omega^2_1$, $\omega^2_2$, $\omega^2_3$, $a_1$ $a_2$ and $a_3$.} 
\label{fig:sim5_trace_plots}
\end{figure}

\begin{figure}
\centering
\subfigure[Posterior of $\mu_1$.]{ \label{fig:sim5_mu1}
\includegraphics[width=4.5cm,height=5cm]{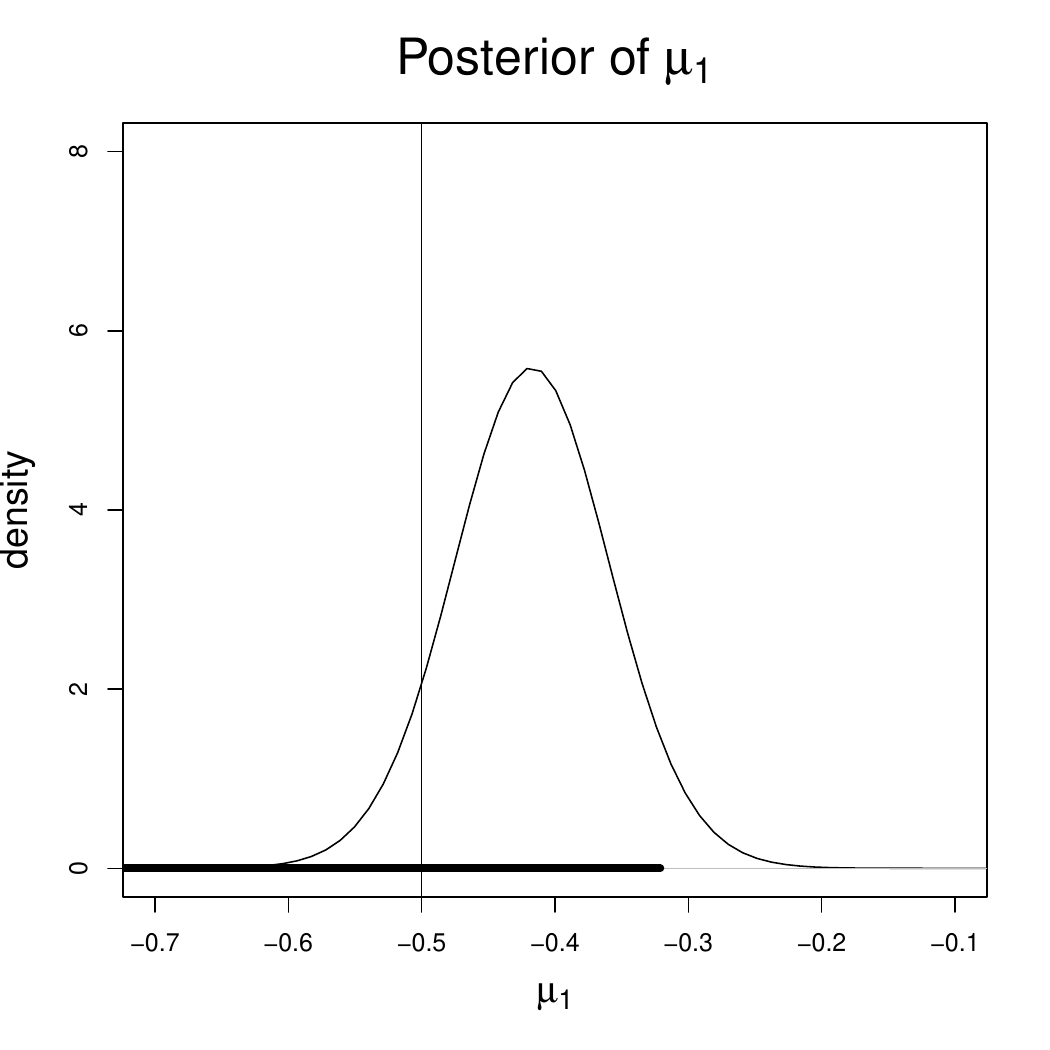}}
\hspace{2mm}
\subfigure[Posterior of $\mu_2$.]{ \label{fig:sim5_mu2}
\includegraphics[width=4.5cm,height=5cm]{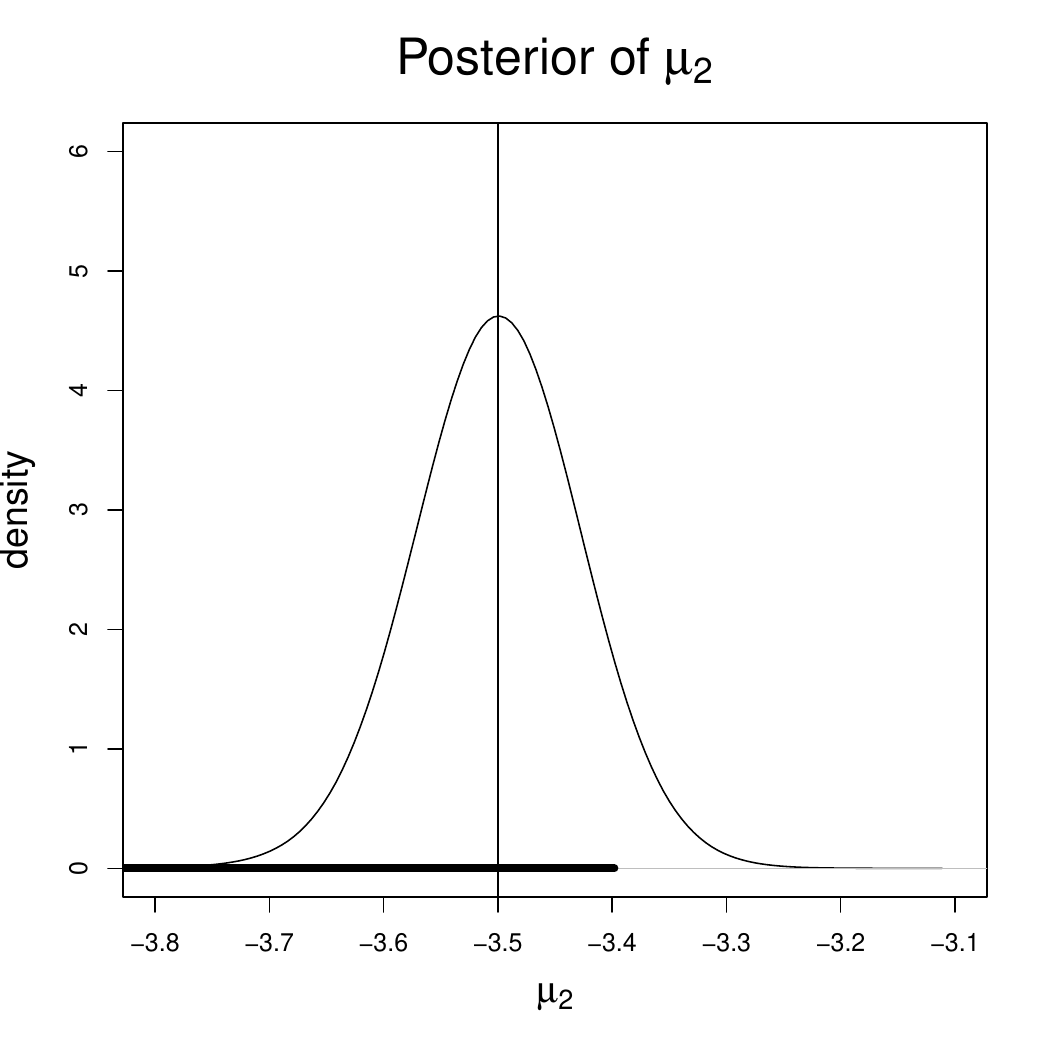}}
\hspace{2mm}
\subfigure[Posterior of $\mu_3$.]{ \label{fig:sim5_mu3}
\includegraphics[width=4.5cm,height=5cm]{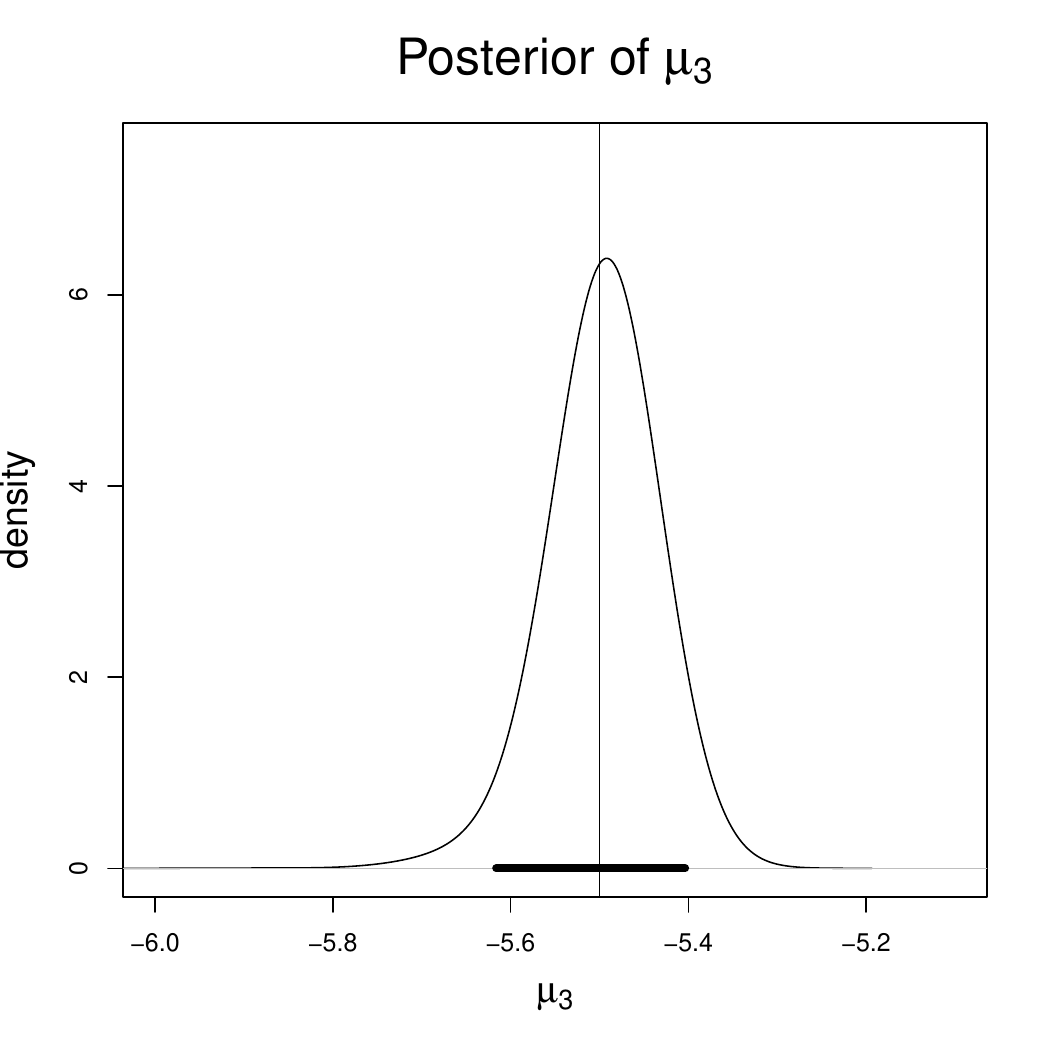}}\\
\vspace{2mm}
\subfigure[Posterior of $\omega^2_1$.]{ \label{fig:sim5_omegasq1}
\includegraphics[width=4.5cm,height=5cm]{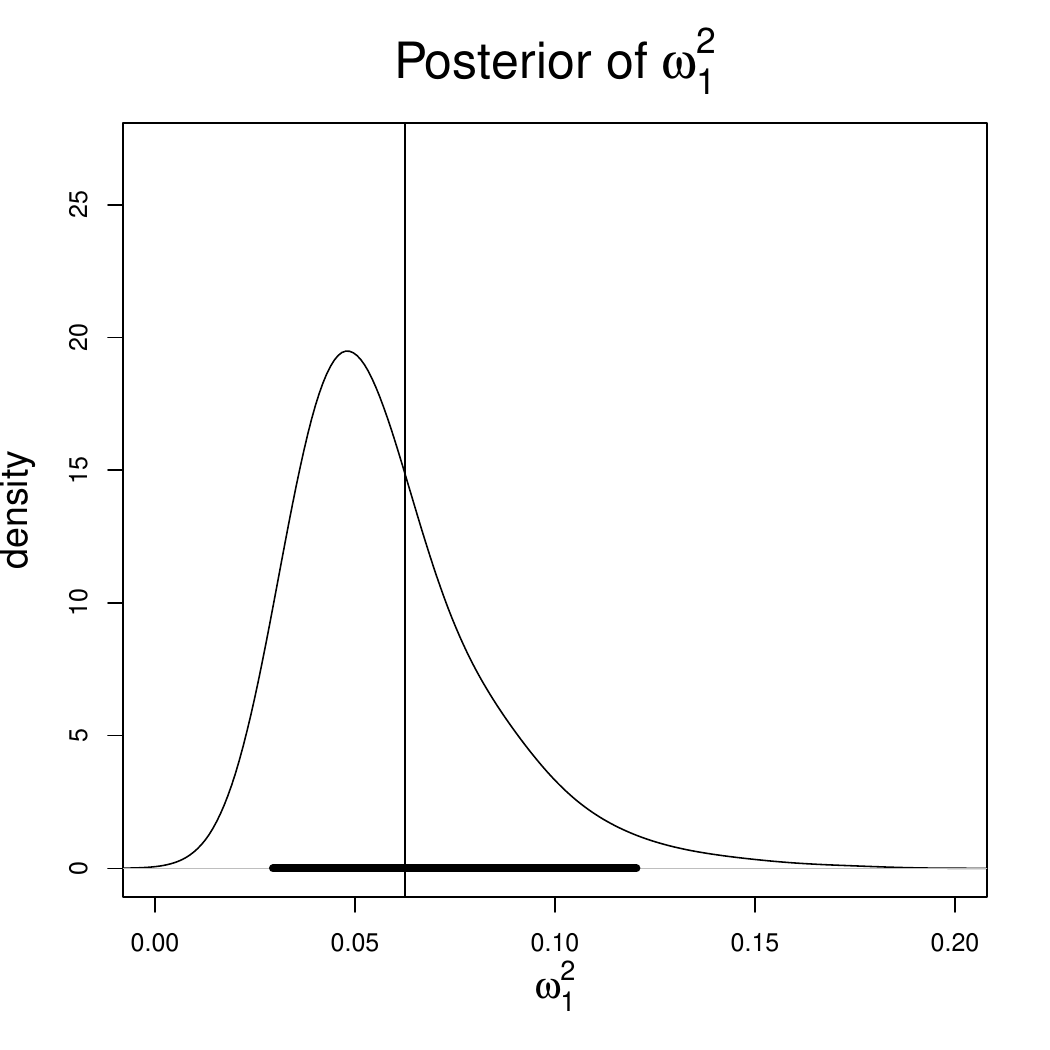}}
\hspace{2mm}
\subfigure[Posterior of $\omega^2_2$.]{ \label{fig:sim5_omegasq2}
\includegraphics[width=4.5cm,height=5cm]{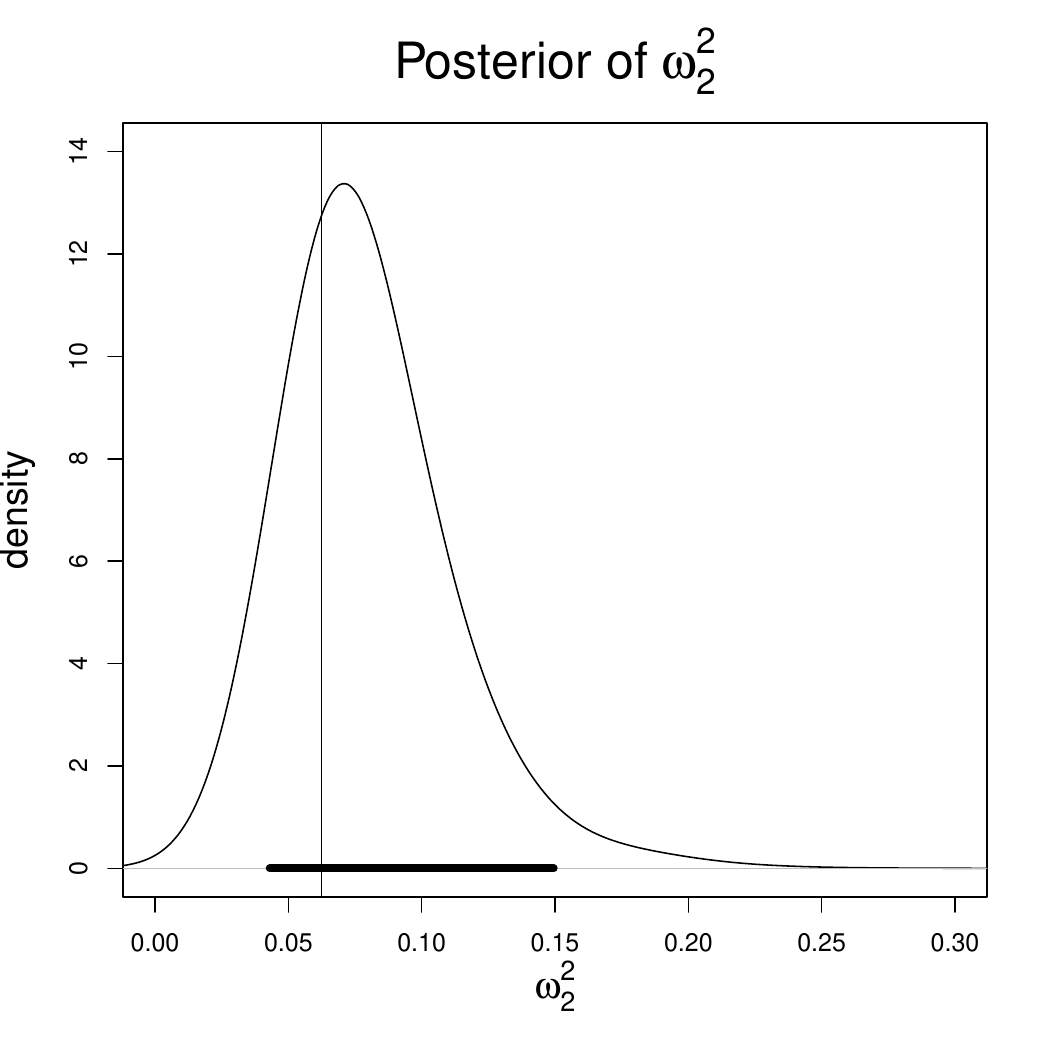}}
\hspace{2mm}
\subfigure[Posterior of $\omega^2_3$.]{ \label{fig:sim5_omegasq3}
\includegraphics[width=4.5cm,height=5cm]{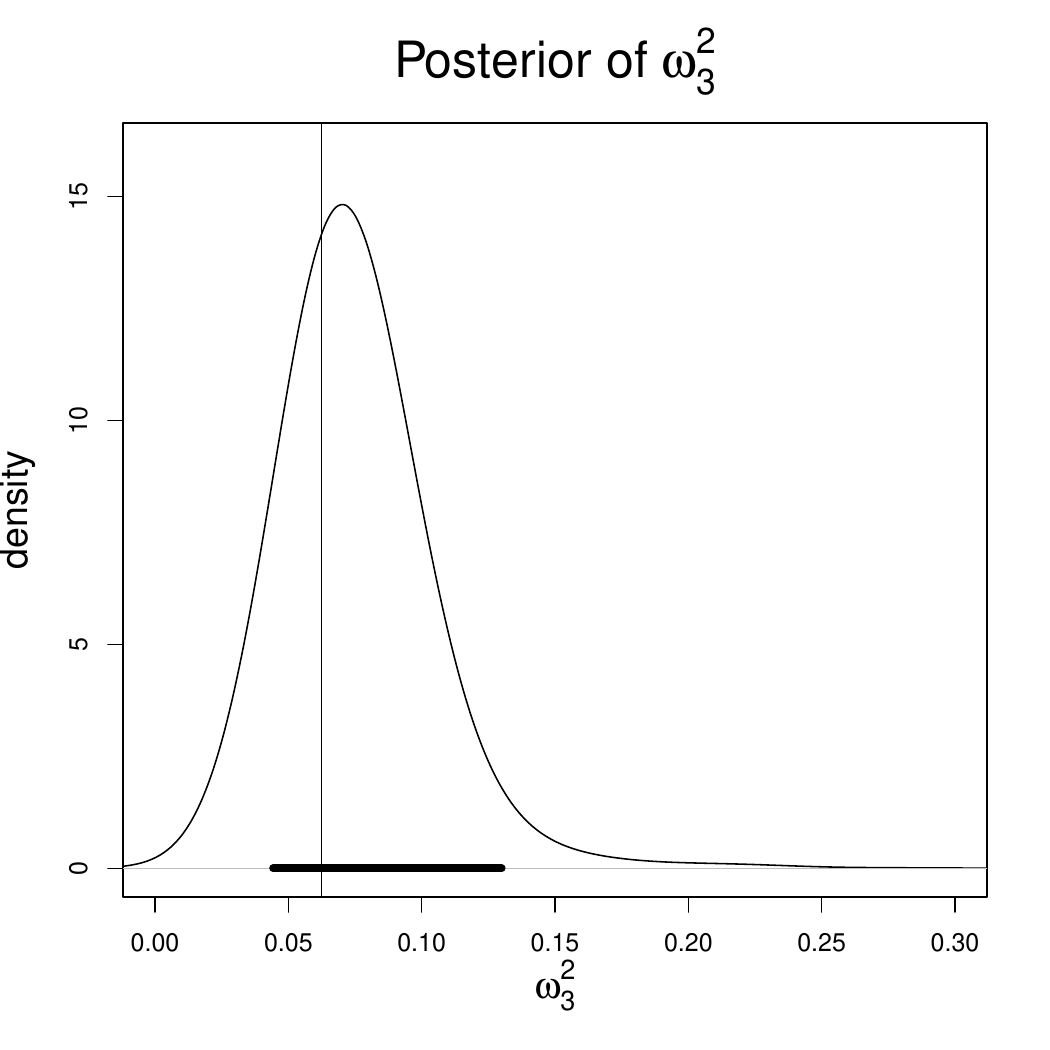}}\\
\vspace{2mm}
\subfigure[Posterior of $a_1$.]{ \label{fig:sim5_p1}
\includegraphics[width=4.5cm,height=5cm]{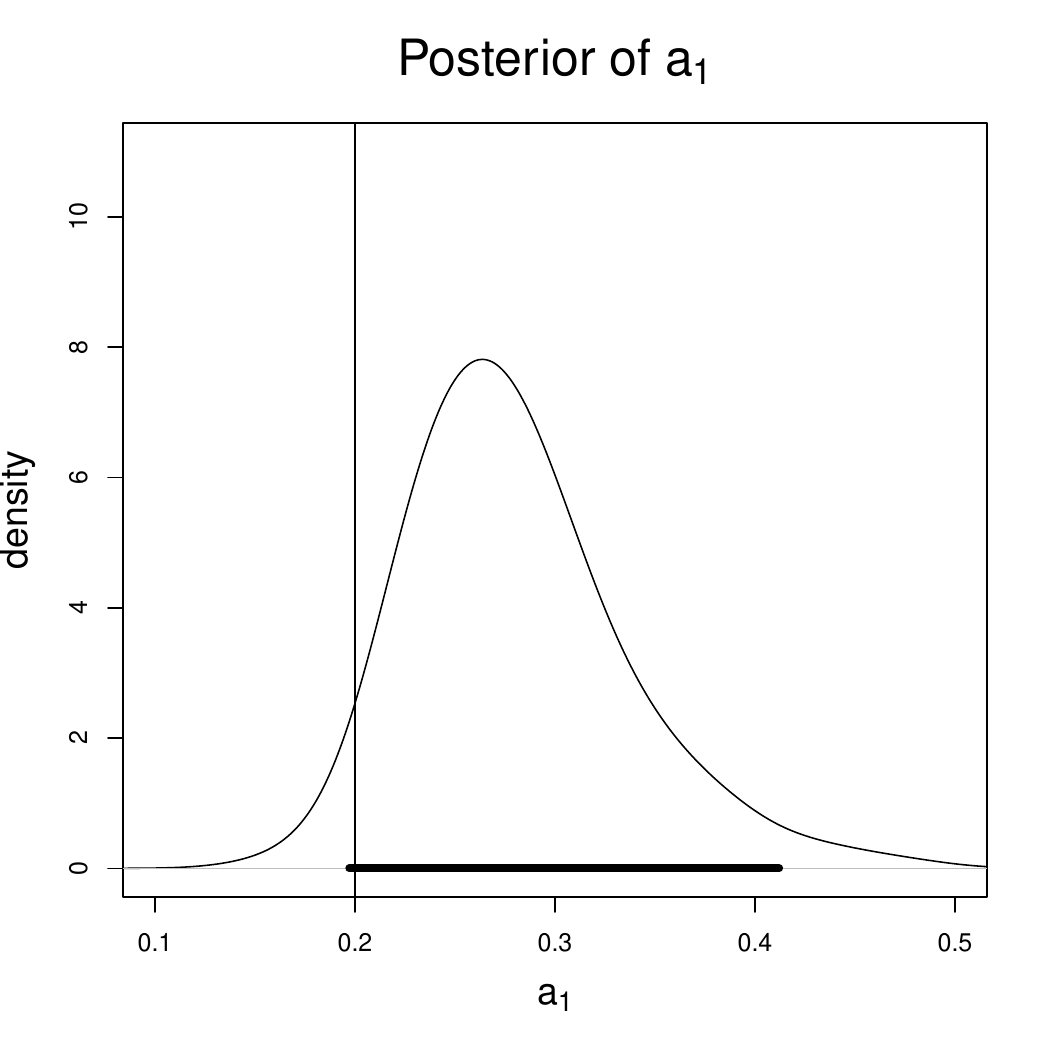}}
\hspace{2mm}
\subfigure[Posterior of $a_2$.]{ \label{fig:sim5_p2}
\includegraphics[width=4.5cm,height=5cm]{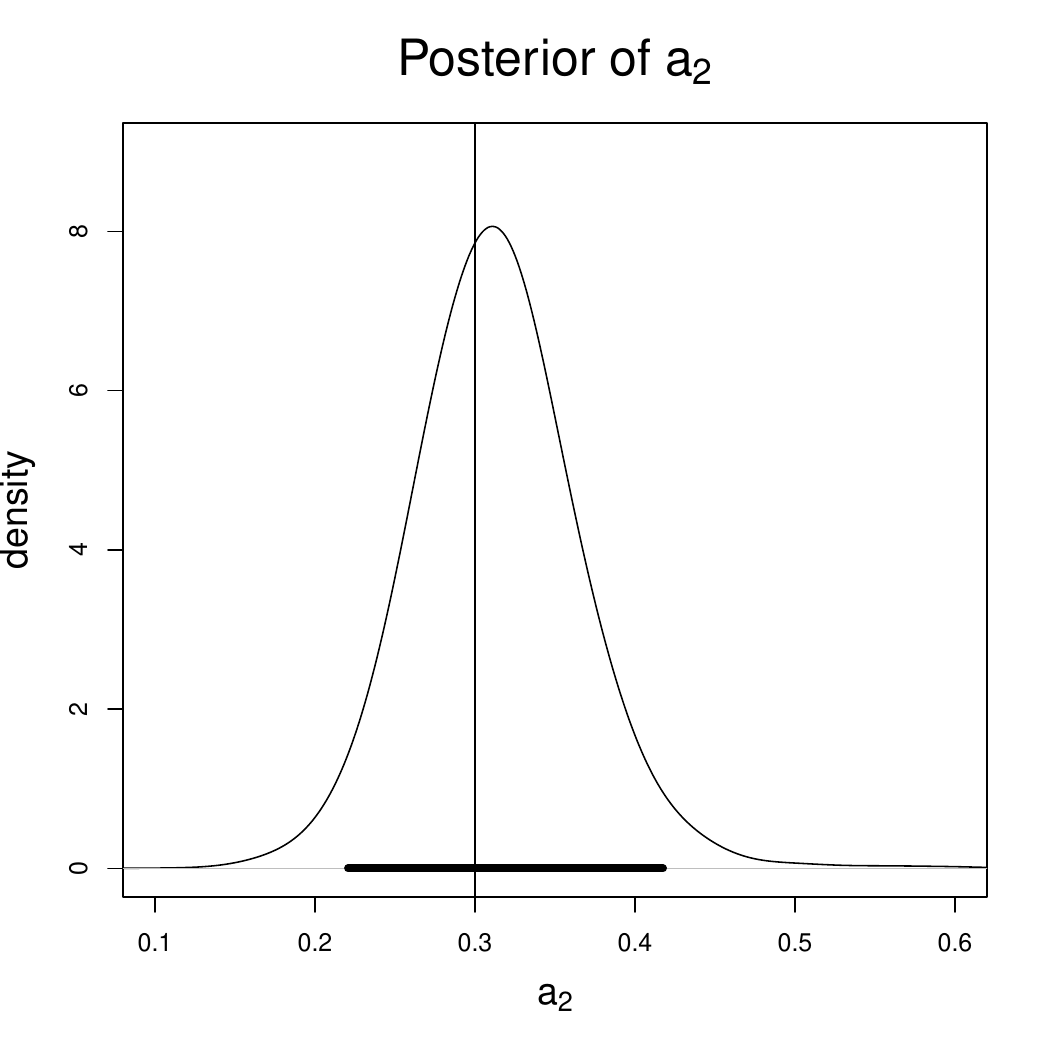}}
\hspace{2mm}
\subfigure[Posterior of $a_3$.]{ \label{fig:sim5_p3}
\includegraphics[width=4.5cm,height=5cm]{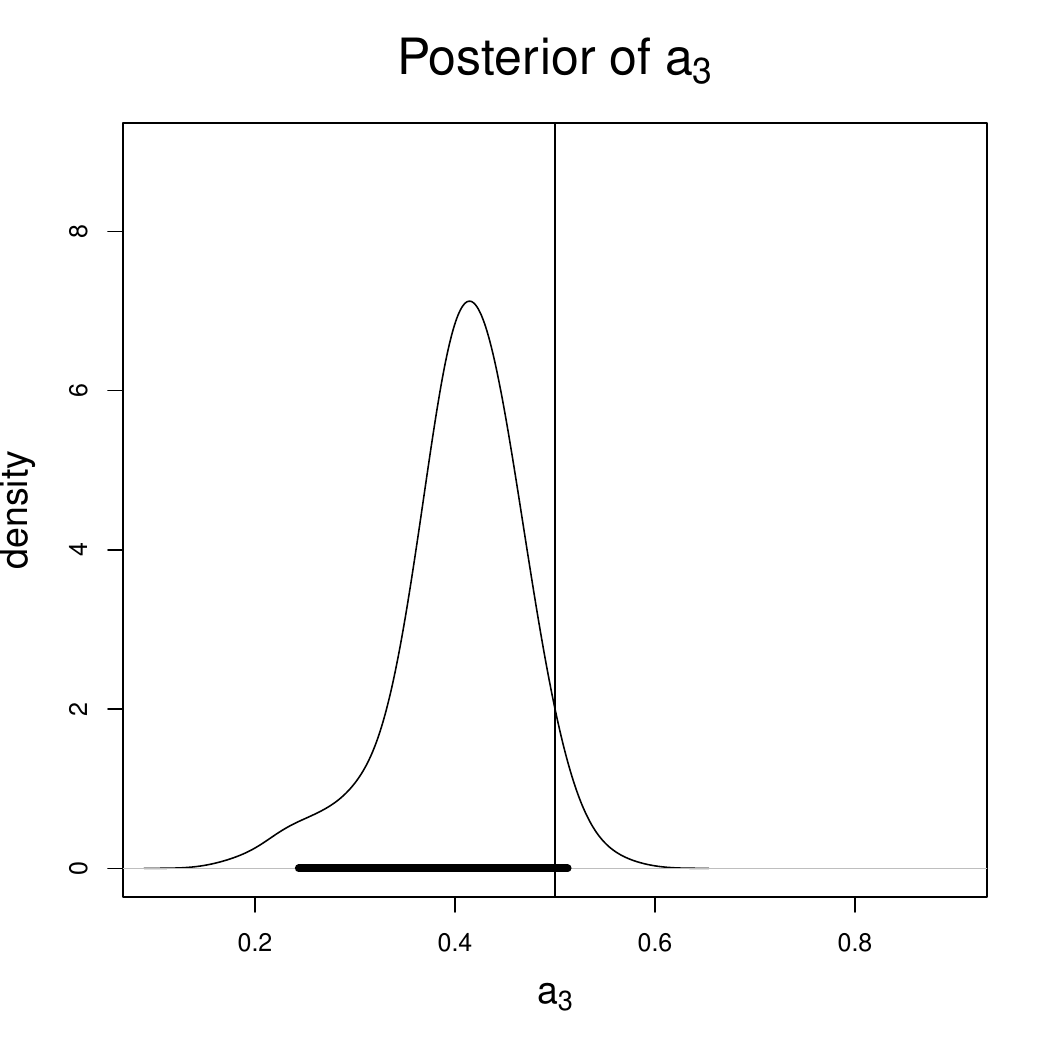}}
\caption{{\bf TTMCMC for $SDE_2$ and $\pi_3$:} Posteriors of $M$, $\mu_1$, $\mu_2$, $\mu_3$, $\omega^2_1$, $\omega^2_2$, $\omega^2_3$, $a_1$ $a_2$ and $a_3$. 
The vertical lines stand for the true values, while the thick horizontal lines denote the 95\% credible intervals.} 
\label{fig:sim5_posterior_plots}
\end{figure}

Figures \ref{fig:sim5_trace_plots} and \ref{fig:sim5_posterior_plots}, the relevant posterior plots associated with $SDE_2$ and $\pi_3$ indicate excellent
mixing as before and satisfactory performance. 
The overall acceptance rate here is about $2\%$.

In this case, $M=3$ and $M=4$ received the posterior probabilities $0.8228$ and $0.1418$, respectively, while all the other values received the zero posterior
probability.

\subsection{Sixth simulation study: $SDE_2$ with $\pi_4$} 
\label{subsec:sim6}

\begin{figure}
\centering
\subfigure[Trace plot of $M$.]{ \label{fig:sim6_trace_comp}
\includegraphics[width=4.5cm,height=5cm]{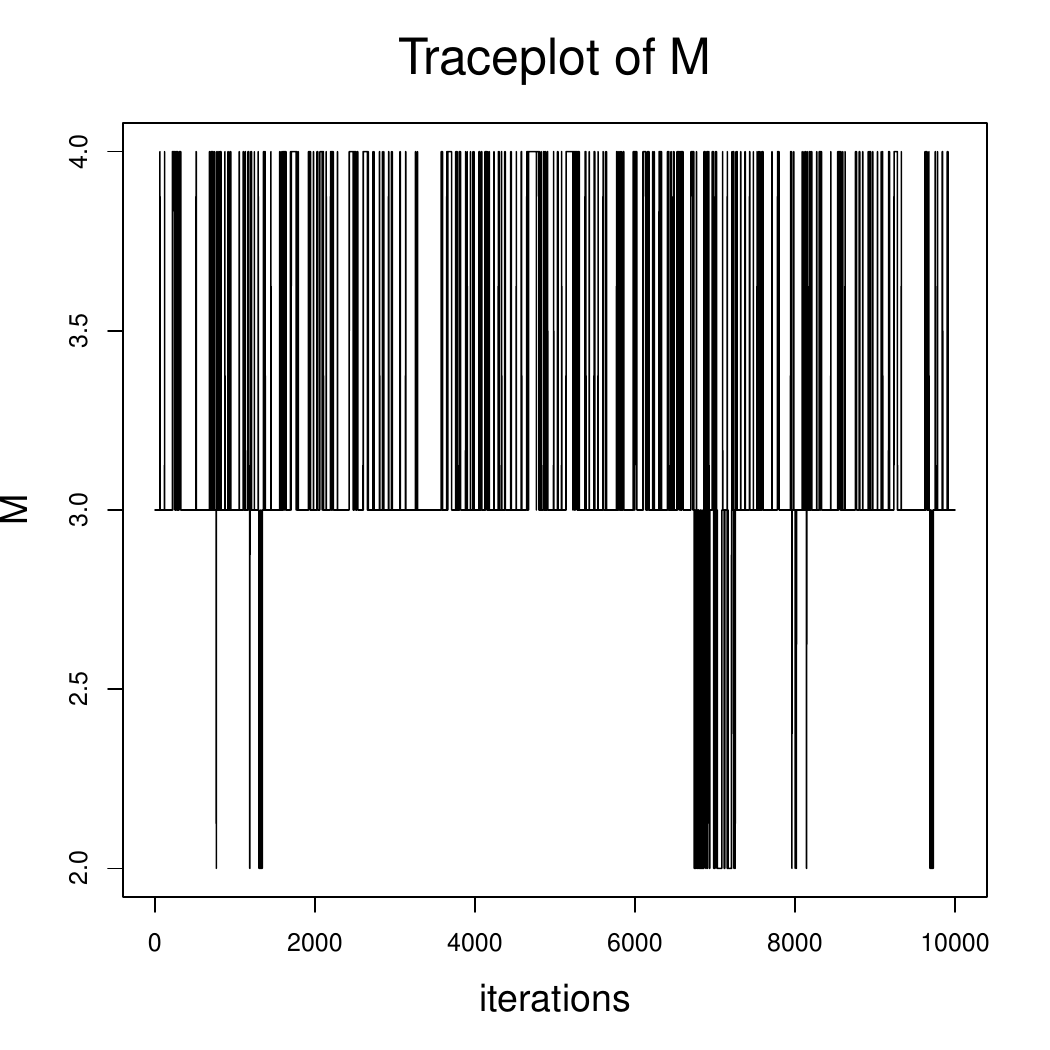}}
\hspace{2mm}
\subfigure[Trace plot of $\mu_1$.]{ \label{fig:sim6_trace_mu1}
\includegraphics[width=4.5cm,height=5cm]{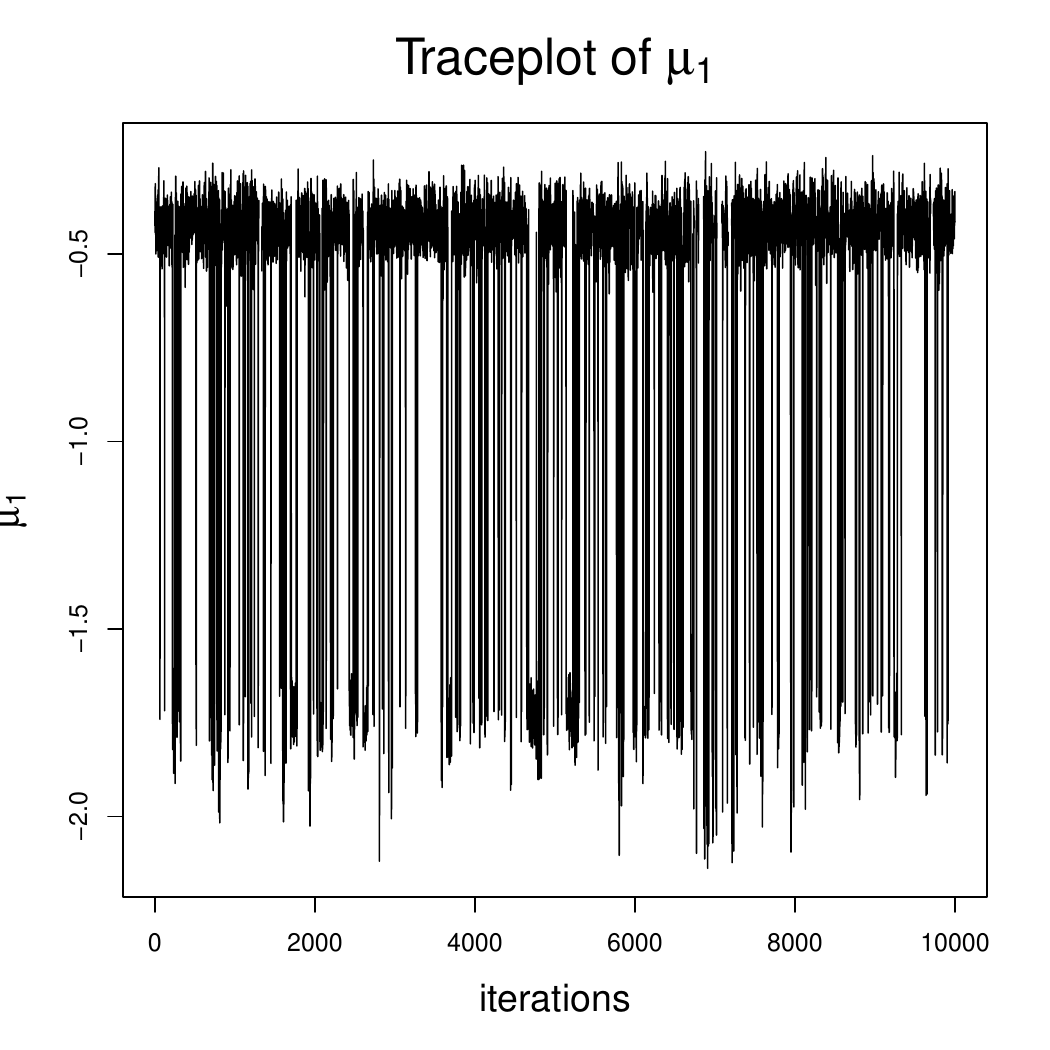}}
\hspace{2mm}
\subfigure[Trace plot of $\mu_2$.]{ \label{fig:sim6_trace_mu2}
\includegraphics[width=4.5cm,height=5cm]{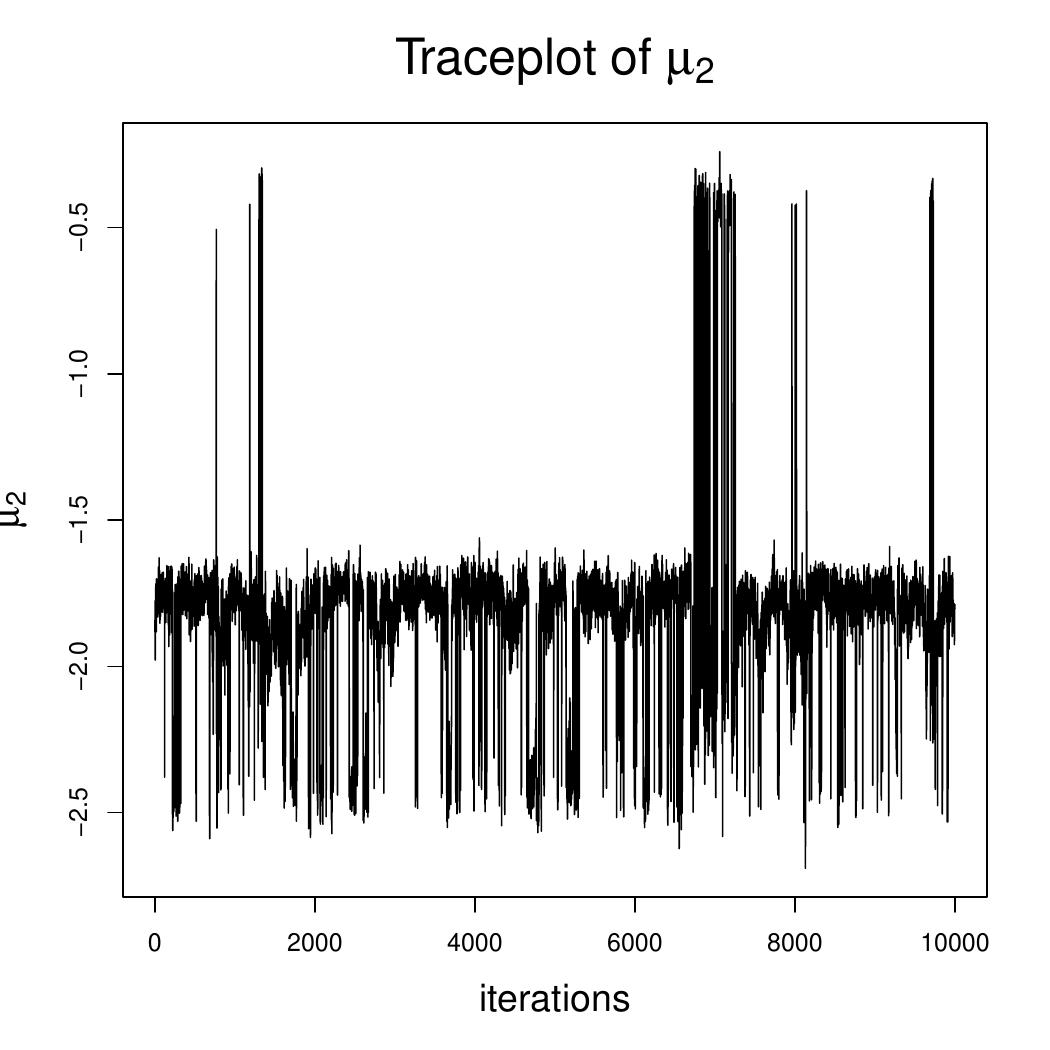}}\\
\vspace{2mm}
\subfigure[Trace plot of $\mu_3$.]{ \label{fig:sim6_trace_mu3}
\includegraphics[width=4.5cm,height=5cm]{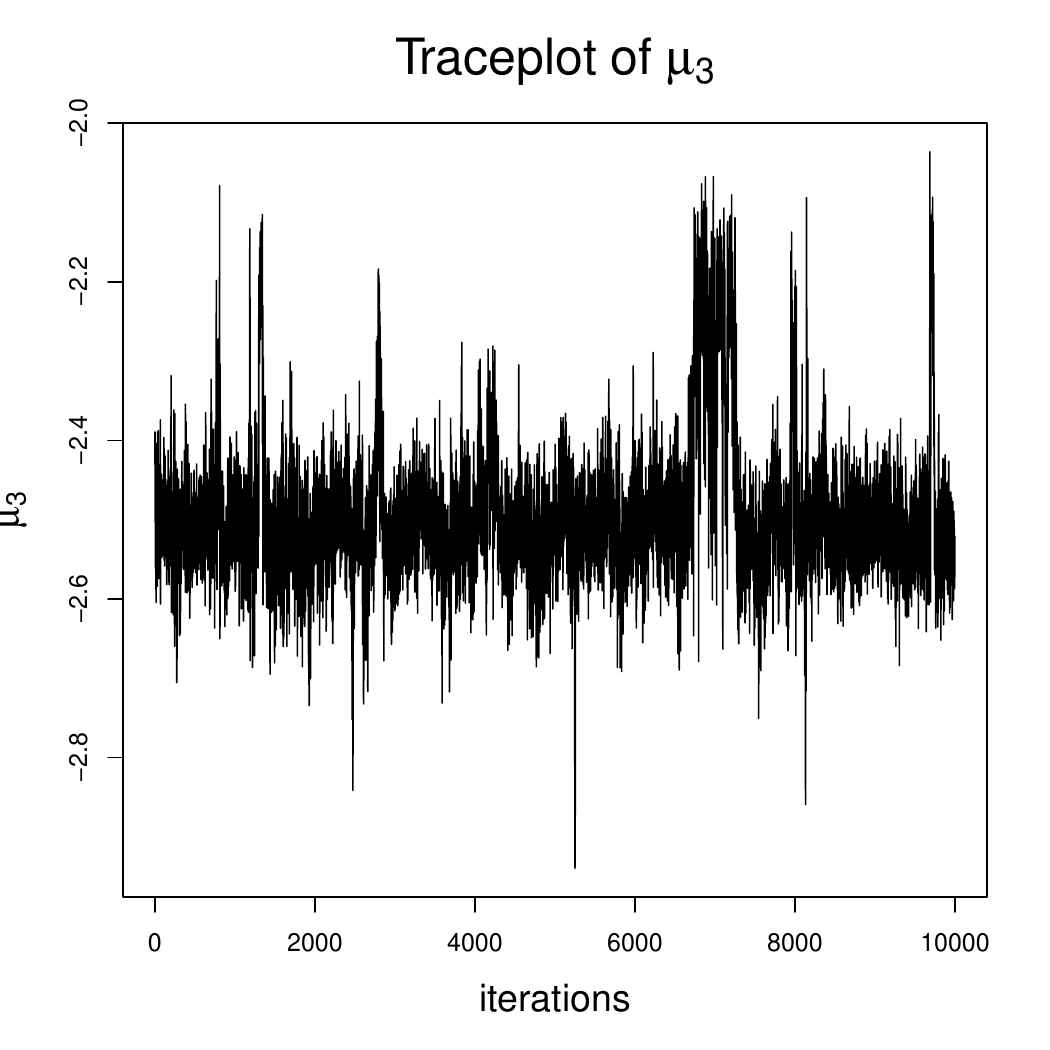}}
\hspace{2mm}
\subfigure[Trace plot of $\omega^2_1$.]{ \label{fig:sim6_trace_omegasq1}
\includegraphics[width=4.5cm,height=5cm]{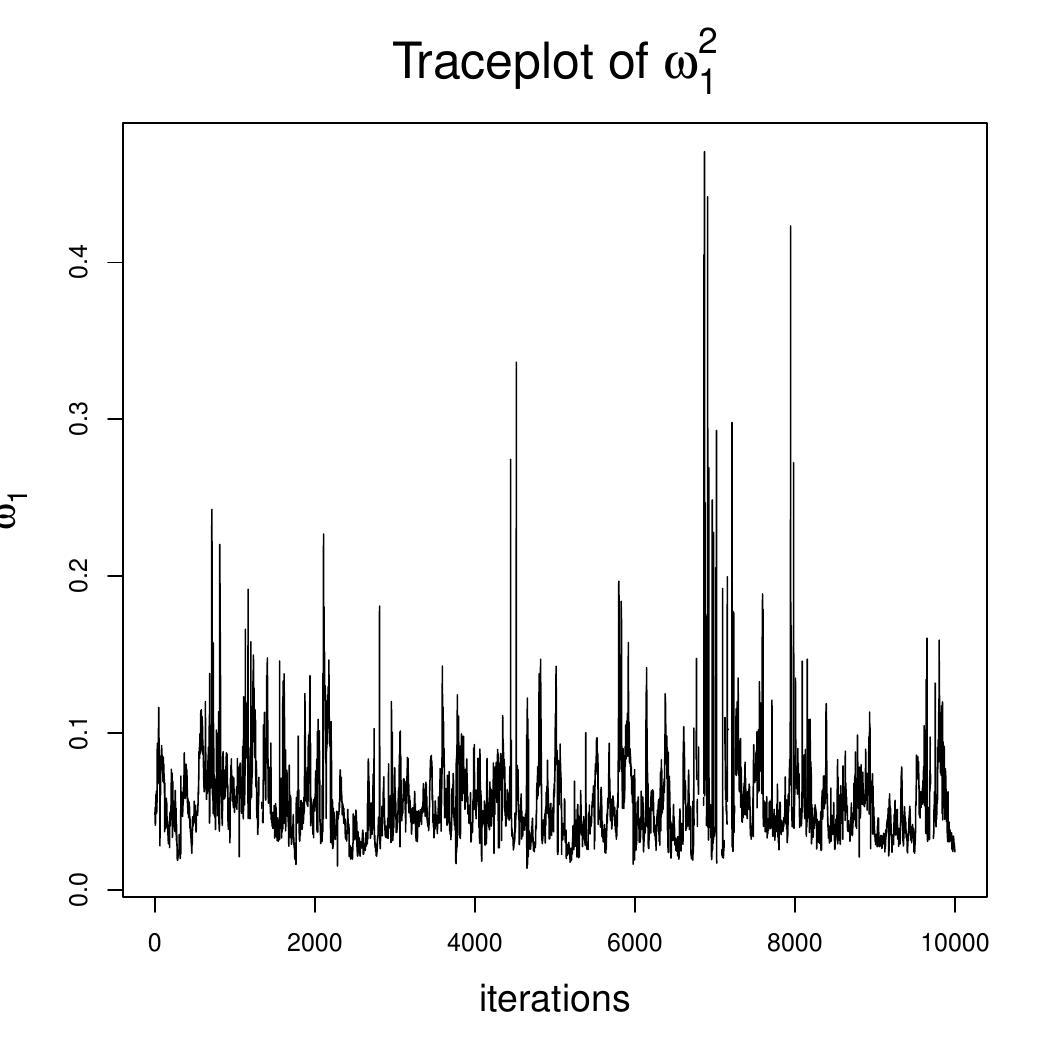}}
\hspace{2mm}
\subfigure[Trace plot of $\omega^2_2$.]{ \label{fig:sim6_trace_omegasq2}
\includegraphics[width=4.5cm,height=5cm]{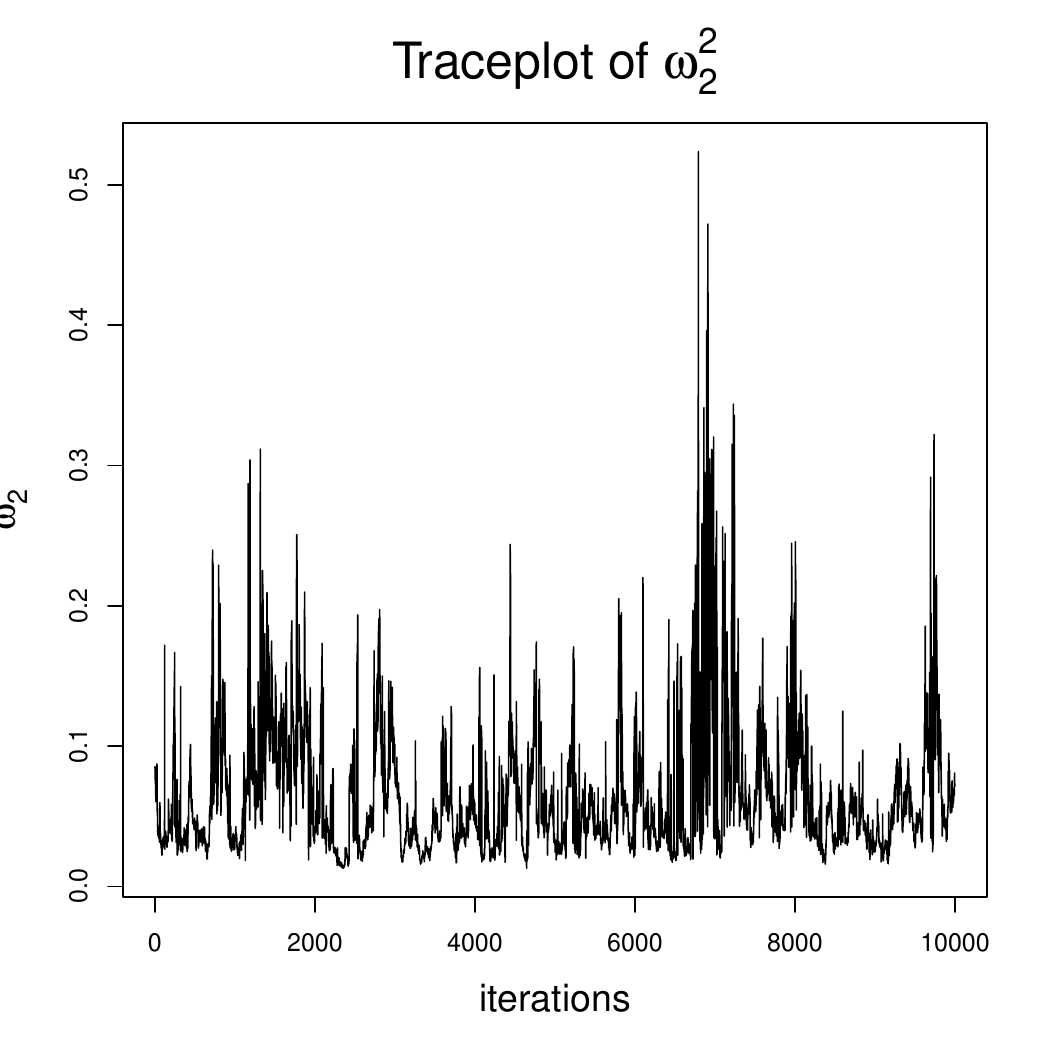}}\\
\vspace{2mm}
\subfigure[Trace plot of $\omega^2_3$.]{ \label{fig:sim6_trace_omegasq3}
\includegraphics[width=4.5cm,height=5cm]{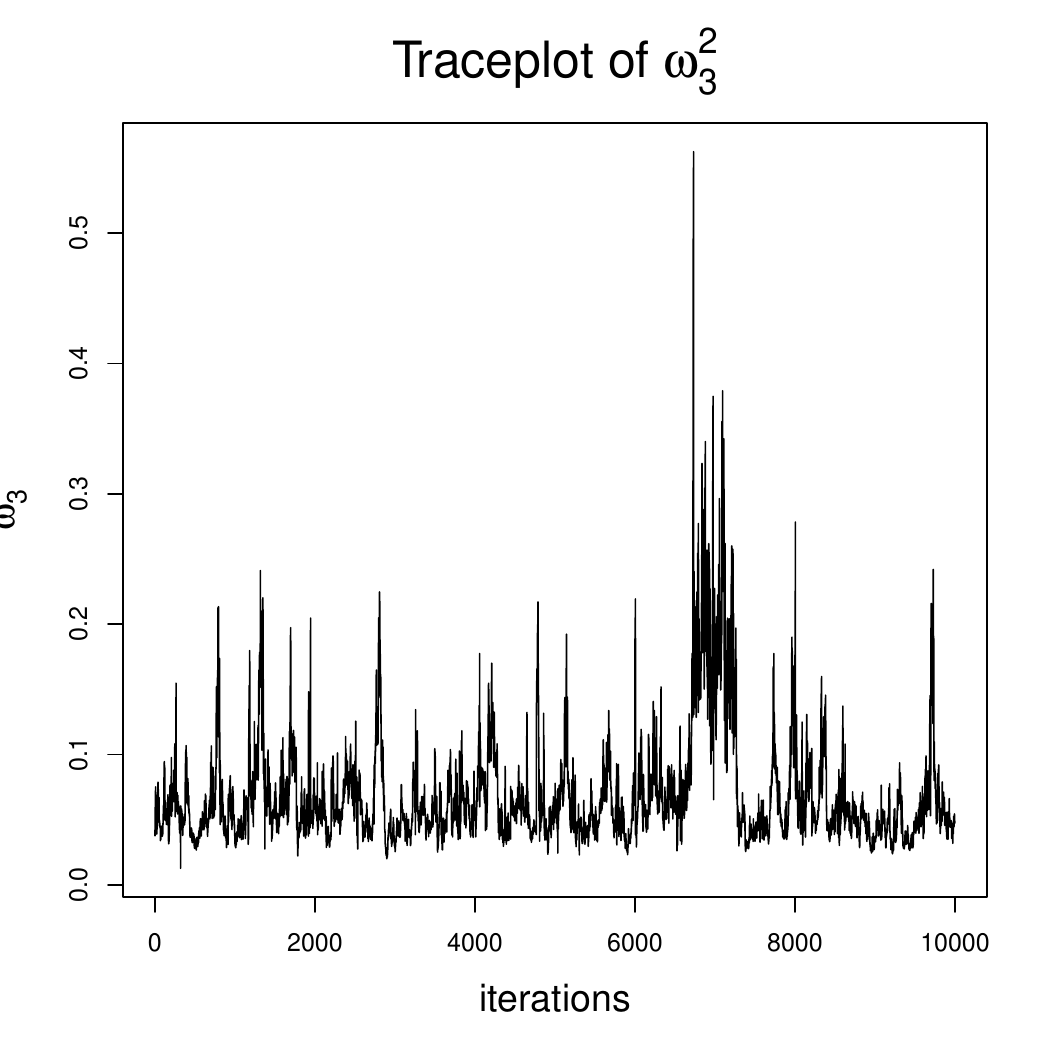}}
\hspace{2mm}
\subfigure[Trace plot of $a_1$.]{ \label{fig:sim6_trace_p1}
\includegraphics[width=4.5cm,height=5cm]{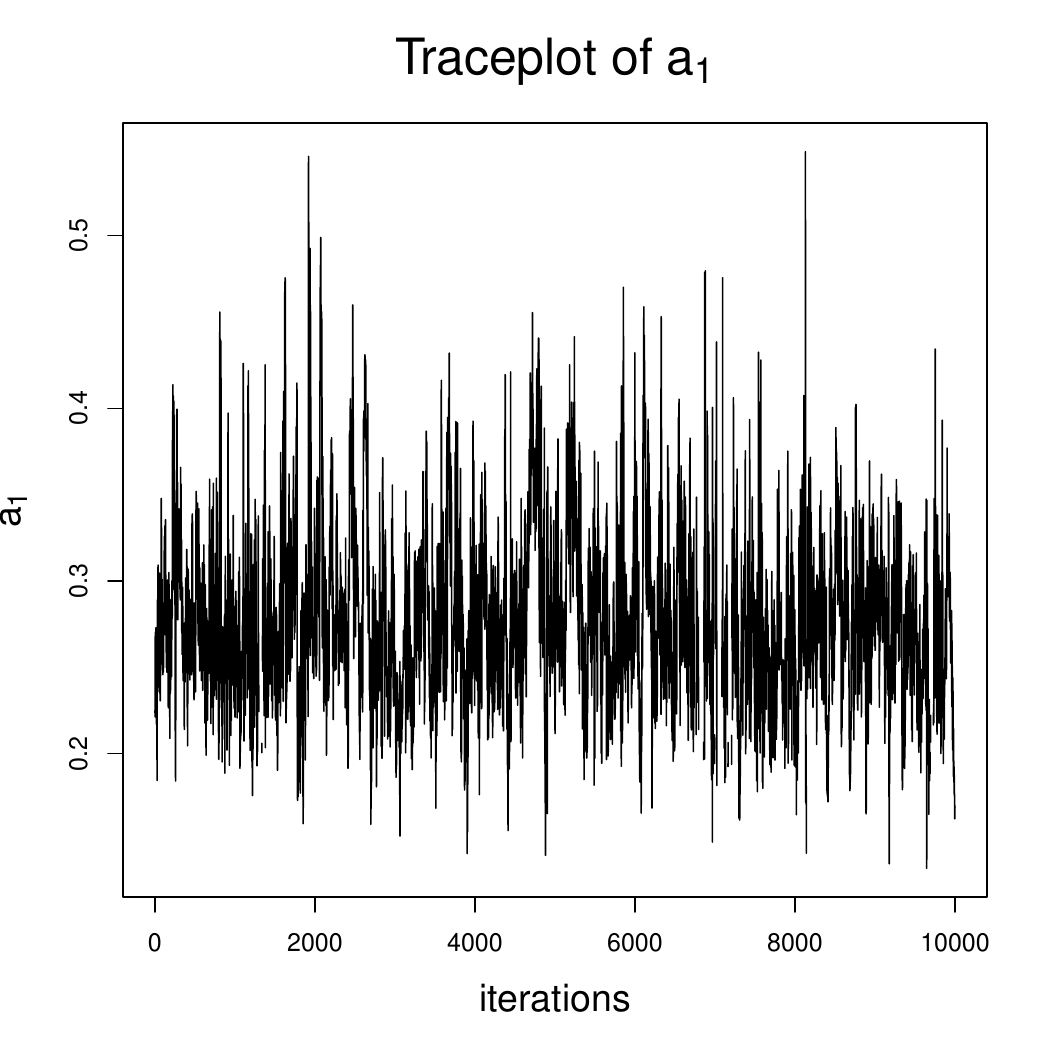}}
\hspace{2mm}
\subfigure[Trace plot of $a_2$.]{ \label{fig:sim6_trace_p2}
\includegraphics[width=4.5cm,height=5cm]{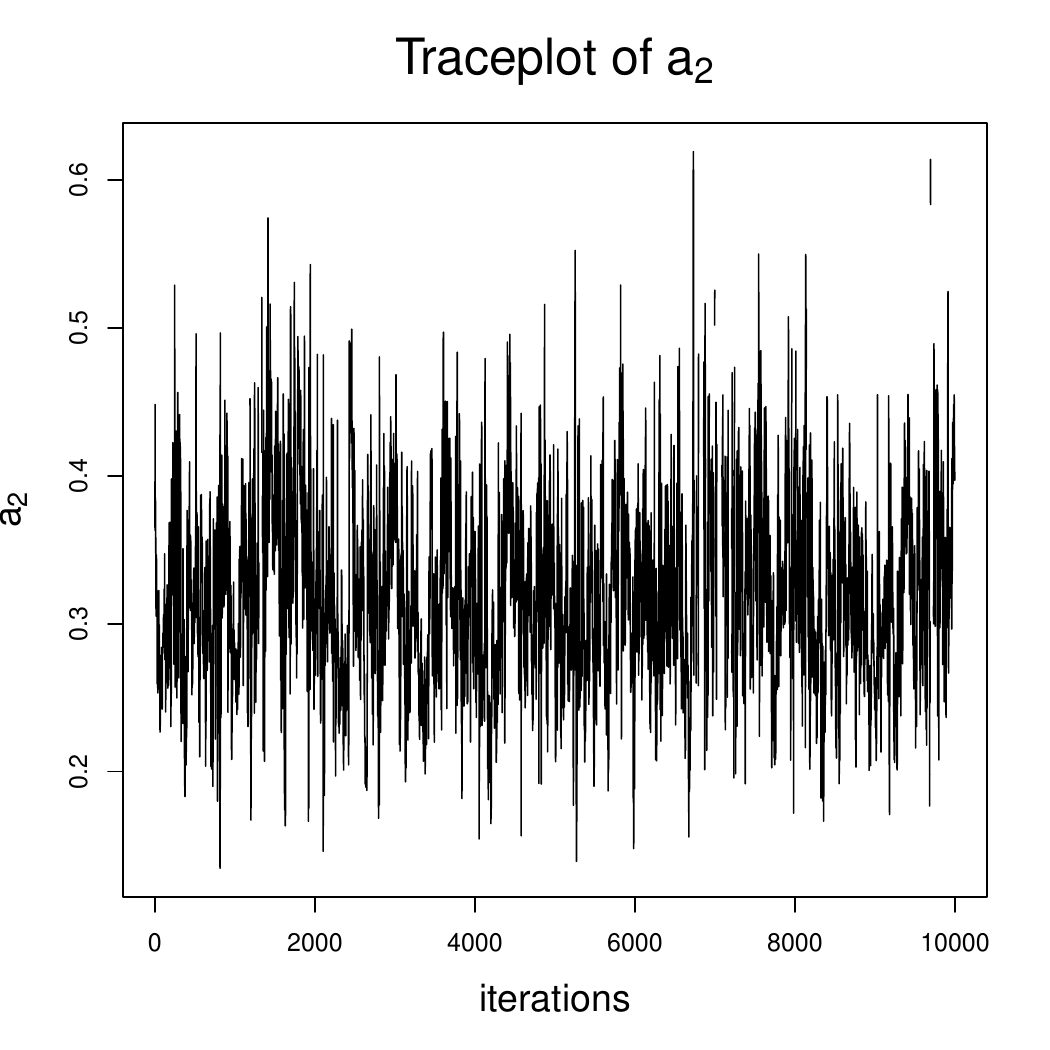}}\\
\vspace{2mm}
\subfigure[Trace plot of $a_3$.]{ \label{fig:sim6_trace_p3}
\includegraphics[width=4.5cm,height=5cm]{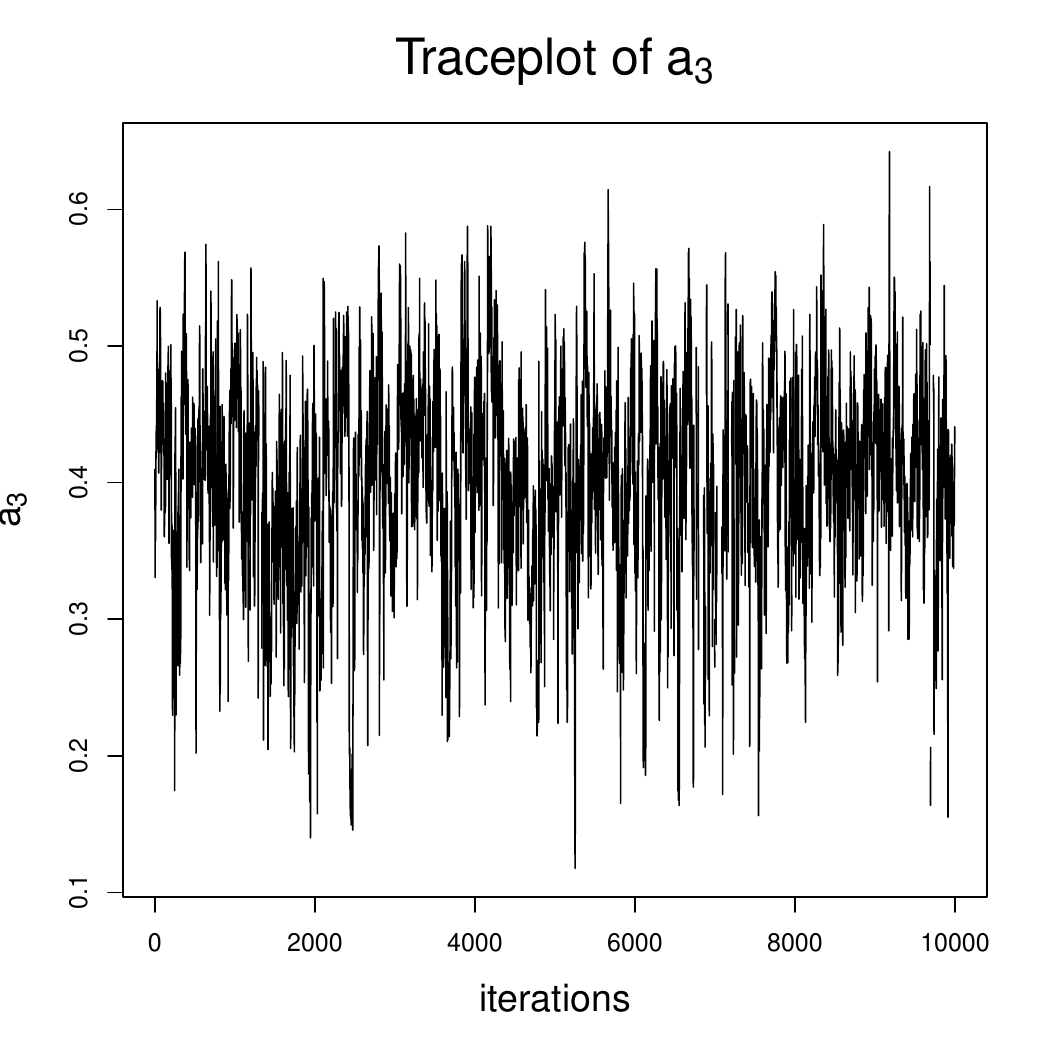}}
\caption{{\bf TTMCMC for $SDE_2$ and $\pi_4$:} Trace plots of $M$, $\mu_1$, $\mu_2$, $\nu_3$, $\omega^2_1$, $\omega^2_2$, $\omega^2_3$, $a_1$ $a_2$ and $a_3$.} 
\label{fig:sim6_trace_plots}
\end{figure}

\begin{figure}
\centering
\subfigure[Posterior of $\mu_1$.]{ \label{fig:sim6_mu1}
\includegraphics[width=4.5cm,height=5cm]{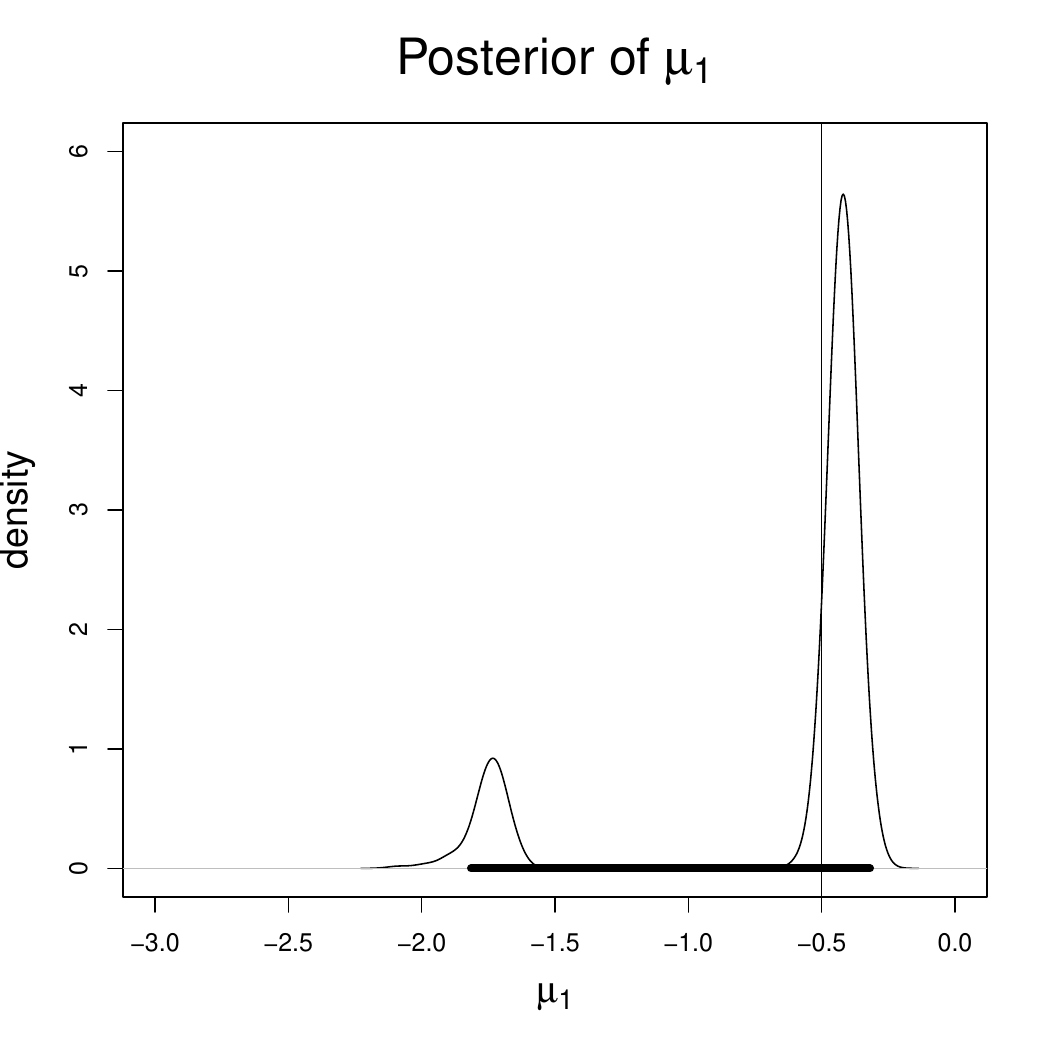}}
\hspace{2mm}
\subfigure[Posterior of $\mu_2$.]{ \label{fig:sim6_mu2}
\includegraphics[width=4.5cm,height=5cm]{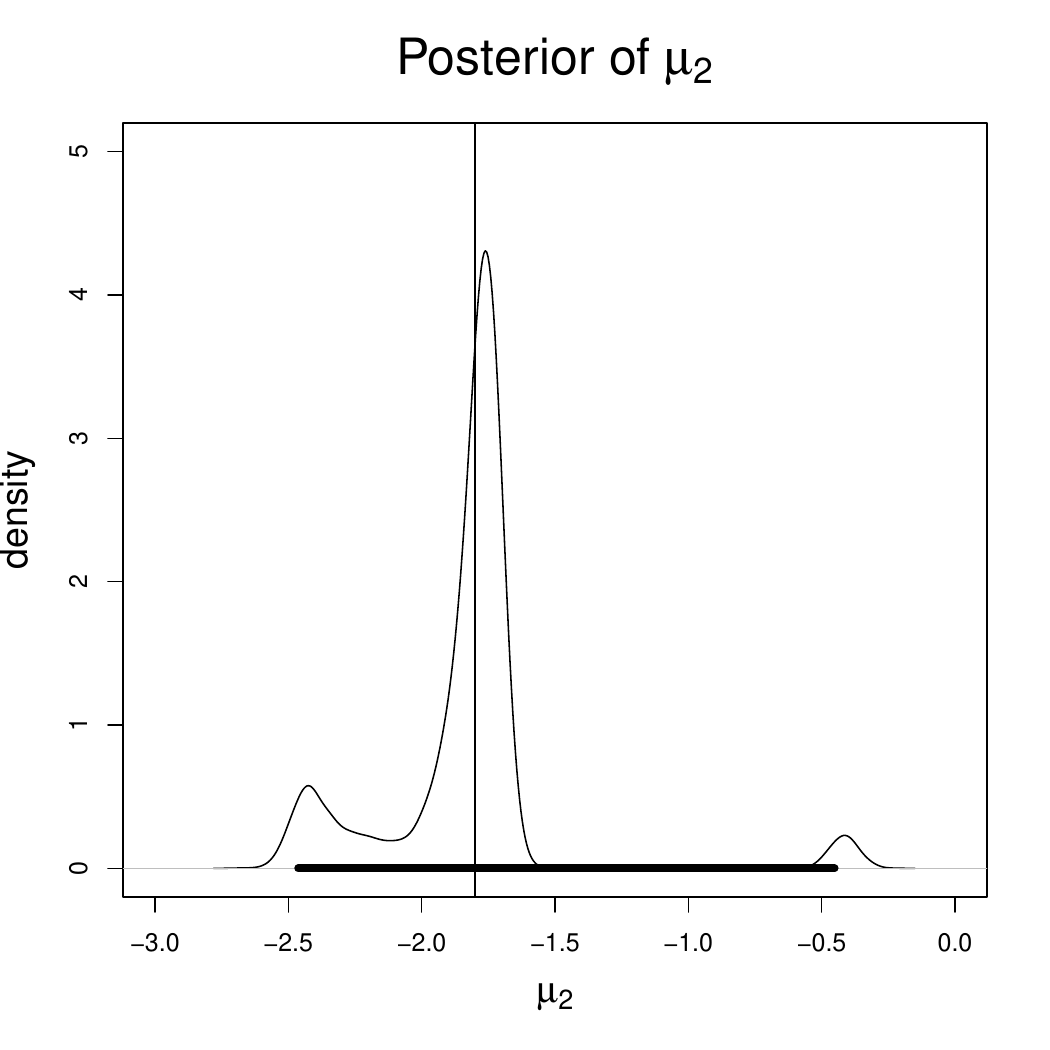}}
\hspace{2mm}
\subfigure[Posterior of $\mu_3$.]{ \label{fig:sim6_mu3}
\includegraphics[width=4.5cm,height=5cm]{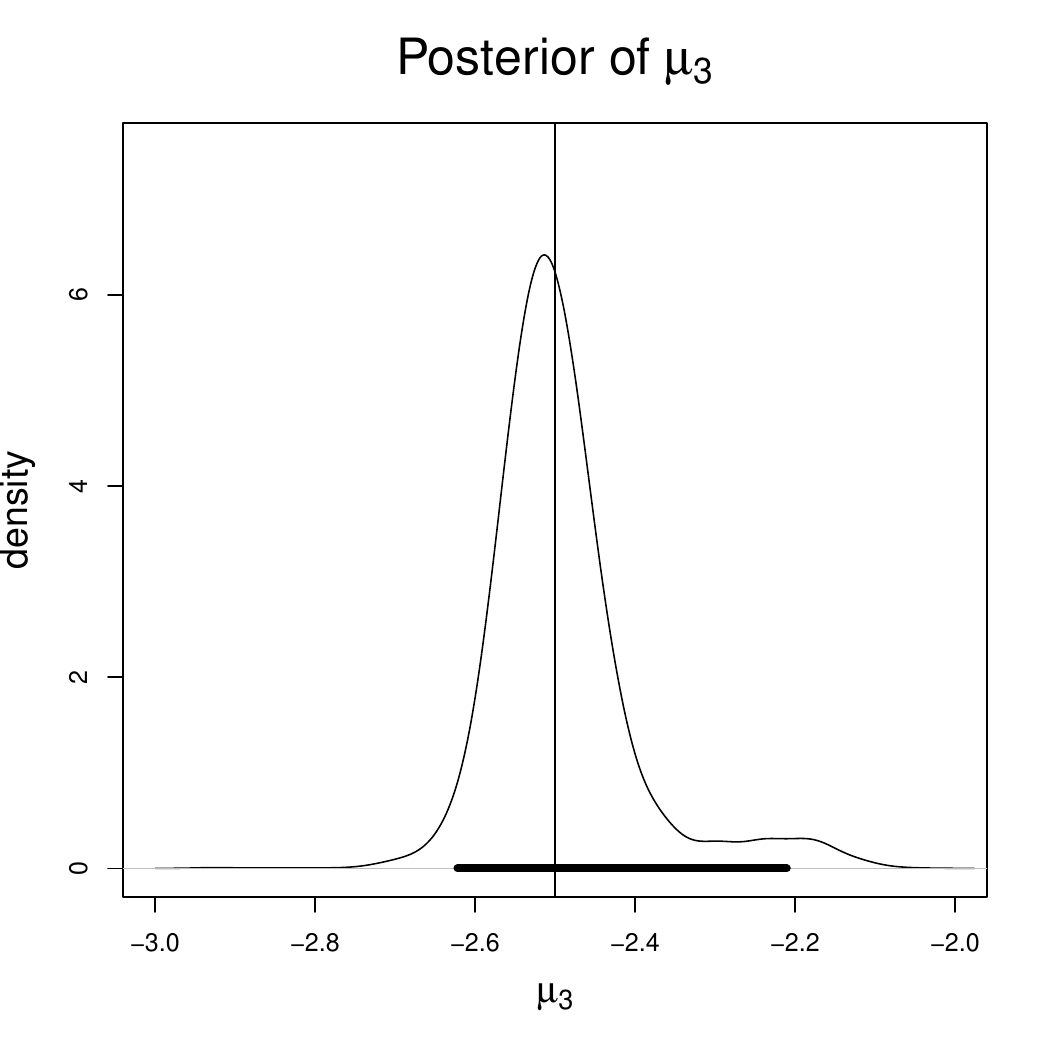}}\\
\vspace{2mm}
\subfigure[Posterior of $\omega^2_1$.]{ \label{fig:sim6_omegasq1}
\includegraphics[width=4.5cm,height=5cm]{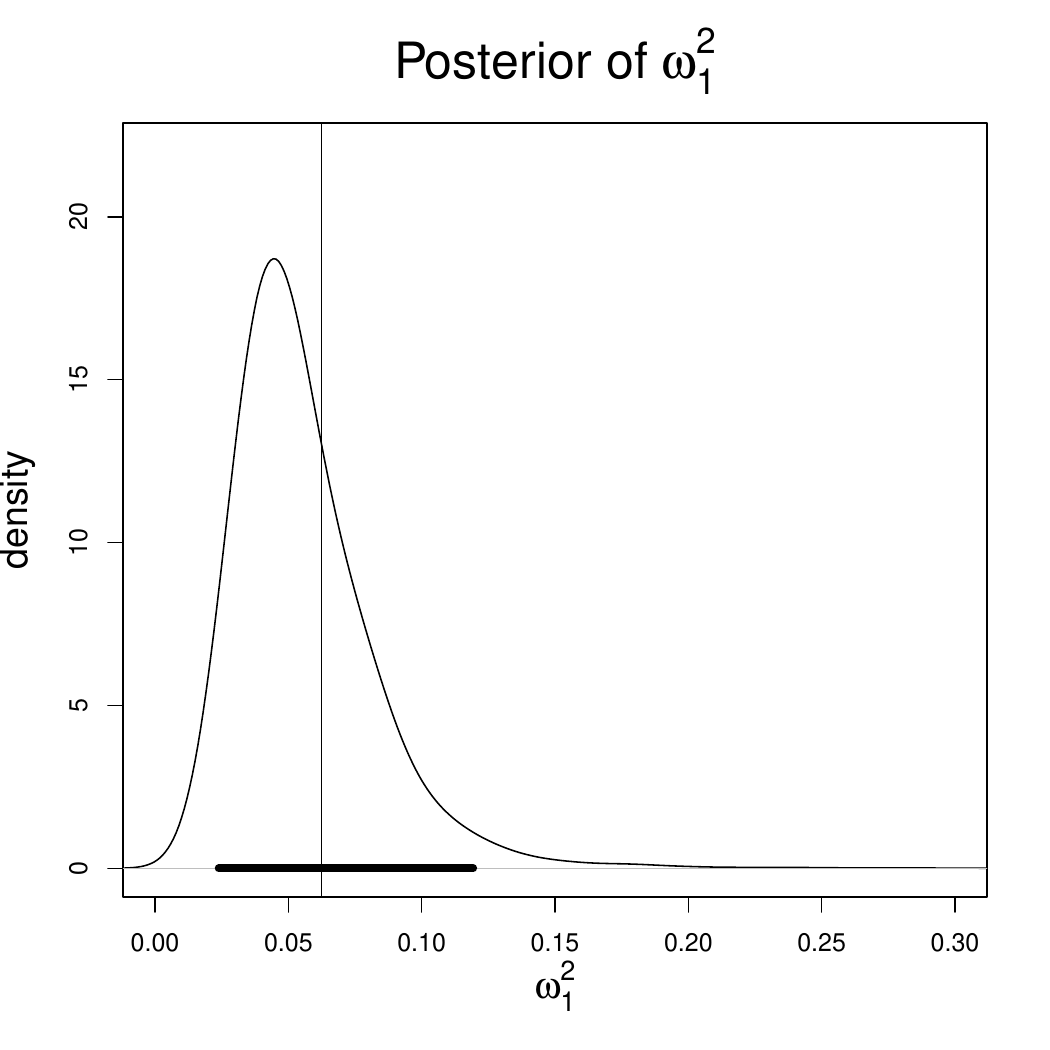}}
\hspace{2mm}
\subfigure[Posterior of $\omega^2_2$.]{ \label{fig:sim6_omegasq2}
\includegraphics[width=4.5cm,height=5cm]{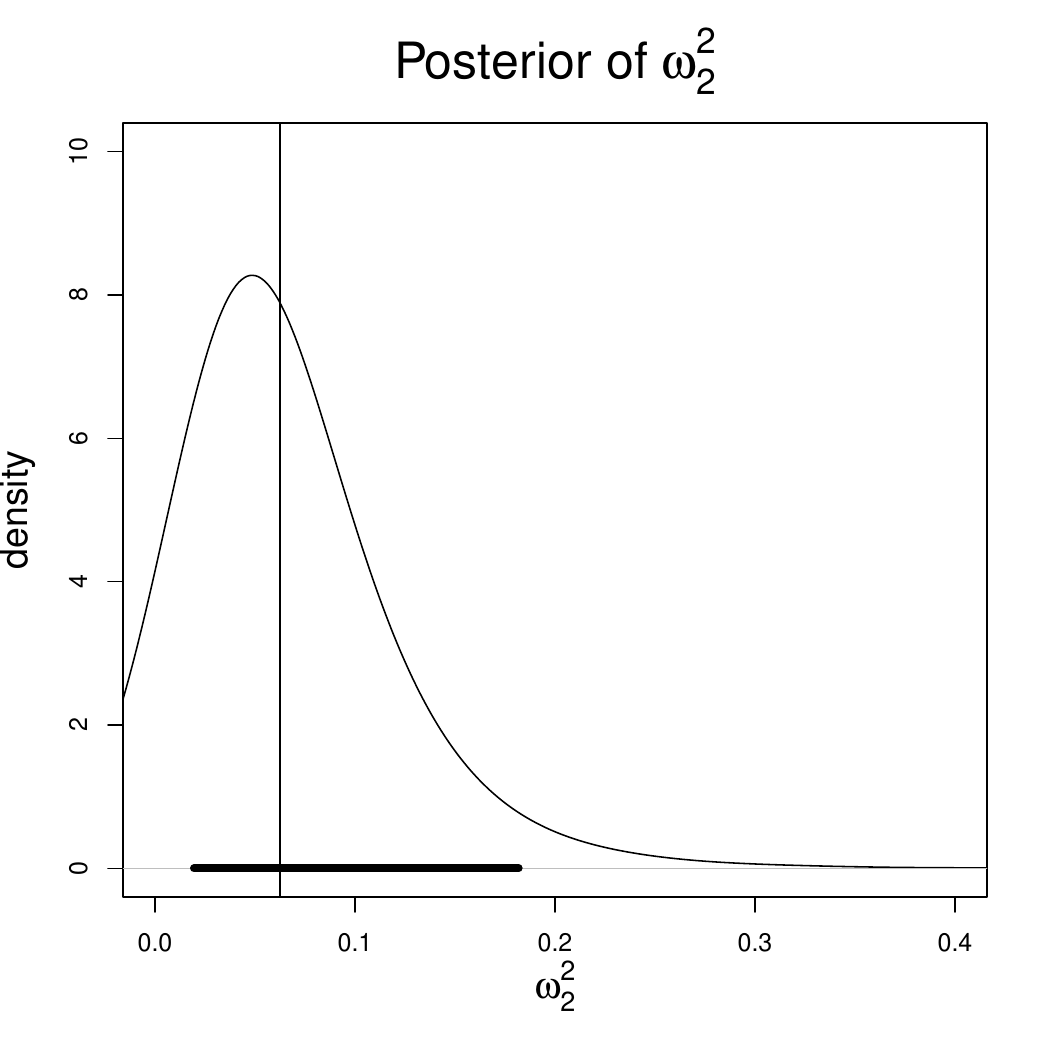}}
\hspace{2mm}
\subfigure[Posterior of $\omega^2_3$.]{ \label{fig:sim6_omegasq3}
\includegraphics[width=4.5cm,height=5cm]{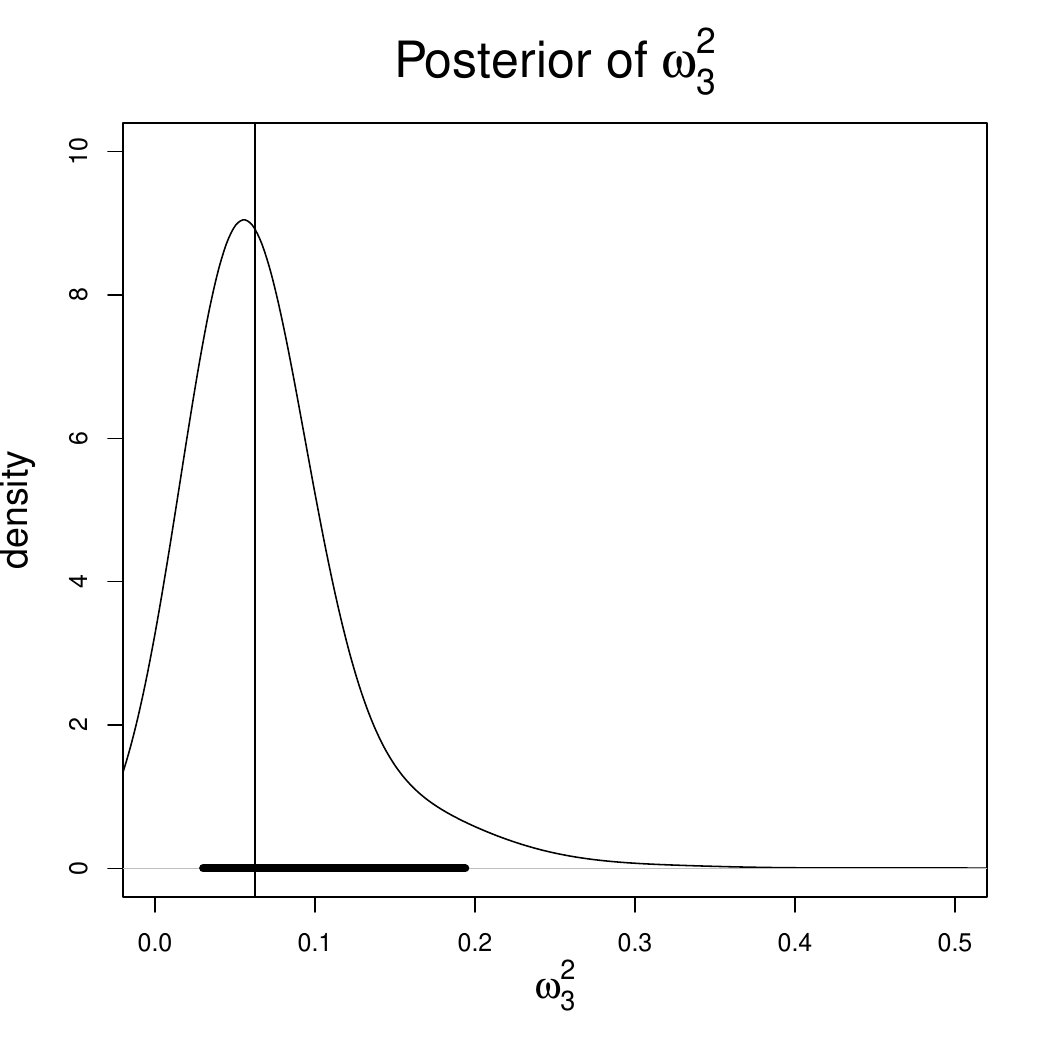}}\\
\vspace{2mm}
\subfigure[Posterior of $a_1$.]{ \label{fig:sim6_p1}
\includegraphics[width=4.5cm,height=5cm]{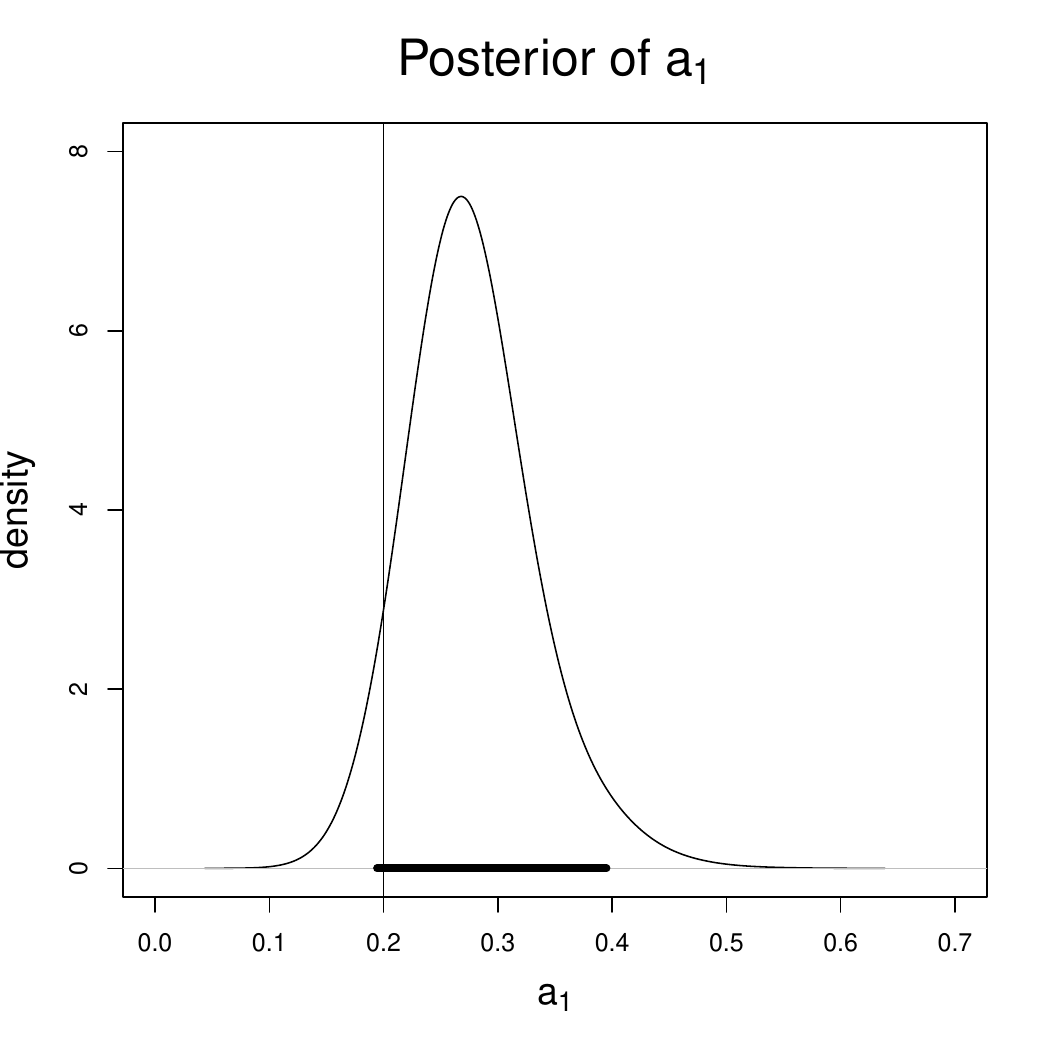}}
\hspace{2mm}
\subfigure[Posterior of $a_2$.]{ \label{fig:sim6_p2}
\includegraphics[width=4.5cm,height=5cm]{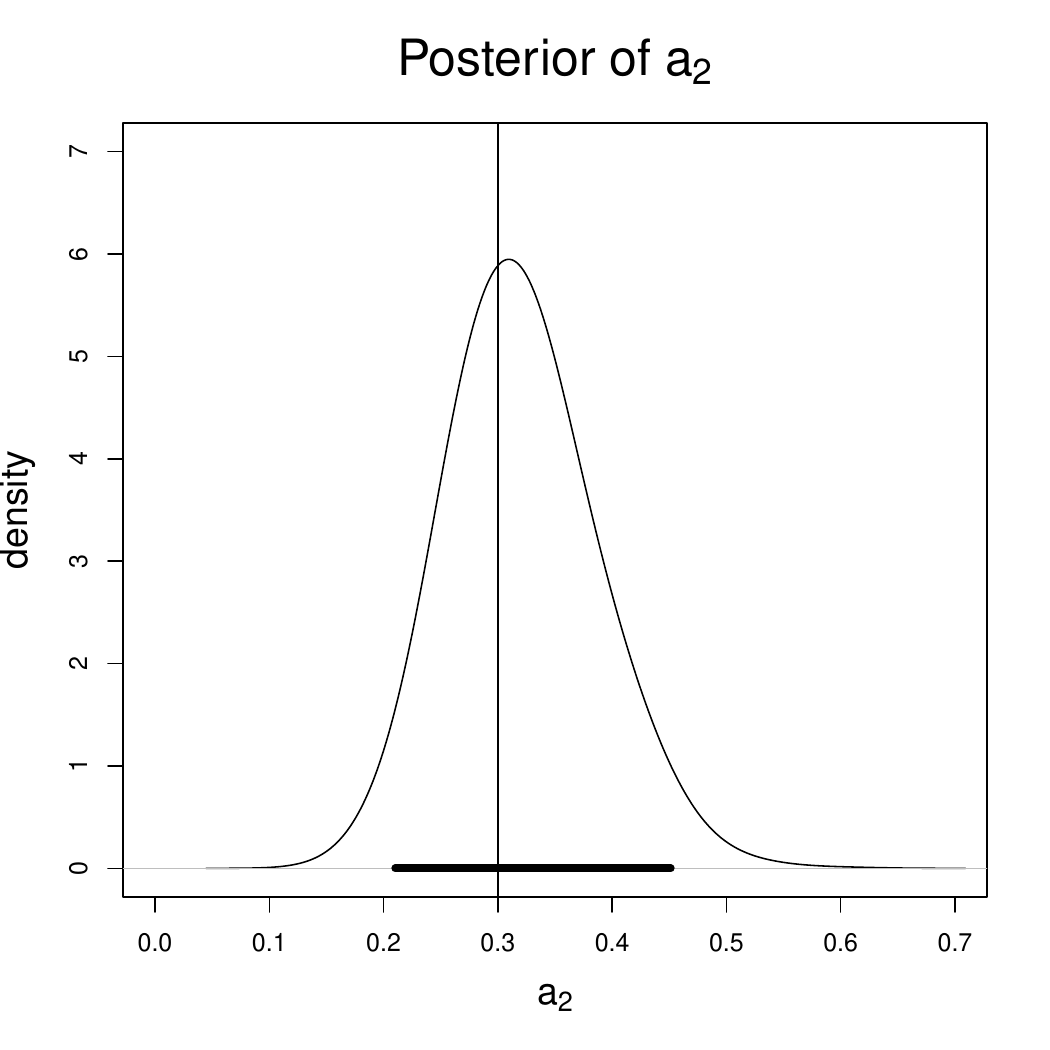}}
\hspace{2mm}
\subfigure[Posterior of $a_3$.]{ \label{fig:sim6_p3}
\includegraphics[width=4.5cm,height=5cm]{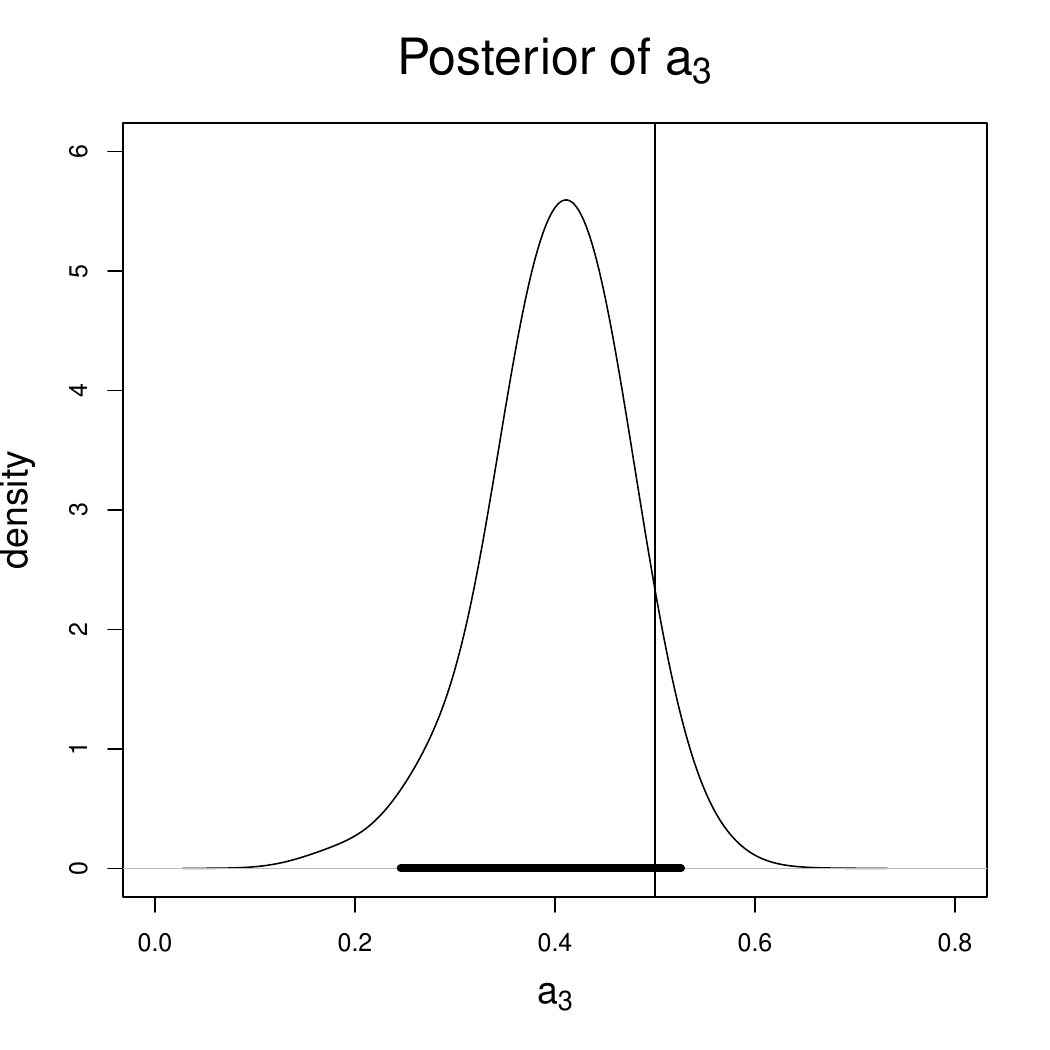}}
\caption{{\bf TTMCMC for $SDE_2$ and $\pi_4$:} Posteriors of $M$, $\mu_1$, $\mu_2$, $\mu_3$, $\omega^2_1$, $\omega^2_2$, $\omega^2_3$, $a_1$ $a_2$ and $a_3$. 
The vertical lines stand for the true values, while the thick horizontal lines denote the 95\% credible intervals.} 
\label{fig:sim6_posterior_plots}
\end{figure}

In this case, the mixture components are not well-separated, and as demonstrated by \ctn{Maud16}, selection of $M$ by $BIC$ leads to unsatisfactory results.
Our Bayesian approach, however, overcomes this challenge by acknowledging variable dimensionality and by employing TTMCMC.
Figures \ref{fig:sim6_trace_plots} and \ref{fig:sim6_posterior_plots}, the relevant posterior plots associated with $SDE_2$ and $\pi_4$ indicate excellent
mixing as before and satisfactory performance. 
The overall acceptance rate here is about $2.5\%$.

In this case, $M=2$, $M=3$ and $M=4$ received the posterior probabilities $0.0328$, $0.8138$ and $0.1534$, respectively, 
while all the other values received the zero posterior probability. In other words, our Bayesian approach successfully identifies the true number of components,
unlike $BIC$.

\subsection{Seventh simulation study: $SDE_2$ with $\pi_5$} 
\label{subsec:sim7}

\begin{figure}
\centering
\subfigure[Trace plot of $M$.]{ \label{fig:sim7_trace_comp}
\includegraphics[width=4.5cm,height=5cm]{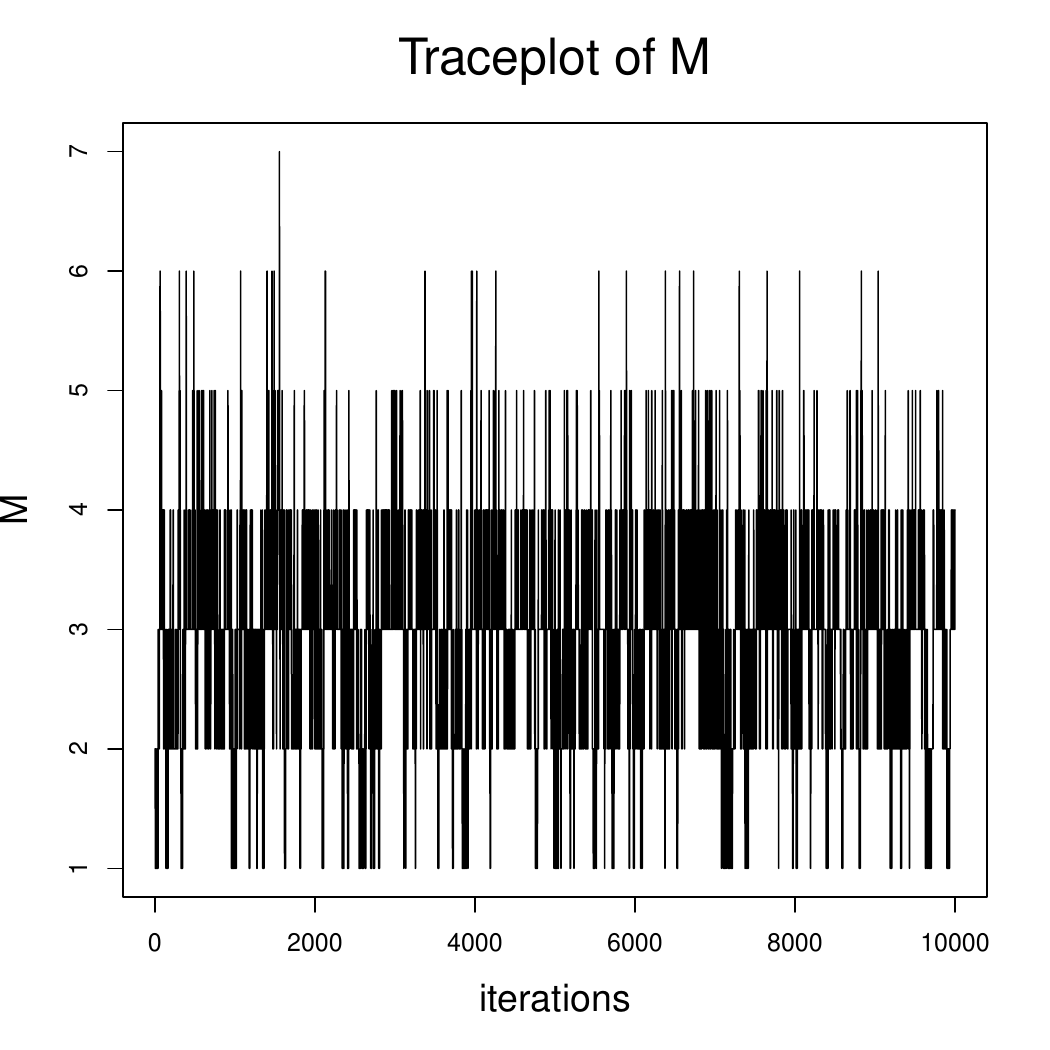}}
\hspace{2mm}
\subfigure[Trace plot of $\mu_1$.]{ \label{fig:sim7_trace_mu1}
\includegraphics[width=4.5cm,height=5cm]{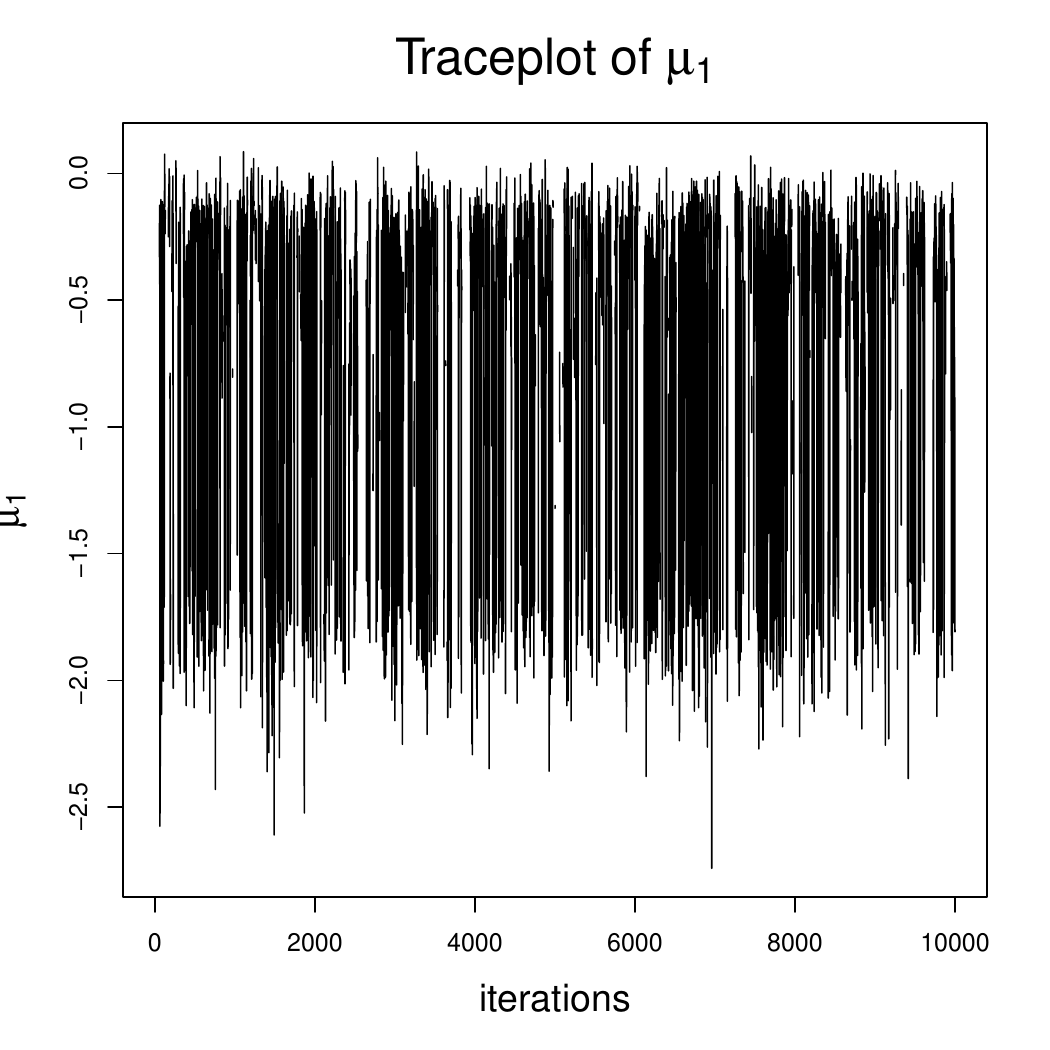}}
\hspace{2mm}
\subfigure[Trace plot of $\mu_2$.]{ \label{fig:sim7_trace_mu2}
\includegraphics[width=4.5cm,height=5cm]{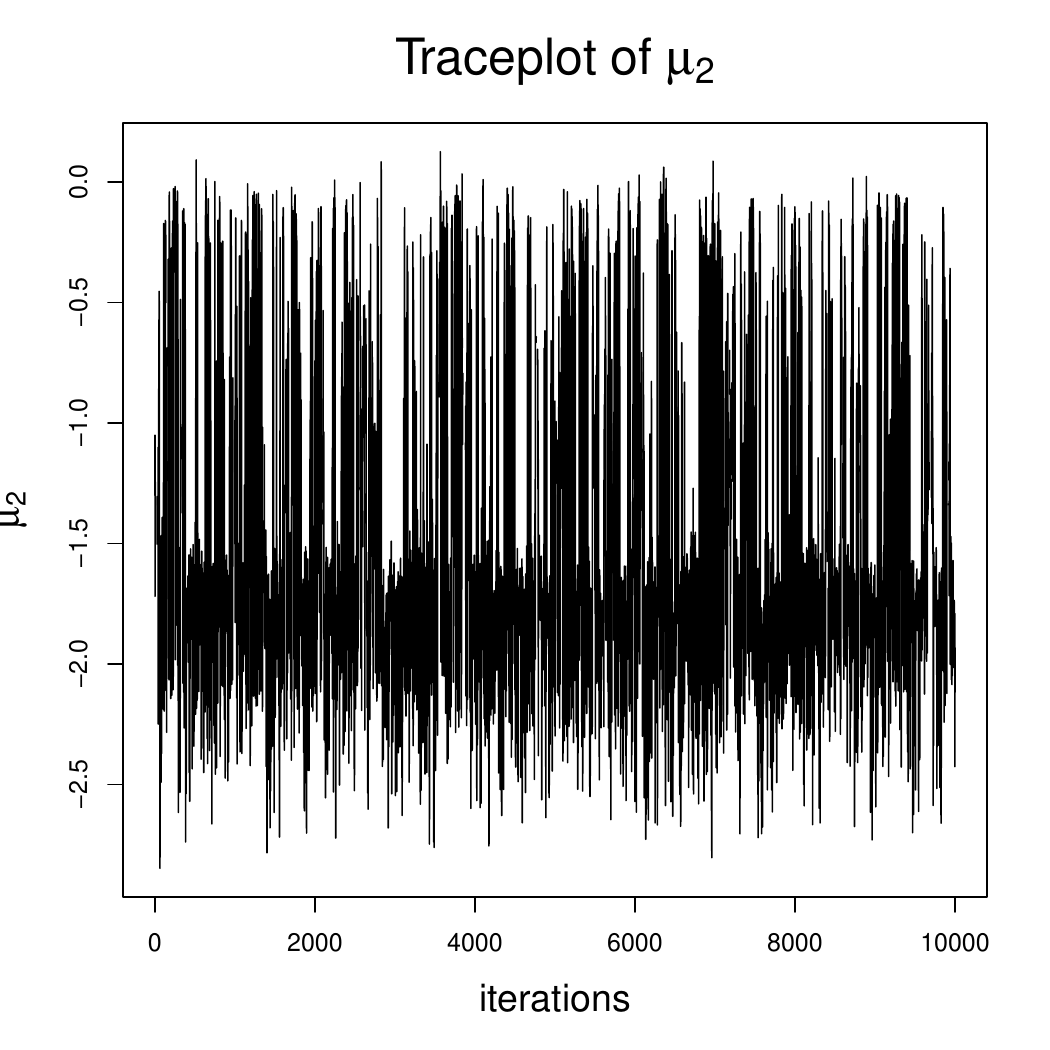}}\\
\vspace{2mm}
\subfigure[Trace plot of $\mu_3$.]{ \label{fig:sim7_trace_mu3}
\includegraphics[width=4.5cm,height=5cm]{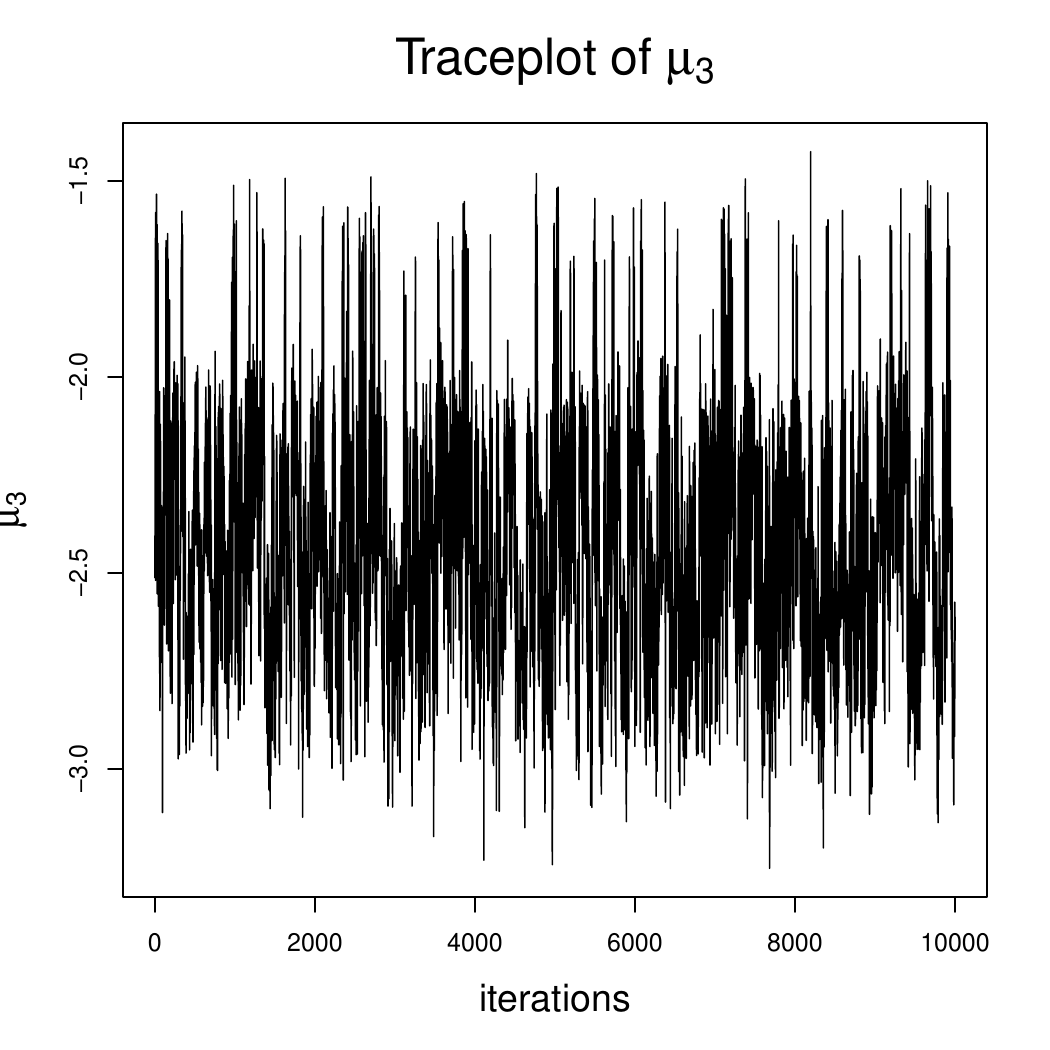}}
\hspace{2mm}
\subfigure[Trace plot of $\omega^2_1$.]{ \label{fig:sim7_trace_omegasq1}
\includegraphics[width=4.5cm,height=5cm]{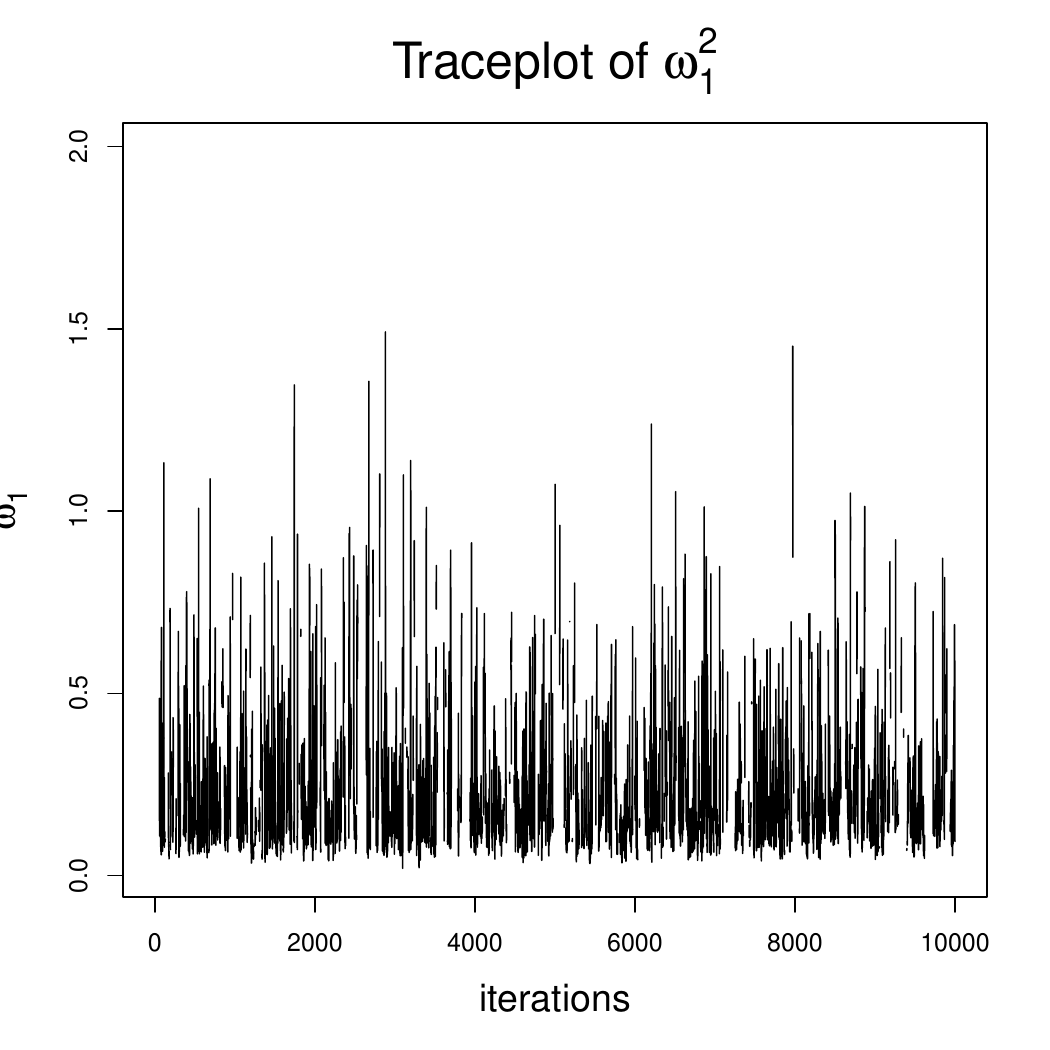}}
\hspace{2mm}
\subfigure[Trace plot of $\omega^2_2$.]{ \label{fig:sim7_trace_omegasq2}
\includegraphics[width=4.5cm,height=5cm]{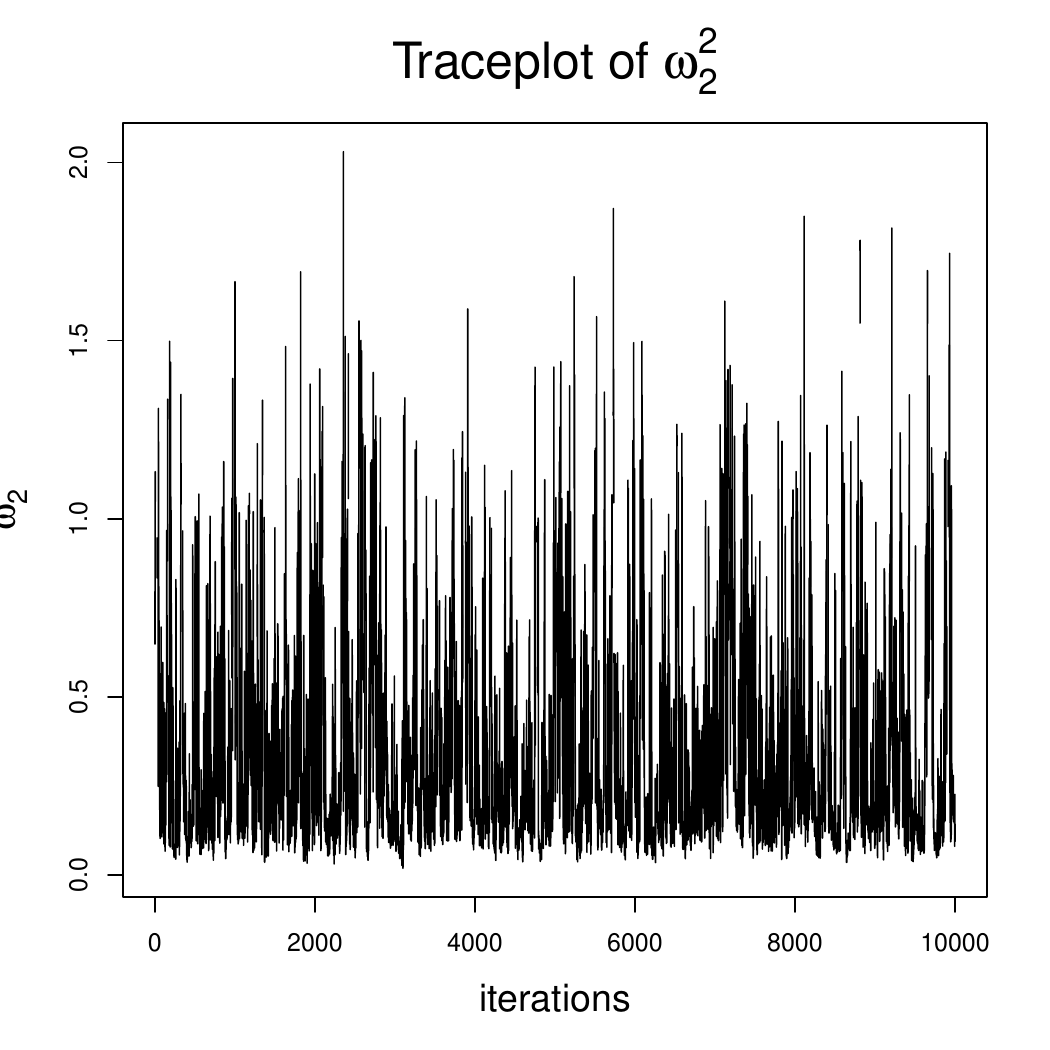}}\\
\vspace{2mm}
\subfigure[Trace plot of $\omega^2_3$.]{ \label{fig:sim7_trace_omegasq3}
\includegraphics[width=4.5cm,height=5cm]{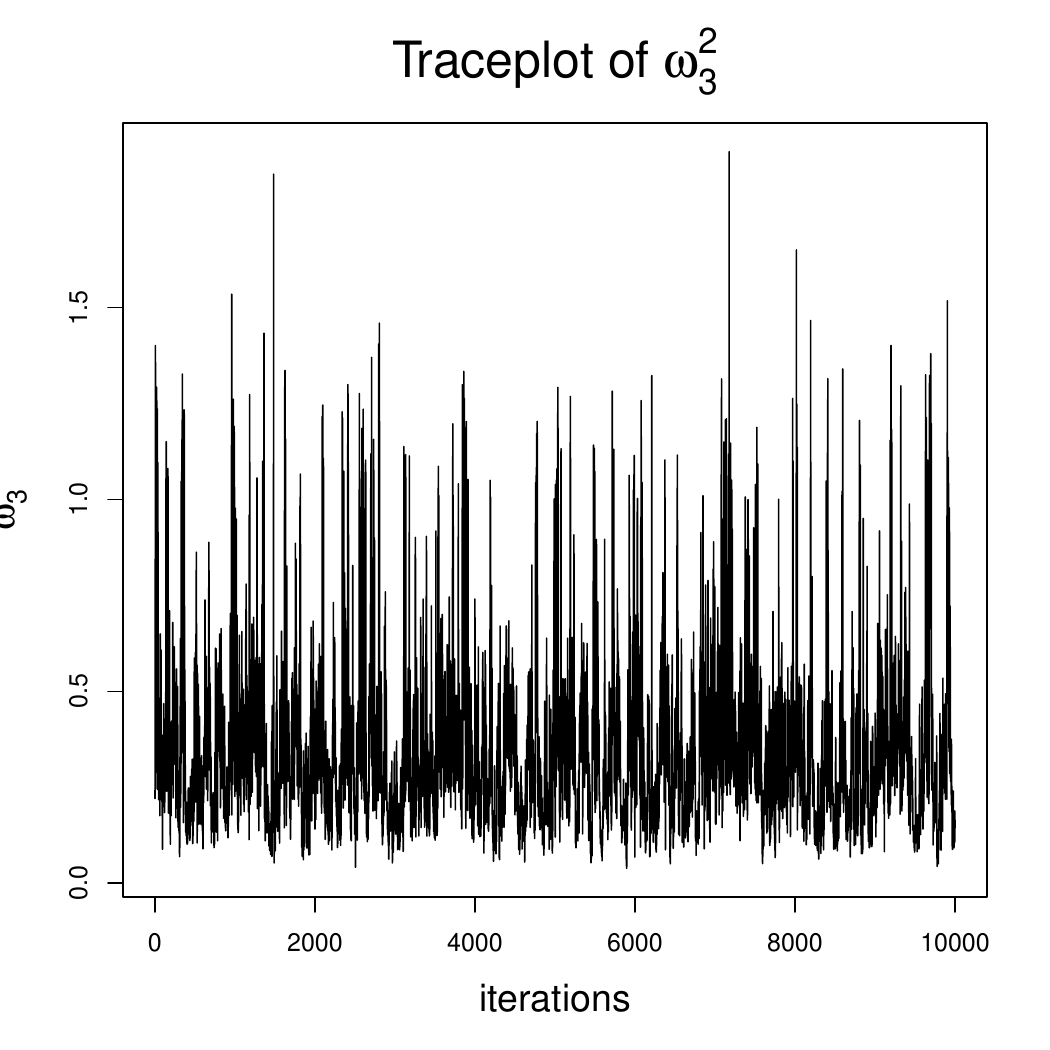}}
\hspace{2mm}
\subfigure[Trace plot of $a_1$.]{ \label{fig:sim7_trace_p1}
\includegraphics[width=4.5cm,height=5cm]{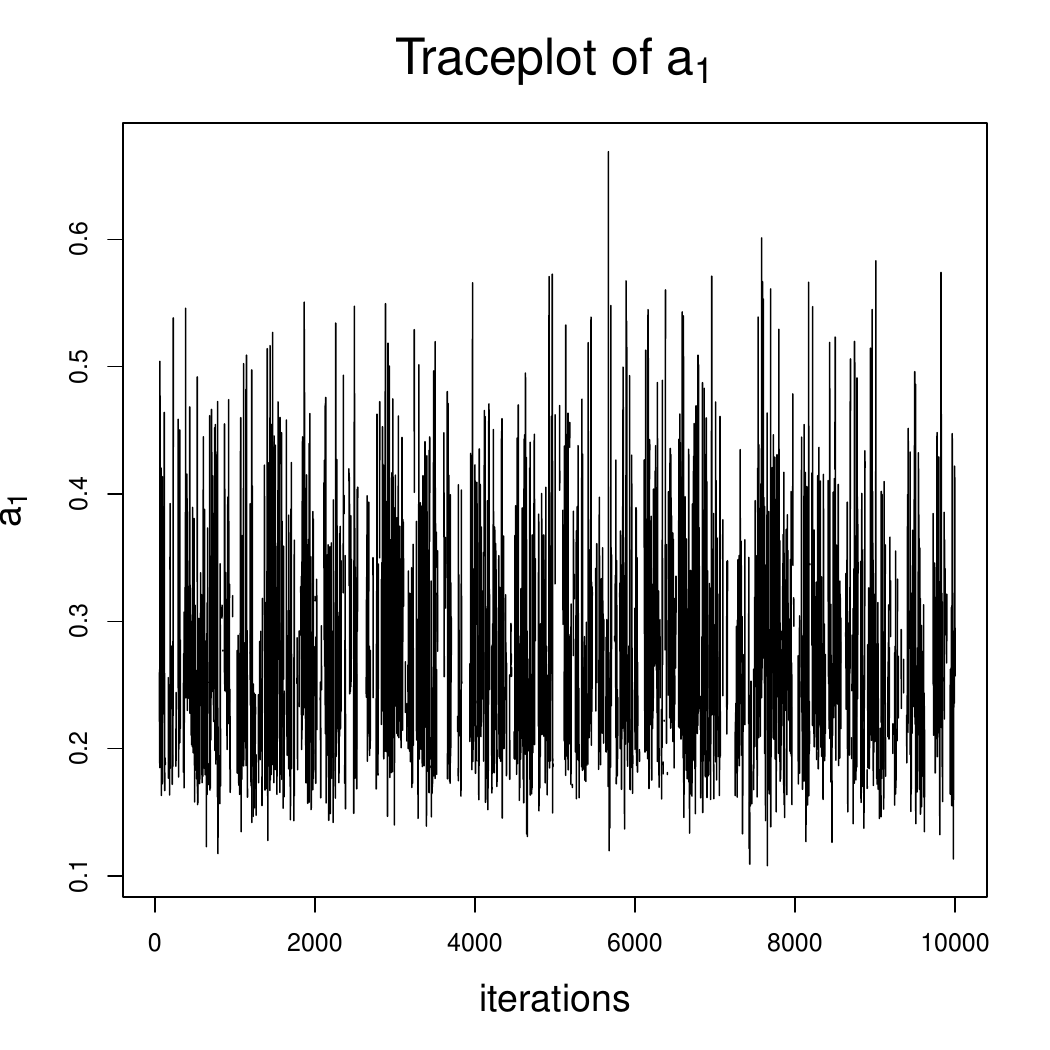}}
\hspace{2mm}
\subfigure[Trace plot of $a_2$.]{ \label{fig:sim7_trace_p2}
\includegraphics[width=4.5cm,height=5cm]{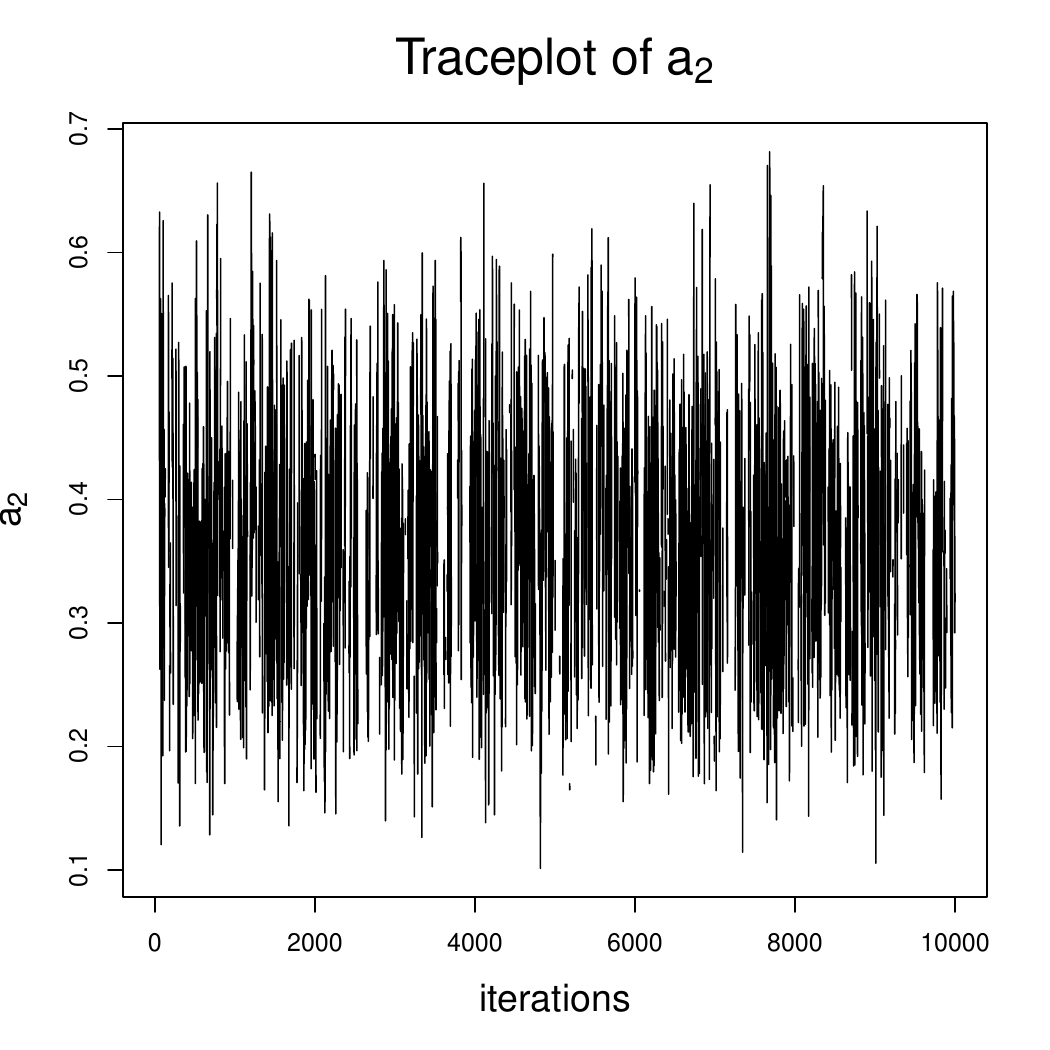}}\\
\vspace{2mm}
\subfigure[Trace plot of $a_3$.]{ \label{fig:sim7_trace_p3}
\includegraphics[width=4.5cm,height=5cm]{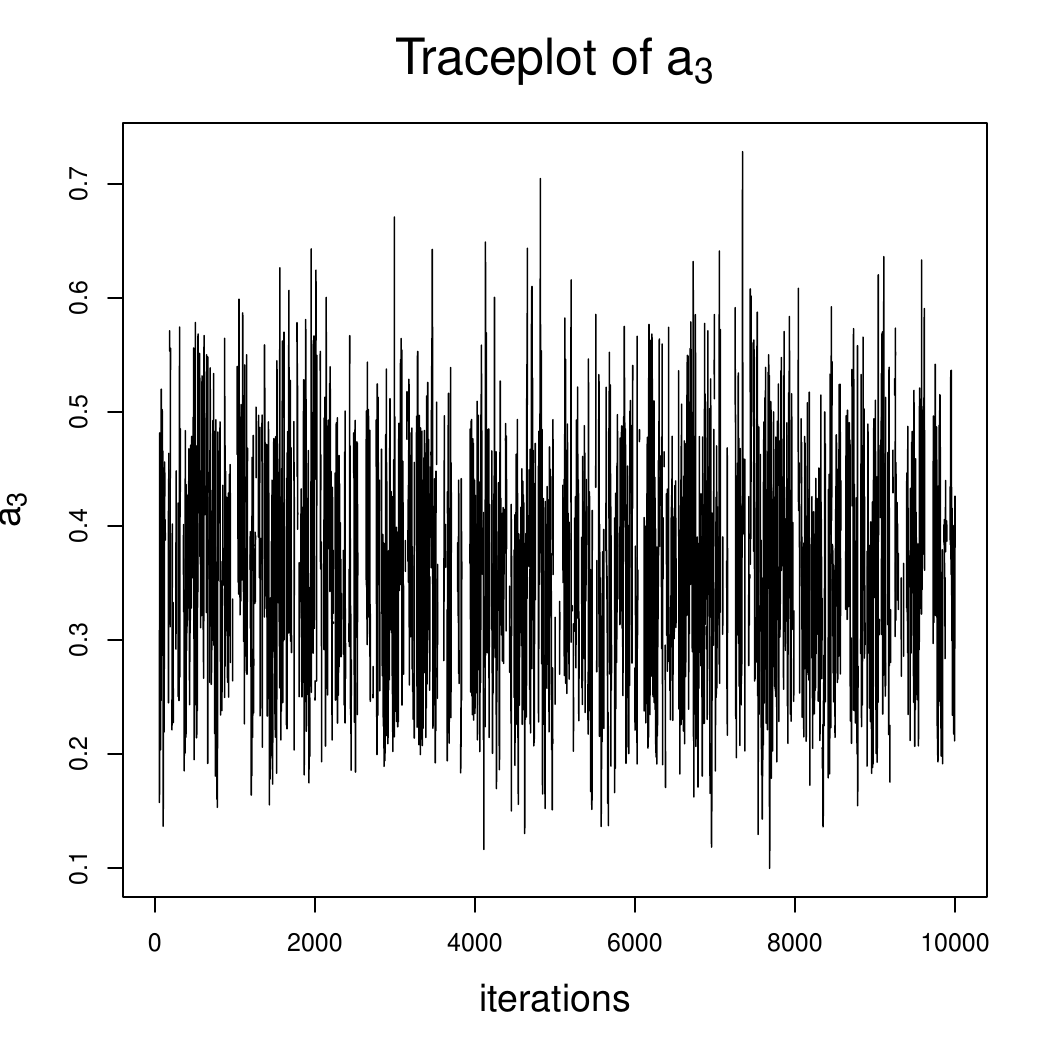}}
\caption{{\bf TTMCMC for $SDE_2$ and $\pi_5$:} Trace plots of $M$, $\mu_1$, $\mu_2$, $\nu_3$, $\omega^2_1$, $\omega^2_2$, $\omega^2_3$, $a_1$ $a_2$ and $a_3$.} 
\label{fig:sim7_trace_plots}
\end{figure}

\begin{figure}
\centering
\subfigure[Posterior of $\mu_1$.]{ \label{fig:sim7_mu1}
\includegraphics[width=4.5cm,height=5cm]{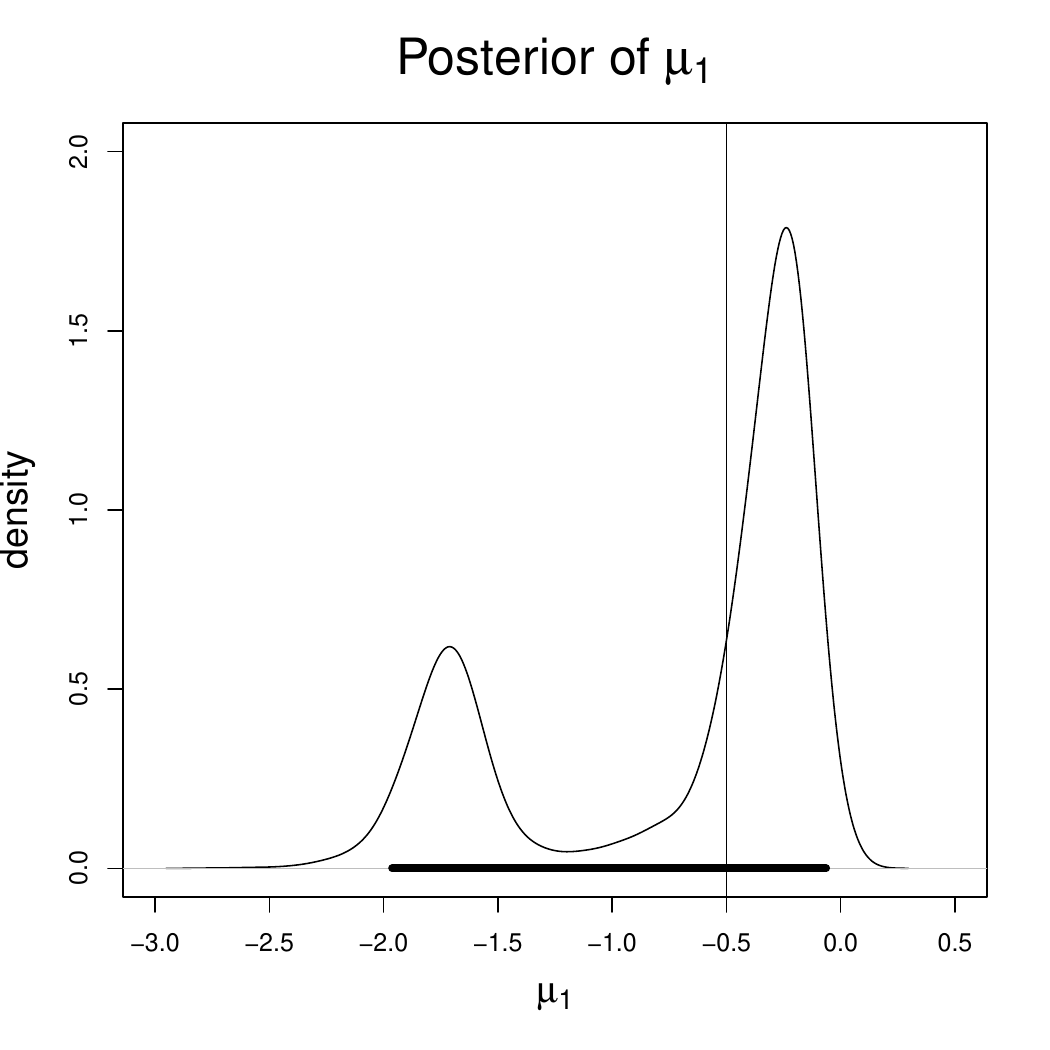}}
\hspace{2mm}
\subfigure[Posterior of $\mu_2$.]{ \label{fig:sim7_mu2}
\includegraphics[width=4.5cm,height=5cm]{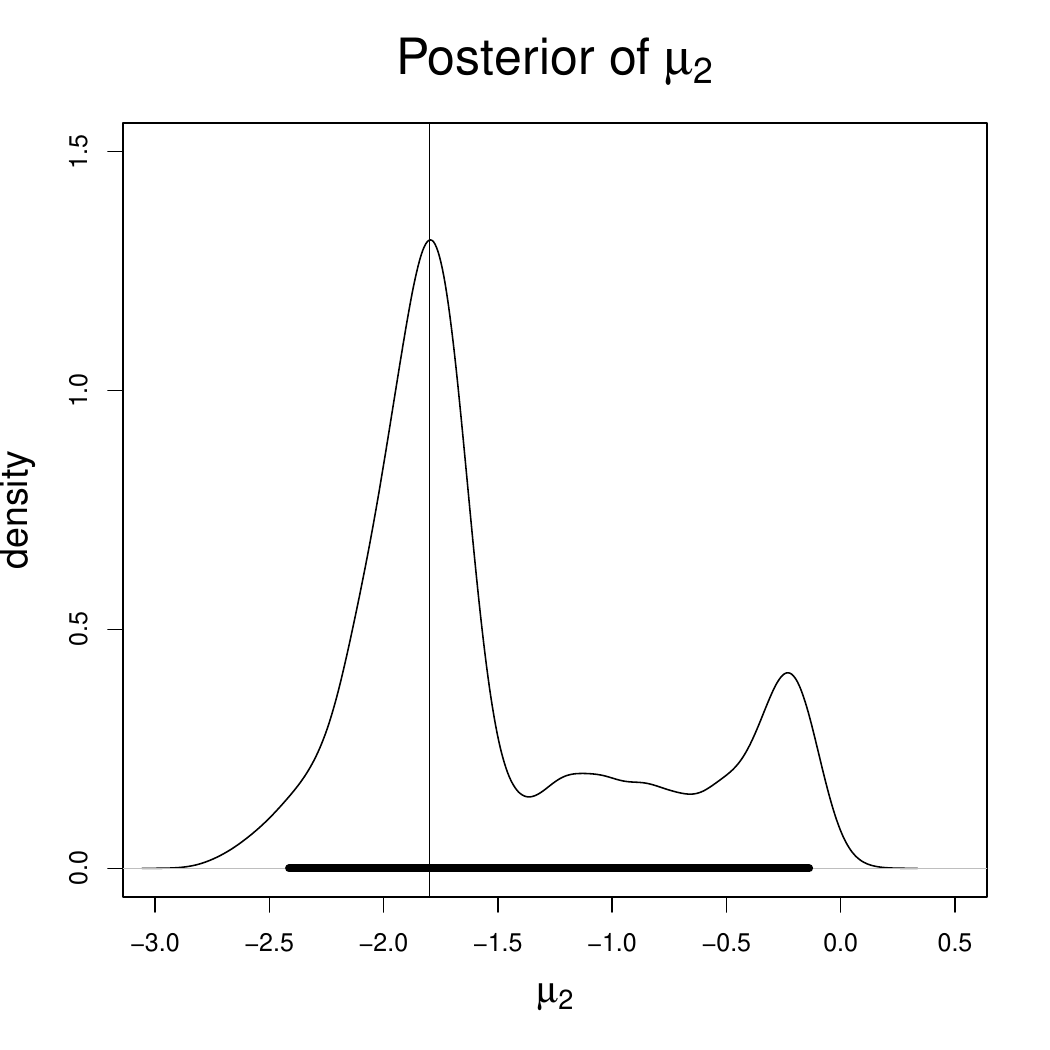}}
\hspace{2mm}
\subfigure[Posterior of $\mu_3$.]{ \label{fig:sim7_mu3}
\includegraphics[width=4.5cm,height=5cm]{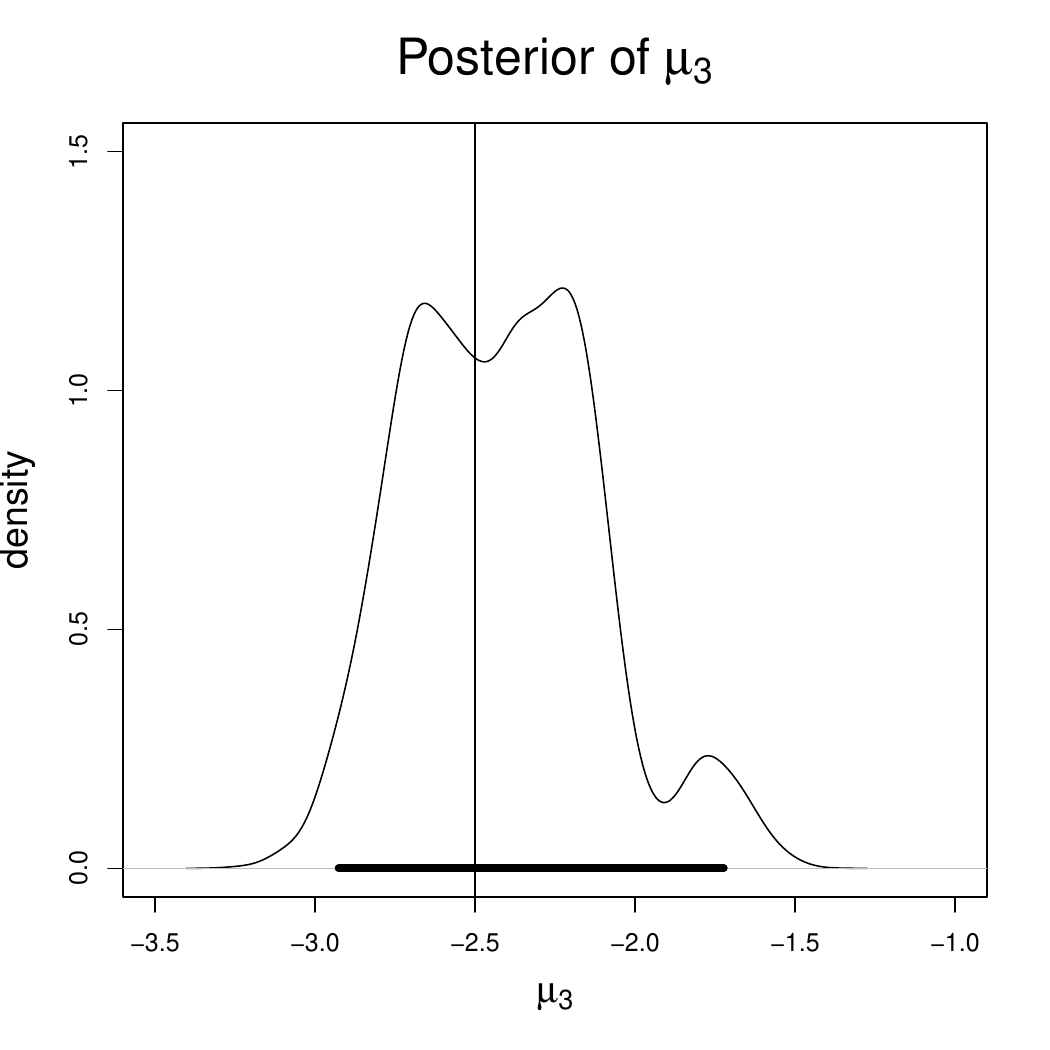}}\\
\vspace{2mm}
\subfigure[Posterior of $\omega^2_1$.]{ \label{fig:sim7_omegasq1}
\includegraphics[width=4.5cm,height=5cm]{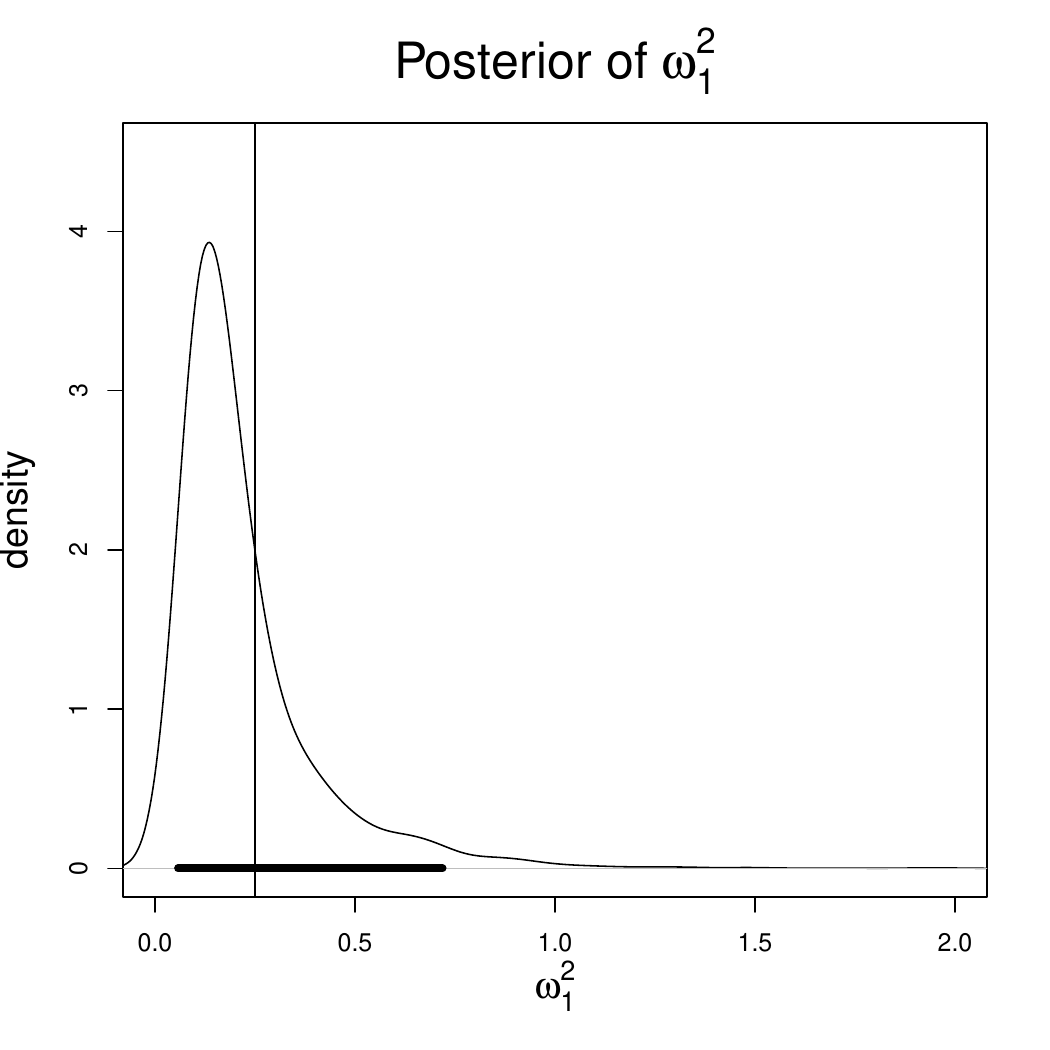}}
\hspace{2mm}
\subfigure[Posterior of $\omega^2_2$.]{ \label{fig:sim7_omegasq2}
\includegraphics[width=4.5cm,height=5cm]{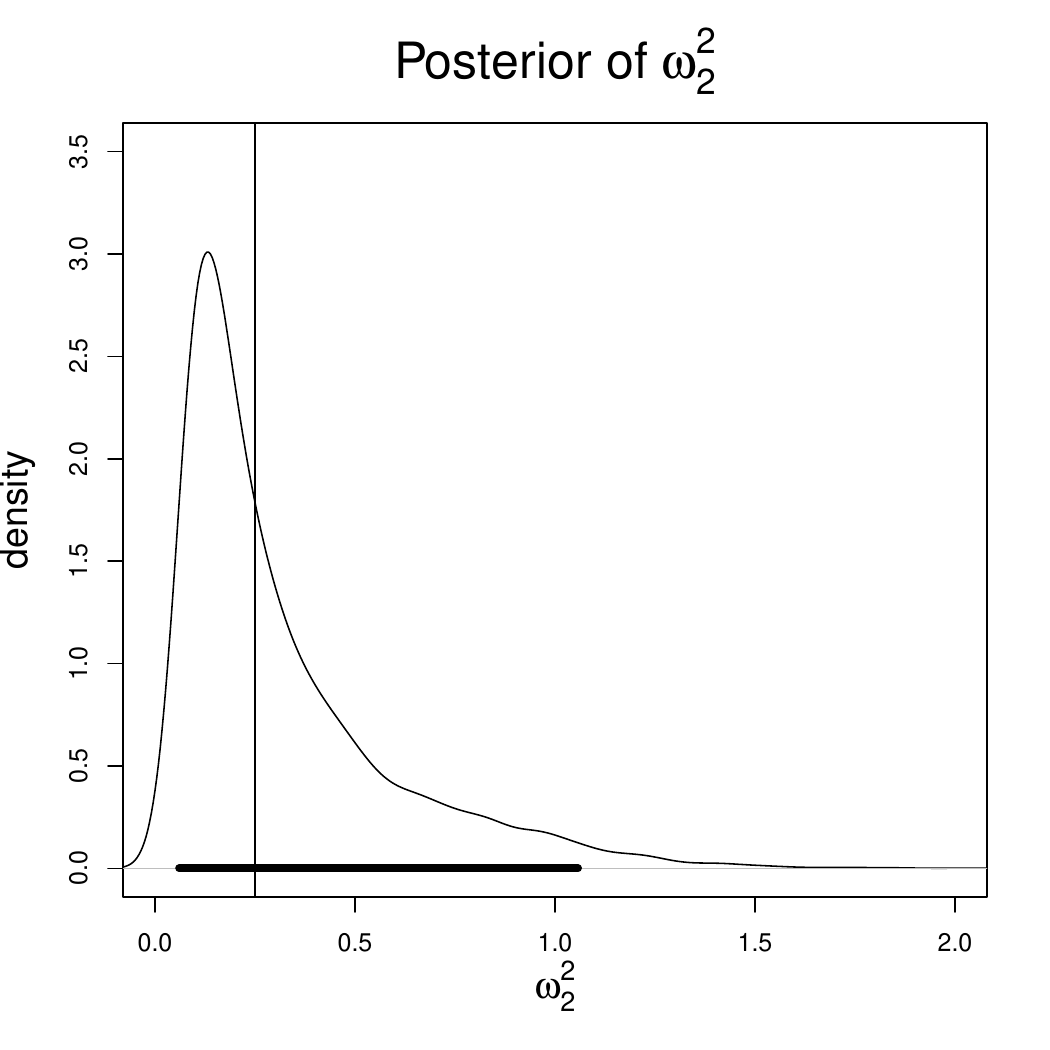}}
\hspace{2mm}
\subfigure[Posterior of $\omega^2_3$.]{ \label{fig:sim7_omegasq3}
\includegraphics[width=4.5cm,height=5cm]{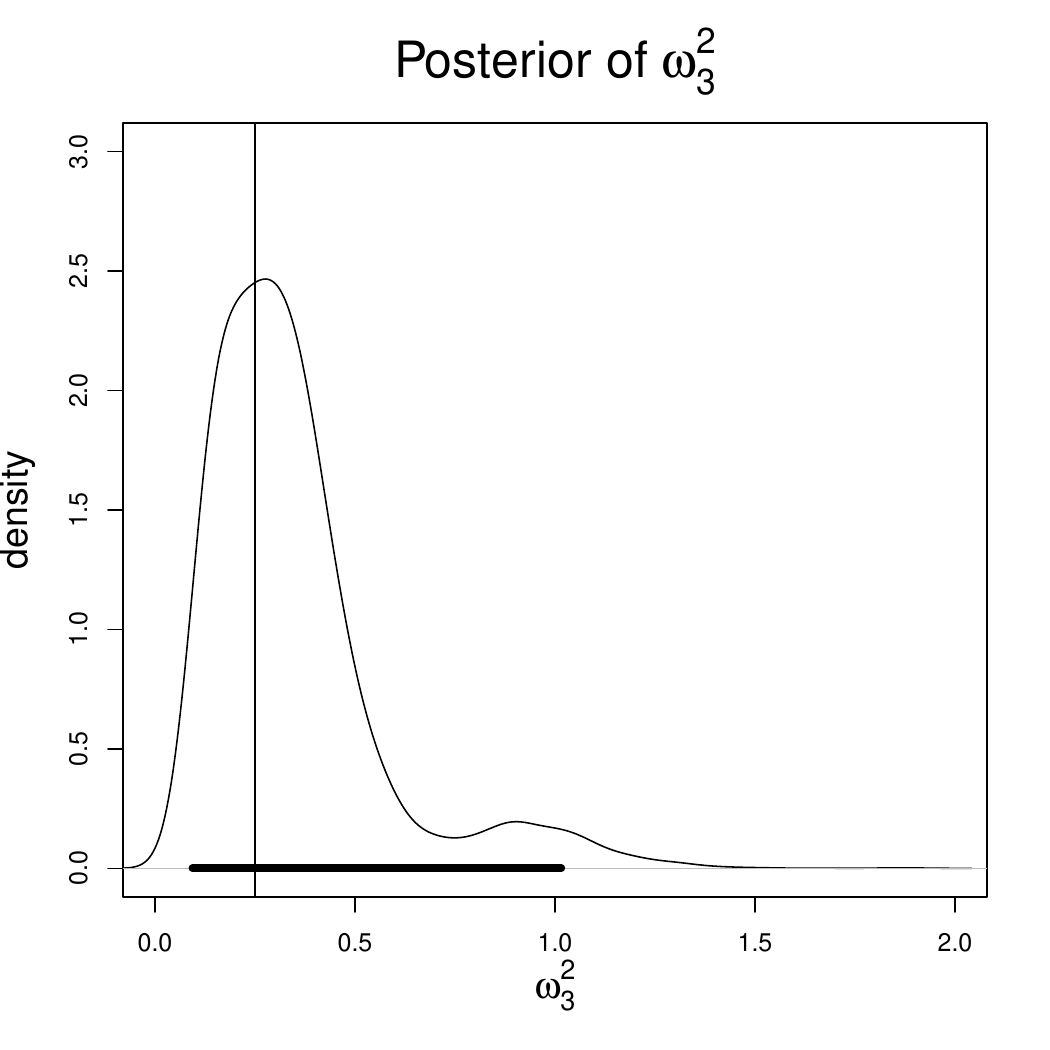}}\\
\vspace{2mm}
\subfigure[Posterior of $a_1$.]{ \label{fig:sim7_p1}
\includegraphics[width=4.5cm,height=5cm]{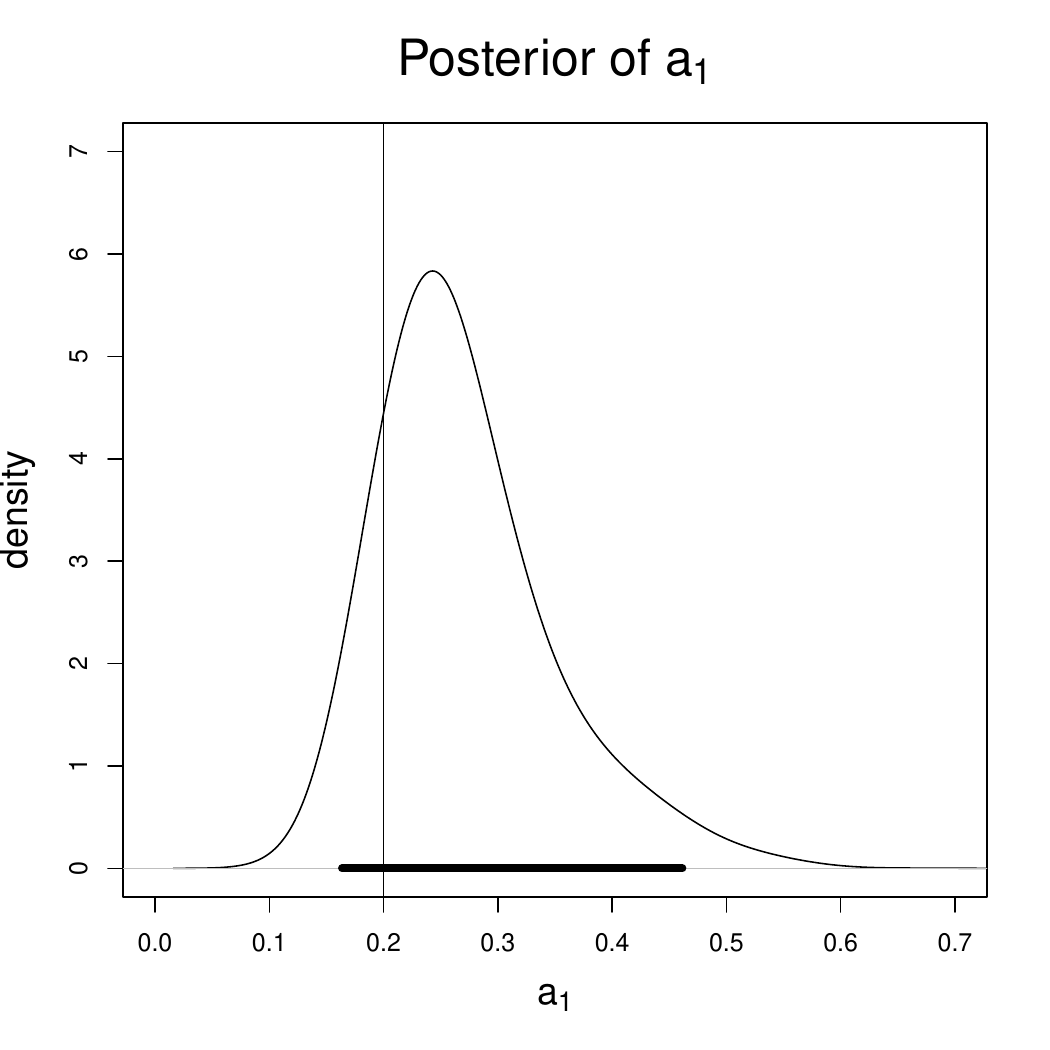}}
\hspace{2mm}
\subfigure[Posterior of $a_2$.]{ \label{fig:sim7_p2}
\includegraphics[width=4.5cm,height=5cm]{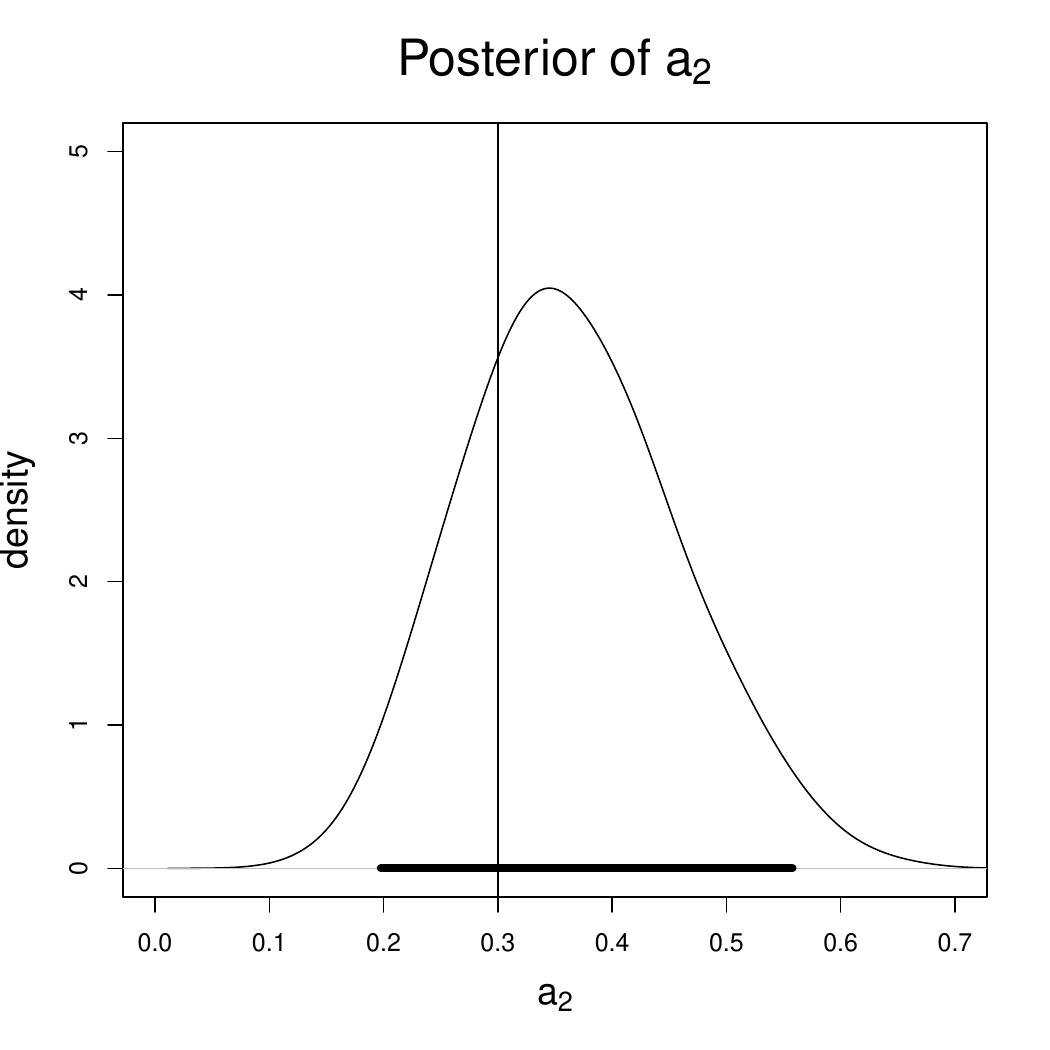}}
\hspace{2mm}
\subfigure[Posterior of $a_3$.]{ \label{fig:sim7_p3}
\includegraphics[width=4.5cm,height=5cm]{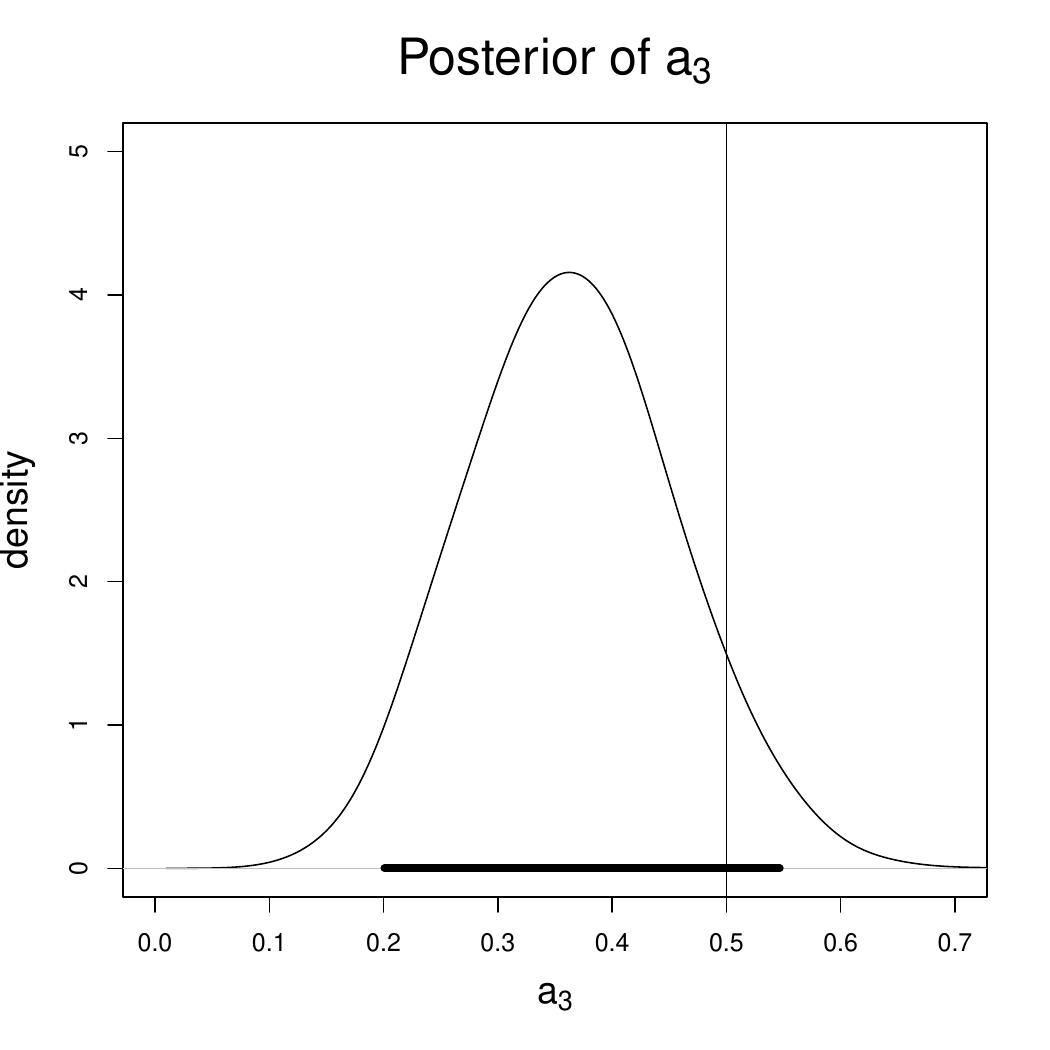}}
\caption{{\bf TTMCMC for $SDE_2$ and $\pi_5$:} Posteriors of $M$, $\mu_1$, $\mu_2$, $\mu_3$, $\omega^2_1$, $\omega^2_2$, $\omega^2_3$, $a_1$ $a_2$ and $a_3$. 
The vertical lines stand for the true values, while the thick horizontal lines denote the 95\% credible intervals.} 
\label{fig:sim7_posterior_plots}
\end{figure}

In this example, the mixture components are severely ill-separated, and Table 4 of \ctn{Maud16} shows that the $BIC$ approach in this case selects 
$1$ or $2$ mixture components, instead of the correct value $M_0=3$. However,
our Bayesian approach based on TTMCMC again successfully captures the correct number of components. Indeed, $M=1,2,3,4$ received the posterior
probabilities $0.0641$, $0.295$, $0.4668$ and $0.1505$, respectively. 
Figures \ref{fig:sim7_trace_plots} and \ref{fig:sim7_posterior_plots} demonstrates quite adequate performance of our approach.
The overall acceptance rate here is about $5.9\%$.

\subsection{Eighth simulation study: $SDE_3$ with $\pi_1$} 
\label{subsec:sim8}

\begin{figure}
\centering
\subfigure[Trace plot of $M$.]{ \label{fig:sim8_trace_comp}
\includegraphics[width=7cm,height=5cm]{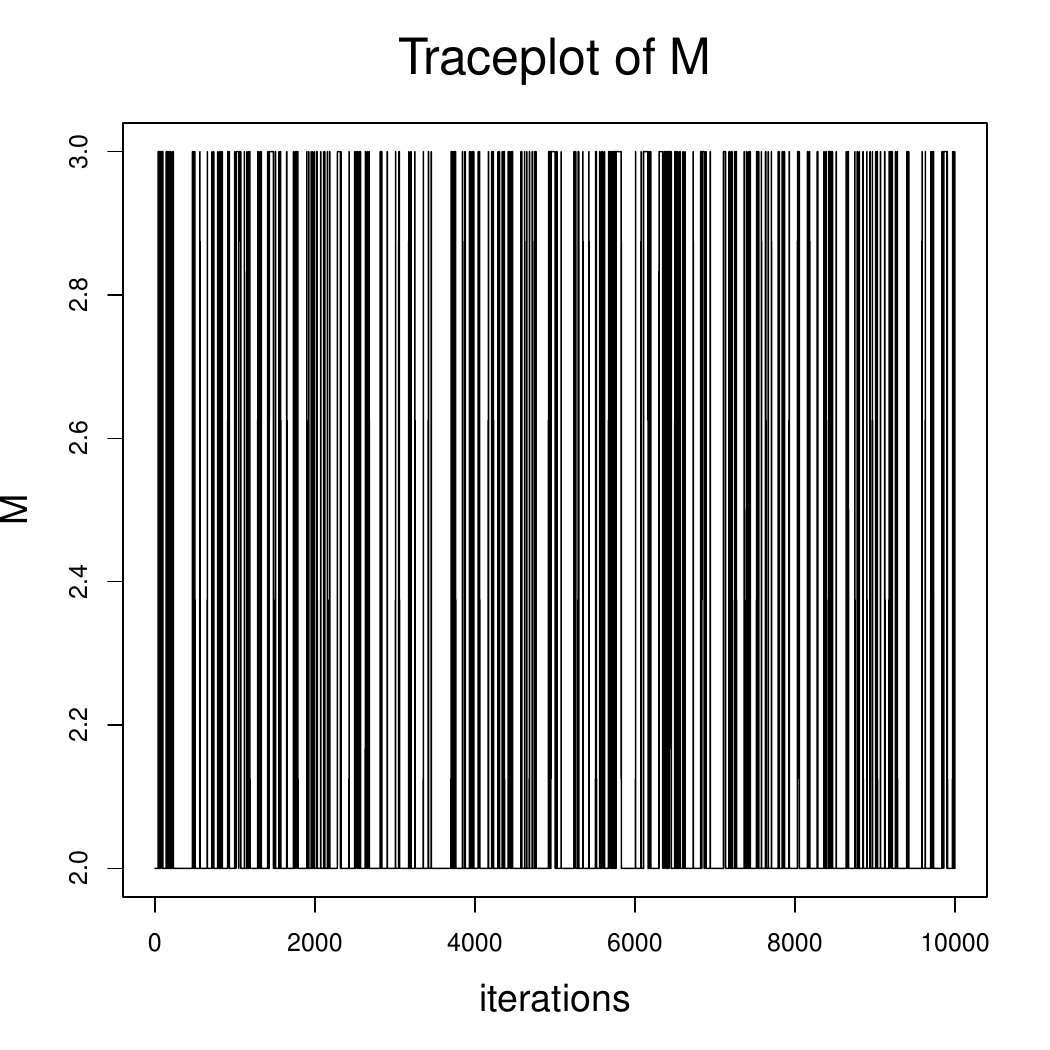}}
\hspace{2mm}
\subfigure[Trace plot of $\mu_1$.]{ \label{fig:sim8_trace_mu1}
\includegraphics[width=7cm,height=5cm]{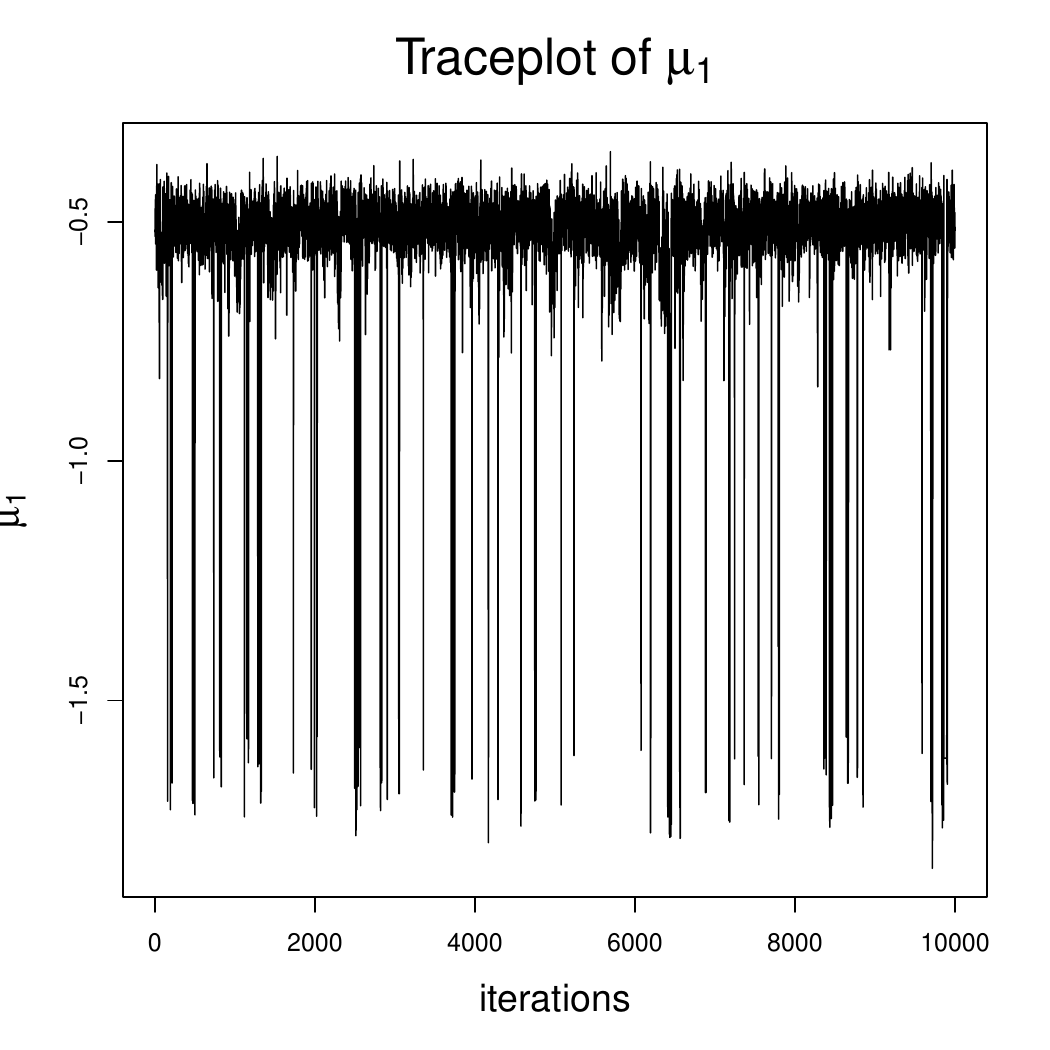}}\\
\vspace{2mm}
\subfigure[Trace plot of $\mu_2$.]{ \label{fig:sim8_trace_mu2}
\includegraphics[width=7cm,height=5cm]{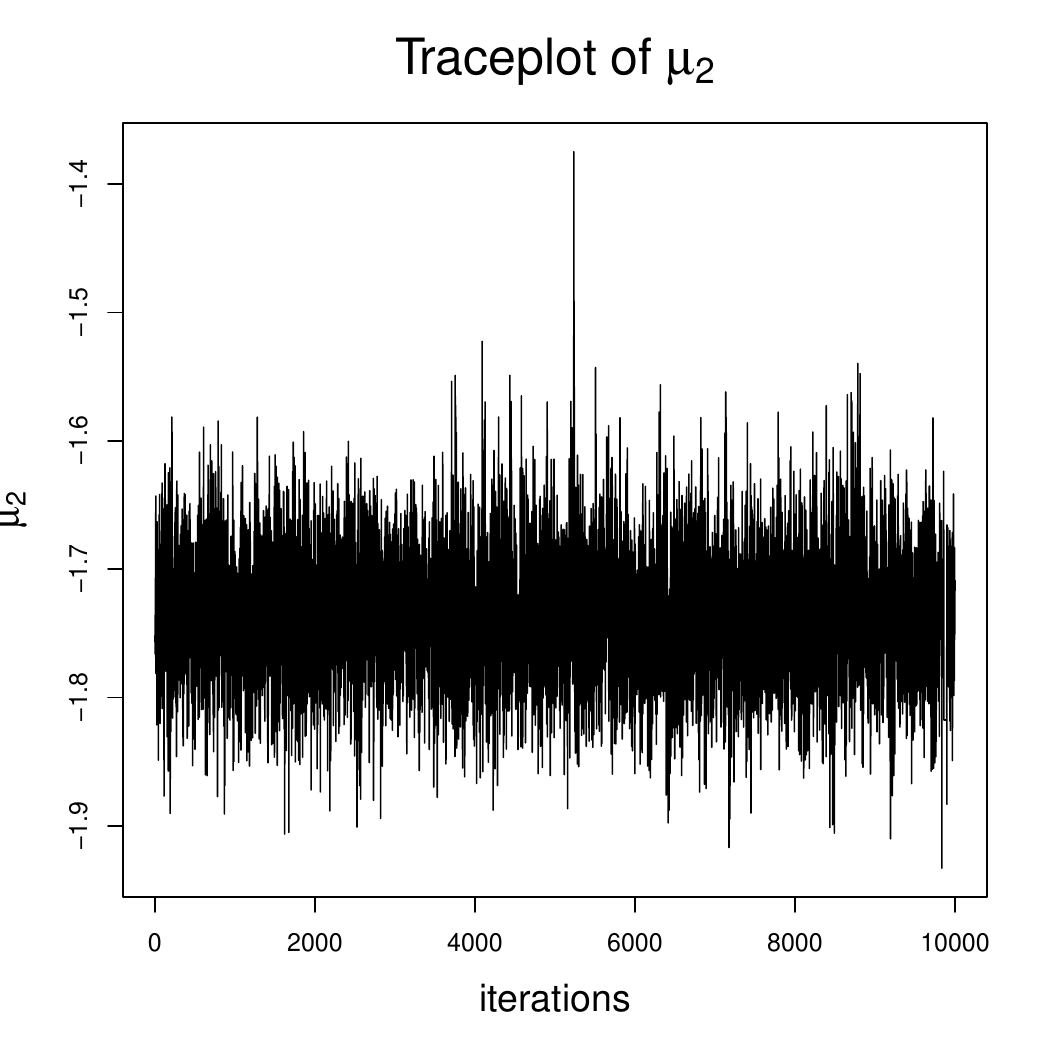}}
\vspace{2mm}
\subfigure[Trace plot of $\omega^2_1$.]{ \label{fig:sim8_trace_omegasq1}
\includegraphics[width=7cm,height=5cm]{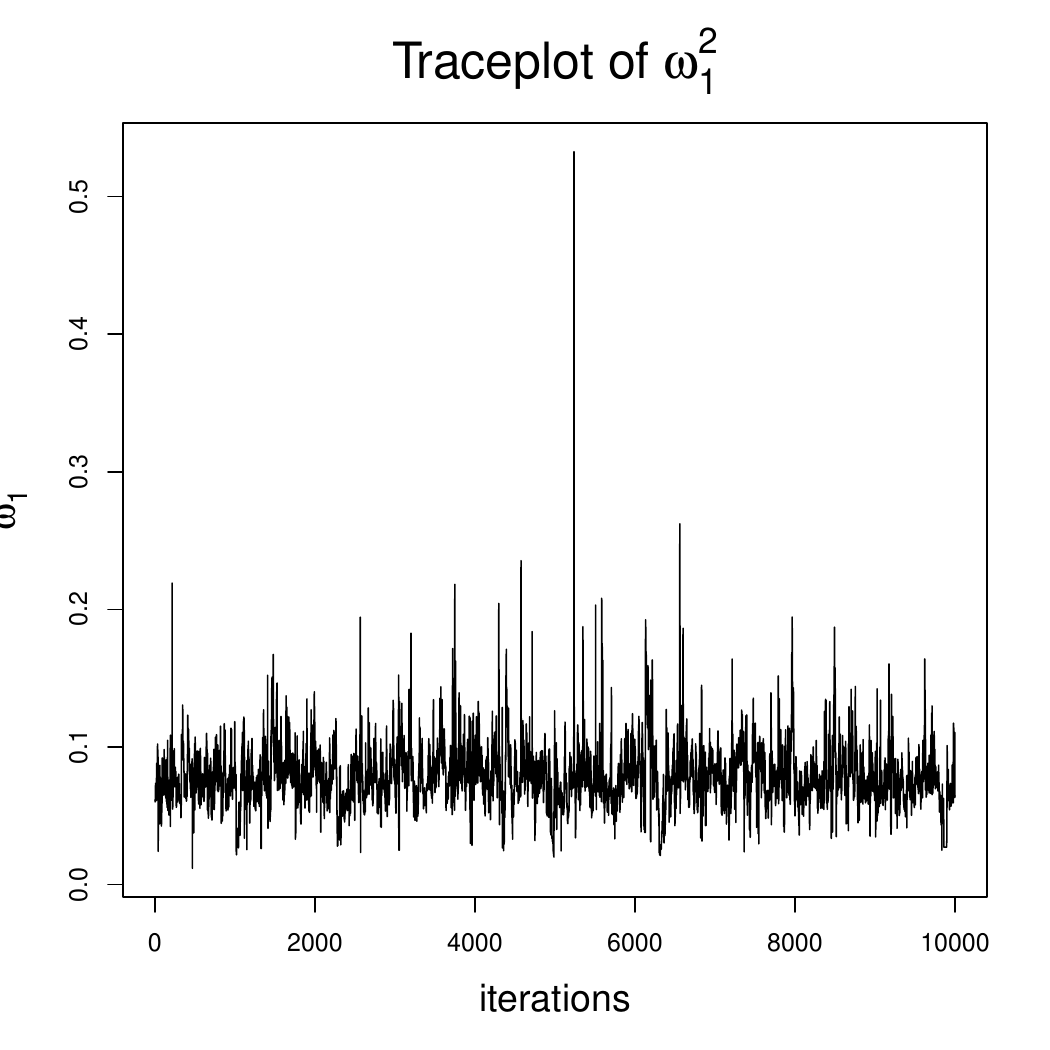}}\\
\vspace{2mm}
\subfigure[Trace plot of $\omega^2_2$.]{ \label{fig:sim8_trace_omegasq2}
\includegraphics[width=7cm,height=5cm]{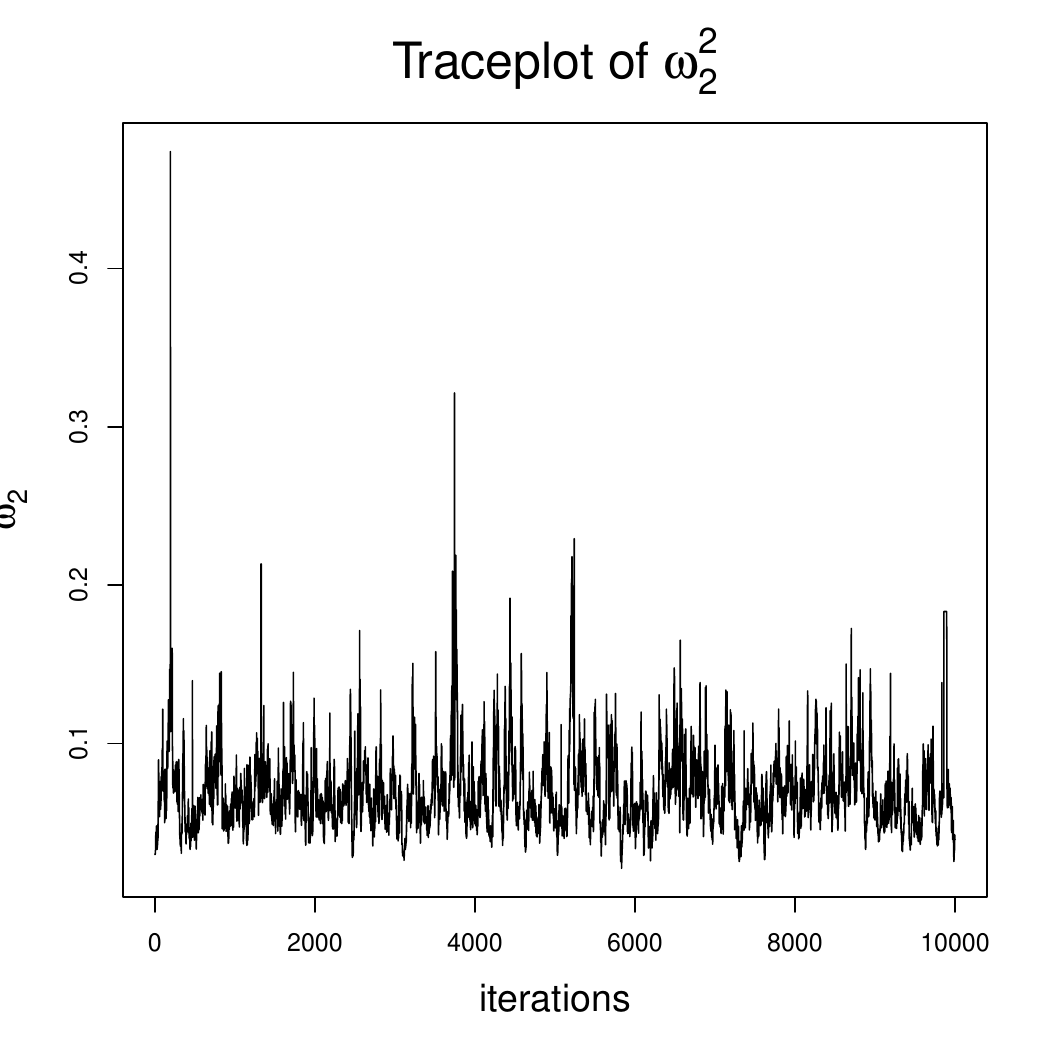}}
\hspace{2mm}
\subfigure[Trace plot of $a_1$.]{ \label{fig:sim8_trace_p1}
\includegraphics[width=7cm,height=5cm]{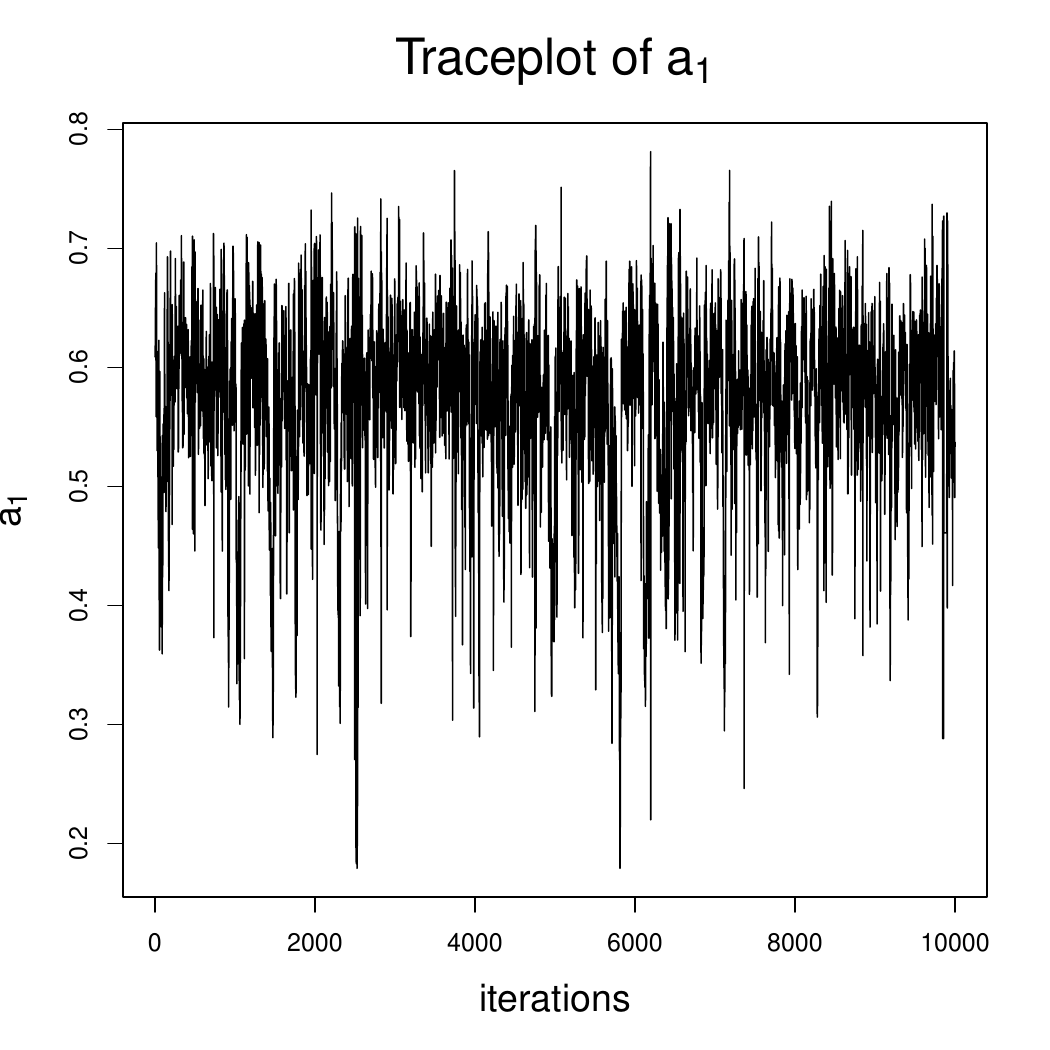}}\\
\vspace{2mm}
\subfigure[Trace plot of $a_2$.]{ \label{fig:sim8_trace_p2}
\includegraphics[width=7cm,height=5cm]{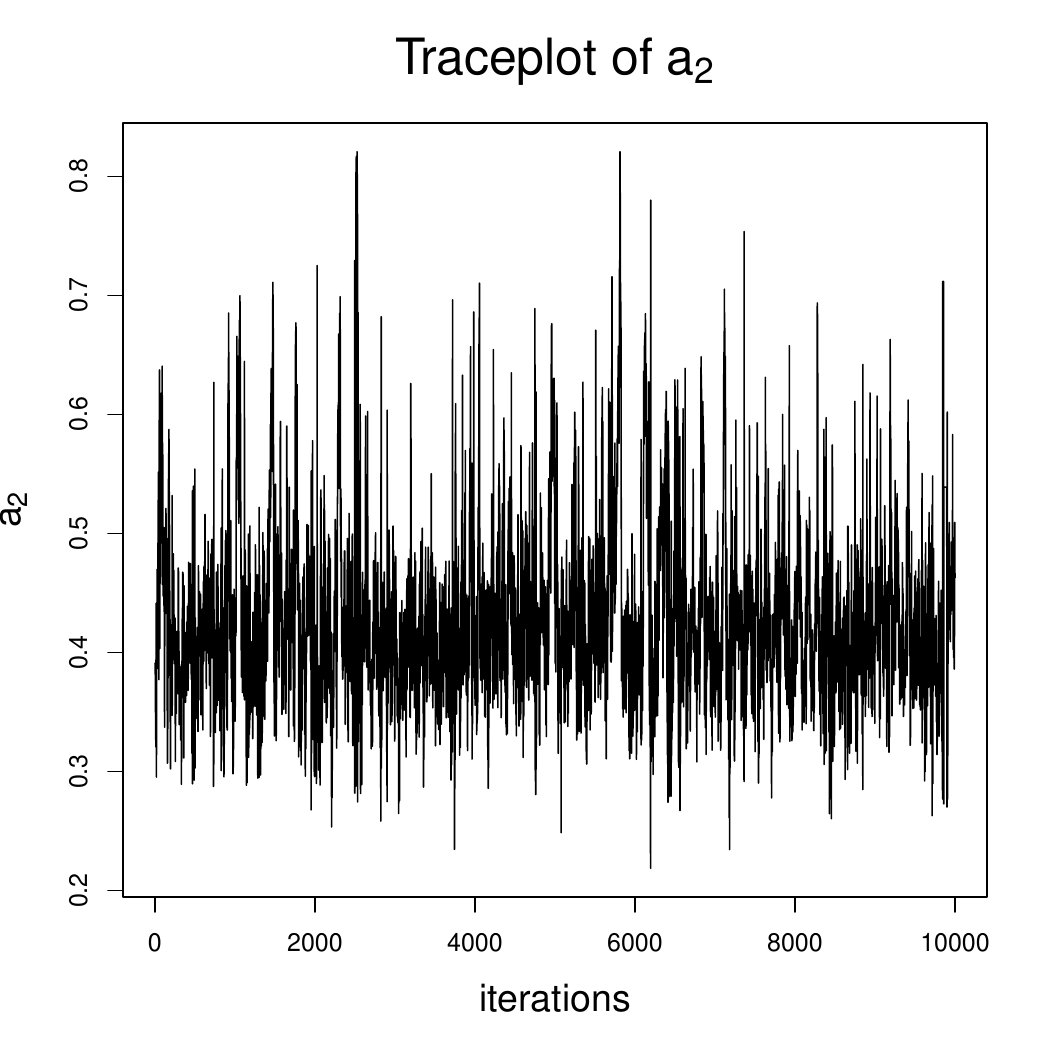}}
\caption{{\bf TTMCMC for $SDE_3$ and $\pi_1$:} Trace plots of $M$, $\mu_1$, $\mu_2$, $\omega^2_1$, $\omega^2_2$, $a_1$ and $a_2$.} 
\label{fig:sim8_trace_plots}
\end{figure}

\begin{figure}
\centering
\subfigure[Posterior of $\mu_1$.]{ \label{fig:sim8_mu1}
\includegraphics[width=7cm,height=6cm]{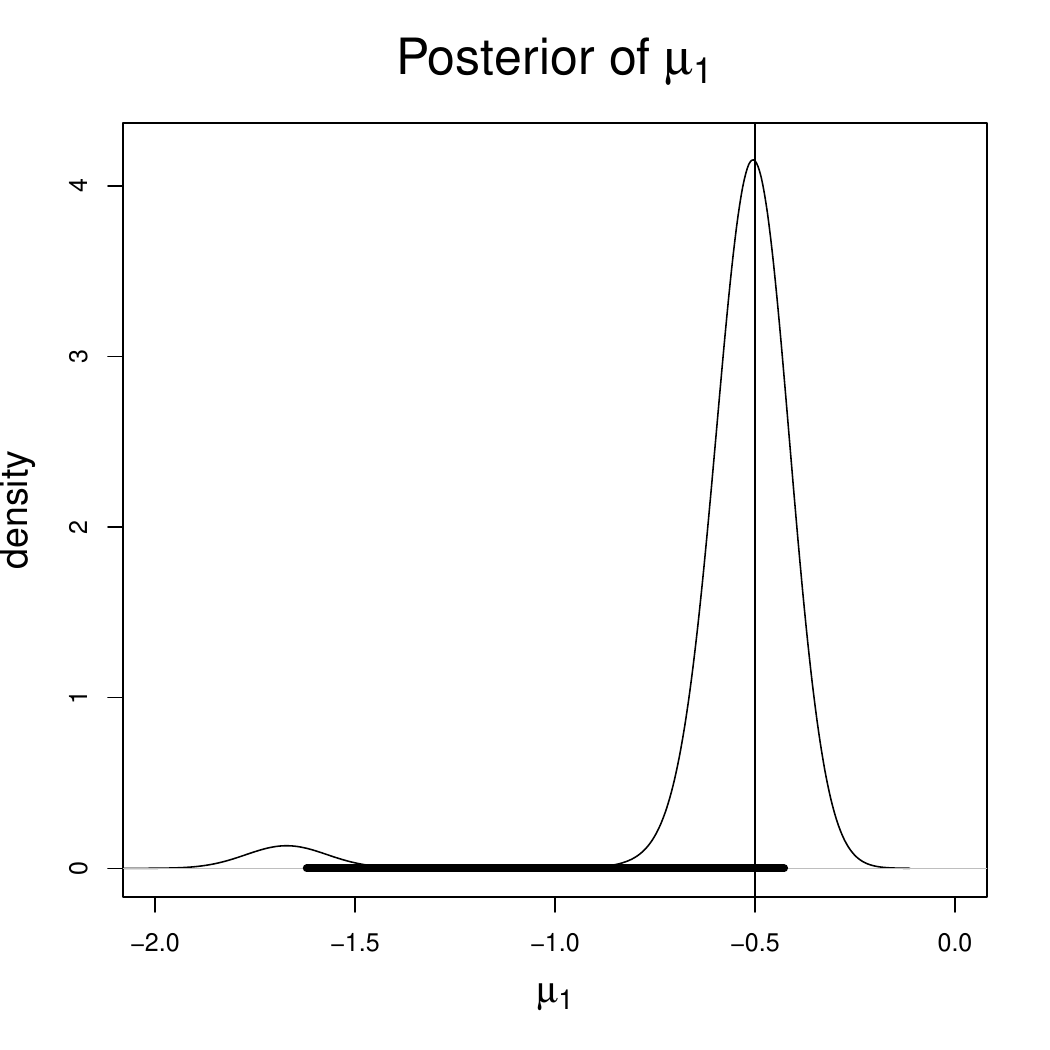}}
\hspace{2mm}
\subfigure[Posterior of $\mu_2$.]{ \label{fig:sim8_mu2}
\includegraphics[width=7cm,height=6cm]{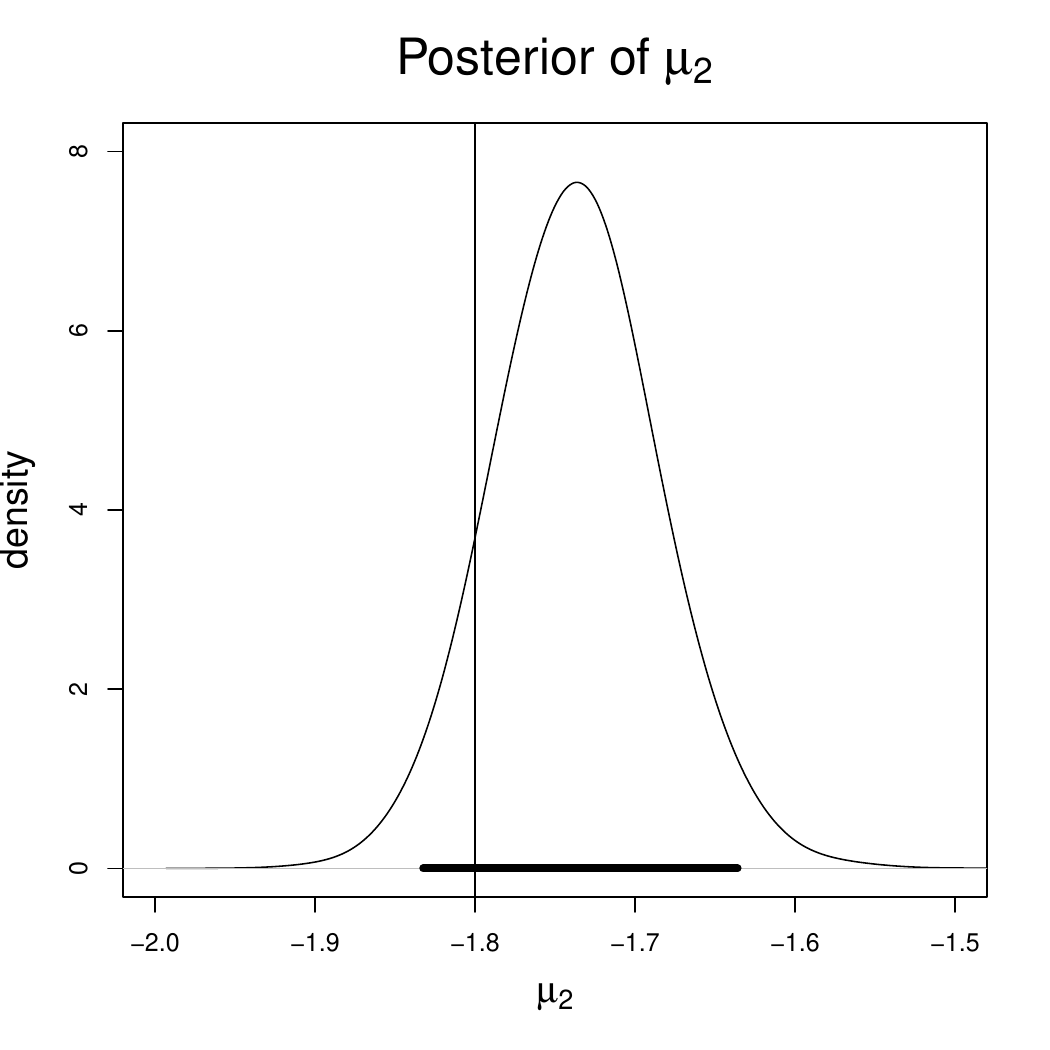}}\\
\vspace{2mm}
\subfigure[Posterior of $\omega^2_1$.]{ \label{fig:sim8_omegasq1}
\includegraphics[width=7cm,height=6cm]{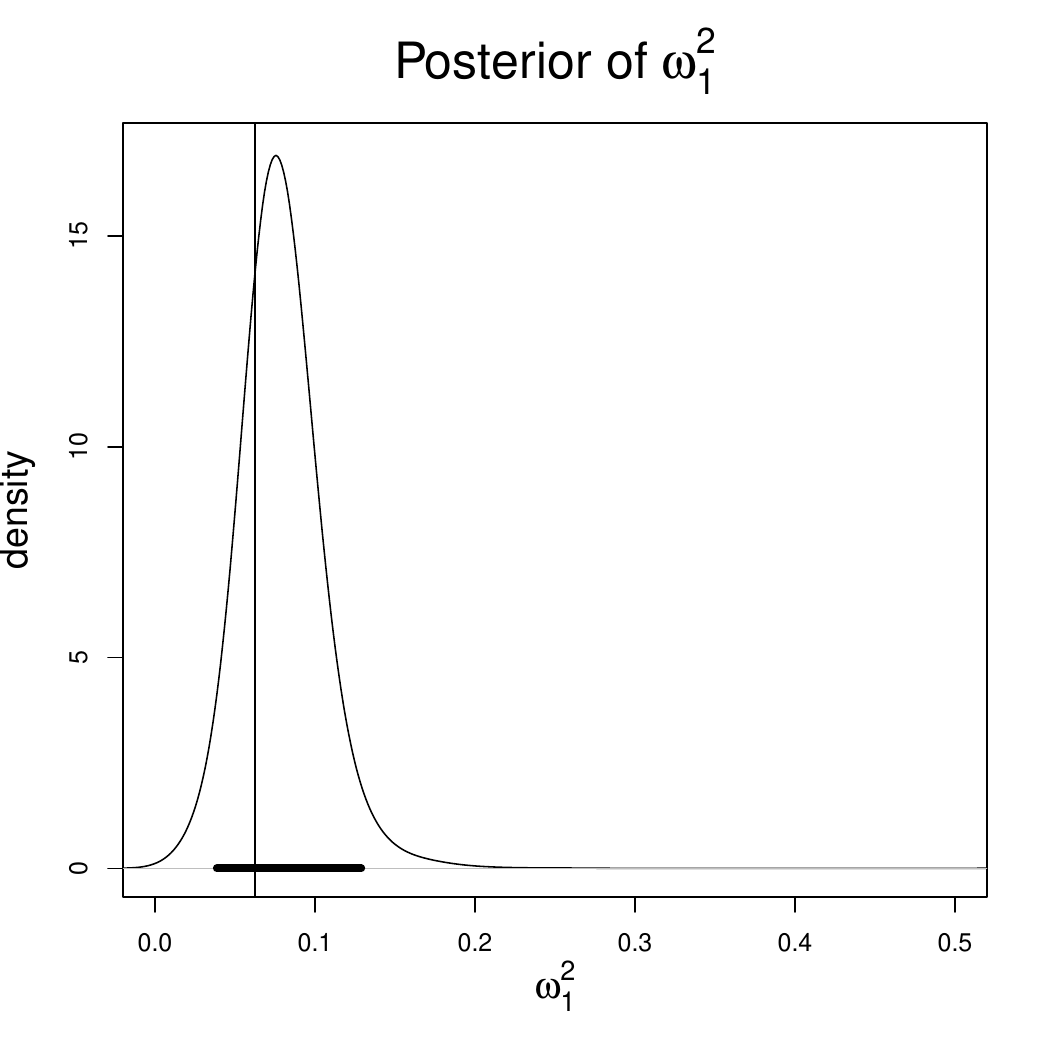}}
\hspace{2mm}
\subfigure[Posterior of $\omega^2_2$.]{ \label{fig:sim8_omegasq2}
\includegraphics[width=7cm,height=6cm]{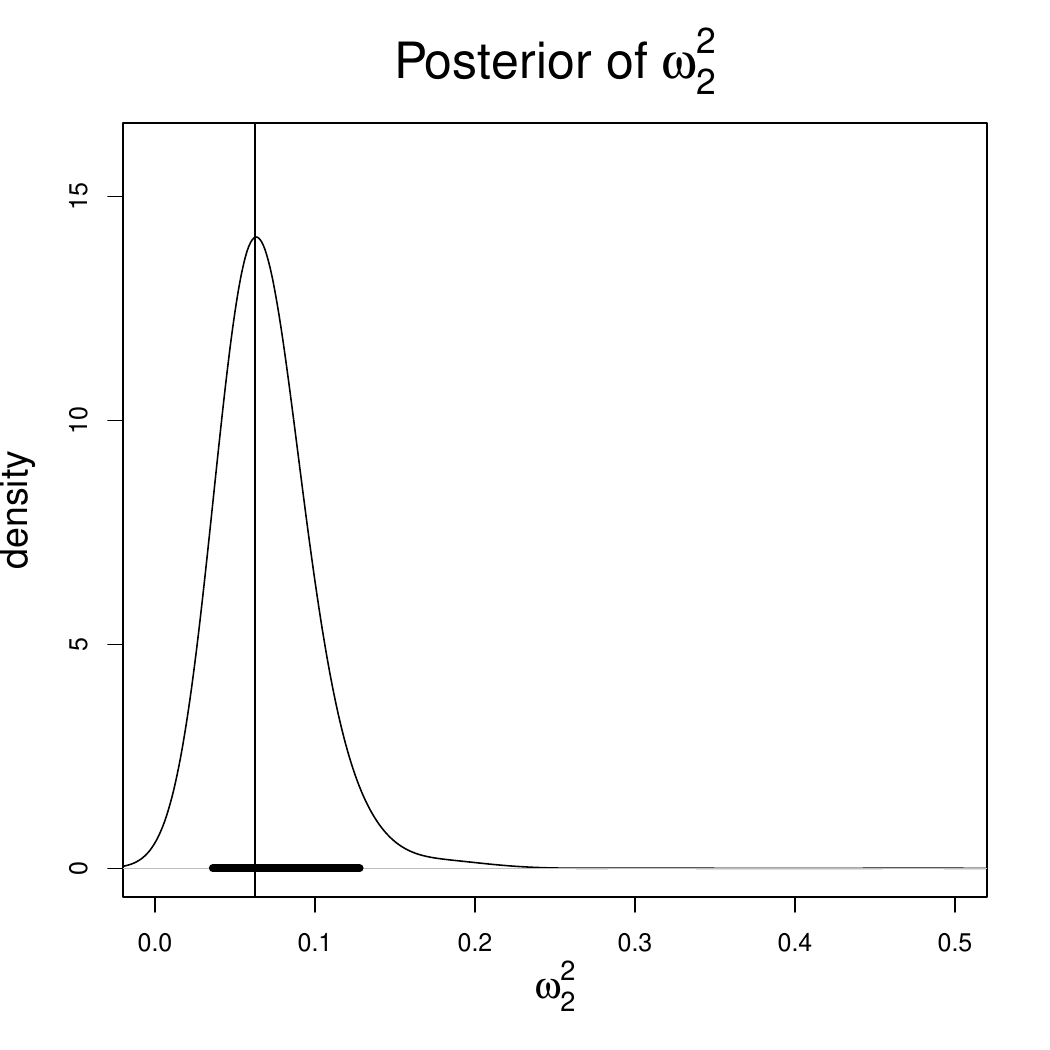}}\\
\vspace{2mm}
\subfigure[Posterior of $a_1$.]{ \label{fig:sim8_p1}
\includegraphics[width=7cm,height=6cm]{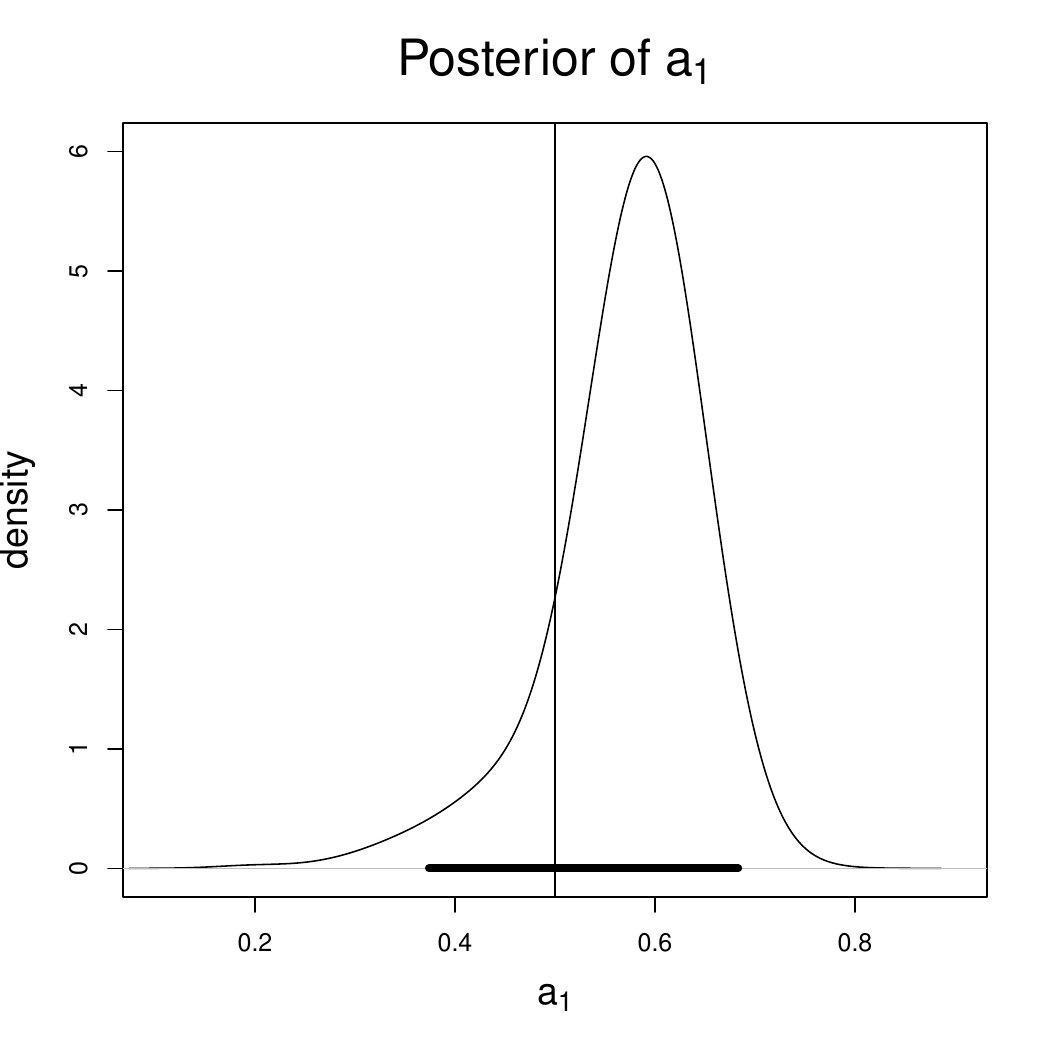}}
\hspace{2mm}
\subfigure[Posterior of $a_2$.]{ \label{fig:sim8_p2}
\includegraphics[width=7cm,height=6cm]{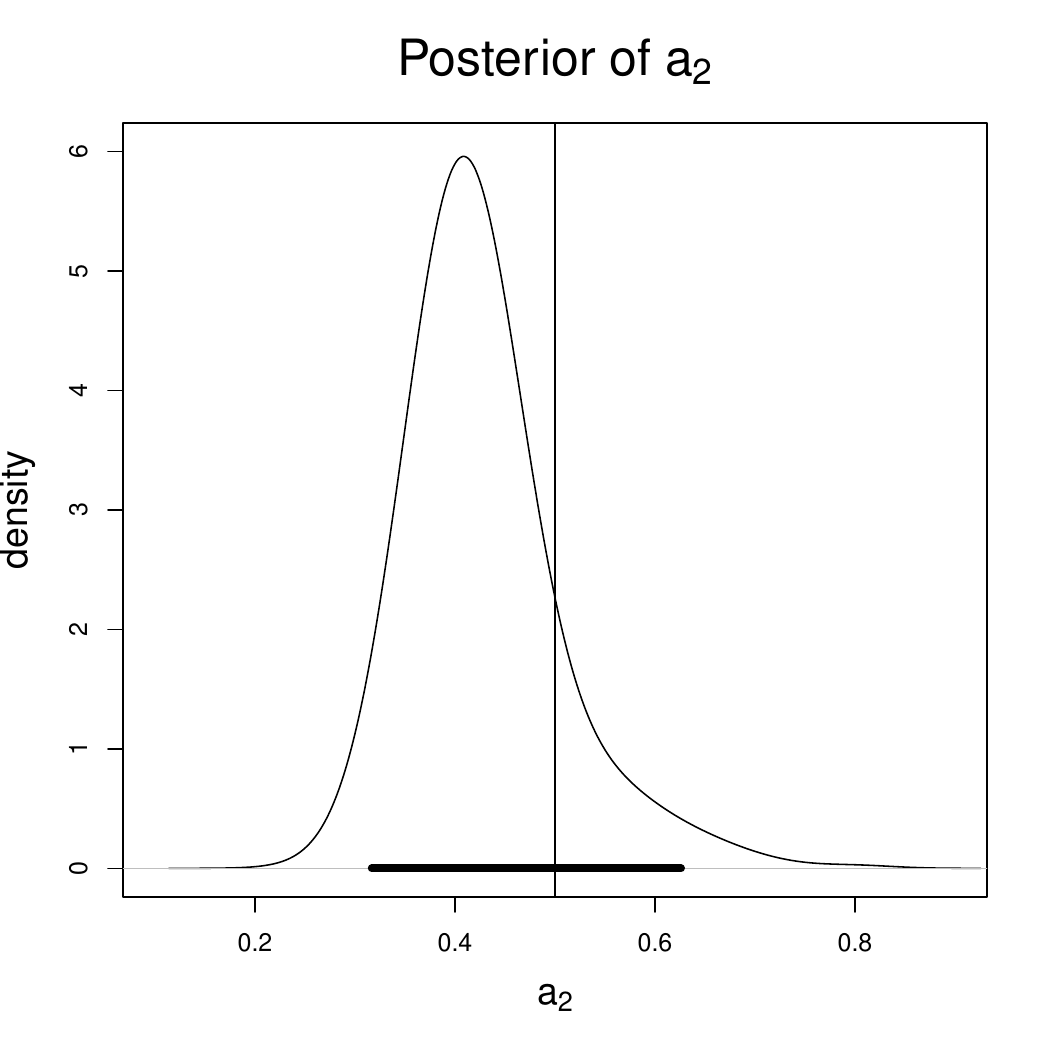}}
\caption{{\bf TTMCMC for $SDE_3$ and $\pi_1$:} Posteriors of $M$, $\mu_1$, $\mu_2$, $\omega^2_1$, $\omega^2_2$, $a_1$ and $a_2$. The vertical lines stand
for the true values, while the thick horizontal lines denote the 95\% credible intervals.} 
\label{fig:sim8_posterior_plots}
\end{figure}

In this case as well, excellent mixing is vindicated by the trace plots shown in Figure \ref{fig:sim8_trace_plots}. Furthermore,
as shown in Figure \ref{fig:sim8_posterior_plots}, all the true values well-captured by our Bayesian approach. 
The posterior probabilities of $M=2$ and $M=3$ turned out to be $0.8629$ and $0.1371$, so that all other values of $M$ have zero posterior mass. 
The overall acceptance rate of our algorithm in this case is about $2.5\%$.

\subsection{Ninth simulation study: $SDE_3$ with $\pi_2$} 
\label{subsec:sim9}

\begin{figure}
\centering
\subfigure[Trace plot of $M$.]{ \label{fig:sim9_trace_comp}
\includegraphics[width=7cm,height=5cm]{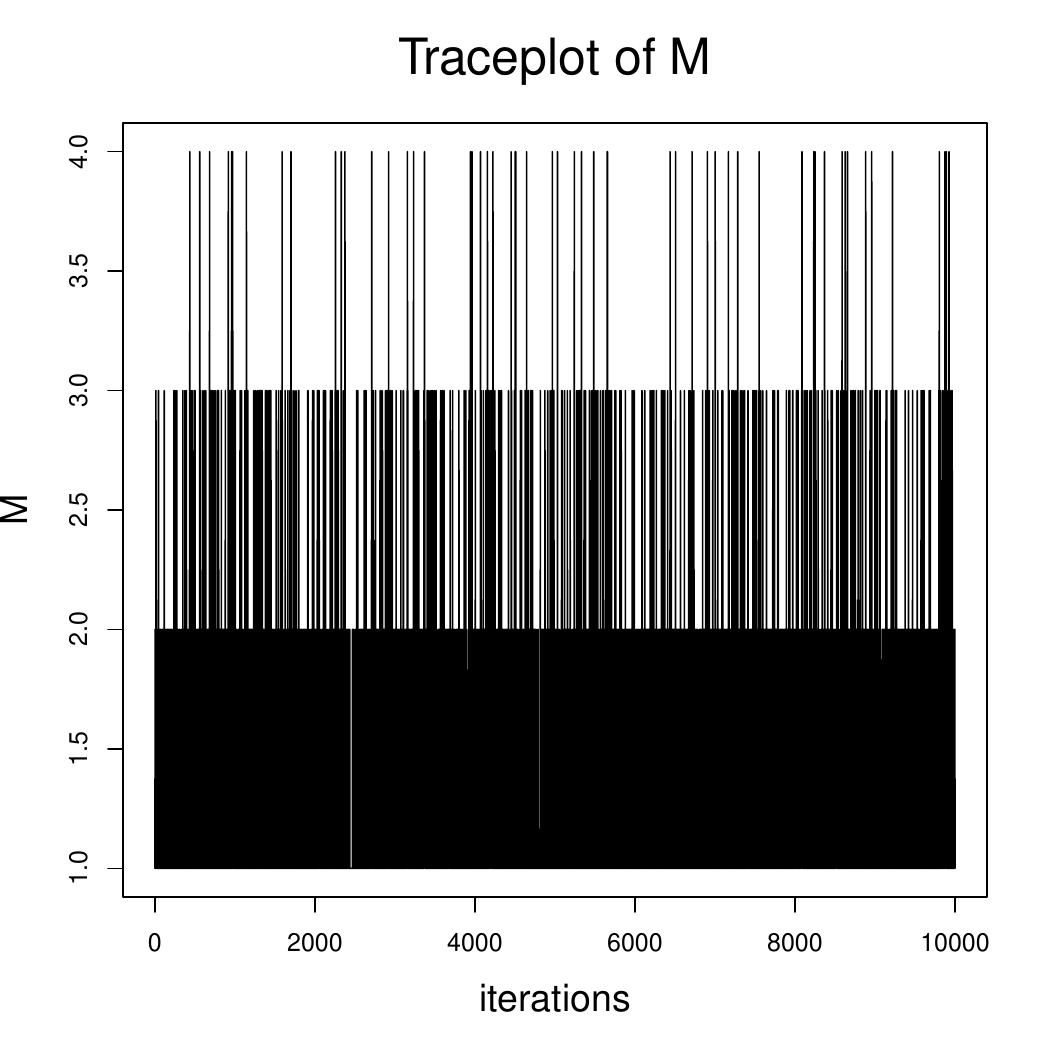}}
\hspace{2mm}
\subfigure[Trace plot of $\mu_1$.]{ \label{fig:sim9_trace_mu1}
\includegraphics[width=7cm,height=5cm]{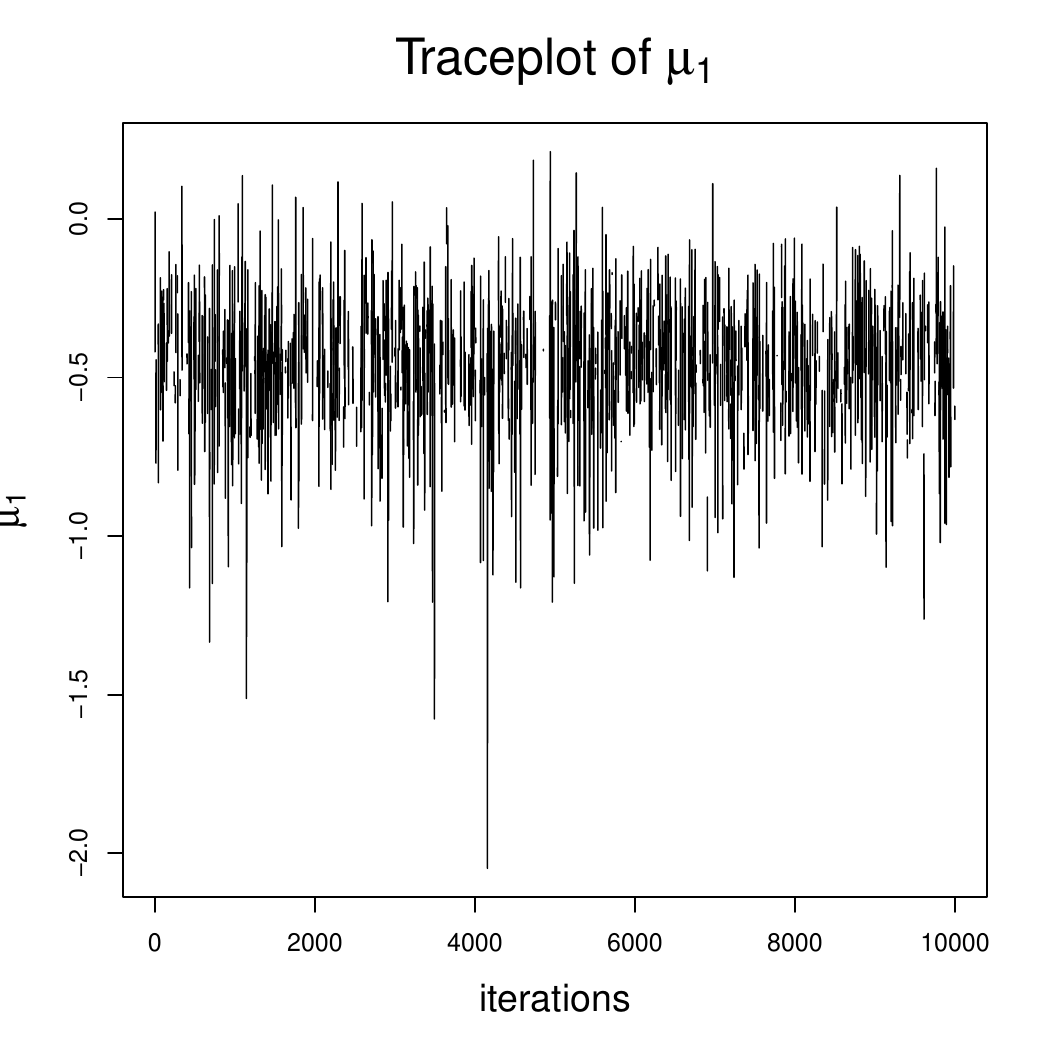}}\\
\vspace{2mm}
\subfigure[Trace plot of $\mu_2$.]{ \label{fig:sim9_trace_mu2}
\includegraphics[width=7cm,height=5cm]{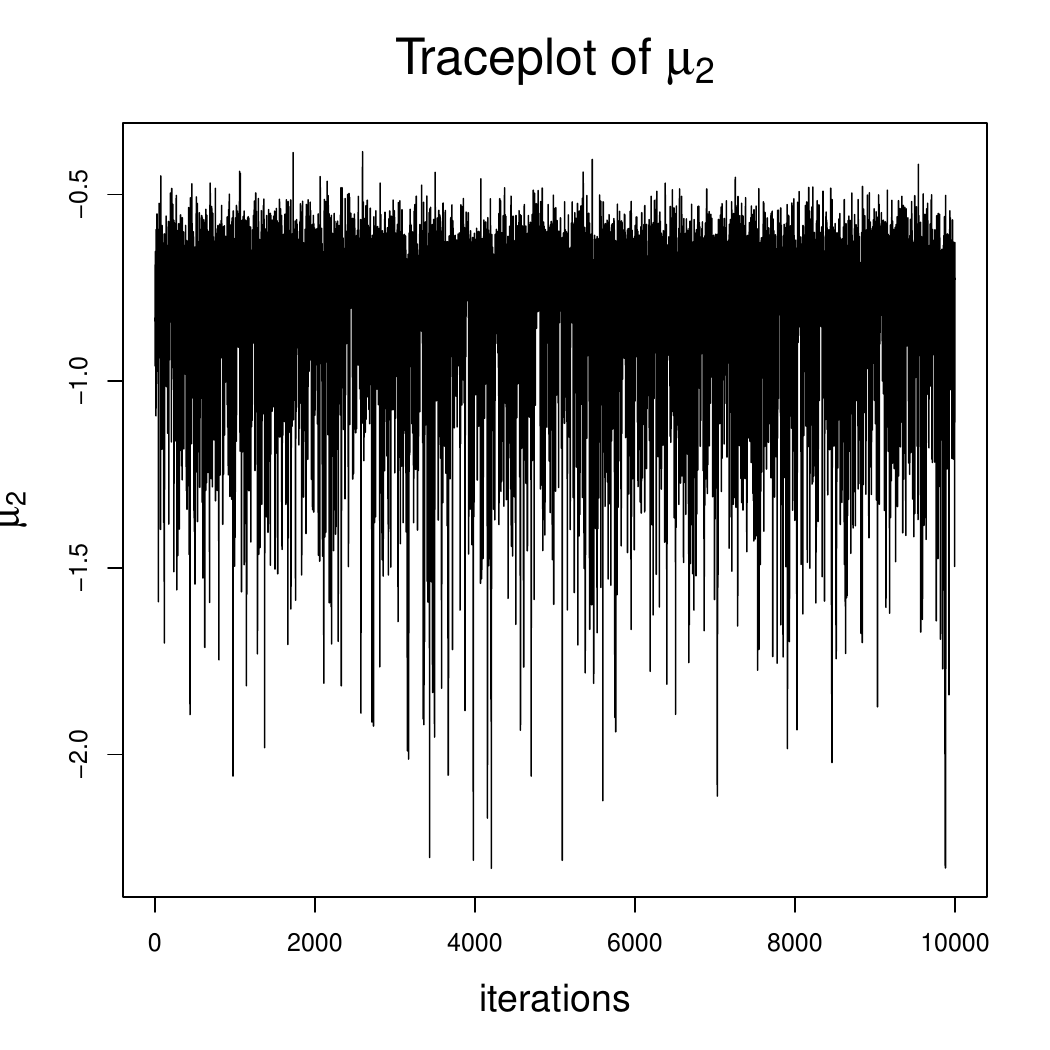}}
\vspace{2mm}
\subfigure[Trace plot of $\omega^2_1$.]{ \label{fig:sim9_trace_omegasq1}
\includegraphics[width=7cm,height=5cm]{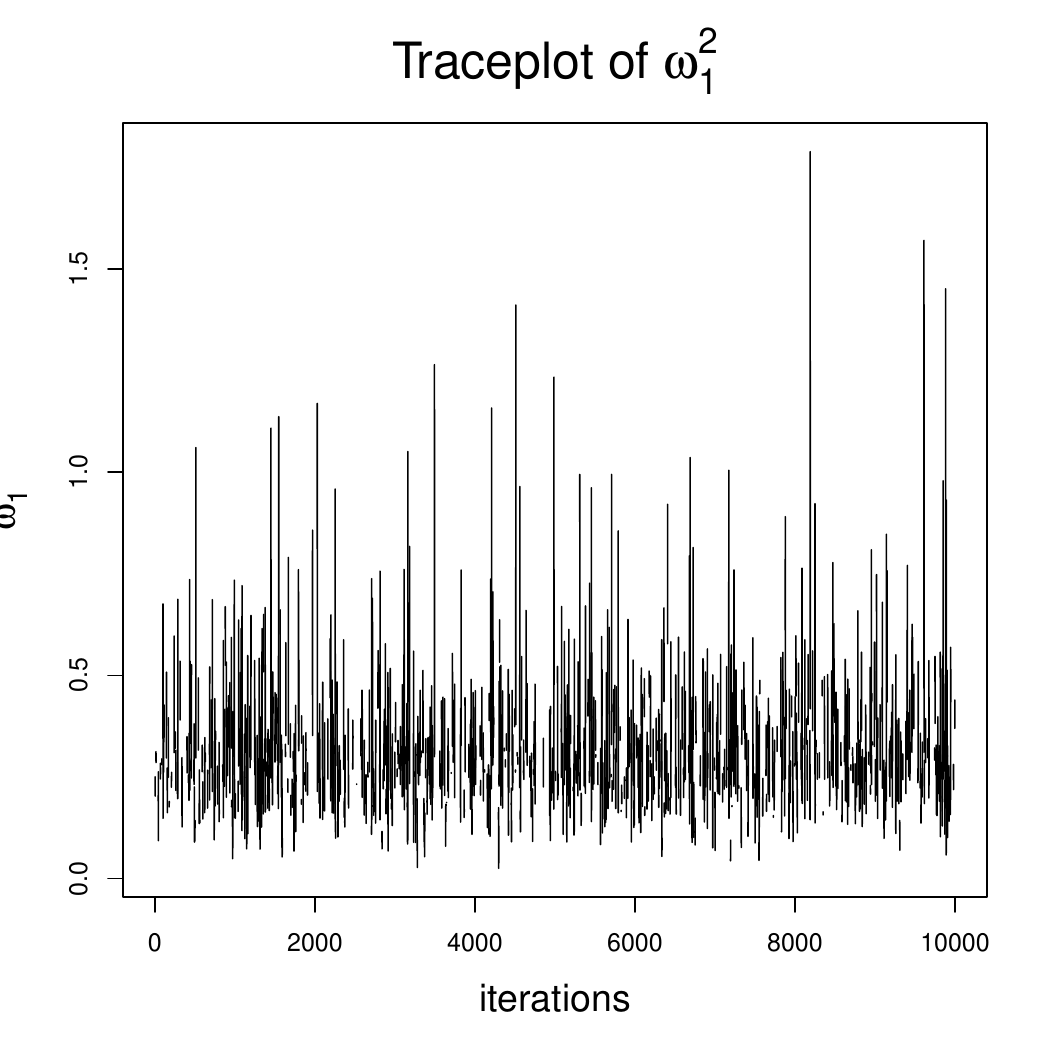}}\\
\vspace{2mm}
\subfigure[Trace plot of $\omega^2_2$.]{ \label{fig:sim9_trace_omegasq2}
\includegraphics[width=7cm,height=5cm]{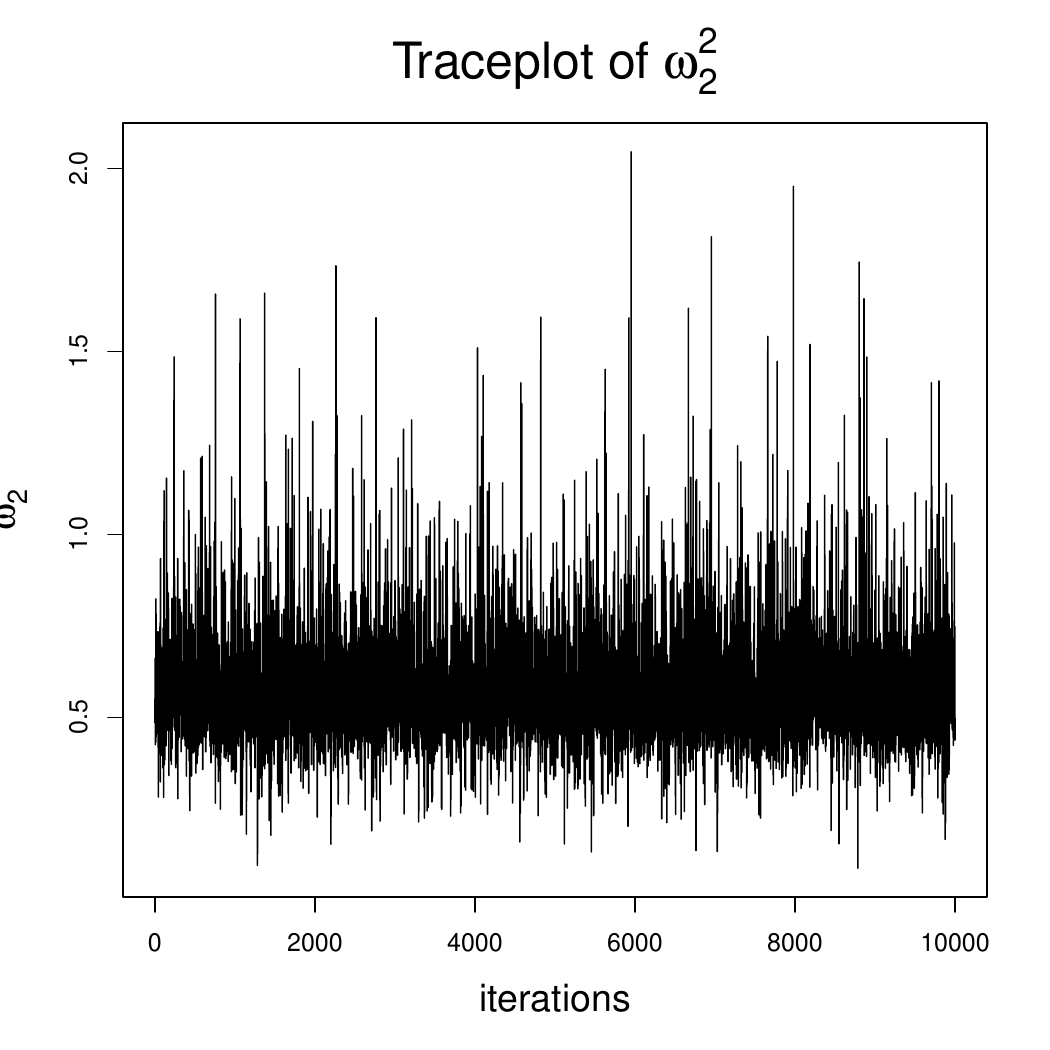}}
\hspace{2mm}
\subfigure[Trace plot of $a_1$.]{ \label{fig:sim9_trace_p1}
\includegraphics[width=7cm,height=5cm]{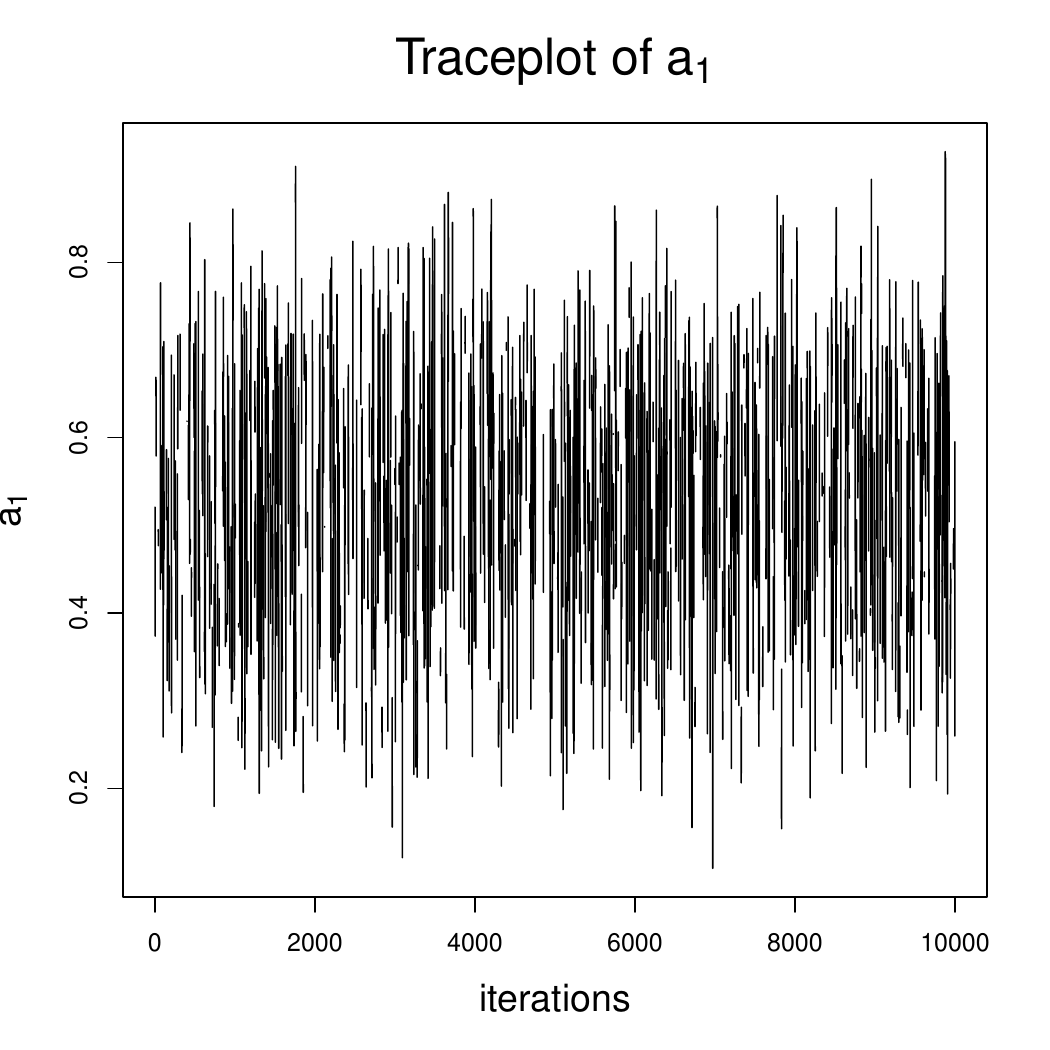}}\\
\vspace{2mm}
\subfigure[Trace plot of $a_2$.]{ \label{fig:sim9_trace_p2}
\includegraphics[width=7cm,height=5cm]{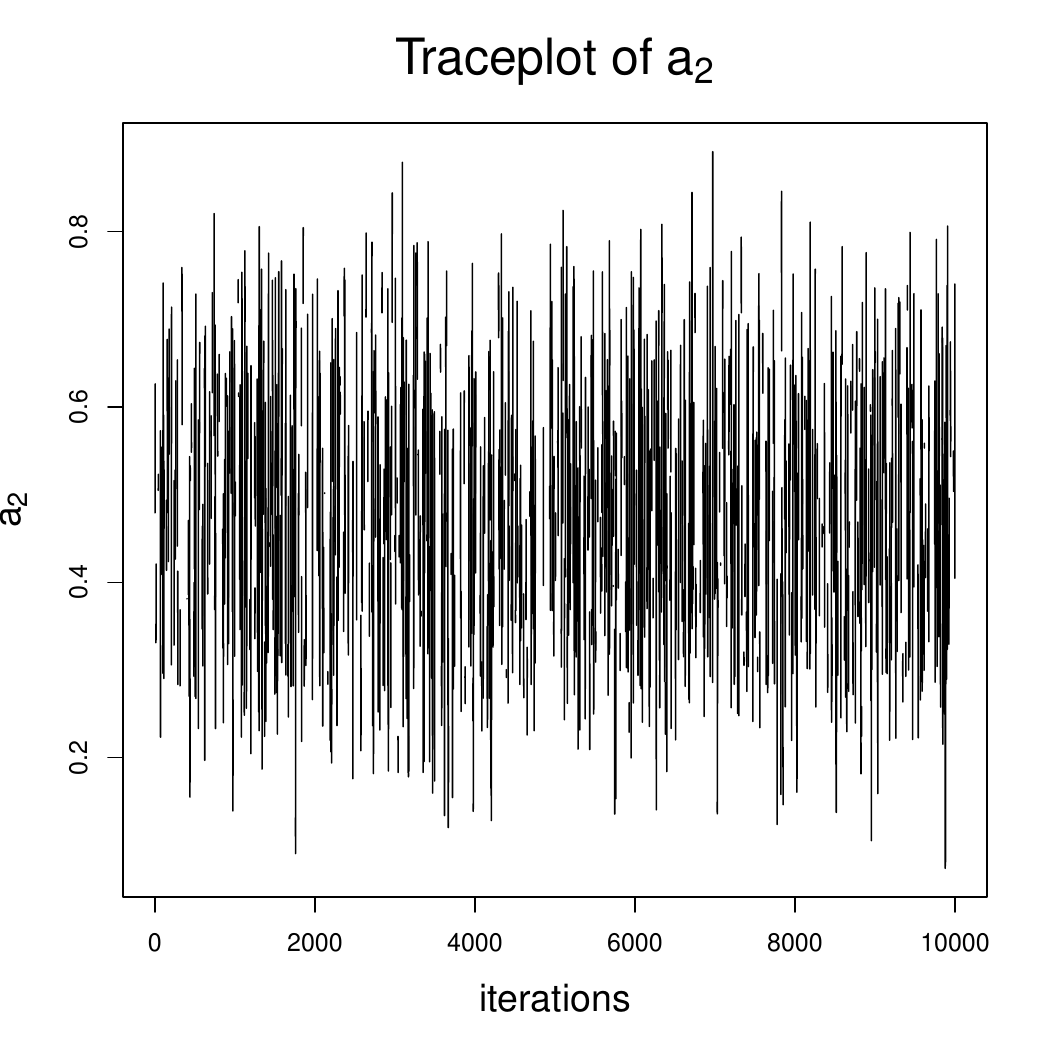}}
\caption{{\bf TTMCMC for $SDE_3$ and $\pi_2$ with $n=100$:} Trace plots of $M$, $\mu_1$, $\mu_2$, $\omega^2_1$, $\omega^2_2$, $a_1$ and $a_2$.} 
\label{fig:sim9_trace_plots}
\end{figure}

\begin{figure}
\centering
\subfigure[Posterior of $\mu_1$.]{ \label{fig:sim9_mu1}
\includegraphics[width=7cm,height=6cm]{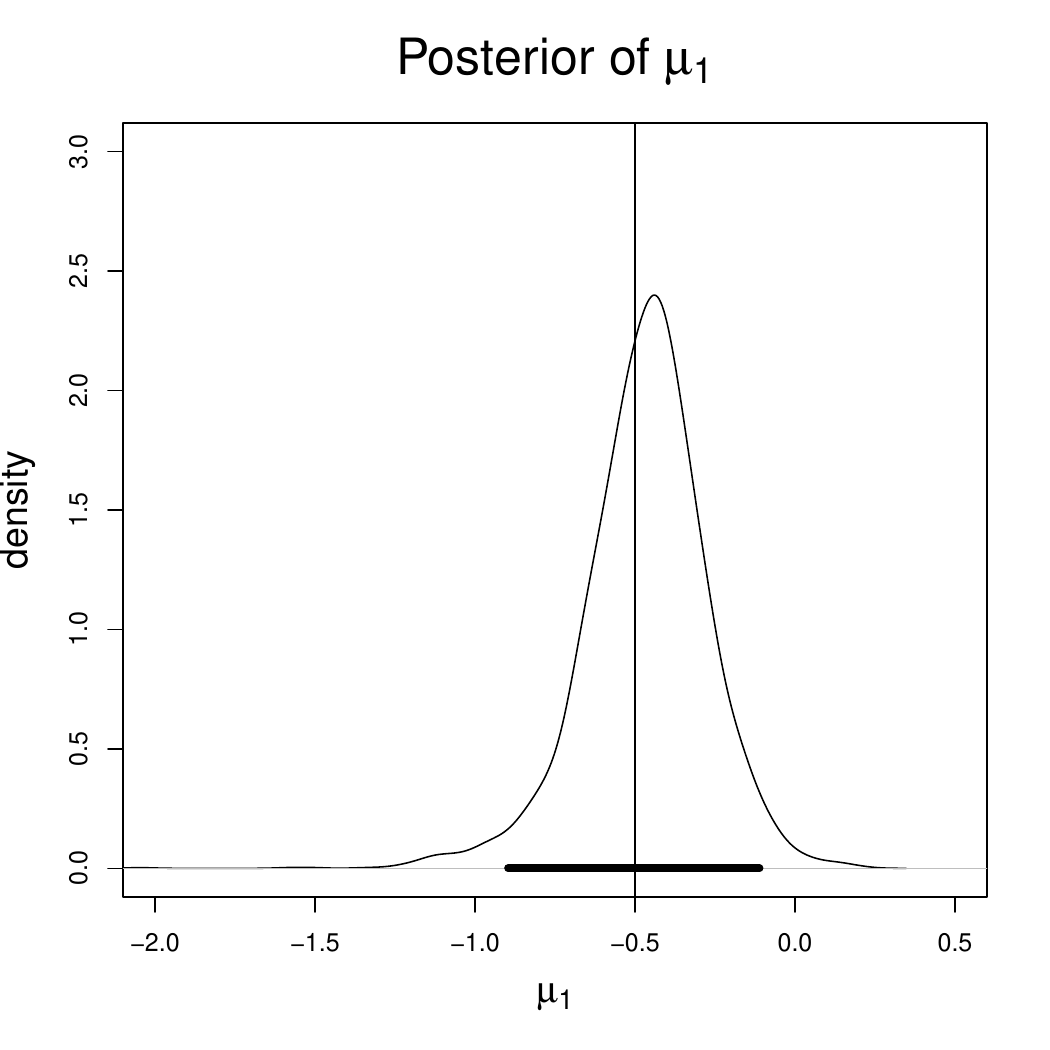}}
\hspace{2mm}
\subfigure[Posterior of $\mu_2$.]{ \label{fig:sim9_mu2}
\includegraphics[width=7cm,height=6cm]{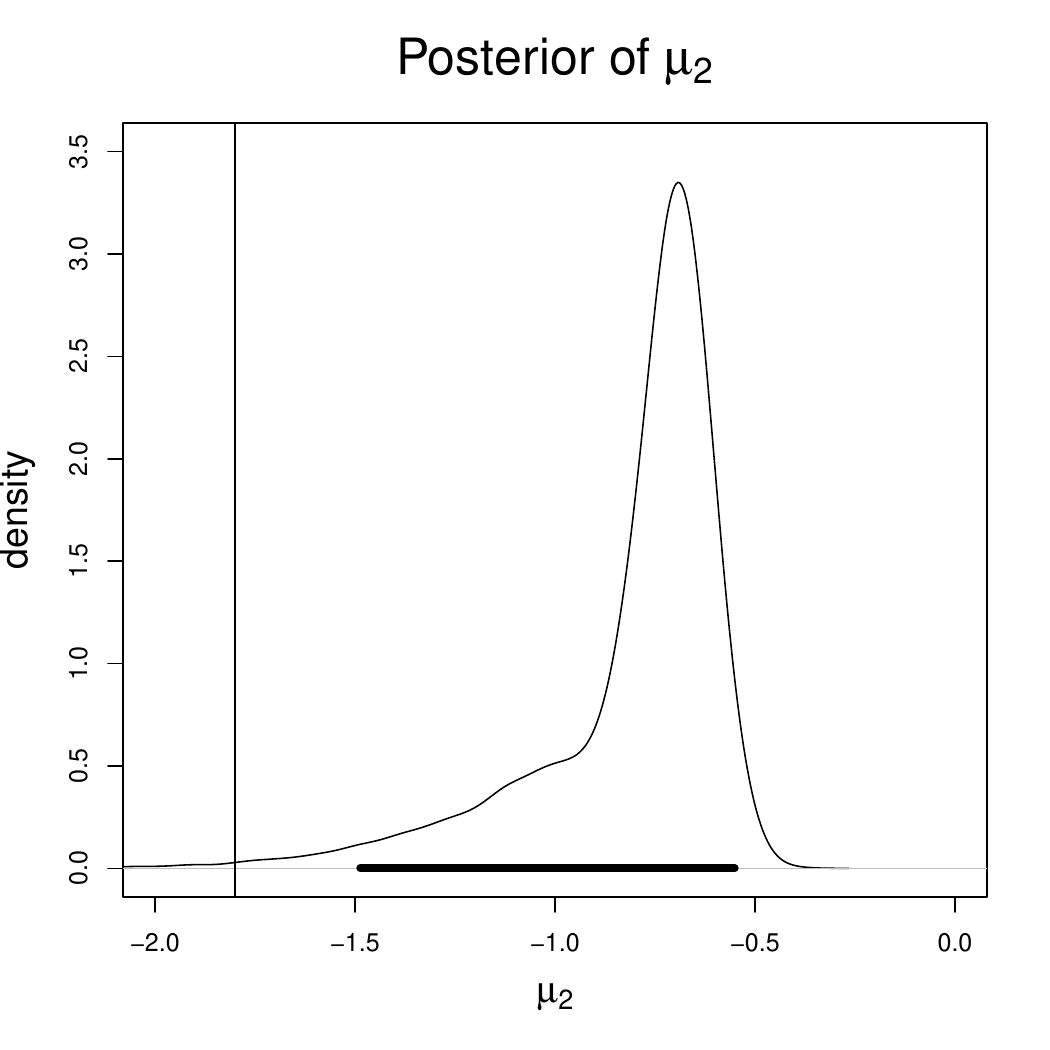}}\\
\vspace{2mm}
\subfigure[Posterior of $\omega^2_1$.]{ \label{fig:sim9_omegasq1}
\includegraphics[width=7cm,height=6cm]{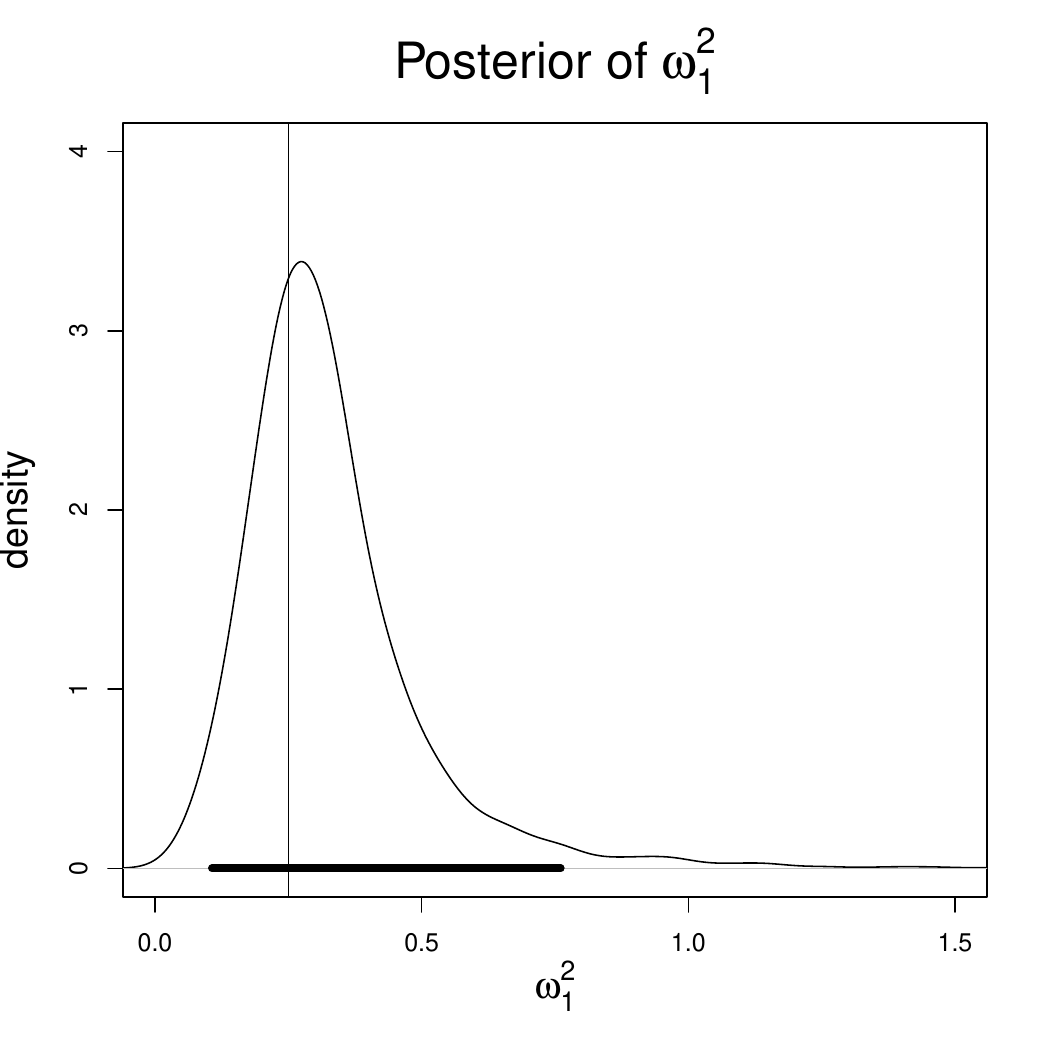}}
\hspace{2mm}
\subfigure[Posterior of $\omega^2_2$.]{ \label{fig:sim9_omegasq2}
\includegraphics[width=7cm,height=6cm]{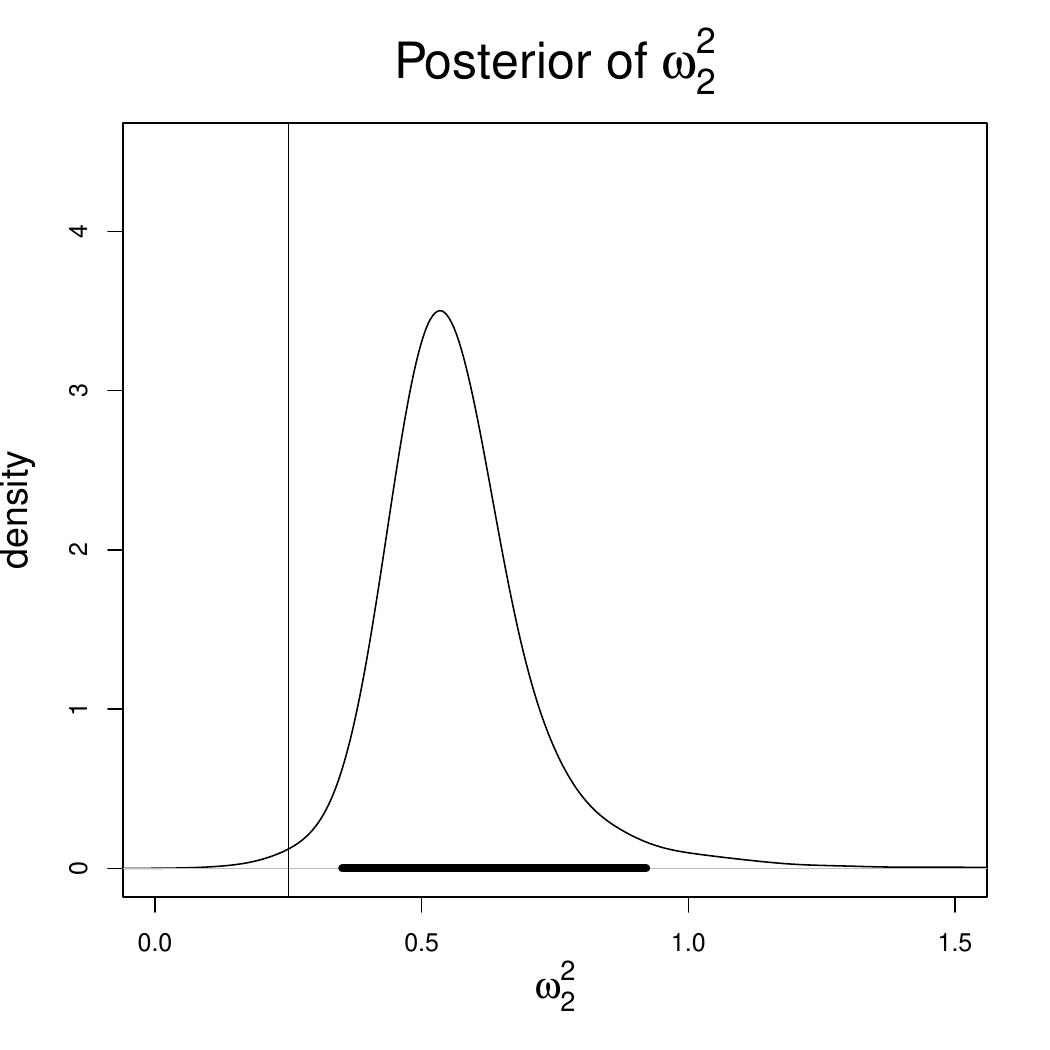}}\\
\vspace{2mm}
\subfigure[Posterior of $a_1$.]{ \label{fig:sim9_p1}
\includegraphics[width=7cm,height=6cm]{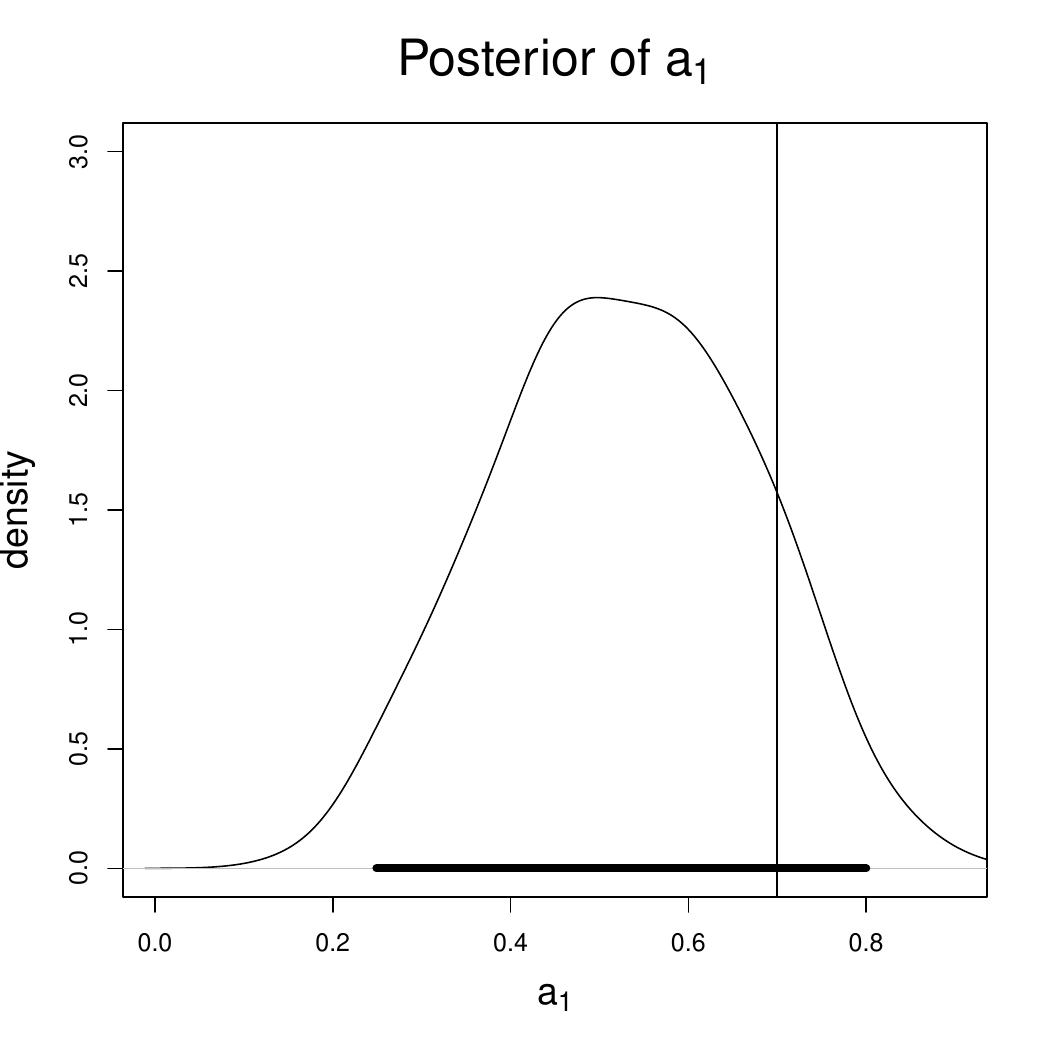}}
\hspace{2mm}
\subfigure[Posterior of $a_2$.]{ \label{fig:sim9_p2}
\includegraphics[width=7cm,height=6cm]{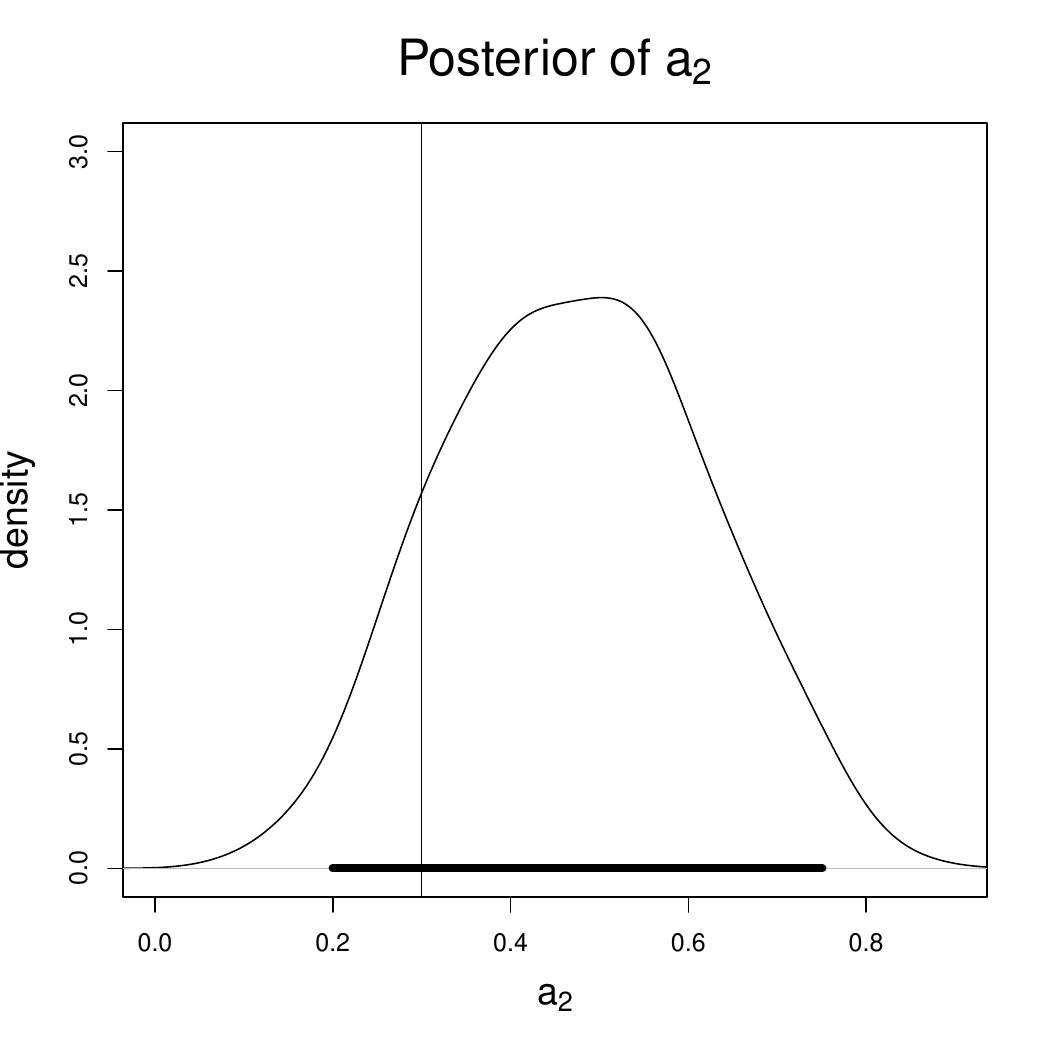}}
\caption{{\bf TTMCMC for $SDE_3$ and $\pi_2$ with $n=100$:} Posteriors of $M$, $\mu_1$, $\mu_2$, $\omega^2_1$, $\omega^2_2$, $a_1$ and $a_2$. The vertical lines stand
for the true values, while the thick horizontal lines denote the 95\% credible intervals.} 
\label{fig:sim9_posterior_plots}
\end{figure}

\begin{figure}
\centering
\subfigure[Histogram of $\hat\phi_i$'s with $n=100$. The vertical lines denote the true values.]{ \label{fig:hist1}
\includegraphics[width=7cm,height=6cm]{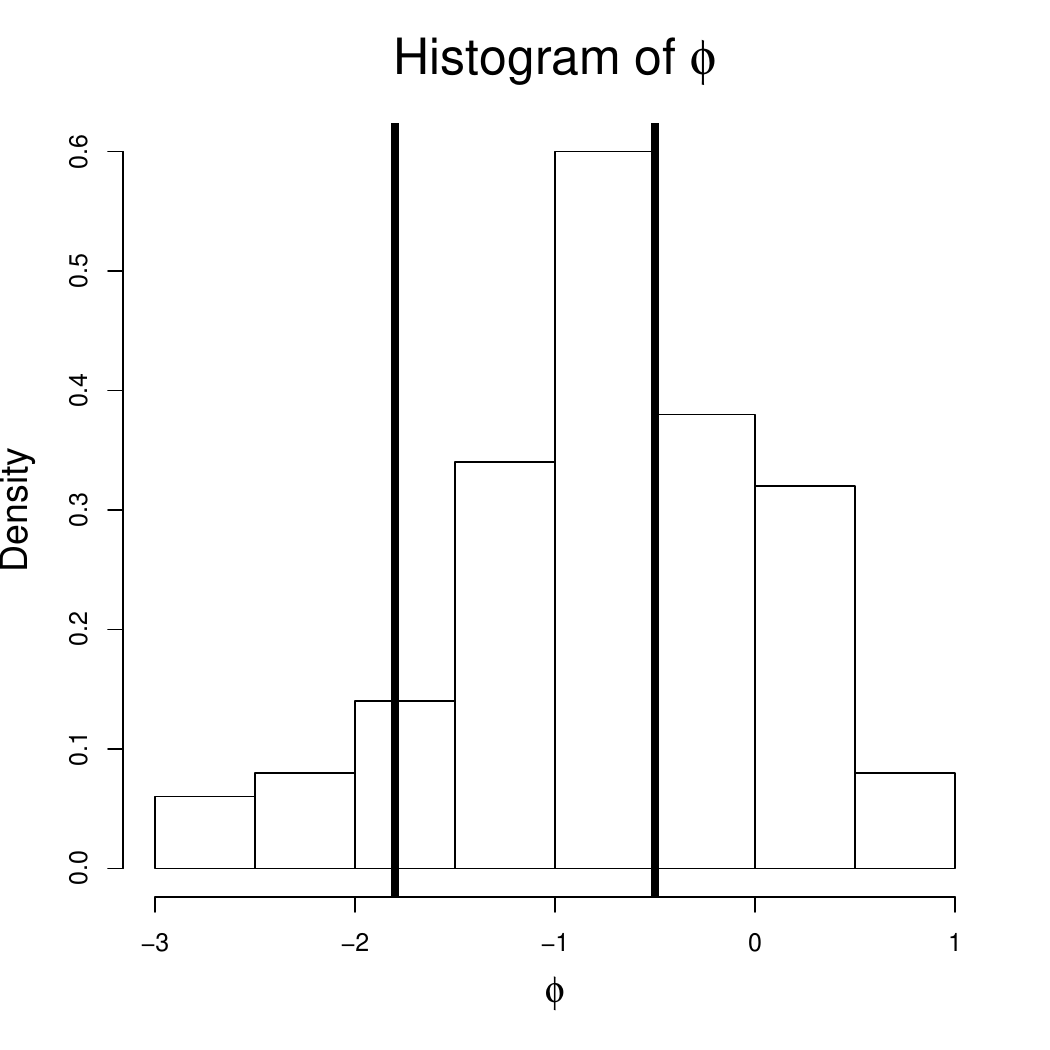}}
\hspace{2mm}
\subfigure[Histogram of $\hat\phi_i$'s, with $n=200$. The vertical lines denote the true values.]{ \label{fig:hist2}
\includegraphics[width=7cm,height=6cm]{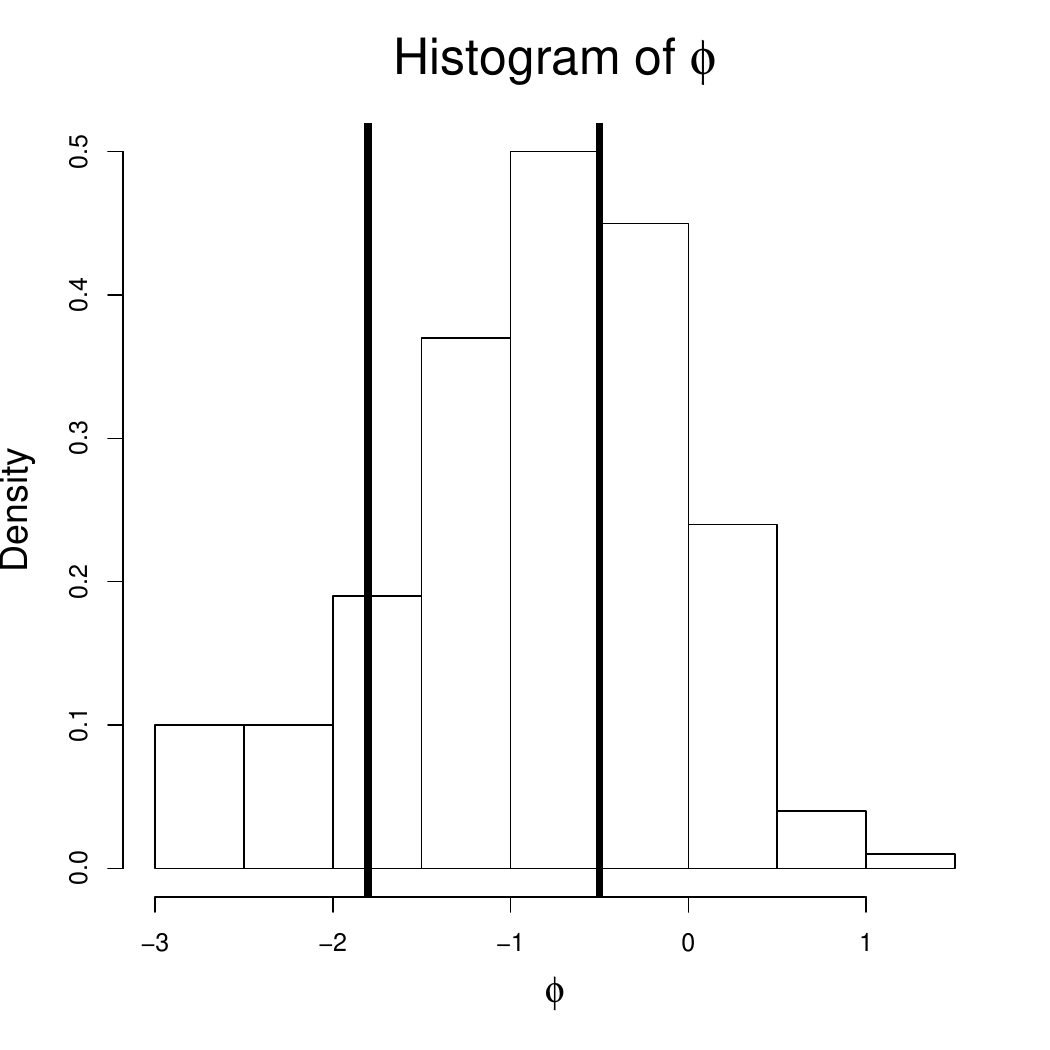}}
\end{figure}

In this setup, as usual, excellent mixing is demonstrated by the trace plots shown in Figure \ref{fig:sim9_trace_plots}. 
However, Figure \ref{fig:sim9_posterior_plots} shows that the true values of $\mu_2$ and $\omega^2_2$ do not fall within the respective 95\% credible intervals.
Also, $M=1,2,3,4$ receive posterior mass $0.6840$, $0.2657$, $0.0444$ and $0.0059$, respectively, showing that the true number of components, $M_0=2$, does not
receive the maximum posterior probability. An examination of the histogram of $\hat\phi_i=U_i/V_i$, shown in Figure \ref{fig:hist1}, where the dark vertical lines
represent the vales $-0.5$ and $-1.8$, reveals that it is indeed very difficult to identify $\mu_2=-1.8$, which is very ill-supported by the data. 
Consequently, the joint prior and posterior dependence between $\mu_2$ and $\omega^2_2$ seems to be responsible for pulling the significant posterior mass of 
$\omega^2_2$ away from the true value.

\begin{figure}
\centering
\subfigure[Trace plot of $M$.]{ \label{fig:sim10_trace_comp}
\includegraphics[width=7cm,height=5cm]{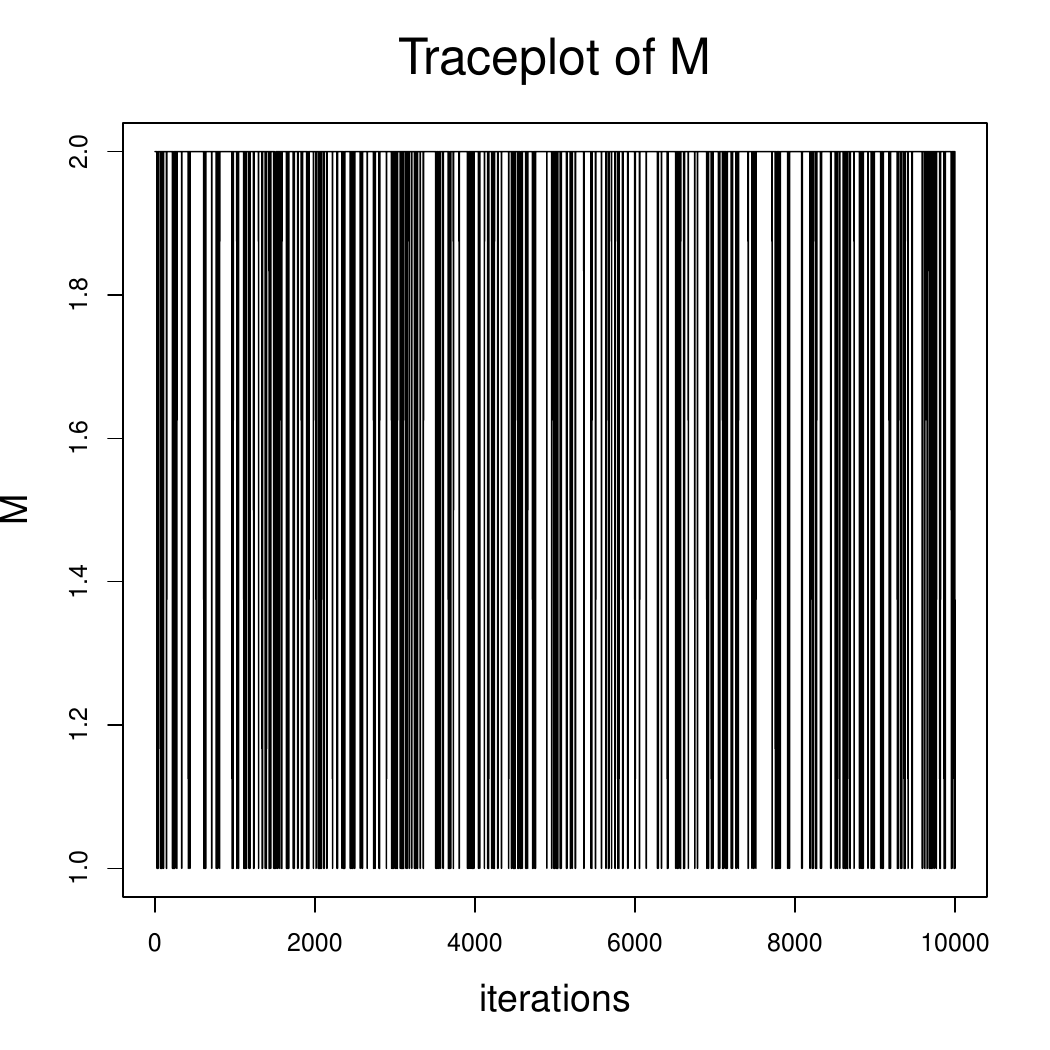}}
\hspace{2mm}
\subfigure[Trace plot of $\mu_1$.]{ \label{fig:sim10_trace_mu1}
\includegraphics[width=7cm,height=5cm]{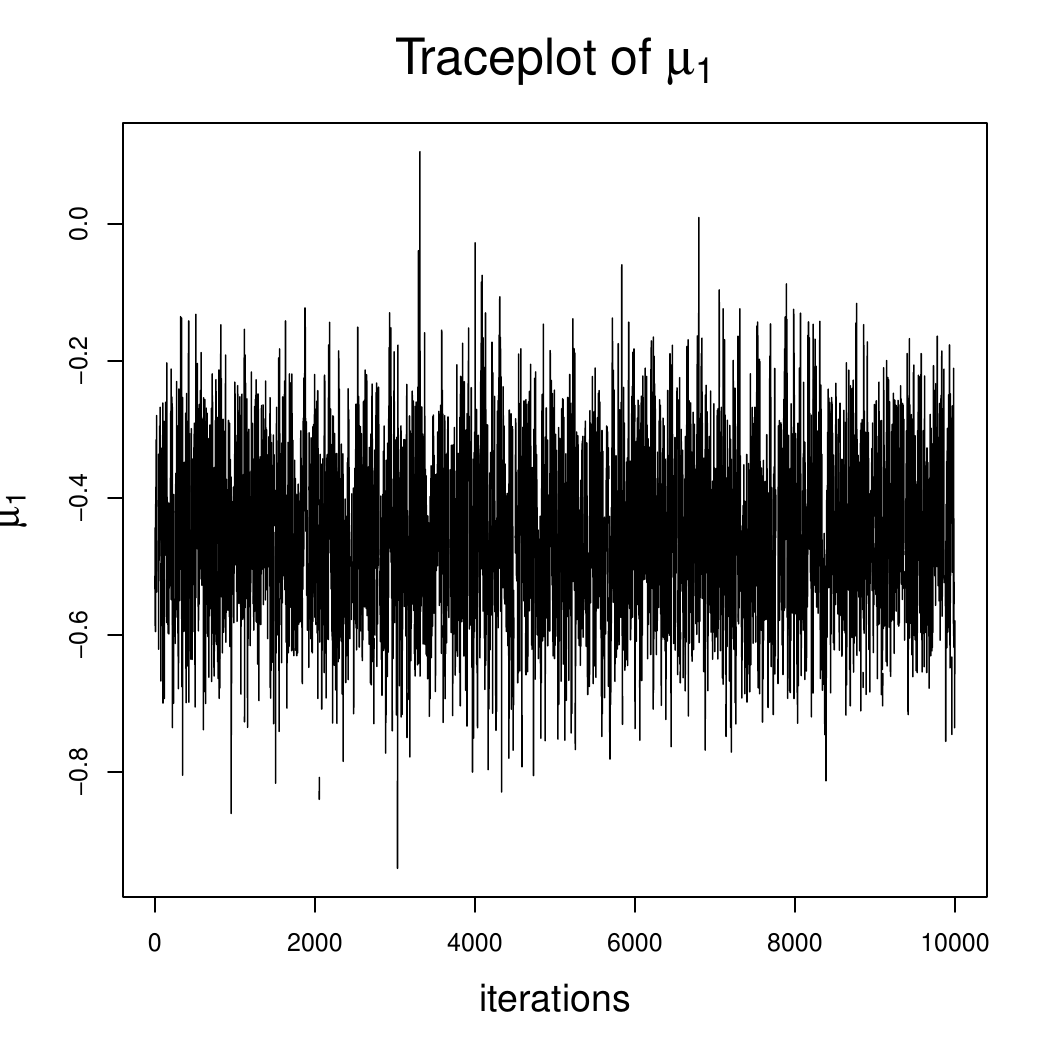}}\\
\vspace{2mm}
\subfigure[Trace plot of $\mu_2$.]{ \label{fig:sim10_trace_mu2}
\includegraphics[width=7cm,height=5cm]{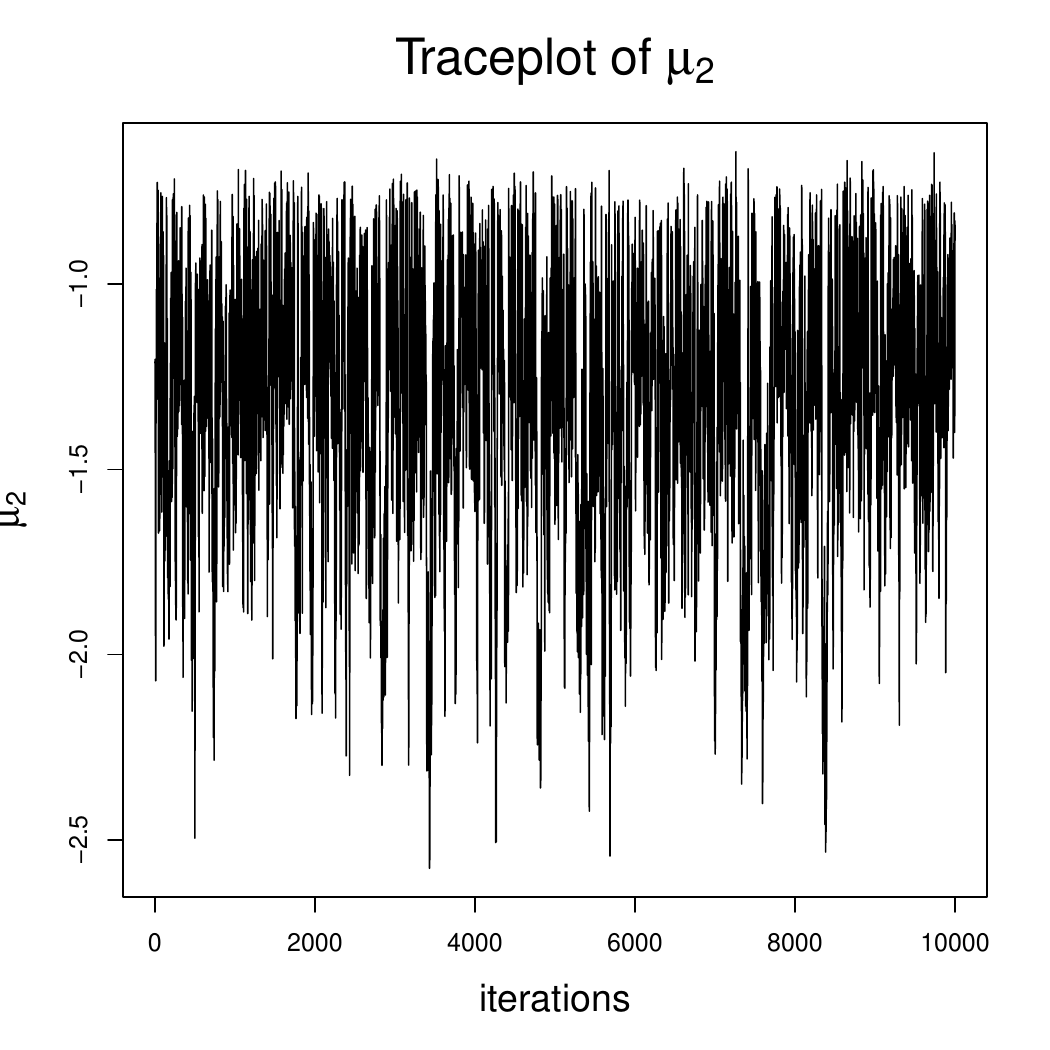}}
\vspace{2mm}
\subfigure[Trace plot of $\omega^2_1$.]{ \label{fig:sim10_trace_omegasq1}
\includegraphics[width=7cm,height=5cm]{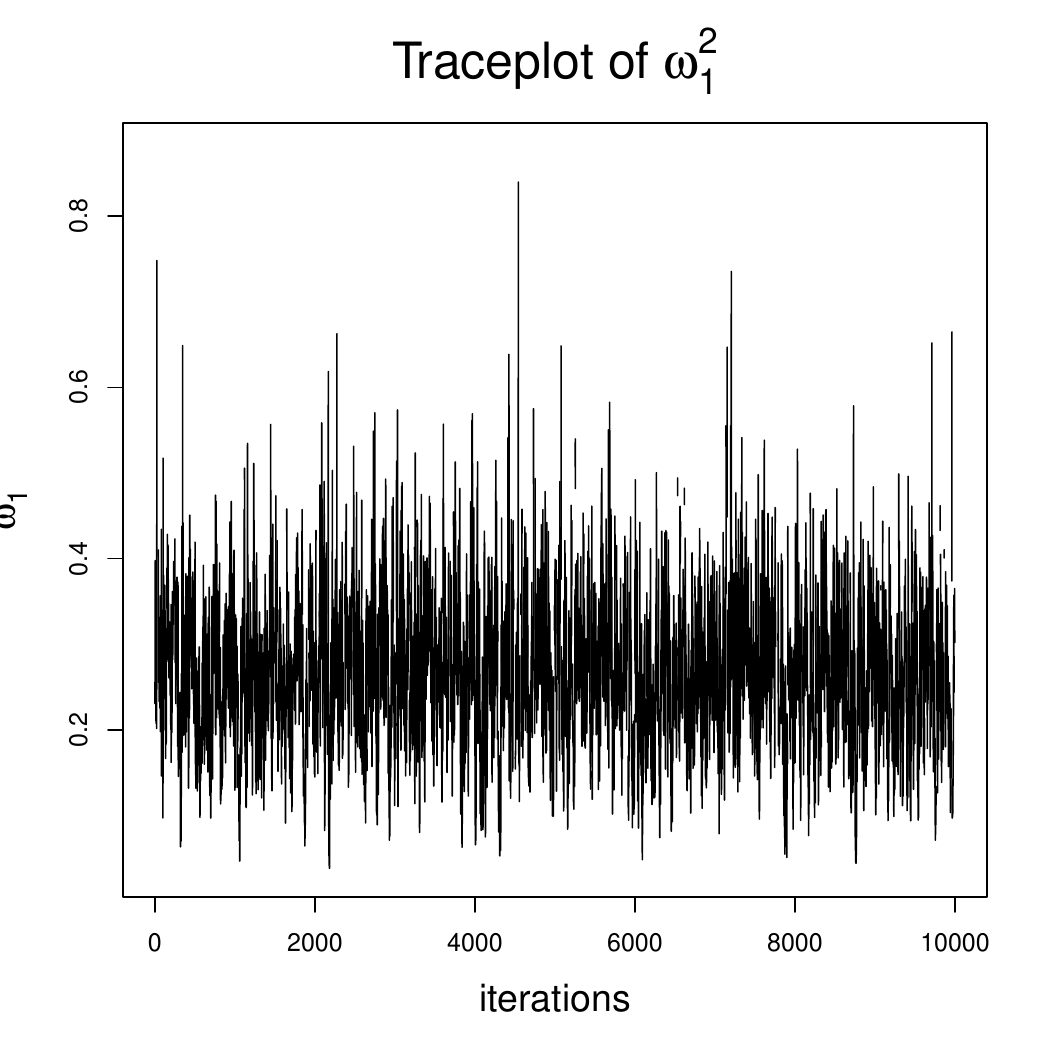}}\\
\vspace{2mm}
\subfigure[Trace plot of $\omega^2_2$.]{ \label{fig:sim10_trace_omegasq2}
\includegraphics[width=7cm,height=5cm]{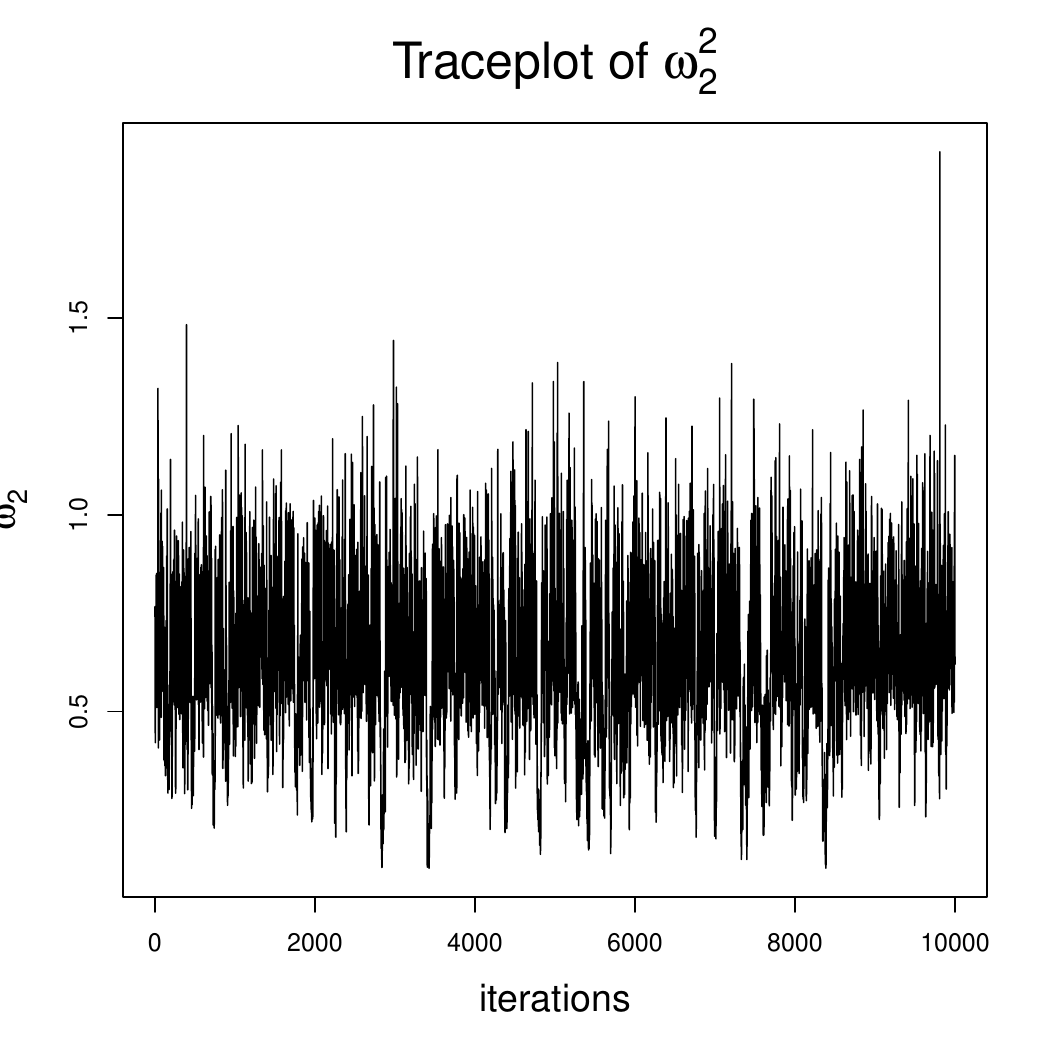}}
\hspace{2mm}
\subfigure[Trace plot of $a_1$.]{ \label{fig:sim10_trace_p1}
\includegraphics[width=7cm,height=5cm]{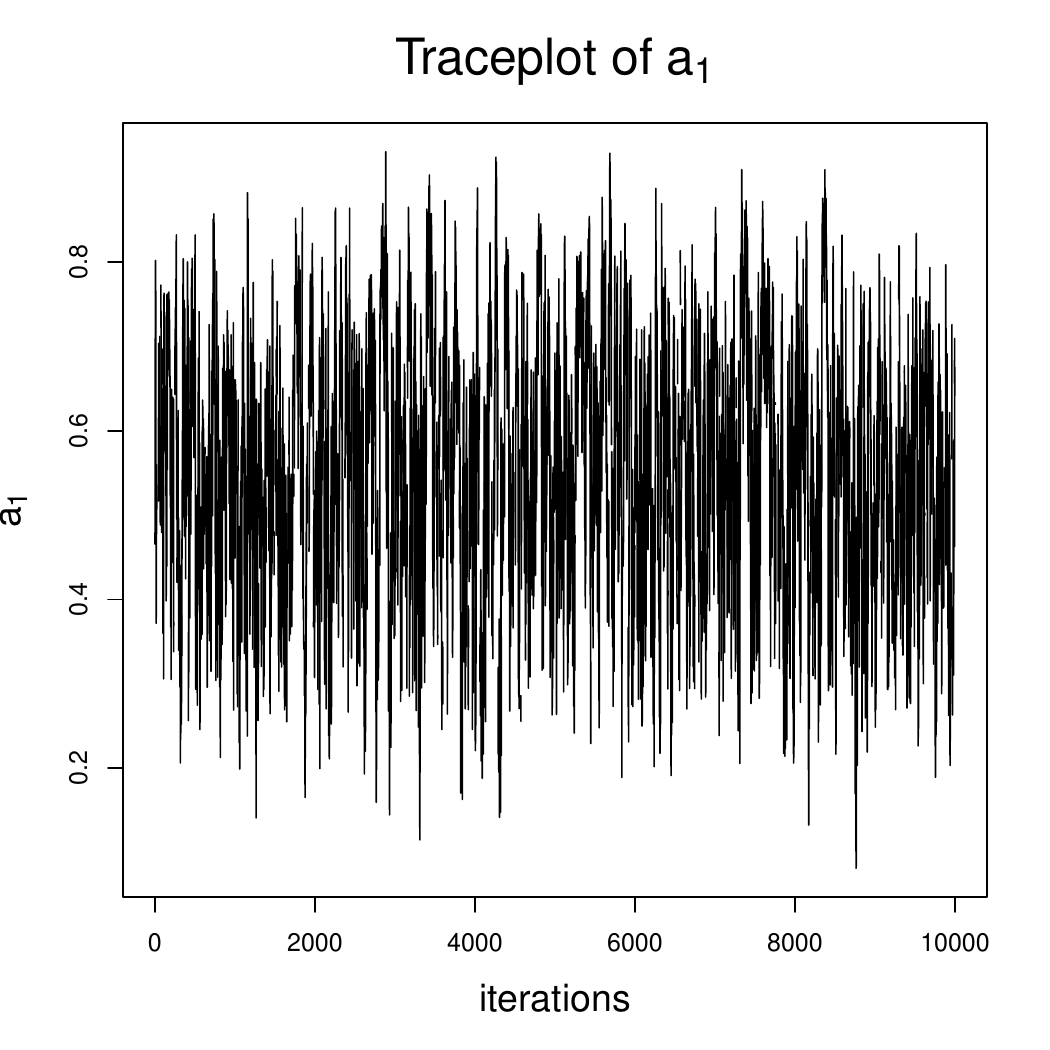}}\\
\vspace{2mm}
\subfigure[Trace plot of $a_2$.]{ \label{fig:sim10_trace_p2}
\includegraphics[width=7cm,height=5cm]{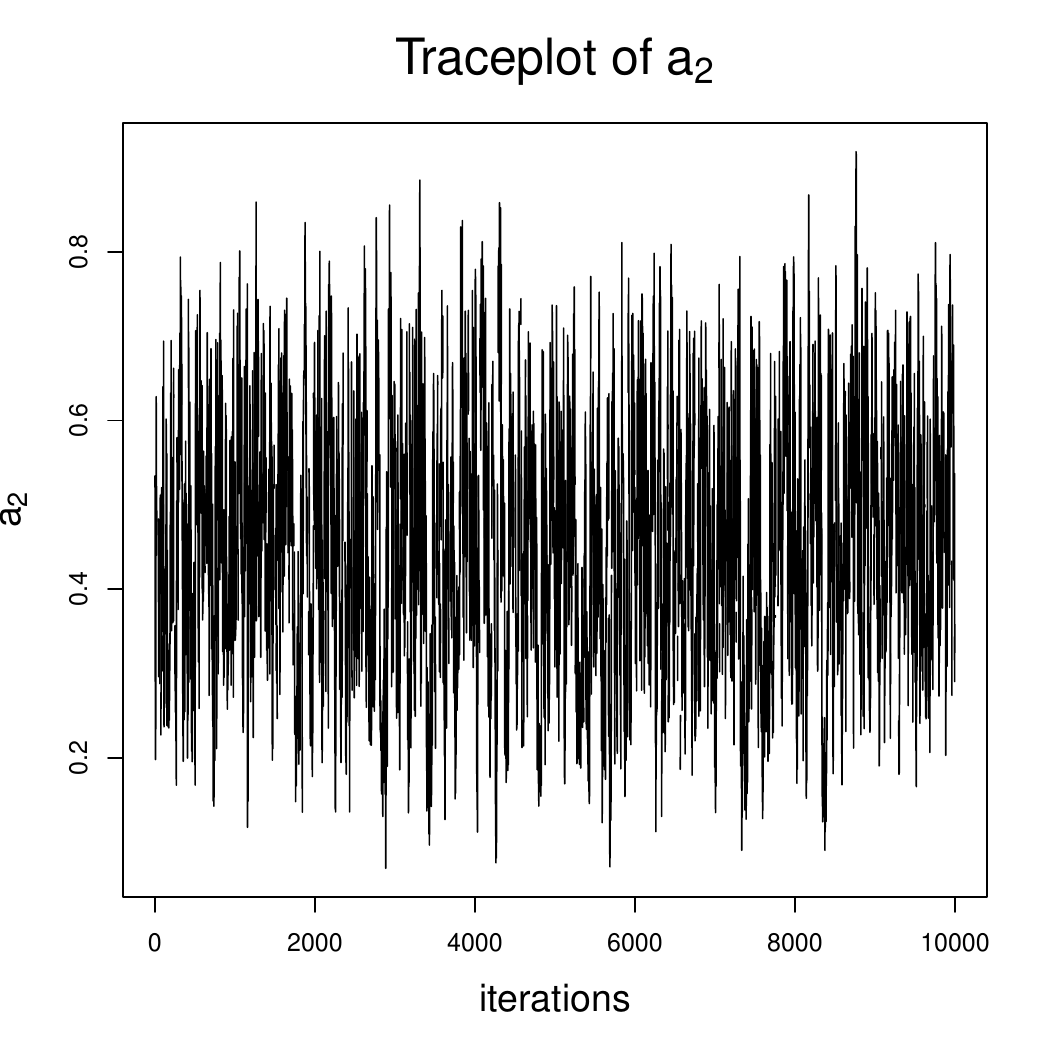}}
\caption{{\bf TTMCMC for $SDE_3$ and $\pi_2$ with $n=200$:} Trace plots of $M$, $\mu_1$, $\mu_2$, $\omega^2_1$, $\omega^2_2$, $a_1$ and $a_2$.} 
\label{fig:sim10_trace_plots}
\end{figure}

\begin{figure}
\centering
\subfigure[Posterior of $\mu_1$.]{ \label{fig:sim10_mu1}
\includegraphics[width=7cm,height=6cm]{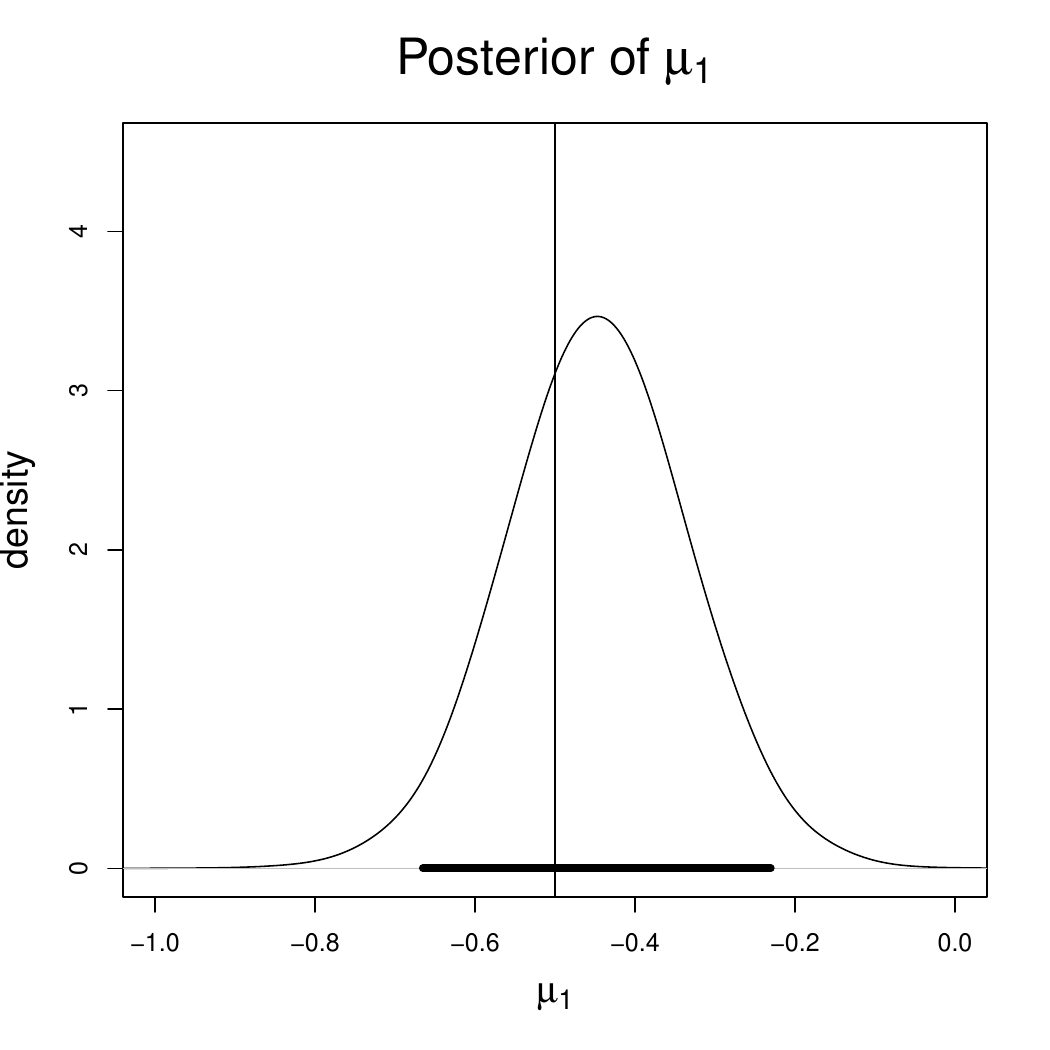}}
\hspace{2mm}
\subfigure[Posterior of $\mu_2$.]{ \label{fig:sim10_mu2}
\includegraphics[width=7cm,height=6cm]{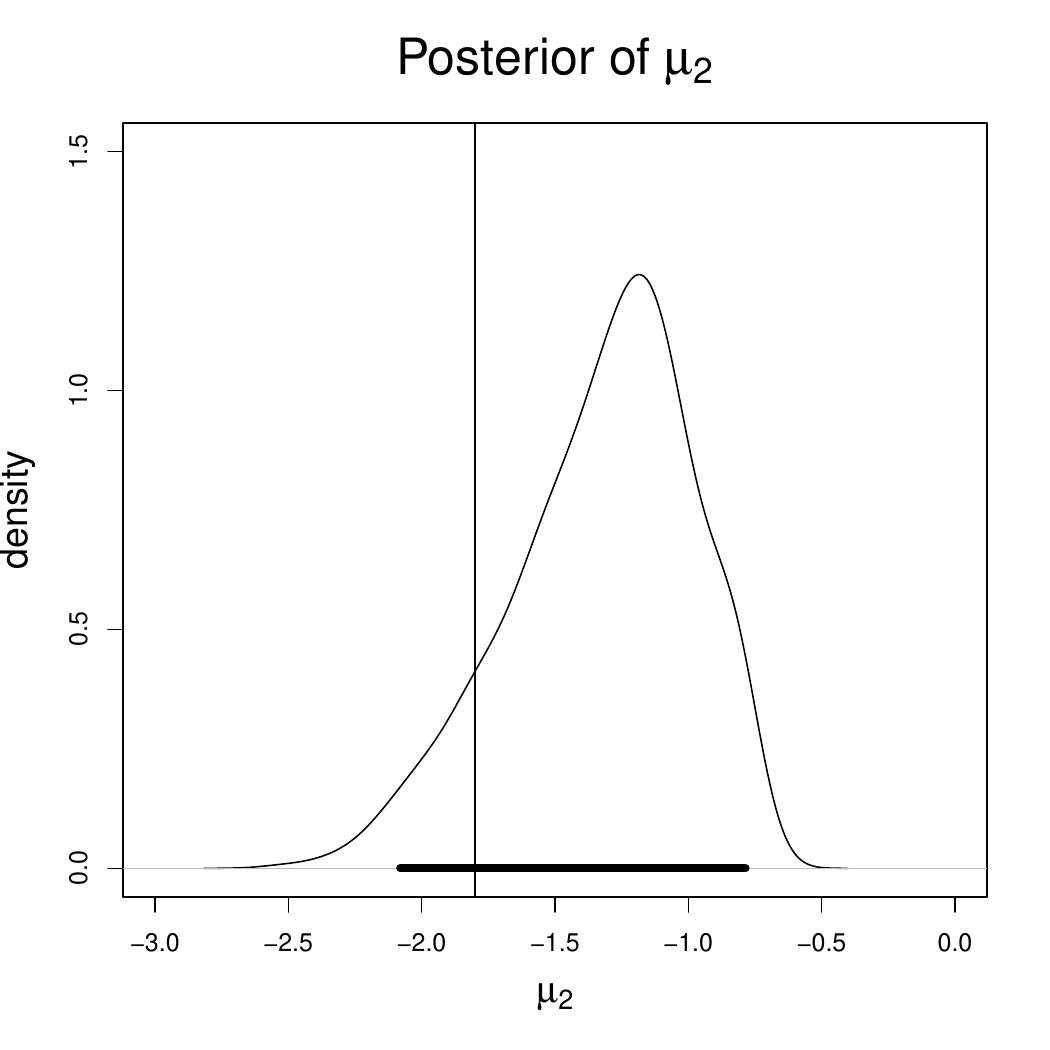}}\\
\vspace{2mm}
\subfigure[Posterior of $\omega^2_1$.]{ \label{fig:sim10_omegasq1}
\includegraphics[width=7cm,height=6cm]{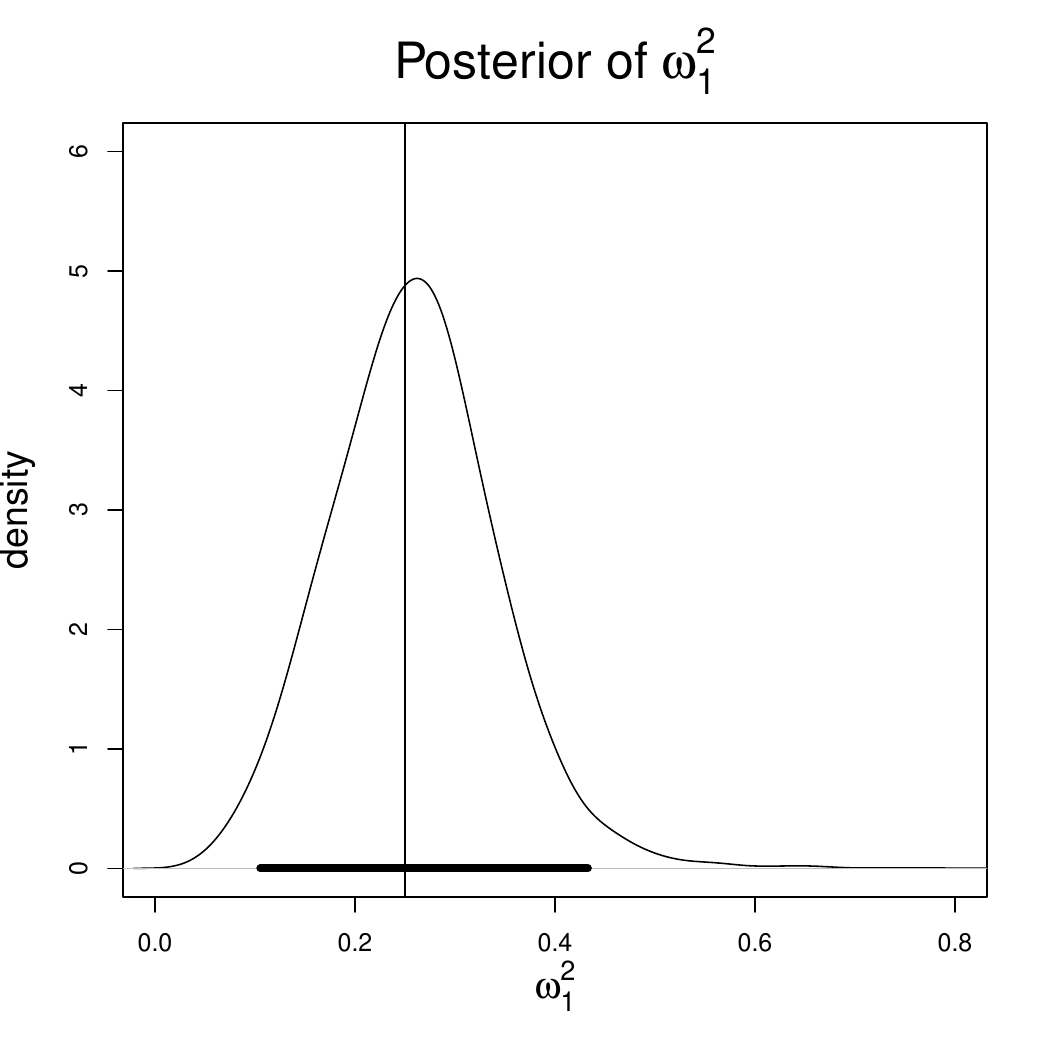}}
\hspace{2mm}
\subfigure[Posterior of $\omega^2_2$.]{ \label{fig:sim10_omegasq2}
\includegraphics[width=7cm,height=6cm]{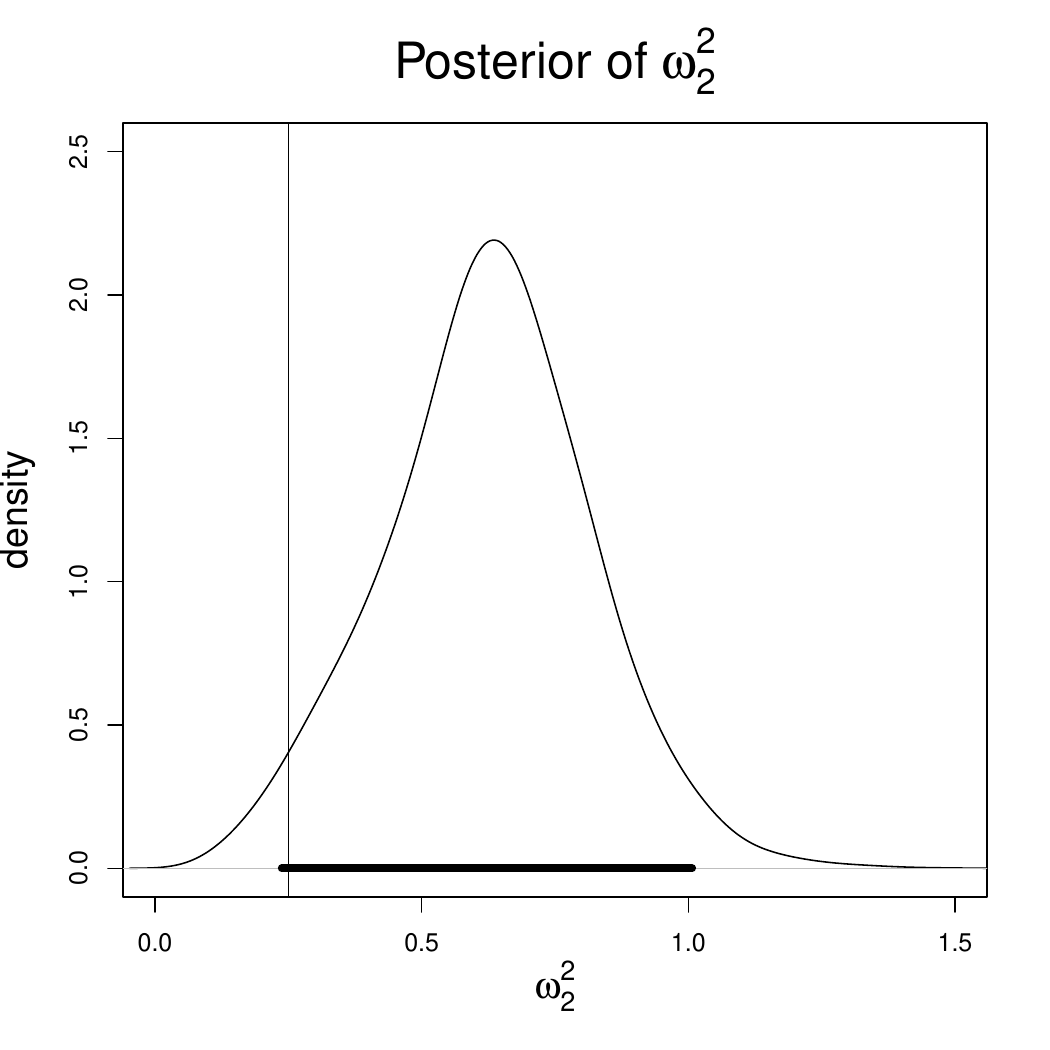}}\\
\vspace{2mm}
\subfigure[Posterior of $a_1$.]{ \label{fig:sim10_p1}
\includegraphics[width=7cm,height=6cm]{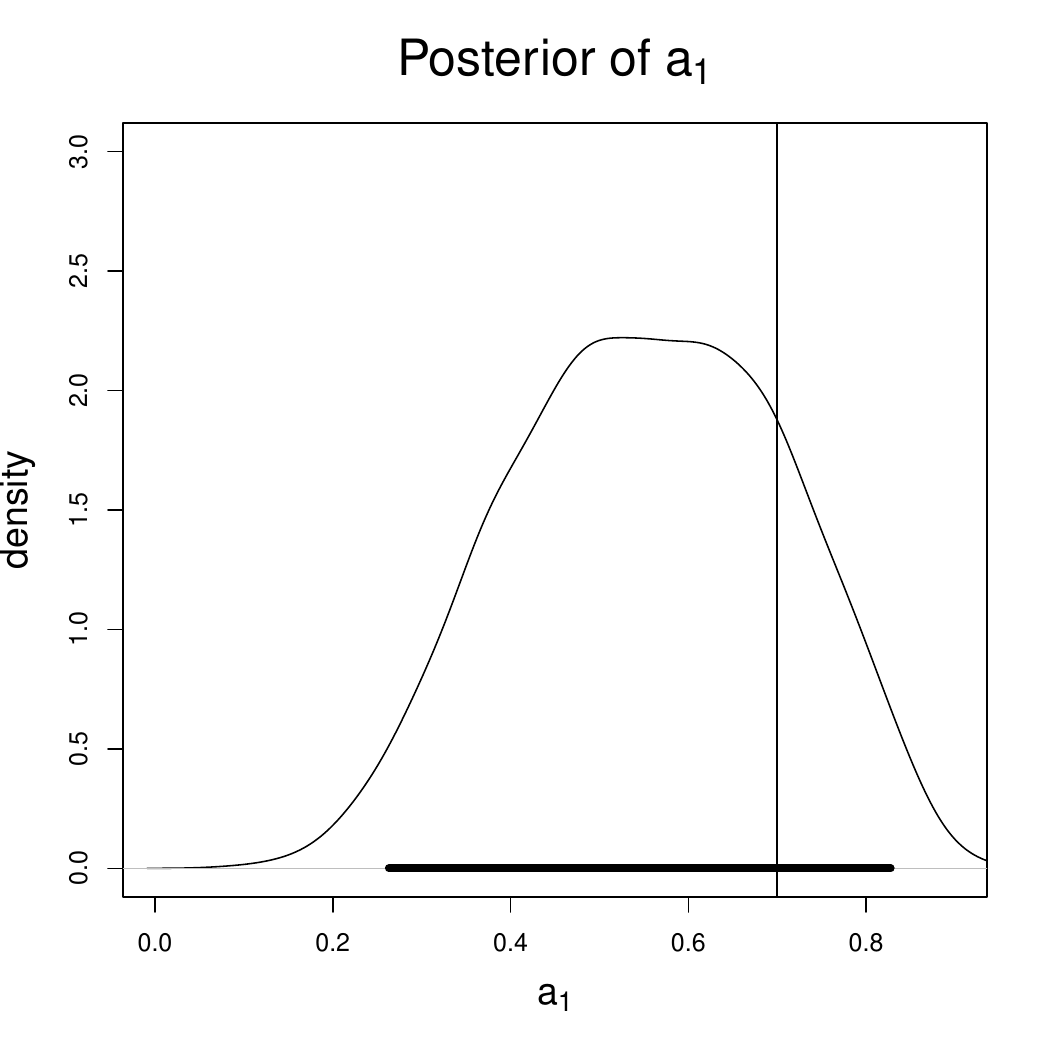}}
\hspace{2mm}
\subfigure[Posterior of $a_2$.]{ \label{fig:sim10_p2}
\includegraphics[width=7cm,height=6cm]{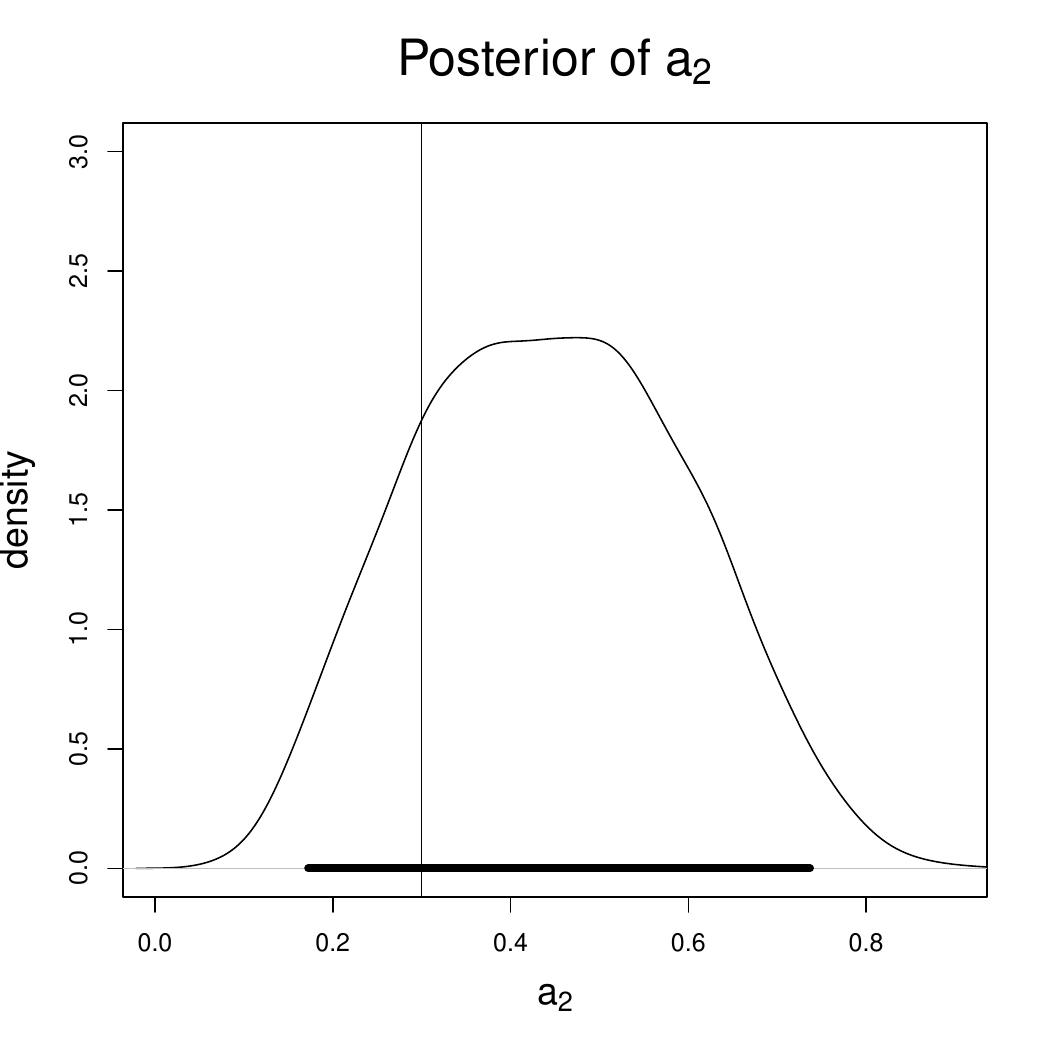}}
\caption{{\bf TTMCMC for $SDE_3$ and $\pi_2$ with $n=200$:} Posteriors of $M$, $\mu_1$, $\mu_2$, $\omega^2_1$, $\omega^2_2$, $a_1$ and $a_2$. The vertical lines stand
for the true values, while the thick horizontal lines denote the 95\% credible intervals.} 
\label{fig:sim10_posterior_plots}
\end{figure}

We repeated the experiment with $n=200$. The new histogram of $\hat\phi_i$, depicted in Figure \ref{fig:hist2}, seems to be more informative about $\mu_2=-1.8$
than Figure \ref{fig:hist2}. That the modified posterior probabilities of $M=1$ and $M=2$ are $0.0956$ and $0.9044$, respectively, with all other values of $M$
having zero posterior probability confirms this. Figure \ref{fig:sim10_trace_plots} suggests excellent mixing as before and Figure \ref{fig:sim10_posterior_plots} 
shows that all the true values are now within the 95\% credible intervals of the modified posterior.  

The overall acceptance rate of our algorithm in this case is about $5.3\%$.

\normalsize

\bibliographystyle{natbib}
\bibliography{irmcmc}

\end{document}